\documentclass[11 pt,a4paper,leqno]{article}
\usepackage[all]{xy}
\usepackage {makeidx}
\usepackage[french,english]{babel}

\makeindex

\RequirePackage{amsmath}
\RequirePackage{amssymb}
\RequirePackage{ifthen}
\RequirePackage{theorem}
\RequirePackage{pslatex}

\usepackage{amsmath}
\usepackage{amssymb}
\usepackage{ifthen}
\usepackage{theorem}
\usepackage{pslatex}

\title{\sc The space of monodromy data for the Jimbo-Sakai family of
$q$-difference equations}

\author{Yousuke Ohyama
  \footnote{Department of Mathematical Sciences, Tokushima University,
    2-1 Minamijyousanjima-cho, Tokushima 770-8506, Japan;
    E-mail:ohyama@tokushima-u.ac.jp},
Jean-Pierre Ramis
\footnote{Institut de France (Acad\'emie des Sciences) and
  Institut de Math\'ematiques de Toulouse,
  CNRS UMR 5219, Universit\'e Paul Sabatier (Toulouse 3),
  118 route de Narbonne, 31062 Toulouse CEDEX 9, France;
  E-mail: ramis.jean-pierre@wanadoo.fr},
Jacques Sauloy
\footnote{Toulouse; E-mail:jacquessauloy@gmail.com}}

\theorembodyfont{\sl} 
\newtheorem{thm}{Theorem}[section]
\newtheorem{lem}[thm]{Lemma}
\newtheorem{prop}[thm]{Proposition}
\newtheorem{cor}[thm]{Corollary}

\theorembodyfont{\rm}
\newtheorem{rmk}[thm]{Remark}
\newtheorem{exa}[thm]{Example}

\newtheorem{defn}[thm]{Definition}
\newtheorem{conj}[thm]{Conjecture}
\newtheorem{fact}[thm]{Fact}

\numberwithin{equation}{thm}

\def\Pr {\textsl{Proof. - }}
\def\Finpr{\hfill $\Box$ \\}
\def\Finprcourt{\hfill $\Box$}

\def\N{{\mathbf N}}                      
\def\Z{{\mathbf Z}}                      
\def\Q{{\mathbf Q}}                      
\def\C{{\mathbf C}}                      

\def\P{{\mathbf{P}}}                     
\def\PC{{\P^1(\C)}}                      
\def\PCg{{\P^1(\C)^\bullet}}              

\def\Id{{\text{Id}}}                     
\def\Ker{{\text{Ker}}}                   
\def\Sp{{\text{Sp}}}                     
\def\Tr{{\text{Tr}}}                     
\def\Im{{\text{Im}}}                     
\def\ii{{\text{i}}}                      
\def\lmod{\left|}                        
\def\rmod{\right|}                       

\def\Vect{{\text{Vect}}}                 
\def\card{{\text{card}}}                 
\def\Diag{{\text{Diag}}}                 

\def\Sr{{\mathbf{S}}}                    
\def\Do{{\overset{\circ}{\text{D}}}}     

\def\M{{\mathcal{M}}}                    
\def\div{{\text{div}}}                    

\def\GL{{\text{GL}}}                     
\def\GLn{{\text{GL}_{n}}}                 
\def\Mat{{\text{Mat}}}                   
\def\Matn{{\text{Mat}_{n}}}               
\def\GLnc{\GLn(\C)}                      
\def\Matnc{\Matn(\C)}                    
\def\Dnc{{\text{D}_n(\C)}}               
\def\Ddc{{\text{D}_2(\C)}}               
\def\SLdc{{\text{SL}_2(\C)}}             
\def\SLc{{\text{SL}}_{2}(\C)}
\def\slc{{{sl}}_{2}(\C)}

\def\Hom{{\text{Hom}}}                   

\def\ie{{\emph{i.e.}}}
\def\eg{{\emph{e.g.}}}
\def\tq{\mid}
\def\Tq{\bigm|}
\def\TQ{\Bigm|}
\def\cf{{\emph{cf.}}}

\def\Rg{{\C[x]}}
\def\Kg{{\C(x)}}
\def\Ra{{\C\{x\}}}                              
\def\Ka{{\C(\{x\})}}                            
\def\Rf{{\C[[x]]}}                              
\def\Kf{{\C((x))}}                              
\def\Rwg{{\mathcal{O}(\C^{*})}}
\def\Kwg{{\mathcal{M}(\C^{*})}}
\def\Rw{{\mathcal{O}(\C^{*},0)}}
\def\Kw{{\mathcal{M}(\C^{*},0)}}

\def\Rai{{\C\{1/x\}}}                            
\def\Kai{{\C(\{1/x\})}}                          
\def\Rfi{{\C[[1/x]]}}                            
\def\Kfi{{\C((1/x))}}                            

\def\thq{{\theta_{q}}}                            
\def\thqp{{\theta_{q}'}}                          
\def\sq{\sigma_q}                                 
\def\Eq{{\mathbf{E}_{q}}}                         
\def\Eqp{{\mathbf{E}_{q}^*}}                      
\def\Eqg{{\mathbf{E}_{q}^\bullet}}          

\def\X{{\mathcal{X}}}
\def\Y{{\mathcal{Y}}}
\def\Cq{{\mathcal{C}_q}}
\def\Cqg{{\mathcal{C}_q^\bullet}}

\def\x{{\underline{x}}}

\def\E{{\mathcal{E}}}                   
\def\F{{\mathcal{F}}}               
\def\G{{\mathcal{G}}}

\def\Rep{{\text{Rep}}}

\def\Sym{{\mathbf{S^\bullet}}}

\begin{document}

\maketitle


\hfill \emph{To the memory of our friend Hiroshi Umemura}


\begin{abstract}
We formulate a geometric Riemann-Hilbert correspondence that applies to the
derivation by Jimbo and Sakai of equation $q$-PVI from ``isomonodromy'' 
conditions. This is a step within work in progress towards the application
of $q$-isomonodromy and $q$-isoStokes to $q$-Painlev\'e.
\end{abstract}

\selectlanguage{french}

\begin{abstract}
Nous formulons une correspondance de Riemann-Hilbert g\'eom\'etrique qui s'applique
\`a la d\'erivation par Jimbo et Sakai de l'\'equation $q$-PVI \`a partir de conditions
\og d'isomonodromie \fg. C'est une \'etape d'un travail en cours en vue de
l'application de la $q$-isomonodromie et des $q$-isoStokes \`a $q$-Painlev\'e.
\end{abstract}

\selectlanguage{english}

\tableofcontents




\section{Introduction}


\subsection{Position of the problem}

This paper is a first step towards a formulation of a Riemann-Hilbert correspondence
for the $q$-Painlev\'e equations, the $q$-analogs of the classical differential Painlev\'e
equations. In this very first step we limit ourselves to the case of the Jimbo-Sakai
$q$-PVI equation \cite{JimboSakai} for fixed generic values of the ``local parameters". \\

For the $q$-Painlev\'e equations, according to the pioneering work of Hidetaka Sakai
\cite{Sak}, we have an exhaustive information \emph{``on the left hand side"} of the $q$-analog
of the Riemann-Hilbert map: the $q$-analogs of the \emph{Okamoto spaces of initial
conditions} are open rational surfaces obtained by blowing up {\small $\P^2(\C)$} in nine
points and removing some lines. But \emph{``on the right hand side"} the $q$-analogs of the
\emph{character varieties} with their structure of open cubic surface are unknown. 
Our initial aim was to fill this gap in the $q$-PVI case. \\

We consider, just as Birkhoff did, the family of $q$-difference systems $\sigma_q X=AX$
\[
\label{birksyst}
A = A_0 + \cdots + x^n A_n \in \Matn(\Rg), ~~ A_0,\ldots,A_n \in \Matnc, ~~ A_0,A_n \in \GLnc.
\]
Following Birkhoff \cite{Birkhoff1}, we associate to a system (\ref{birksyst}) a matrix
$M$ (a variant of Birkhoff connection matrix) representing in some sense (this will be
commented in more detail in \ref{subsection:firstage} and \ref{subsection:secondage})
some kind of $q$-analog of the monodromy data for differential equations. This map induces
\emph{an isomorphism}, the Riemann-Hilbert-Birkhoff correspondence,
between the systems modulo rational gauge equivalence on one side and the ``matrices of
monodromy data'' modulo a natural equivalence on the other side. \\

The Jimbo-Sakai family studied in \cite{JimboSakai} is associated to a subspace of
the space:
\[
\left\{A_0 + x A_1 + x^2 A_2 \in \Mat_2(\Rg) \Tq
A_0,A_2 \in \GL_2(\C), ~~ A_1 \in \Mat_2(\C)\right\}.
\]
The subspace is restrained by conditions on the ``local data", \ie\ the conjugacy classes
of $A_0$ and $A_2$ (actually their spectra, for they are assumed to be semi-simple).
Sakai gave a direct description of the space of equations
$\sq X = \left(A_0 + x A_1 + x^2 A_2\right) X$
as an open rational surface; this is what we consider as the ``left hand side''
of the Riemann-Hilbert correspondence. The ``right hand side'' is our space $\F$ of
``monodromy data'' modulo equivalence. We call it the \emph{space of monodromy data} or
\emph{$q$-character variety}. Then we give a first geometric description of $\F$ as an algebraic
surface. In particular we give an embedding of $\F$ into $\left(\P^1(\C)\right)^4$. \\

In the second part of our paper we introduce a new tool: we call it \emph{the Mano
decomposition} and use it to get a more precise description of the algebraic variety
$\F$. \\

This extremely useful process was inspired to us by the paper \cite{Mano} of
Toshiyuki Mano. The equations that appear in the Jimbo-Sakai family  can, in some
sense, be \emph{split} into $q$-hypergeometric components and the corresponding
monodromy matrix $M$ can be \emph{split} into the monodromy matrices of these
$q$-hypergeometric components. So, Mano decomposition can be understood as providing
a splitting of the \emph{global} monodromy around the four intermediate singularities
into \emph{local} monodromies around two pairs of singularities. \\

Mano decompositions allow us to describe parameterizations (the $q$-pants parameterizations)
of $F$. They are $q$-analogs of the classical parametrizations of the Fricke cubic surface
(the character variety of PVI) associated to a pant decomposition \cite{ILT}. \\
 
We tried to ``identify" $\F$ among the ``classical" surfaces. We only got some partial
informations\footnote{In particular about elliptic fibrations.} allowing some guesses:
in particular $\F$ could be a Zariski open subset of a K3 surface.


\subsection{Isomonodromy and Painlev\'e equations}
\label{subsection:IsomonodromyandPainleve}

The purpose of this part is to recall briefly some classical results on Painlev\'e
equations and linear representations. It is necessary for the detailed description
of some (not so evident) analogies between the differential case and the $q$-difference
case that we will present below. \\

Theory of the Painlev\'e differential equations has developed through two very different
lines. One is the classification of second order algebraic ordinary differential equations  
which satisfy the Painlev\'e property (\ie\ the movable singularities are poles).
The other one is a deformation theory of linear ordinary differential equations:
one asks to move the coefficients of the equation without changing its monodromy data
(or more generally its generalized monodromy data). At the very beginning of the XX-th
century, P. Painlev\'e \cite{Pain} and B. Gambier initiated the first line and R. Fuchs
\cite{Fu} initiated the second. In the present paper everything is in the spirit of
this second line.


\subsubsection{Representations of the free group of rank $3$ into $\SLc$.
Character varieties}
\label{subsubiso}

We present the character varieties in elementary purely algebraic terms (no differential
equations here). At the end of the paragraph we will introduce some topology: fundamental
groups of punctured spheres. \\

We denote $\Gamma_3:= \langle u_0,u_t,u_1 \rangle$ the free group of rank $3$ generated by
the letters\footnote{The motivation for the indices $0,t,1,\infty$ will appear
at the end of this paragraph.} $u_0,~u_t,~u_1$. It is identified with the free group 
$\langle u_0,u_t,u_1,u_\infty\vert\,  u_0u_tu_1u_\infty=1 \rangle$ 
generated by $u_0,~u_t,~u_1,~u_\infty$ up to the relation $u_0 u_t u_1 u_\infty = 1$. \\

Let $\rho: \Gamma_3\rightarrow \SLc$ be a linear representation.
We set $M_l:=\rho(u_l)$ ($l=0,t,1,\infty$). We denote 
$e_l$ and $e_l^{-1}$ ($l=0,t,1,\infty$) the eigenvalues of $M_l$.
The representation $\rho$ can be identified with 
$(M_0,M_t,M_1)\in \left(\SLc\right)^3$. Therefore the set of such representations 
$\Hom\left(\Gamma_3,\SLc\right)$ modulo the adjoint action of ${\text{SL}}_{2}(\C)$
can be identified with 
$\left(\SLc\right)^3/\SLc$ (the set of triples of matrices up to overall conjugation)~:
\[
\Hom(\Gamma_3,\SLc)/\SLc = \left(\SLc\right)^3/\SLc;
\]
$\left(\SLc\right)^3$ is a complex affine variety of dimension $9$. \\

To a representation $\rho:\Gamma_3\rightarrow \SLc$ we associate its seven
\emph{Fricke coordinates} (or trace coordinates), the four ``parameters"~:
\[
a_l:=\Tr~M_l=e_l+e_l^{-1}, \quad l=0,t,1,\infty
\]
and the three ~``variables":
\[
X_0=\Tr~M_1M_t, \quad  X_t=\Tr~M_1M_0, \quad  
X_1=\Tr~M_tM_0.
\]
These seven coordinates satisfy the \emph{Fricke relation}
$F(X,a)=0$ (cf. \cite{Mag}), where~:
\begin{equation}
\label{equacubicsurf1}
F(X,a) := F\left((X_0,X_t,X_1);(a_0,a_t,a_1,a_\infty)\right) :=
X_0X_tX_1+X_0^2+X_t^2+X_1^2-A_0X_0-A_tX_t-A_1X_1
+A_\infty,
\end{equation}
with~:
\begin{equation}
\label{desaauxA}
A_i := a_ia_\infty+a_ja_k, ~~\text{for} ~~ i=0,t,1, ~~ \text{and}~
A_\infty := a_0a_ta_1a_\infty+a_0^2+a_t^2+a_1^2+a_\infty^2-4.
\end{equation}

The seven Fricke coordinates of $\rho$ are clearly invariant by equivalence of
representations. Then, using the seven Fricke coordinates, we get an algebraic map
from $\left(\SLc\right)^3/\SLc$ to $\C^7$. The image is the six dimensional quartic
hypersurface of $\C^7$ defined by the equation $F(X,a)=0$. \\

We fix the parameter $a$ and denote ${\mathcal S}(a)$ or ${\mathcal S}_{A_0,A_t,A_1,A_\infty}$
or\footnote{For reasons that will appear in the next paragraph.} ${\mathcal S}_{VI}(a)$
the cubic surface of $\C^3$ defined by the equation $F(X,a)=0$. We call this surface
the \emph{character variety} of PVI. \\

By a theorem of Fricke, Klein and Vogt \cite{Mag,Gold} 
the equivalence class of an \emph{irreducible} representation is completely determined
by its seven Fricke coordinates. \\

We denote $\overline{\mathcal S}(a)$ the projective completion\footnote{As an abstract
algebraic surface it is a del Pezzo surface of degree $3$.} of
${\mathcal S}(a)$ in $\mathbf{P}^3(\C)$. The family 
$\{\overline{\mathcal S}(a)\}_{a\in \C^4}$ contains all  \emph{smooth} projective cubic
surfaces (up to linear transformations). The list of projective cubic surfaces was
given by Schl\"afli \cite{Schla} over a century ago. For this list we refer to
\cite{BW}: \cf\ table $4$, page $255$. There are $20$ families of singular projective
cubic surfaces. An excellent reference is \cite[\S 3, page 11)]{IwaU}. \\

The surface ${\mathcal S}(a)$ is \emph{simply connected} \cite{CL}. It can be smooth
or have singular points according to the values of $a$. The number of singular points
is at most $4$. Singular points of ${\mathcal S}(a)$ appear from semi-stable
representations which are of two kinds~:
\begin{itemize}
\item 
Either $M_l=\pm I_2$ \big(that is $\rho(u_l)$ belongs to the center of
${\text{SL}}_{2}(\C)$\big) for some $l=0,t,1,\infty$, hence $e_l=\pm 1$ and 
$a_l=\pm 2$. This case is called \emph{the resonant case}.
\item
Or the representation is reducible. 
This condition can be translated into an algebraic condition on $a$,
\cf\ \cite{Iwa}, \cite{Kli} page 22, \cite{Mazz}; we have~:
\begin{equation}
\label{reduct3}
e_0\, e_t^{\pm 1}e_1^{\pm 1}e_\infty^{\pm 1}=1
\end{equation}
for some triple of signs.
\end{itemize}

\bigskip

An example of a singular cubic surface with $4$ singular points is the Cayley
cubic \cite{CL}. We get it for $(A_0,A_t,A_1,A_\infty)=(0,0,0,-4)$ (this is true
either if $a=(0,0,0,0)$ or if $a=(\pm 2,\pm 2,\pm 2,\pm 2)$ with product $-16$):
\begin{equation}
\label{equacubicsurf2}
X_0X_tX_1+X_0^2+X_t^2+X_1^2-4=0.
\end{equation}


We denote~:
$
F_{X_i}:=\frac{\partial F(X,a)}{\partial X_i}=X_jX_k+2X_i-A_i.
$
The character variety 
${\mathcal S}_{VI}(a)={\mathcal S}_{A_0,A_t,A_1,A_\infty}$
is equipped with a ``natural" algebraic symplectic form (Poincar\'e residue)~:
\begin{equation}
\label{symplectic}
\omega_{VI,a}:=\frac{dX_t\wedge dX_0}{2i\pi F_{X_1}}
=\frac{dX_1\wedge dX_t}{2i\pi F_{X_0}}
=\frac{dX_1\wedge dX_t}{2i\pi F_{X_t}}
\end{equation}
We have\footnote{The motivation for the choice of the factor
$-\frac{1}{2i\pi}$ will appear in the next paragraph, cf. footnote\ref{symplectomorphism}.}
$dF \wedge \omega_{VI,a}=-\frac{1}{2i\pi} dX_0\wedge dX_t\wedge dX_1$.
The Poisson bracket associated to $-2i\pi\, \omega_{VI,a}$ is the \emph{Goldman bracket}
defined by~: $\{X_i,X_j\}=F_{X_k}$, and circular permutations. \\

Let $S_4^2$ be the four punctured sphere. Its fundamental group $\pi_1(S_4^2)$ is
isomorphic to a free group of rank $3$: we can choose as generators the homotopy
classes of three simple loops turning around three punctures. \\

Therefore we can apply the preceding results to the study of  equivalence classes
of representations of $\pi_1(S_4^2)$ into $\SLc$. It is a purely topological matter
and the choice of the punctures is indifferent up to an homeomorphism. But in the
following we will need the complex structure: $S^2 = \PC$. Then, starting from $4$
arbitrary punctures, up to a M\"{o}bius transformation,  we can choose as punctures
$0,t,1,\infty$ for some value of $t$. This explains our initial notation. \\

For $t \in \PC \setminus \{0,1,\infty\}$ we set~:
\[
\widetilde{\Rep}_t:= 
\Hom\left(\pi_1\left(\PC\setminus \{0,t,1,\infty\}\right),\SLc\right)/\SLc.
\]
For small changes\footnote{More precisely if $t$ remains in an open disc of the
$3$-punctured sphere $\PC \setminus \{0,1,\infty\}$.} of $t$, the group
$\pi_1\left(\PC\setminus \{0,t,1,\infty\}\right)$ remains constant, more precisely
there exist canonical isomorphisms~:
\[
\pi_1\left(\PC\setminus \{0,t_1,1,\infty\}\right)\rightarrow
\pi_1\left(\PC\setminus \{0,t_2,1,\infty\}\right),
\]
therefore there are canonical isomorphisms
$\widetilde{\Rep}_{t_2}\rightarrow \widetilde{\Rep}_{t_1}$. Geometrically this says
that the space of representations 
$\widetilde{\Rep}:=
\{\widetilde{\Rep}_t\}_{t\in \mathbf{P}^1(\C)\setminus \{0,1,\infty\}}$
can be interpreted as ``a local system of varieties" parameterized by
$t\in \PC\setminus \{0,1,\infty\}$:
the fibration $\widetilde{{\Rep}}\rightarrow \PC\setminus \{0,1,\infty\}$ (whose
fiber over $t$ is $\widetilde{\Rep}_t$) has a natural flat Ehresmann connection
on it \cite{Boa}.

\begin{rmk}
There are nice relations involving the coordinates of the gradient of $F$ and
some determinants (cf. \cite{ILT}, 3.9, page 10)~:
\begin{equation}
\label{detfricke}
F_{X_1}^2=(X_0X_t+2 X_1-A_1)^2=
\begin{vmatrix}
2 & -a_0 & -a_1 & X_0 \\
- a_0 & 2 & X_t & -a_\infty \\
-a_1 & X_t & 2 &-a_1 \\
X_0 & -a_\infty & -a_1 & 2
\end{vmatrix}
\end{equation}
and the circular permutations. Each relation is equivalent to $F=0$.
\end{rmk}


\subsubsection{Isomonodromy and PVI}

We recall briefly some basics about the sixth Painlev\'e equation and its relation with
isomonodromic families of linear Fuchsian differential equation. For more details, \cf\
\cite{CL}\footnote{We used the excellent presentation of \cite{Kli}.}. 

\bigskip

The sixth Painlev\'e equation is~:

\begin{equation}
\label{equap6}
      {\text (PVI)} \quad \frac{d^2y}{dt^2} =
      \frac12 \left(\frac{1}{y}+\frac{1}{y-1} +
      \frac{1}{y-t}\right)\left(\frac{dy}{dt}\right)^2 -
      \left(\frac{1}{t}+\frac{1}{t-1}+\frac{1}{y-t}\right)\frac{dy}{dt}
\end{equation}
\[
+\frac{y(y-1)(y-t)}{t^2(t-1)^2}\left( \alpha+ \beta\frac{t}{y^2}+
\gamma\frac{t-1}{(y-1)^2}+\delta\frac{t(t-1)}{(y-t)^2}\right);
\]
$\alpha,~\beta,~\gamma,~\delta\in \C$ are the parameters. \\

The generic solution of PVI has essential singularities and/or branch points in
the points $0,~1,\infty$. These points are called \emph{fixed singularities}. The other
singularities, the \emph{moving singularities} (so called because they depend on the
initial conditions) are \emph{poles}: it is the Painlev\'e property. A solution of PVI
can be analytically continued to a meromorphic function on the universal covering of 
$\PC\setminus \{0,1,\infty\}$. For generic values of the integration constants and of
the parameters $\alpha,~\beta,~\gamma,~\delta$, it cannot be expressed via elementary or 
classical transcendental functions\footnote{It is the irreducibility property of PVI,
  \cf\ for example \cite{CL} 1.8, page $2937$.}. For this reason, Painlev\'e called
these functions: ``transcendantes nouvelles" (new transcendental functions). \\

In modern formulation, solutions of PVI parameterize isomonodromic deformations (in $t$)
of rank two meromorphic connections over the Riemann sphere having simple poles at
the $4$ points $0,~t,1,\infty$. \\

We consider traceless $2\times 2$ linear differential systems with four fuchsian
singularities on the Riemann sphere $\PC$ (parameterized by a complex variable $t$)~:

\begin{equation}
\label{systemp6}
\frac{dY}{dz}=A(z;t)Y, \quad \quad 
A(z;t):=\frac{A_0(t)}{z}+\frac{A_t(t)}{z-t}+\frac{A_1(t)}{z-1}
\end{equation}
with the residue matrices $A_l(t)\in \slc$ ($l=0,t,1$) having 
$\pm \frac{\theta_l}{2}$ as eigenvalues (independantly of $t$).
We set $\theta:=(\theta_0,\theta_t,\theta_1,\theta_\infty)$: it encodes
(through a \emph{transcendental} mapping) the \emph{local monodromy data}. \\
 
Choosing a germ of a fundamental matrix solution $\Phi(z,t)$ of the above system near some 
nonsingular point $z_0$, one has a linear monodromy representation (anti-homomorphism)~:
\[
\rho: \pi_1\left(\PC\setminus \{0,t,1,\infty\};z_0\right)
\rightarrow \SLc
\]
such that the analytic continuation of $\Phi$ along a loop $\gamma$ (at $z_0$) defines
another fundamental matrix solution $\Phi\, \rho (\gamma)$. The equivalence class of 
$\rho$ in $\SLc$ is independant of the choice of the fundamental solution $\Phi$.
The system \eqref{systemp6} is said \emph{isomonodromic} if this conjugation class
is locally constant with respect to $t$, or equivalently if the matrices $A_l$ ($l=0,t,1$)
depends on $t$ in such a way that the monodromy of a fundamental solution $\Phi(z:t)$
does not change for small deformations of $t$. \\

A meromorphic connection can be interpreted as an equivalence class of systems modulo
rational equivalence (gauge transformation).
If two systems $\frac{dY}{dz}=A(z;t)Y$ and $\frac{dY}{dz}=B(z;t)Y$,
satisfying the conditions \eqref{systemp6}, are rationally equivalent on $\PC$,
that is if there exists a rational matrix $P$ such that $B=P^{-1}AP-P^{-1}\frac{dP}{dz}$,
then the two corresponding monodromy representation are equivalent. The isomonodromy
property is invariant by a rational equivalence. We can speak of isomonodromic
deformations of connections. \\

Schlesinger \cite{Schle} found that the isomonodromy condition is equivalent to having
the linear differential equation\footnote{We need a condition on $Y(\infty,t)$ to fix
$B(z,t)$, see \cite[p. 432]{SatoMiwaJimbo}.}:
\begin{equation}
\label{systemp6s}
\frac{dY}{dt}=B(z,t)Y, ~~\text{with}~~ B(z,t):=-\frac{A_t(t)}{z-t}Y.
\end{equation}

We define the \emph{Schlesinger system} as the system \eqref{systemp6} and
\eqref{systemp6s}~:
\[
\frac{dY}{dz}=A(z,t)Y, \quad \frac{dY}{dt}=B(z,t)Y,
\]
Then the isomonodromy of the system \eqref{systemp6} is equivalent to the
\emph{complete integrability condition}  (also called \emph{zero curvature
condition}) of the Schlesinger system~:
\begin{equation}
\label{LaxPair}
\frac{\partial B}{\partial z}-\frac{\partial A}{\partial t}=[A,B].
\end{equation}
Expliciting this condition, we see that the isomonodromicity of the system
\eqref{systemp6} is expressed by the following equations (called the Schlesinger
equations) on $(A_0,A_t,A_1)$:

\begin{equation}
\label{eqschles}
\frac{dA_0}{dt}=\frac{[A_t,A_0]}{t}, \quad 
\frac{dA_t}{dt}=\frac{[A_0,A_t]}{t}+\frac{[A_1,A_t]}{t-1}, \quad
\frac{dA_1}{dt}=\frac{[A_t,A_1]}{t-1} \cdot
\end{equation}
These equations correspond equivalently to the integrability of the logarithmic
connection in variables $(z,t)$~:
\[
\nabla:=d-\left(A_0(t)d \log z+A_t(t)d \log (z-t)+A_1(t)d \log (z-1)\right).
\]
on the trivial rank two vector bundle on $\PC$. \\

We suppose now that the Schlesinger equations are satisfied by the matrix $A$ of
the system \eqref{systemp6} and (following \cite{JMUII})  we will derive the non
linear second order PVI for some values of the parameter (under some genericity
condition on the local monodromy exponents $\pm\theta_l/2$). \\

We set $A_\infty:=-A_0-A_t-A_1$ and we suppose that the matrices $A_l$ ($l=0,t,1,\infty$)
are semi-simple. The eigenvalues of the $A_l$ ($l=0,t,1,\infty$) are independant of $t$
and we denote them by $e_l, e_l^{-1}$. We suppose $e_l\neq \pm 1$ or equivalently
$\pm\theta_l\notin \pi\Z$ (non-resonance conditions). \\

From Schlesinger equations we get $\displaystyle \frac{dA_\infty}{dt}=0$, therefore,
up to a constant gauge transformation, we can suppose
$A_\infty=\begin{pmatrix}
\theta_\infty & 0 \\
0 & -\theta_\infty\end{pmatrix}$. \\
 
We denote $[A]_{ij}$ the $(i,j)$ entry of the matrix of the differential system
\eqref{systemp6}. We suppose that the system is \emph{irreducible}. Then $[A]_{12}$
is not identically $0$. We have $A_0+A_t+A_1=-A_\infty$, therefore $[A_0+A_t+A_1]_{12}=0$.
Hence $z(z-t)(z-1)[A]_{12}$ is linear in $z$ and it admits a unique zero at the point
$z=q(t)$, where~:
\[
q(t)=-\frac{t[A_0]_{12}}{t[A_t]_{12}+[A_1]_{12}} \cdot
\]

The point $q(t)$ is \emph{an apparent singularity} of the second order linear ODE
satisfied by the first component $y$ of any solution $Y$ of the system  \eqref{systemp6}.
We denote~:
\[
p(t):=[A\left(q(t),t)\right)]_{11}+ \frac{\theta_0}{2q}+\frac{\theta_t}{2(q-t)}
+\frac{\theta_1}{2(q-1)},
\]
then the Schlesinger system \eqref{eqschles} is equivalent to the Hamiltonian system
of PVI whose (non autonomous) Hamiltonian is~:
\[
H_{VI}(q,p,t):=\Tr \left[\left(\frac{A_0(t)}{t}+\frac{A_1(t)}{t-1}\right)A_t(t)\right]-\frac{\theta_0\theta_t}{2t}
-\frac{\theta_t\theta_1}{2(t-1)} \cdot
\]
(\cf\ \cite{JMUII}). \\

Now we can write the Hamiltonian system in
PVI form with the following values for the parameters~:
\[
\alpha=(\theta_\infty-1)^2 \quad \beta=-\theta_0^2 , \quad 
\gamma= \theta_1^2,\quad \delta=1-\theta_t^2.
\]

The Riemann-Hilbert correspondence RH is given by the monodromy map between
the space of linear systems \eqref{systemp6} with prescribed poles and local
exponents $\pm \theta_l/2$, modulo $\SLc$-gauge transformations, on one side
(the source or ``left hand side"), and the space of monodromy representations with
prescribed local exponents modulo conjugation in $\SLc$ on the other side
(the target or ``right hand side"). \\

The relation with the notations introduced in \ref{subsubiso} is:
\[
e_l=e^{i\pi \theta_l}, \quad a_l=\Tr~M_l=2\cos \pi\theta_l.
\]
The Riemann-Hilbert correspondence can be  translated into a correspondence between
solutions of PVI and equivalence classes of monodromy representations. \\

We recall that an analytic complex vector field on a complex manifold (resp. the associate
flow) is called \emph{complete} if complex solutions (flow curves) \emph{exist for all
complex time}. The very naive phase space\footnote{A solution is defined by its initial
values $y(t_0)$ and $y'(t_0)$.} of the system associated to PVI is
$\left(\PC\setminus \{0,t,1,\infty\}\right)\times \C^2$.
It is not a good phase space because the solutions have poles: the Painlev\'e flow
is not complete. Using a series of blowing-ups K. Okamoto introduced a good space
of initial conditions $\mathcal{M}_{t_0}(\theta)$ at any point $t_0\in \C$ \cite{Ok1,Ok2}.
It is a convenient semi compactification of the naive phase space $\C^2$, an open
rational surface\footnote{A $8$ point blow-up of the Hirzebruch surface $\Sigma_2$ minus
an anti-canonical divisor.}. 
This surface is endowed with an algebraic symplectic structure given by the extension of
the standard symplectic form\footnote{The pole divisor of this extension is the
anticanonical divisor of a compactification of the Okamoto variety: the vertical
leaves. The vertical leaves configuration is described by a Dynkin diagram: today
a ``good list" of the Painlev\'e equations is labelled by such diagrams.} $dp\wedge dq$.
The Okamoto variety of initial conditions at $t_0$ can be identified with the moduli space
of meromorphic connections over the Riemann sphere\footnote{In the non resonant case.
In the resonant case, that is if one of the $\theta_l$ is an integer, then
$\mathcal{M}_{t_0}(\theta)$ is the moduli space of \emph{parabolic} connections
\cite{IKSstab}.}  having simple poles at the four points $0,~t_0,1,\infty$ with
local exponents $\{\pm \theta_l\}_{l=0,t,1,\infty}$. \\

For $\theta$ fixed, we have a fiber bundle
$\mathcal{M}(\theta)\rightarrow \PC\setminus \{0,1,\infty\}$: the fiber above $t_0$ is
$\mathcal{M}_{t_0}(\theta)$. \\

The naive Painlev\'e foliation extends to this fiber bundle.
This extension is transverse to the fibers and we get a complete (symplectic) flow,
the Painlev\'e flow. For all $t_0,t_1\neq 0,1,\infty$ this flow induces an analytic
symplectic diffeomorphism $\mathcal{M}_{t_0}(\theta)\rightarrow \mathcal{M}_{t_1}(\theta)$.
We get also \emph{analytic} maps (Riemann-Hilbert maps)~:
\[
\text{RH}: \mathcal{M}_{t}(\theta) \rightarrow \mathbf{S}_{P_{VI}}.
\]
Such a map can be interpreted as an analytic map:
\[
\text{RH}: \mathcal{M}_{t}(\theta)\rightarrow 
\mathcal{S}_{A_0A_tA_1A_\infty},
\]
where (using \eqref{desaauxA})~:
\begin{equation}
\label{desthetaaux A1}
A_i=4\left(\cos \theta_i\cos \theta_\infty+\cos \theta_j\cos \theta_k \right),
\end{equation}
where $(i,j,k)$ is a permutation of $(0,t,1)$.
and
\begin{equation}
\label{desthetaaux A2}
A_\infty
=16(\cos \theta_0\cos \theta_t\cos \theta_1\cos \theta_\infty)+
4(\cos^2 \theta_0+\cos^2 \theta_t+\cos^2 \theta_1 +\cos^2 \theta_\infty-1).
\end{equation}
This map is always proper. If the cubic surface 
$\mathcal{S}_{A_0A_tA_1A_\infty}$ is smooth, then this map is an analytic symplectic
isomorphism\footnote{The pull back by RH of the symplectic form $\omega_{VI,a}$
is the standard symplectic form on the Okamoto variety of initial conditions
\cite{Kli},  Proposition 4.3 \label{symplectomorphism}}. In the singular case
RH is a proper map, more precisely it realizes an analytic minimal resolution of
singularities of $\mathcal{S}_{A_0A_tA_1A_\infty}$. Along the irreducible components of
the exceptional divisor, PVI restricts to a Riccati equation\footnote{The singular
points of type $A_1$ , $A_2$ , $A_3$ , $D_4$ on the cubic surface yield $1$, $2$, $3$
and $4$ exceptional Ricatti curves \cite{vdPSaito}.}.

\begin{rmk}
``Pulling back'' the fiber bundle
$\widetilde{{\Rep}}\rightarrow \PC\setminus \{0,1,\infty\}$ and its connection by
the Riemann-Hilbert map (\ie\ keeping the base and changing the fibers through RH)
yields the fiber bundle $\mathcal{M}\rightarrow \PC\setminus \{0,1,\infty\}$ with
its PVI connection. This allows one to give an important interpretation of the
non-linear monodromy of PVI using a braid group \cite{Iwa}, \cite{CL}, page 2.
\end{rmk}


\subsubsection{Iso-irregularity and the Painlev\'e equations. Wild character varieties}

This paper is limited to fuchsian $q$-difference equations. The Lax pairs of irregular
$q$-Painlev\'e equations are listed up by Murata \cite{Murata} (see subsection
\ref{subsubsection:irregularqde}). We look here at the irregular differential
case briefly and think in terms of $q$-analogies. \\

In order to classify \emph{irregular} connections, the monodromy is no longer
sufficient, it is necessary to introduce \emph{generalized monodromy data} (formal
monodromy, Stokes multipliers \emph{and} links).
Martinet-Ramis \cite{MRAcceleration} have constructed a local \emph{wild fundamental 
group}, so that germs of connections with irregular singularities can be interpreted 
as finite dimensional representations of this group. This construction uses
\emph{multisummability} of divergent series as an essential
ingredient\footnote{``Generically" Ramis $k$-summability is sufficient.
 It is the case for the Painlev\'e equations.}. In the global case a group
is no longer sufficient, it is necessary to introduce a wild \emph{groupoid}
\cite{vdPSaito}. 
The generalized monodromy data are in some sense representations of this wild groupoid.
The quotient of the set of these representations by the natural equivalence relation is
a \emph{wild character variety} \cite{vdPSaito,Boa3,Boa4}. \\

In all the cases PI, PII, PIII, PIV, PV, the wild character varieties are, as in
the PVI case, cubic surfaces. The interested reader will find a list of equations
of these surfaces in \cite{vdPSaito} (pages $19$-$20$) 
and (in a nice form) in \cite{CMR} (Table 1, page 2). \\

We consider some\footnote{\cf\ for explicit conditions \cite{JMUI,JMUII,Ok3}.} 
traceless $2\times 2$ linear differential systems
with at most $3$ singularities (one at least being irregular) on the Riemann sphere 
$\PC$ (parameterized by a complex variable $t$)~:

\begin{equation}
\label{systemp}
\frac{dY}{dz}=A(z;t)Y
\end{equation}
The system \eqref{systemp} parameterized by $t$ is said \emph{iso-irregular} if 
the conjugation class of generalized monodromy data (or wild monodromy representation)
is locally constant with respect to $t$, or equivalently if the matrices $A_l$ ($l=0,t,1$)
depends on $t$ in such a way that the generalized monodromy data of a fundamental solution 
$\Phi(z:t)$ does not change for small deformations of $t$. \\

As in the fuchsian case, the iso-irregularity condition is equivalent to an integrability
condition (Schlesinger equation) and therefore it is possible in each case to express it
as a Painlev\'e equation, along similar lines \cite{JMUI,JMUII}. \\

Ren\'e Garnier was the first to interpret a Painlev\'e equation as an iso-irregular
deformation of a linear equation \cite{Garnier}. He did not define the generalized
monodromy data using Stokes phenomena (as Birkhoff did some years before \cite{Birkhoff1}),
he defined them by an interesting process of \emph{confluence} of classical monodromy data
of fuchsian equations.


\subsection{The first age of $q$-monodromy; $q$-PVI according to Jimbo and Sakai}
\label{subsection:firstage}


\subsubsection{Riemann-Hilbert correspondence for $q$-difference equations}
\label{subsubsection:RHCforqDE}

In his celebrated 1913 article \cite{Birkhoff1} ``The generalized Riemann problem
for linear differential equations and the allied problems for linear difference
and $q$-difference equations'', Birkhoff looks for ``transcendental invariants''
in order to classify rational fuchsian $q$-difference equations or systems.
The systems have the form:
$$
Y(q x) = A(x) Y(x),
$$
$x$ a complex variable, $q$ a complex number such that $0 < \lmod q \rmod < 1$
and $A(x)$ an invertible $n \times n$ matrix of rational functions (so the unknown
$Y$ is a vector of functions). \\

Rational equivalence is induced by gauge transformations $Z = Q Y$,
$Q \in \GLn(\Kg)$, so that $Z$ is a candidate solution of $Z(qx) = B(x) Z(x)$,
where:
$$
B(x) := Q(qx) A(x) Q(x)^{-1} \text{~is declared rationally equivalent to~} A(x).
$$

The problem of classification is not changed if $A$ and $B$ are replaced respectively
by $fA$ and $fB$ with $f$ any scalar function, so one can as well (and Birkhoff does)
assume that $A$ is polynomial:
$$
A = A_0 + \cdots + A_\mu x^\mu, \quad A_0,\ldots,A_\mu \in \Matnc.
$$
Birkhoff moreover assumes that $A_0,A_\mu \in \GLnc$ (this means in essence that
$0$ and $\infty$ are regular singularities) and implicitly\footnote{Almost no
assumption or definition is explicit in \cite{Birkhoff1}, and many conclusions
are not either. However \emph{the heart of the matter} is dealt with.} that 
$A_0,A_\mu$ are semi-simple and ``non-resonant'' (such details will be explained
section \ref{section:A Birkhoff type classification theorem}). Their eigenvalues,
seen as elements of $\C^* \pmod{q^\Z}$, are considered as \emph{exponents} at $0$
and $\infty$ and should encode the local monodromies there. \\

Fuchs-Frobenius type algorithms yield local fundamental solutions $\Y^{(0)}$ and
$\Y^{(\infty)}$, made up of multivalued functions. \emph{Birkhoff connection matrix}
is then $P := (\Y^{(\infty)})^{-1} \Y^{(0)}$. The main results of Birkhoff (in the part
devoted to $q$-difference equations) are then that:
\begin{enumerate}
\item The local exponents being fixed, $P$ classifies $A$ up to rational
equivalence.
\item $P$ can be characterized by $\mu n^2 + 1$ ``characteristic constants'',
the transcendental invariants looked for.
\end{enumerate}
The second statement comes from the fact that $P$ is almost $q$-invariant (the
defect comes from the multivaluedness) so its elements can almost be identified to
elliptic functions and those are very much controlled by their zeroes and poles.
More precisely, each coefficient of $P$ has $\mu$ zeroes by which it is determined
up to a constant; this altogether yields $(\mu + 1) n^2$ degrees of freedom, but
taking in account gauge freedom reduces this dimension to $\mu n^2 + 1$, see the
Master: \cite[\S 20]{Birkhoff1}. (Also see, in the case $\mu = n = 2$, remark
\ref{rmk:Birkhoffmu=n=2} at the beginning of \ref{subsection:Setting, general facts}.) \\

Birkhoff's paper has some drawbacks:
\begin{itemize}
\item Contrary to the case of differential equations, multivaluedness can (and
should) be avoided.
\item The problem is solved only under generically true assumptions.
\item Irregular equations are not considered.
\end{itemize}

As for the first two drawbacks, see \ref{subsubsection:uniformsolutions} below.
As for the third one, Birkhoff himself with his student Guenther made a decisive
step in \cite{BirkhoffGuenther}, but the sequel had to wait for seven decades, see
\ref{subsubsection:irregularqde}. However, the main question from our point of
view is: in what sense does Birkhoff connection matrix encode monodromy ?


\subsubsection{$q$-analogues of Painlev\'e equations}

The search for $q$-analogues of classical special functions has been a flourishing
industry in the best part of twentieth century. Some physicists have been specially
(!) interested in discrete analogues of Painlev\'e functions, see \eg\ \cite{RGH}.
One way to specify them was \emph{confinement of singularities}, invented by
Grammatikos, Ramani and Papageorgiou see \cite{GRP}. It seems that it can be
considered as a sensible discrete analogue of the Painlev\'e property, which
was the guiding criterion of Painlev\'e himself. However, this did not lead
to the discovery of a $q$-analog to PVI. \\

In \cite{JimboSakai}, Jimbo and Sakai adapted the isomonodromy approach to the
$q$-difference setting. They considered a family of order $2$ degree $2$ systems:
$$
Y(qx) = A_t(x) Y(x), \quad A_t(x) = A_0(t) + A_1(t) x + A_2(t) x^2, \quad
\forall t \;,\; A_0(t), A_2(t) \in \GL_2(\C), A_1(t) \in \Mat_2(\C),
$$
with conditions similar to those imposed by Birkhoff in \cite{Birkhoff1}. \\

Then they imposed that the family has constant local data, \ie\ that $A_0(t)$
and $A_2(t)$ have constant eigenvalues\footnote{Actually there is a subtle twist
in the case of $A_0$, but this does not matter here; also, the singularities (zeroes
of $\det A(x)$) are subject to some similar condition.}. To express isomonodromy,
they consider Birkhoff connection matrix as depending on $t$ and, in a bold step,
assume that it is $q$-constant:
$$
\forall t,x \;,\; P(qt,x) = P(t,x).
$$
They deduce a ``Lax pair'', some kind of integrability condition analogous to
\eqref{LaxPair}: it is the system \eqref{equation:qLaxPair} herebelow. From this
they derive a nonlinear $q$-difference equation they consider as the adequate
analogue of PVI. Applying the usual test to support such a claim\footnote{This
is a standard process in the history of $q$-analogues. For a very detailed
(and unusually rigorous) such study, see \cite{JSAIF} which tackles the case
of $q$-hypergeometric equations under the name of ``confluence''. We shall
not further delve into these matters in the present work.}, they go to the
``continuous limit'' $q \to 1$ and show how to recover classical PVI equation. \\
  
The succesful attack of Jimbo and Sakai was very influential. Any attempt at
a theory of monodromy for $q$-difference equations should use it as a touchstone.
For the sake of completeness, we now give a description of their model in their
own notations.

\paragraph{Connection preserving deformation and  $q$-$P_{\mathrm{VI}}$.}

We review here the $q$-analogue of the sixth Painlev\'e equation obtained by Jimbo
and Sakai in \cite{JimboSakai}. \\

We denote $y=y(t), z=z(t)$, $\bar{y}=y(qt), \bar{z}=z(qt)$. We take $a_1$, $a_2$, $a_3$,
$a_4$, $b_1$, $b_2$, $b_3$, $b_4$ as complex parameters of the equation. The $q$-analogue
$q$-$P_{\mathrm{VI}}$ of the sixth Painlev\'e equation considered here is: 
$$
\dfrac{y\bar{y}}{a_3a_4}
= \dfrac{(\bar{z}-b_1t)(\bar{z}-b_2t)}{(\bar{z}-b_3)(\bar{z}-b_4)},\quad
\frac{z\bar{z}}{b_3b_4}
 = \dfrac{(y-a_1t)(y-a_2t)}{(y-a_3)(y-a_4)}
\quad \text{with}
\frac{b_1b_2}{b_3b_4}=q\frac{a_1a_2}{a_3a_4} \cdot
$$
Jimbo and Sakai derive $q$-$P_{\mathrm{VI}}$ from \emph{connection preserving deformation}
of a fuchsian linear $q$-difference equation with rank two and order two. What is called
here ``connection preserving deformation'' is a compatibility condition of the following
two $q$-difference equations:
\begin{equation}
\label{equation:qLaxPair}
\begin{cases}
Y(qx,t) = A(x,t)Y(x,t), \\
Y(x,qt) = B(x,t)Y(x,t).
\end{cases}
\end{equation}
The system above \eqref{equation:qLaxPair} can be best understood by writing
$Y_t(x) := Y(x,t)$ and $A_t(x) := A(x,t)$. The first relation says that we have
a family (parameterized by $t$) of $q$-difference equations; the second relation
states a gauge equivalence of $Y_t$ with $Y_{qt}$. The system \eqref{equation:qLaxPair}
is compatible if and only if:
$$
A(x, q t) B(x, t)=B(q x, t) A(x, t). 
$$
By the compatibility condition, the Birkhoff connection matrix $P(t)$ is a
``quasi-constant'', \ie\ $P(qt)=P(t)$: this is Jimbo and Sakai interpretation
of $q$-isomonodromy. \\

Jimbo and Sakai set  the rank two matrices $A(x,t)$ and $B(x,t)$ as
$$
A(x,t) = A_0(t)+xA_1(t)+x^2A_2
$$
and 
$$
B(x,t) = \dfrac{x}{(x-a_1 qt)(x-a_2qt)}(xI+B_0(t)) \cdot
$$
We assume that $A_2=\Diag(\sigma_1,\sigma_2)$ and that the eigenvalues of $A_0(t)$ are
$\rho_1=\theta_1 t, \rho_2= \theta_2 t$. We set
$\det A(x,t)=\sigma_1\sigma_2(x-x_1)(x-x_2 )(x-x_3)(x-x_4)$. We define the parameters
$a_j$ and $b_k$ by:
$$
x_1=a_1t,  \quad x_2=a_2t, \quad x_3= a_3, \quad x_4=a_4, \quad
b_1 = \dfrac{a_1a_2}{\rho_1}, \quad
b_2 = \dfrac{a_1a_2}{\rho_2}, \quad
b_3 = \dfrac{1}{\sigma_1q}, \quad
b_4 = \dfrac{1}{\sigma_2} \cdot
$$
We take  variables  $y=y(t)$, $z_i=z_i(t)\ (i=1,2)$ such that 
$$
A_{12}(y,t)=0, \quad A_{11}(y,t)=\sigma_1 z_1,\quad 
A_{22}(y,t)= z_2,
$$
We set a variable $z$ in such a way that:
$$
z_2=\sigma_1\sigma_2qz(y-a_3).
$$
Then we obtain $q$-$P_{\mathrm{VI}}$ by the compatibility condition
$A(x, q t) B(x, t)=B(q x, t) A(x, t)$.


\subsubsection{Does Birkhoff connection matrix encode monodromy ?}

Since in some sense Birkhoff connection matrix connects solutions at $0$ and
at $\infty$, it is comparable to connection matrices of the classical theory
(those related to analytic continuation of solution along pathes, see for instance
\cite{IKSY}). So it was generally felt that it should relate to the monodromy of
the $q$-difference system if \emph{that} could be defined somehow. However it was not
clear what was the topology underlying it. So for some time the confirmations
of the monodromy interpretation were indirect. Actually they came from Galois
theory. \\

Without going in any detail, let us say that differential Galois theory, as created
by Picard and Vessiot, attaches to a differential equation or system with rational
coefficients a linear algebraic group $G$. This is related to the monodromy group
$M$ in the following ways:
\begin{enumerate}
\item In all cases, $M$ naturally embeds into $G$: $M \subset G$.
\item In case of a fuchsian differential equation, $G$ is the Zariski closure of $M$
(Schlesinger density theorem).
\end{enumerate}

In \cite{Etingof}, Etingof proved a $q$-analogue of Schlesinger density theorem
with Birkhoff connection matrix in the role of monodromy in the following way.
He assumes that the rational system $Y(qx) = A(x) Y(x)$, $A \in \GLn(\Kg)$ is
such that $A(0) = A(\infty) = I_n$, the identity matrix; this means in essence
that $0$ and $\infty$ are not merely regular singularities, as in Birkhoff's
paper, but \emph{ordinary} points (indeed, it is the case when the equation
can be solved with power series, without the need of special transcendental
functions). In this case, the connection matrix
$P(x)$ as built by Birkhoff is uniform over $\C^*$ and truly elliptic. On the
other hand, in the mean time (since Picard and Vessiot), differential Galois
theory had been extended to difference equations over fields more general than
the complex numbers, so that there is a linear algebraic group $G$ attached to
the equation. Then the \emph{values} $P(a)^{-1} P(b)$, where defined, generate
a subgroup $M$ of $G$ which is Zariski-dense in $G$. So it would be a natural
conjecture that the $q$-analogue of the monodromy group is the group generated
by all the $P(a)^{-1} P(b)$. \\

Van der Put and Singer then extended in \cite{vdPSinger} this result to the case
of general fuchsian systems. However, difficulties appear that are not present
in the classical case of differential equations. First, the natural field of
constants in the $q$-different setting is the field of meromorphic functions
over $\C^*$ that are $q$-invariant: $f(q x) = f(x)$. This field can be identified
with a field of elliptic functions (see the end of \ref{subsection:q-notations}).
It is not algebraically closed, which is a severe drawback for Picard-Vessiot
theory. Second, natural solutions to basic $q$-difference equations have bad
multiplicative properties. For instance, writing $e_c$ a non trivial solution
of the constant scalar equation $f(qx) = c f(x)$ ($c \in \C^*$), it is not
possible to impose that $e_c e_d = e_{cd}$ or even that $e_c e_d/e_{cd} \in \C^*$.
To overcome these difficulties, van der Put and Singer introduced symbolic
solutions. The theory then develops nicely, in particular (to stick to our
monodromy-headed point of view\footnote{Note however that the theory expounded
in \cite{vdPSinger} has many more advantages, including a tannakian interpretation
and a description of the universal Galois group for fuchsian equations.}) it does
contain a Schlesinger type density theorem for general fuchsian equations. \\

So in some sense it was established that Birkhoff connection matrix has something
to do with monodromy. However the transcendental point of view of Birkhoff seemed
partly abandonned.


\subsection{The second age of $q$-monodromy}
\label{subsection:secondage}

Since the end of the last century (and millenium), mainly under the influence of
the second author, transcendental methods in the theory of $q$-difference equations
(including Galois theory) have relied on the use of theta functions. This is related
to the fact that $\Eq := \C^*/q^\Z$, as a Riemann surface, can be seen as an elliptic
curve; and solutions of $q$-difference systems as sections of holomorphic vector bundles
over $\Eq$.


\subsubsection{Uniform solutions to $q$-difference equations}
\label{subsubsection:uniformsolutions}

In the work of Praagman on formal classification of difference and $q$-difference
operators \cite{Praagman} appears an argument based on the fact that every
holomorphic vector bundle on the elliptic curve $\Eq$ is meromorphically trivial.
An easy consequence of this fact is that any rational $q$-difference system:
\begin{equation}
\label{eqn:qDEstandard}
Y(qx) = A(x) Y(x), \quad A \in \GLn(\Kg),
\end{equation}
admits a ``full complement'' (that is a system of maximal possible rank) of
solutions meromorphic over $\C^*$. Therefore, contrary to the case of differential
equations \emph{it is not necessary to use multivalued functions}. From our point
of view (Riemann-Hilbert correspondence), this means that what we consider as
monodromy should not be related on ambiguity of analytic continuation. \\

In \cite{JSAIF}, the third author gave a concrete content to this result by solving
explicitly fuchsian systems in a way similar to the Fuchs-Frobenius method for
differential equations. This was applied to the Riemann-Hilbert correspondence
and also to Galois theory in \cite{JSGAL}. In the latter paper, the Galois group
was defined by tannakian means and a Schlesinger density theorem similar to those
quoted above was proved. Moreover, as a bonus rewarding the use of ``true'' (not
symbolic) functions, a very precise meaning could be given to the degeneracy
(``continuous limit''), when $q \to 1$, towards monodromy and towards differential Galois
theory. In particular, the values $P(a)^{-1} P(b)$ degenerate, when $q \to 1$, into
monodromy matrices of the differential system. \\

However some undue complications in the computations led to the idea that Birkhoff
connection matrix mixes, in some sense, local monodromies at $0$ and $\infty$ with
monodromy at the ``intermediate singularities'' (those in $\C^*$). We find it
relevant to explain this point in some detail, because the way we define and use
monodromy in the present work is directly related to it. \\

We suppose that $A(0), A(\infty) \in \GLnc$, which, as already noted, means that
the above system is fuchsian at $0$ and $\infty$ (it can indeed be characterized
by the fact that solutions satisfy some kind of moderate growth condition, \cite{JSAIF}).
Then there exist constant invertible matrices\footnote{Generically $A^{(0)} = A(0)$
and $A^{(\infty)} = A(\infty)$, but this is not the case if there are ``resonancies''.}
$A^{(0)}, A^{(\infty)} \in \GLnc$ such that:
$$
A(x) = M^{(0)}(qx) A^{(0)} \left(M^{(0)}(x)\right)^{-1} \text{~~and~~}
A^{(\infty)}(x) = M^{(\infty)}(qx) A^{(\infty)} \left(M^{(\infty)}(x)\right)^{-1},
$$
where $M^{(0)} \in \GLn\left(\C(\{x\})\right)$ and
$M^{(\infty)} \in \GLn\left(\C(\{1/x\})\right)$.
This implies that one can look at fundamental solutions of \eqref{eqn:qDEstandard}
in the form:
$$
\Y^{(0)} = M^{(0)} e_{A^{(0)}} \text{~~and~~}
\Y^{(\infty)} = M^{(\infty)} e_{A^{(\infty)}},
$$
where $e_{A^{(0)}}$, $e_{A^{(\infty)}}$ are respectively solutions of the systems 
\emph{with constant coefficients}
$$
Y(qx) = A^{(0)} Y(x), \quad \text{resp.} \quad Y(qx) = A^{(\infty)} Y(x).
$$
Birkhoff (along with his predecessors) 
solves those systems using multivalued functions such as $x^{\ln c/\ln q}$ (where
$c$ is an eigenvalue of $A^{(0)}$, resp. $A^{(\infty)}$); van der Put and Singer
use a symbol $e_c$; and Sauloy uses $\thq(x)/\thq(cx)$ (the theta function $\thq$
will be precisely defined later)\footnote{Functions described here suffice in
the generic case that $A^{(0)}, A^{(\infty)}$ are semi-simple. Otherwise, one also
introduces ``$q$-logarithms'', see \ref{subsubsection:qlog}.}. Then Birkhoff connection
matrix writes:
$$
P := (\Y^{(\infty)})^{-1} \Y^{(0)} = (e_{A^{(\infty)}})^{-1} M e_{A^{(0)}}, \text{~~where~~}
M := (M^{(\infty)})^{-1} M^{(0)}.
$$
It comes out that $e_{A^{(0)}}, e_{A^{(\infty)}}$ really encode the local monodromies
at $0$ and $\infty$ and the corresponding local Galois groups can be directly
computed from them\footnote{The local Galois groups were independently found
by Baranovsky and Ginzburg \cite{BaranovskyGinzburg} in the context of loop
groups.}. And it has been verified in many contexts that $M$ indeed encodes
the monodromy at intermediate singularities. In this paper we define a space
of monodromy data for the Jimbo-Sakai family using $M$ instead of $P$.

\begin{rmk}
Birkhoff matrix $P$ still plays an important role, since it directly relates to
solutions. For instance, in \cite{OhyamaConnectionHG}, the first author computes
it for basic hypergeometric equations; also see \cite{JRHG2}, where Roques uses
it to study Galois groups.
\end{rmk}


\subsubsection{Irregular $q$-difference equations and other $q$-Painlev\'e
 equations: Murata's list}
\label{subsubsection:irregularqde}

As we said before, Birkhoff and Guenther had led a first attack at irregular
$q$-difference equations in \cite{BirkhoffGuenther}, but that part of the theory
remained dormant for quite a long time. In \cite{RSZ}, the second and third authors
along with Changgui Zhang defined a $q$-analog of Stokes phenomenon and applied it
to Riemann-Hilbert correspondence for irregular $q$-difference equations. This was
further used for Galois theory in \cite{RS3}. \\

On the other hand, Murata, in \cite{Murata}, extended the work of Jimbo and Sakai
to various degeneracies of $q$-PVI related to families of irregular equations. \\

It is natural to envision an application of the tools of \cite{RSZ,RS3} to extend
the methods and results of the present paper to Murata's list. A first attempt was
sketched in Anton Eloy's thesis \cite{Eloy}. We hope to pursue this goal in a near
future.


\subsubsection{Families, moduli}

In all versions of Riemann-Hilbert correspondence for $q$-difference equations
during the second age, moduli problems and behaviour of continuous families were
not properly adressed. In the present work, we fix the local data: obviously
this should give rise to a fibering of some global space of monodromy data over
a space of local monodromy data. An attempt at this appears in already quoted
Eloy's thesis, but most of the work is yet to be done. We also hope to pursue
this goal in a near future.


\subsubsection{``Intermediate'' singularities}
\label{subsubsection:intermediatesingularities}

One of the successes of classical Riemann-Hilbert theory lies in the ability
to decompose global phenomena into local ones, in particular, to define local
monodromies, local Galois groups, etc. From the beginning, it has seemed very
difficult to do something similar for $q$-difference equations. \\

One aspect of the problem is that the obvious singularities other than $0$ and
$\infty$, \ie\ the poles of $A(x)$ and those of $A(x)^{-1}$ in $\C^*$, are not
really local: they are moved under the action of the dilatation operator
$x \mapsto q x$. Therefore it seems that they should be replaced either by the
corresponding $q$-spirals (discrete spirals of the form $a q^\Z$); or by the
corresponding points in $\Eq$. \\

In \cite{JSGAL}, reduction of the global Galois group to local contributions was
accomplished only in the trivial case of an abelian Galois group\footnote{In that
case, it boils down to ``class field theory over $\Eq$'', as described in Serre's
book \cite{SerreGACC}.}. But the first significant progress in this direction
(understanding local contributions) was accomplished much later by
Roques in \cite{JRHG2}. He used the Lie algebra instead of the Galois group to
take in account the local contribution of the connection matrix (via the residue
of its logarithmic derivative) at the only intermediate singularity of a ``basic''
hypergeometric (\ie\ $q$-hypergeometric) equation. In a somewhat different vein
(sheaf theoretic approach), Roques and the third author gave in \cite{RoquesSauloy}
a cohomological interpretation of the rigidity index defined by Sakai and Yamaguchi
in \cite{SakaiYamaguchi}. There, the local contributions of intermediate
singularities to an Euler characteristics can be measured. \\

In the present paper, a new technique is developped under the name of ``Mano
decomposition'' (as it has its roots in Mano's paper \cite{Mano}) which in some
sense allows us to localize the monodromy at \emph{pairs of points}. We put great
hopes in this process for future progress.


\subsection{Contents of this paper}
\label{subsection:contents}

Section \ref{section:tools} is devoted to general notations and conventions, along
with some basic tools for dealing with $q$-difference equations. \\

In section \ref{section:A Birkhoff type classification theorem}, we state and prove
a variant of Riemann-Hilbert correspondence from \cite{Birkhoff1} but using the
matrix $M$ described above in \ref{subsubsection:uniformsolutions}. We apply it
first to a criterion of reducibility, second to the case of ``hypergeometric''
systems (actually a slightly more general class allowing for some degeneracies). \\

In section \ref{section:JimboSakaiFamilyI} we define the Jimbo-Sakai family, thus
formalising the objects studied in \cite{JimboSakai}. We introduce the space $\F$
of its monodromy data, defined as the space of rational equivalence classes of
such equations but translated through our Riemann-Hilbert correspondence; we do
this for fixed local monodromy data (exponents at $0$ and $\infty$) and also fixed
singular set. Then we give a first geometric description of $\F$ as an algebraic
surface. This part (subsections \ref{subsection:Falgsurf} to \ref{rough4surf})
has the character of a preliminary exploration, collecting as much information
as possible in order to be later able to identify our cubic surface, which will
be done to some extent conjecturally in sections \ref{section:JimboSakaiFamilyII}
and \ref{section:Geometrysurgeryandpants}. Thus for instance we give a close look
to incidence relations in \ref{subsection:Spacesoffunctions again}. \\
    
Section \ref{section:ManoDecomposition} deals with a new process inspired by the
paper \cite{Mano} of Mano. This allows to decompose the monodromy matrix $M$ of a
system in the Jimbo-Sakai family into the product $M = P Q$ of two hypergeometric
monodromy matrices, while distributing the four singularities of $M$ among $P$ and
$Q$. The proof is very detailed because it involves some new objects, techniques
and tools which we hope will be handy in the future. The main results in this part
are the existence theorem \ref{thm:existenceMD}, and the gauge freedom and normal forms
(propositions \ref{prop:gaugefreedom} and \ref{prop:NormalformsforC}). \\
      
In section \ref{section:JimboSakaiFamilyII} we apply Mano decomposition to obtain
a more precise description of the space $\F$ as an algebraic \emph{fibered} surface.
We do that under the same assumptions as Jimbo and Sakai, plus some more that are
generically true and that seem reasonable; actually, they are essentially the same
as those underlying similar works in the classical case of differential equations
related to Painlev\'e and isomonodromy. \\
        
Section \ref{section:Geometrysurgeryandpants} is devoted to a larger picture and
tries to formulate analogies between the character varieties and their dynamics in
the differential and in the $q$-difference case. By nature, it is partly conjectural. \\

In conclusive section \ref{section:conclusion} we describe some interesting open problems
and perspectives.

\subsection*{Acknowledgements}

The first author would like to express his very great appreciation to the
Institut de Math\'ematiques de Toulouse for their hospitality during his visit on
February 2018. He was partly supported by JSPS KAKENHI Grant Number JP19K03566. \\
The second author thanks Arnaud Beauville for some help with Kummer surfaces, and
Martin Klimes and Emmanuel Paul for illuminating discussions about the relations
between local reducibility of connections and lines on cubic surfaces. \\
Most of the contribution of the third author was done while he was a
member of the Institut de Math\'ematiques de Toulouse.



\section{Tools}
\label{section:tools}


\subsection{General notations}

Here are some standard notations of general use:

\begin{itemize}
\item{$\Kg$ is the field of rational fractions over $\C$.}
\item{$\Ka$, the field of meromorphic germs at $0$ (or Laurent\footnote{We shall
sometimes - as here - understand Laurent power series to have bounded below exponents,
whence the form $\sum\limits_{n \geq n_0}$ for some $n_0 \in \Z$; and sometimes not. The
context should make it clear.} convergent power series), is the quotient field
of $\Ra$, the ring of holomorphic germs at $0$ (\ie\ convergent power series).}
\item{$\Kf$, the field of Laurent formal power series, is the quotient field
of $\Rf$, the ring of formal power series.}
\item{$\Rai,\Kai,\Rfi,\Kfi$ are similarly defined replacing $x$ by $1/x$.}
\item{$\M(\Omega)$ is the ring of meromorphic functions on the open subset $\Omega$
of a Riemann surface (thus a field if $\Omega$ is a domain). Most of the time, $\Omega$
will be a domain of the Riemann sphere $\Sr$ or 
of the elliptic curve $\Eq$ defined further below.}
\item{$\Sr$, the Riemann sphere and its open subsets
$\C = \Sr \setminus \{\infty\}$ and
$\C_\infty := \Sr \setminus \{0\}$.}
\item{$\Matn$, $\Mat_{m,n}$, $\GLn$ are spaces of square, resp. rectangular 
matrices and the linear group; $\Dnc \subset \GLnc$ is the subgroup of diagonal
invertible matrices.}
\item{$\Diag(a_1,\ldots,a_n)$, $\Sp(A)$ respectively denote a diagonal matrix 
and the spectrum of an arbitrary matrix $A$. Most of the time we consider the
spectrum as a \emph{multiset}, \ie\ its elements have multiplicities.}
\end{itemize}


\subsection{$q$-notations}
\label{subsection:q-notations}

Here are some notations related to $q$ but of general interest:

\begin{itemize}
\item{$q$ is a complex number such that $0 < \lmod q \rmod < 1$.}
\item{$\sq$ is the $q$-dilatation operator $f(z) \mapsto f(qz)$.}
\item{$\Cq := \{z \in \C \tq \lmod q \rmod < \lmod z \rmod \leq 1\}$, the
\emph{fundamental annulus}.}
\item{For $x \in \C^*$, we write $R(x) \in \Cq$ its unique representative
modulo $q^\Z$.}
\item{$\Eq$ is $\C^*/q^\Z$ either seen as a group or, more frequently,
as a Riemann surface (a complex torus, or ``elliptic curve''). Indeed, the
compositum of the canonical projection $\C^* \rightarrow \Eq$ with the map
$z \mapsto e^{2 \ii \pi z}$ is a covering map between Riemann surfaces and also
a group morphism with kernel $\Z + \Z \tau$, where $q = e^{2 \ii \pi \tau}$,
whence an identification of $\Eq$ with $\C/(\Z + \Z \tau)$.}
\item{The canonical projection $\pi: \C^* \rightarrow \Eq$ is also denoted
$a \mapsto \overline{a}$. It is bijective from $\Cq$ to $\Eq$.}
\item{We write $[a;q] := a q^\Z$ the discrete logarithmic
$q$-spiral $\pi^{-1}(\overline{a})$.}
\item{$(a;q)_n = \prod\limits_{0 \leq i < n} (1 - a q^i)$ and
$(a;q)_\infty := \prod\limits_{n \geq 0} (1 - a q^n)$ are the $q$-Pochhammer
symbols.}
\item{These notations are ``collectivized'' as follows:
$$
[a_1,\ldots,a_m;q] := \bigcup_{i=1}^m [a_i;q], \quad
(a_1,\ldots,a_m;q)_n := \prod_{i=1}^m (a_i;q)_n, \quad
(a_1,\ldots,a_m;q)_\infty := \prod_{i=1}^m (a_i;q)_\infty.
$$
}
\end{itemize}

The operator $\sq$ acts naturally on the field $\Kwg$, the subfield of
``constants'': 
$$
\Kwg^{\sq} := \{f \in \Kwg \tq \sq f = f\}
$$ 
has a natural identification with the field of elliptic functions $\M(\Eq)$;
any $f \in \M(\Eq)$ can at will be seen as a meromorphic function on $\Eq$;
as a meromorphic function on $\C^*$ such that $f(qx) = f(x)$; or as a
meromorphic function on $\C$ with $(\Z + \Z \tau)$-periodicity.

\paragraph{A convention for notations of congruences.}

Since in all the text most congruences in $\C^*$ are modulo $q^\Z$, we shall
systematically (when $a,b \in \C^*$) write $a \equiv b$ for $a \equiv b \pmod{q^{\Z}}$.


\subsection{Some functions}
\label{subsection:somefunctions}

The main one is the following theta function:
$$
\thq(x) := \sum_{n \in \Z} q^{n(n-1)/2} x^n.
$$
It is holomorphic over $\C^*$ and satisfies the functional equations:
$$
\thq(qx) = \dfrac{1}{x} \thq(x) = \thq(1/x).
$$
Thanks to \emph{Jacobi's triple product formula}:
$$
\thq(x) = (q;q)_\infty (-x;q)_\infty (-q/x;q)_\infty,
$$
it has simple zeroes over $[-1;q]$ and nowhere else, which we
summarize\footnote{We write $\sum m_i [x_i]$ the divisors on a Riemann
surface $X$, where the $m_i \in \Z$ and the $x_i \in X$; and
$\div_X(f)$ the divisor of a function $f$ on $X$ or of a section of
a line bundle (when this divisor is defined). Note that if $X$ is
non compact, the support of a divisor is not necessarily finite.} by:
$$
\div_{\C^*}(\thq) = \sum_{a \in [-1;q]} [a].
$$
Since this divisor is $\sq$-invariant, or because $\thq$ can be seen as a
section of a line bundle over $\Eq$, we can also write:
$$
\div_\Eq(\thq) = [\overline{-1}].
$$

For every $c \in \C^*$, the function $e_{q,c}(x) := \thq(x/c)/\thq(x)$ is a
non trivial meromorphic solution of the $q$-difference equation $\sq f = c f$
such that $\div_\Eq(e_{q,c}) = [\overline{-c}] - [\overline{-1}]$. (One could
as well use instead the function $\thq(x)/\thq(cx)$.)


\subsection{Some first order equations}
\label{subsection:someequations}


\subsubsection{Equations $\sq f = u f$}

Using Laurent series expansions, one proves easily the following facts:

\begin{itemize}
\item{The equation $\sq f = c x^k f$, $c \in \C^*$, $k \in \Z$, has non
trivial solutions in $\Ka$ if, and only if, $k = 0$ and $c = q^m$, $m \in \Z$;
and then these solutions are the elements of $\C^* x^m$.}
\item{The equation $\sq f = c x^k f$, $c \in \C^*$, $k \in \Z$, has non 
trivial solutions in $\Rwg$ if, and only if $k < 0$ or $k = 0$ and $c = q^m$, 
$m \in \Z$. In the last case, these solutions are the elements of $\C^* x^m$.}
\end{itemize}

The second statement can be completed as follows. Let $c \in \C^*$ and
$k \in \N^*$. Then the solutions of $\sq f = c x^{-k} f$ in $\Rwg$ form
a $\C$-vector space of dimension $k$. Using the theory of theta functions,
one can moreover prove that any non trivial solution can be written:
$$
f = \text{constant} \times \thq(x/x_1) \cdots \thq(x/x_k), 
\quad x_1 \cdots x_k = c.
$$
Thus $\div_\Eq(f)$ is an effective divisor of degree $k$ and evaluation
$\overline{(-1)^k c} \in \Eq$, \ie:
$$
\div_\Eq(f) = [\alpha_1] + \cdots + [\alpha_k], \text{~where~}
\alpha_1,\ldots,\alpha_k \in \Eq \text{~and~}
\alpha_1 + \cdots + \alpha_k = \overline{(-1)^k c}.
$$


\subsubsection{$q$-logarithms}
\label{subsubsection:qlog}

For $E = \Ra, \Rwg, \Rw, \ldots$ there is a so-called ``$q$-De Rham complex'':
$$
0 \longrightarrow \C \longrightarrow E \overset{\sq - 1}{\longrightarrow}
E \longrightarrow \C \longrightarrow 0.
$$
The meaning of the left side is that $q$-constants (\ie\ solutions of
$\sq f = f$) are true constants. The right side map $E \to \C$ sends
$\sum a_n x^n$ to $a_0$. It is related to so-called \emph{$q$-logarithms},
\ie\ solutions of $\sq f - f = 1$ (more generally of $\sq f - f = c \in \C^*$). \\

Any solution $f \in \Kw$ of $\sq f - f = c \in \C^*$ can be uniquely
extended to $\Kwg$ and it has poles, as can be seen for instance by
integration along the boundary of the fundamental annulus $\Cq$
(and using Cauchy formula). The simplest solutions are obtained as
follows; let:
$$
l_q(x) := x \dfrac{\thqp(x)}{\thq(x)} \cdot
$$
Then $\sq l_q - l_q = -1$ and $l_q$ has simple poles over $[-1;q]$ and
nowhere else. More generally (and more precisely), the solutions of
$\sq f - f = c \in \C^*$ having one simple pole modulo $q^\Z$ are the
$-c l_q(x/a) + b$, $a \in \C^*$, $b \in \C$. \\

More generally, we shall repeatedly use the following fact:

\begin{lem}
\label{lem:logarithmiccase}
Let $c \in \C^*$ and set $\phi(x) := \thq\left(\frac{x}{c}\right)$ and 
$\psi(x) := x \phi'(x) = \frac{x}{c} \thqp\left(\frac{x}{c}\right)$.
Then, if $f,g \in \Rwg$ are such that $\sq f = \dfrac{c}{x} f$ and
$\sq g = \dfrac{c}{x} (g - f)$, then we have $f = \alpha \phi$ and
$g = \alpha \psi + \beta \phi$ for some $\alpha,\beta \in \C$.
\end{lem}
\Pr
The function $f/\phi$ must be elliptic with at most simple poles, so it must
be constant: $f = \alpha \phi$. If $f = 0$, a similar argument applies to $g$.
Otherwise, $g/f$ must be a $q$-logarithm with at most a single pole modulo $q^\Z$
and we use the remark preceding the lemma.
\Finprcourt


\subsection{Gauge transformations}
\label{subsection:gaugetransformations}

Let $K$ be any of the fields over which $\sq$ can be defined and let
$A,B \in \GLn(K)$. Then, formally, if $X$ is a column vector solution
of the system:
\begin{equation}
\label{sqX=AX}
\sq X = A X,
\end{equation}
then one gets a solution $Y = F X$ of $\sq Y = B Y$ if:
\begin{equation}
\label{B = F[A]}
B = F[A] := (\sq F) A F^{-1}.
\end{equation}
We shall symbolize this relation by the diagram:
$$
A \overset{F}{\rightarrow} B.
$$
Indeed, $F$ can be seen as a morphism (actually an isomorphism) from $A$ to
$B$ in some category. It is easy to check that $A \overset{I_n}{\rightarrow} A$
and that from $A \overset{F}{\rightarrow} B$ and $B \overset{G}{\rightarrow} C$
one can infer $A \overset{GF}{\rightarrow} C$.


\subsection{Local reduction for fuchsian equations}
\label{subsection:LRFE}

Assumptions, definitions and results are stated here at $0$; the
corresponding facts at $\infty$ are obtained by substituting $1/x$
for $x$. Detailed statements and proofs are given in \cite{JSAIF}. \\

Let $A \in \GLn(\Ka)$ be such that $A(0) \in \GLn(\C)$, meaning that $A(x)$ 
is well defined at $x = 0$ and that its value is invertible. We also say
that $A$ is \emph{regular} at $0$. Thus actually $A \in \GLn(\Ra)$. 

\begin{defn}
\label{defn:usualnonresonant}
We say that $A$ is \emph{non resonant} at $0$ if 
$\Sp\ A(0) \cap q^{\N^*} \Sp\ A(0) = \emptyset$; said otherwise, for every
$c,d \in \Sp\ A(0)$, if $c \equiv d$, then $c = d$. 
\end{defn}

\begin{prop}
\label{prop:usualnonresonant}
(i) Let $A \in \GLn(\Ka)$ be such that $A(0) \in \GLn(\C)$. Then there
exists $F \in \GLn(\Kg)$ such that $B := F[A] \in \GLn(\Ka)$ and $B$ is non
resonant at $0$. \\
(ii) Let $A \in \GLn(\Ka)$ be non resonant at $0$. Then there exists 
a unique $F \in \GLn(\Ra)$ such that $F(0) = I_n$ and $A = F[A(0)]$.
\end{prop}

We shall use the following variant of the second statement:

\begin{cor}
If $A(0) = C R C^{-1}$, $C \in \GLnc$, $R = \Diag(\rho_1,\ldots,\rho_n)$
and if $\rho_i \not\equiv \rho_j$ for $i \neq j$ (so $A(0)$
is at the same time non resonant and semi simple), then there
is a unique $F \in \GLn(\Ra)$ such that $F(0) = C$ and $A = F[R]$.
\end{cor}

\paragraph{Normal forms.}

Combining the two statements in the proposition, we get:

\begin{cor}
Let $A \in \GLn(\Ka)$ be such that $A(0) \in \GLn(\C)$. Then there
exists $F \in \GLn(\Ka)$ and $A^{(0)} \in \GLnc$ such that $A = F[A^{(0)}]$.
Moreover, $A^{(0)}$ can be taken such that $\Sp(A^{(0)}) \subset \Cq$ (the
fundamental annulus); it is then unique up to conjugacy.
\end{cor}

Indeed, more generally, if $A_1,A_2 \in \GLn(\C)$ have all their eigenvalues
in the fundamental annulus, then:
$$
\left(A_2 = F[A_1], F \in \GLn(\Ka)\right) \Longrightarrow F \in \GLn(\C).
$$
(Without the assumption on eigenvalues, one could still deduce that $F$ is
a Laurent polynomial.)


\subsection{Singularities of meromorphic matrices}
\label{subsection:singularitiesofmeromorphicmatrices}

Let $M \in \GLn(\Kwg)$. We call \emph{singularities} of $M$ its poles
as well as the poles of $M^{-1}$. The singular locus is written $\Sigma_M$:
$$
\Sigma_M := \{\text{Poles of~} M\} \cup \{\text{Poles of~} M^{-1}\}.
$$
Writing $\tilde{M} := {}^t \text{com}(M)$ the transpose of the comatrix of 
$M$, we have, by Cramer's relations:
$$
M \tilde{M} = \tilde{M} M = (\det M) I_n,
$$
whence:
$$
\Sigma_M := \{\text{Poles of~} M\} \cup \{\text{Zeroes of~} \det M\}.
$$
In particular, if $M \in \GLn(\Kwg) \cap \Matn(\Rwg)$, then $\Sigma_M$
is the set of zeroes of $\det M$, and we can speak of multiplicity: the
\emph{multiplicity} of a singularity is its multiplicity as a zero of the
non trivial holomorphic function $\det M$.


\subsection{Birkhoff factorisation of analytic matrices}
\label{subsection:BFAN}

The \emph{preliminary theorem} of Birkhoff \cite[p. 266-267]{Birkhoff1}, 
stated in the basic case of a single simple contour, is the following:

\begin{thm}
Let $C$ a simple closed analytic curve on $\Sr$ separating $0$
from $\infty$ and let $D_0 \ni 0$, $D_\infty \ni \infty$ the connected
components of $\Sr \setminus C$. Let $M(x)$ an analytic invertible matrix
in a neighborhood of $C$ (\ie\ $x \mapsto M(x)$ is analytic with values
in $\GLnc$). Then there exists open neighborhoods $V_0$ of $\overline{D_0}$
and $V_\infty$ of $\overline{D_\infty}$ and analytic matrices $M_0$ on $V_0$
and $M_\infty$ on $V_\infty$ such that:
\begin{enumerate}
\item{$M_0 = M_\infty M$ in a neighborhood of $C$ contained in 
$V_0 \cap V_\infty$.}
\item{$M_0$ is regular (\ie\ holomorphic with holomorphic inverse) over $V_0$.}
\item{$M_\infty$ is regular over $V_\infty \setminus \{\infty\}$ and holomorphic
at $\infty$.}
\end{enumerate}
\end{thm}

Note however that the last condition cannot be sharpened, one cannot in general
require that $M_\infty(\infty) \in \GLnc$. It is easy to state variants, where
$M_0$ and $M_\infty$, resp. $M_0$ and $M_0^{-1}$ exchange their positions, etc. \\

We shall apply the theorem with $M \in \GLn(\Kwg)$. Then the relations
$M_0 = M_\infty M$ and $M_\infty = M_0 M^{-1}$ automatically allow for an
invertible meromorphic extension of $M_0$ and $M_\infty$ over $\C^*$, and
we simply write $M_0, M_\infty \in \GLn(\Kwg)$. Moreover, the regularity
conditions in the conclusion of the theorem then say that $M_0$ has the
same singularities as $M$ over $V_\infty \setminus \{\infty\}$ and that
$M_\infty$ as the same singularities as $M$ over $V_0$. In particular:

\begin{cor}
\label{cor:BFAN}
Let $M \in \GLn(\Kwg) \cap \Matn(\Rwg)$ with singular locus 
$\Sigma = \det^{-1}(0)$. We assume that $\det M$ has only simple zeroes.
Let $C$ as in the preliminary theorem (so $C$ does not meet $\Sigma$) and
write $\Sigma_0 := \Sigma \cap D_0$, $\Sigma_\infty := \Sigma \cap D_\infty$.
Then one has a factorisation $M = M_0^{-1} M_\infty$ over $\C^*$, with:
\begin{enumerate}
\item{$M_0$ is regular over $\C \setminus \Sigma_\infty$, $M_0^{-1}$ is
holomorphic over $\C$ and $\det M_0^{-1}$ has simple zeroes over $\Sigma_\infty$.}
\item{$M_\infty$ is regular over $\C^* \setminus \Sigma_0$, holomorphic over
$\Sr \setminus \{0\}$ and $\det M_\infty$ has simple zeroes over $\Sigma_0$.}
\end{enumerate}
\end{cor}


\subsection{Rational classification of fuchsian systems}

As in the theory of differential equations, one of the main problems is
the rational classification of rational systems. We say that
$A,B \in \GLn(\Kg)$ are \emph{globally} or \emph{rationally equivalent}
if there exists a rational gauge transformation $F \in \GLn(\Kg)$ 
such that $B = F[A]$. This is plainly an equivalence relation. \\

Again as in the theory of differential equations, the first step towards
global classification is local classification. The weaker equivalence
relation induced by gauge transformations $F \in \GLn(\Ka)$, resp.
$F \in \GLn(\Kai)$, is called \emph{local analytic\footnote{Since this
work is restricted to fuchsian systems, we have no use for \emph{formal}
classification.} equivalence at $0$}, resp. at $\infty$. \\

As already noted in \ref{subsection:gaugetransformations}, it will be convenient
to denote gauge transformations by diagrams:
$$
B = F[A] := (\sq F) A F^{-1} \text{~is denoted~} A \overset{F}{\longrightarrow} B.
$$
The reason is that there is a more general notion of (rational or local 
analytic) morphism from $A \in \GLn(\Kg)$ to $B \in \GL_p(\Kg)$, defined 
as a rectangular $p \times n$ (rational or local analytic) matrix $F$ such 
that $(\sq F) A = B F$. Gauge transformations then correspond to (rational 
or local analytic) isomorphisms. Thus for instance we can compose gauge
transformations $A \overset{F}{\rightarrow} B$ and 
$B \overset{G}{\rightarrow} C$ to obtain $A \overset{GF}{\rightarrow} C$,
meaning that $(GF)[A] = G[F[A]]$. We also can use commutative diagrams,
invert arrows, etc. 

\begin{defn}
The system $A \in \GLn(\Kg)$ is said to be \emph{strictly fuchsian} at $0$,
resp. at $\infty$, if $A(0) \in \GLnc$, resp. $A(\infty) \in \GLnc$.
It is said to be \emph{fuchsian} at $0$, resp. at $\infty$, if it is 
locally analytically equivalent at $0$, resp. at $\infty$, to a system
which is strictly fuchsian at $0$, resp. at $\infty$.
\end{defn}

The following was proved in \cite[2.1]{JSAIF}:

\begin{prop}
If $A \in \GLn(\Kg)$ is fuchsian at $0$ and at $\infty$, it is rationally
equivalent to a system which is strictly fuchsian at $0$ and at $\infty$.
\end{prop}

Note that for every gauge transformation $F$ and every $f \in \Kg^*$, one has:
$$
F[A] = B \Longrightarrow F[fA] = f B \quad \text{so that} \quad
A \sim B \Longrightarrow fA \sim fB
$$
for any of the above equivalence relations. Thus, for the rational classification
of rational systems, we may and shall restrict to the case that $A$ is a polynomial
matrix which is invertible as a rational matrix: $A \in \GLn(\Kg) \cap \Matn(\Rg)$.
For fuchsian systems, this can be made more precise:

\begin{lem}
Let $B \in \GLn(\Kg)$ be strictly fuchsian at $0$ and at $\infty$ and let 
$f$ the lcm of all the denominators of its coefficients, so that $A := f B$ 
is polynomial:
$$
A = A_0 + x A_1 + \cdots + x^\mu A_\mu \in \Matn(\Rg), \mu \in \N, A_\mu \neq 0.
$$
and the gcd of the coefficients of $A$ is $1$. Then $A_0,A_\mu \in \GLnc$.
\end{lem}
\Pr
Since $A(0) \in \GLnc$, we see that $f(0) \neq 0$, so that $A_0 \in \GLnc$. 
At infinity, $B \equiv B(\infty)$, so $A \equiv C x^N B(\infty)$, where
$C X^N$ is the leading term of $f$. Thus $N = \mu$ and $A_\mu = C B(\infty)$.
\Finprcourt



\section{A Birkhoff type classification theorem}
\label{section:A Birkhoff type classification theorem}

Birkhoff classification theorem in \cite{Birkhoff1} is a form of Riemann-Hilbert
correspondence for $q$-difference equations. For reasons explained in subsection
\ref{subsubsection:Comparison with Birkhoff classification} (see also subsection
\ref{subsubsection:uniformsolutions}), we use a variant of Birkhoff connection matrix
(our matrix $M$ introduced in corollary \ref{cor:traficotedconnectionmatrix}). \\

So from now on, we assume, just as Birkhoff did, that $A$ has the form:
$$
A = A_0 + x A_1 + \cdots + x^\mu A_\mu \in \Matn(\Rg), \mu \in \N, A_0,A_\mu \in \GLnc.
$$
We consider as \emph{local data} the conjugacy classes of $A_0,A_\mu$ (this is
for $0$ and $\infty$) and the zeroes of $\det A(x)$ (this is for so-called
\emph{intermediate singularities}). We do classification for \emph{fixed local data}.
We intend, in a future work, to describe the space of monodromy data as fibered
above a base, the space of possible local data. \\

In a first version of the theorem (theorem \ref{thm:RHBversion1}), we add nonresonancy
assumptions that are generically satisfied (these are the same assumptions as in
\cite{Birkhoff1} and also those taken by Jimbo and Sakai in \cite{JimboSakai}). Then
we give a more general and slightly less precise version (theorem \ref{thm:RHBversion2})
which we shall need in a special case.

\begin{rmk}
Readers interested mainly in character varieties of $q$-Painlev\'e equations should
skip the proofs in this section (they are standard $q$-difference technology) and
concentrate on the constructions and on the statements about them.
\end{rmk}


\subsection{Classification theorem for non resonant systems}

Here we add the following hypotheses:
\begin{itemize}
\item{$A_0$ is non resonant in the following strong sense: 
$$
\Sp A_0 = \{\rho_1,\ldots,\rho_n\} \subset \C^* \text{~and~}
i \neq j \Longrightarrow \rho_i \not\equiv \rho_j.
$$
}
\item{$A_\mu$ is non resonant in the strong sense:
$$
\Sp A_\mu = \{\sigma_1,\ldots,\sigma_n\} \subset \C^* \text{~and~}
i \neq j \Longrightarrow \sigma_i \not\equiv \sigma_j.
$$
}
\end{itemize}

\begin{rmk}
Non resonancy in the ``weak'' sense would allow for multiple eigenvalues (see
definition \ref{defn:usualnonresonant}). This weaker property can always be achieved
up to rational gauge transformation (proposition \ref{prop:usualnonresonant}).
Actually, any fuchsian $A(x)$ is rationally equivalent to some strictly fuchsian
$B$ such that all the eigenvalues of $B(0)$ are in $\Cq$. Strong non resonancy
defined here is equivalent to weak non resonancy plus separability (all eigenvalues
distinct).
\end{rmk}

Obviously, $A_0$ and $A_\mu$ are then semisimple. We shall set:
\begin{align*}
R := \Diag(\rho_1,\ldots,\rho_n), \\
S := \Diag(\sigma_1,\ldots,\sigma_n),
\end{align*}
so that $A_0$ and $R$ are conjugate, and the same for $A_\mu$ and $S$. Note 
that, with those notations, we implicitly fixed an order on the spectra. \\

From the given form $A = A_0 + \cdots + x^\mu A_\mu$, we draw:
$$
\det A(x) = \sigma_1 \cdots \sigma_n (x - x_1) \cdots (x - x_N), \quad
N := n \mu, \quad x_1,\ldots,x_N \in \C^*,$$
subject to Fuchs relation:
$$
x_1 \cdots x_N = (-1)^N \dfrac{\rho_1 \cdots \rho_n}{\sigma_1 \cdots \sigma_n} 
$$
We shall add one more strong non resonancy condition:
\begin{itemize}
\item $k \neq l \Longrightarrow x_k \not\equiv x_l$. \\
\end{itemize}


\subsubsection{Local reductions}
\label{subsubsection:Local reductions}

In this section, we consider $R$, $S$ and $\x := \{x_1,\ldots,x_N\}$ as
fixed and subject to the above strong non resonancy conditions and also
to Fuchs relation. \\

Let $E_{R,S,\x}$ the set of matrices $A = A_0 + \cdots + x^\mu A_\mu$
with all $A_i \in \Matnc$ and such that:
\begin{itemize}
\item{$A_0$ is conjugate to $R$;}
\item{$A_\mu$ is conjugate to $S$;}
\item{$\det A(x)$ vanishes at the $x_i \in \x$.}
\end{itemize}
Together, those conditions imply that 
$\det A(x) = \sigma_1 \cdots \sigma_n (x - x_1) \cdots (x - x_N)$.
(Recall that $\deg \det A = n \mu = N$.) \\

We denote $\sim$ the equivalence relation induced on $E_{R,S,\x}$ by rational
equivalence. \emph{We intend to describe the quotient set $E_{R,S,\x}/\sim$.}
This is the meaning of Birkhoff's interpretation of the Riemann-Hilbert problem.

\begin{lem}
(i) Let $C \in \GLnc$ such that $A_0 = C R C^{-1}$. Then there exists 
a unique $M_0 \in \GLn(\Ra)$ such that $M_0(0) = C$ and $M_0[R] = A$.
(Recall these notations were introduced in \ref{subsection:gaugetransformations}). \\
(ii) Let $D \in \GLnc$ such that $A_\mu = D S D^{-1}$. Then there exists 
a unique $M_\infty \in \GLn(\Rai)$ such that $M_\infty(\infty) = D$ and 
$M_\infty[S x^\mu] = A$. \\
(iii) Let $C'$, $D'$ alternative choices for the conjugating matrices $C$, $D$
and $M_0'$, $M_\infty'$ the resulting gauge transformations as in (i), (ii).
Then there exist constant diagonal $n \times n$ matrices $\Gamma$, $\Delta$
such that $C' = C \Gamma$, $D' = D \Delta$; and then $M_0' = M_0 \Gamma$,
$M_\infty' = M_\infty \Delta$.
\end{lem}
\Pr
Statements (i) and (ii) were proved in \ref{subsection:singularitiesofmeromorphicmatrices}. \\
Proof of (iii): $\Gamma := C^{-1} C'$ commutes with $R$ and $\Delta := D^{-1} D'$
with $S$, so they are diagonal. Then $M_0 \Gamma$ and $M_\infty \Delta$ satisfy
the adequate relations, so by unicity they are respectively equal to $M_0'$, 
$M_\infty'$.
\Finpr

All this can be summarized by the following commutative diagram:
$$
\xymatrix{
R \ar@<0ex>[dd]_{\Gamma} \ar@<0ex>[drr]^{M'_0}  
& & & 
& S x^\mu  \ar@<0ex>[dll]_{M'_\infty} \ar@<0ex>[dd]^{\Delta} \\
& & A & &  \\
R \ar@<0ex>[urr]^{M_0}  
& & & 
& S x^\mu  \ar@<0ex>[ull]_{M_\infty} 
}
$$

\begin{prop}[Properties of $M_0$ and $M_\infty$]
(i) $M_0 \in \GLn(\Ra)$ admits a unique extension $M_0 \in \GLn(\M(\C))$
such that:
\begin{itemize}
\item{$M_0$ has simple poles over $\x q^{-\N}$ (and nowhere else);}
\item{$M_0^{-1}$ is holomorphic all over $\C$.}
\end{itemize}
(ii) $M_\infty \in \GLn(\Rai)$ admits a unique extension 
$M_\infty \in \GLn(\M(\C_\infty))$ such that:
\begin{itemize}
\item{$M_\infty$ is holomorphic all over $\C_\infty$;}
\item{$M_\infty^{-1}$ has simple poles over $\x q^{\N^*}$ (and nowhere else).}
\end{itemize}
\end{prop}
\Pr
(i) We use the arrow $R \overset{M_0}{\longrightarrow} A$, \ie\ the equality
$A = M_0[R] = (\sq M_0) R M_0^{-1}$ first in the clearly equivalent forms:
$M_0 = A^{-1} (\sq M_0) R$ and $M_0^{-1} = R^{-1} (\sq M_0^{-1}) A$. \\
The second relation allows us to extend $M_0^{-1}$, which is initially defined
and holomorphic over some open disk $\Do(0,r)$, $r > 0$, to $\Do(0,q^{-1}r)$,
where $\lmod q^{-1}r \rmod = \lmod q \rmod^{-1} r > r$ since $\lmod q \rmod < 1$.
Iterating, we get a holomorphic extension to $\C$. \\
The first relation shows that on any open disk $\Do(0,r)$, $M_0$ has the same
poles as $\sq M_0$, \ie\ those of $M_0$ over the smaller disk $\Do(0,qr)$; but
one must add the poles of $A^{-1}$ if any. Iterating yields the conclusion. \\
One could also argue using only the determinant of the second relation:
$\frac{\sq \det M_0^{-1}}{\det M_0^{-1}} = \frac{\det A}{\rho_1 \cdots \rho_n} \cdot$ \\
(ii) Similarly, the arrow $S x^\mu \overset{M_\infty}{\rightarrow} A$, \ie\
the equality $A = M_\infty[S x^\mu] = (\sq M_\infty) (S x^\mu) M_\infty^{-1}$
translate into $\sq M_\infty = A M_\infty (S x^\mu)^{-1}$ and the argument goes
on the same lines (here we use $\sq$ to expand disks centered at $\infty$).
\Finprcourt


\subsubsection{Connection matrix and Riemann-Hilbert-Birkhoff correspondence}
\label{subsubsection:CMRHBC}

This section can be best understood with the help of the following commutative diagram:
$$
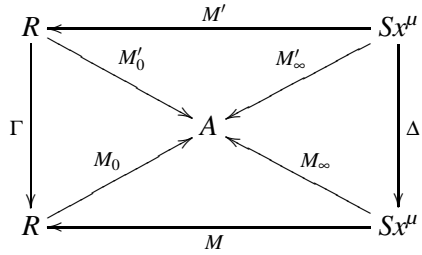
\xymatrix{
R \ar@<0ex>[dd]_{\Gamma} \ar@<0ex>[drr]^{M'_0}  
& & & 
& S x^\mu  \ar@<0ex>[dll]_{M'_\infty} \ar@<0ex>[dd]^{\Delta} \ar@<0ex>[llll]_{M'} \\
& & A & &  \\
R \ar@<0ex>[urr]^{M_0}  
& & & 
& S x^\mu  \ar@<0ex>[ull]_{M_\infty} \ar@<0ex>[llll]^{M} 
}
$$
For instance the south-west and south-east diagonal arrows respectively mean that
$M_0[R] = A$ and that $M_\infty[S x^\mu] = A$, so that:
$$
\left((M_0)^{-1} M_\infty\right)[S x^\mu] = R,
$$
and similarly for $M_0'$ and $M_\infty'$. All this can be read on the diagram.

\begin{cor}
\label{cor:traficotedconnectionmatrix}
Set $M := M_0^{-1} M_\infty \in \GLn(\Kwg)$. Then:
\begin{itemize}
\item{$\sq M = R M (S x^\mu)^{-1}$.}
\item{$M$ is holomorphic all over $\C^*$.}
\item{$M^{-1}$ has simple poles over $[\x;q] = \x q^\Z$ (and nowhere else); 
equivalently, $\det M$ has simple zeroes over $[\x;q]$ (and nowhere else).}
\end{itemize}
\end{cor}

We shall write $F_{R,S,\x}$ the set of such matrices:
$$
F_{R,S,\x} := \left\{M \in \Matn(\Rwg) \TQ \sq M = R M (S x^\mu)^{-1}
\text{~and all zeroes of~} \det M \text{~are simple and lay over~} [\x;q] \right\}.
$$
Note that this set actually depends only of the image of $\x$ in $\Eq$, not on
$\x \subset \C^*$ itself.

\paragraph{Gauge freedom.}
We saw that, $A$ being given, $M_0$ and $M_\infty$ are uniquely determined up
to the right action of the group $\Dnc \subset \GLnc$ of diagonal matrices.
From the relations $M_0 \sim M_0 \Gamma$, $M_\infty \sim M_\infty \Delta$,
$\Gamma, \Delta \in \Dnc$, we deduce that $M \sim \Gamma^{-1} M \Delta$.
We are thus led to introduce the following right action of $\Dnc \times \Dnc$
on $F_{R,S,\x}$:
$$
M^{(\Gamma, \Delta)} := \Gamma^{-1} M \Delta.
$$
The reader may check that this is indeed a right action\footnote{Later in the
text, we shall rather use the left action $M \mapsto \Gamma M \Delta^{-1}$.
Of course, the equivalence classes (orbits) are the same.} and that $F_{R,S,\x}$
is stable under this action. We shall write $M \sim M^{(\Gamma, \Delta)}$ the
corresponding equivalence relation on $F_{R,S,\x}$ and $F_{R,S,\x}/\sim$ the quotient
of $F_{R,S,\x}$ under this action and equivalence relation. As a consequence, we see
that we have constructed a well defined map:
$$
E_{R,S,\x} \, \longrightarrow \, F_{R,S,\x}/\sim,
$$
mapping $A$ to the equivalence class of $M$. \\

From now on, we shall write:
$$
\E_{R,S,\x} := E_{R,S,\x}/\sim \text{~~and~~} \F_{R,S,\x} := F_{R,S,\x}/\sim.
$$
Moreover, if no confusion is to be feared, we shall frequently omit the
indication of local data and abreviate:
$$
\E := \E_{R,S,\x} \text{~~and~~} \F := \F_{R,S,\x}.
$$

\begin{prop}
\label{prop:RHBcorresp}
The above map goes to the quotient and defines a ``Riemann-Hilbert-Birkhoff
correspondence'':
\begin{equation}
\label{equation:RHBcorresp}
\E_{R,S,\x} = E_{R,S,\x}/\sim \;\; \longrightarrow \, \F_{R,S,\x} = F_{R,S,\x}/\sim.
\end{equation}
\end{prop}
\Pr
Let $B = B_0 + \cdots + B_\mu x^\mu \in E_{R,S,\x}$, (so that $B_0 \sim R$,
$B_\mu \sim S$ and $\det B$ has simple zeroes at $\x$) and assume that
$B = F[A]$, $F \in \GLn(\Kg)$. We have a commutative diagram:
$$
\xymatrix{
& & A \ar@<0ex>[dd]^{F} & & \\
R \ar@<0ex>[rru]^{M_0} \ar@<0ex>[rrd]^{N_0}  
& & & & 
S x^\mu  \ar@<0ex>[llu]_{M_\infty} \ar@<0ex>[lld]^{N_\infty}  \\
& & B & & 
}
$$
Let $\Gamma := N_0^{-1} \F M_0 \in \GLn(\Ka)$ and 
$\Delta := N_\infty^{-1} F M_\infty \in \GLn(\Kai)$. 
Then $\Gamma = (\gamma_{i,j})_{1 \leq i,j \leq n}$ is an automorphism of $R$ 
and $\Delta = (\delta_{i,j})_{1 \leq i,j \leq n}$ is an automorphism of $S x^\mu$, 
so that:
$$
\Gamma[R] = R \Longrightarrow \sq \Gamma R = R \Gamma \Longrightarrow
\forall i,j = 1,\ldots,n \;,\; 
\sq \gamma_{i,j} = \dfrac{\rho_i}{\rho_j} \gamma_{i,j}.
$$
Since $\gamma_{i,j} \in \Ka$ and $\rho_i/\rho_j \not\in q^\Z$ for $i \neq j$,
we conclude that $\gamma_{i,j} = 0$ for $i \neq j$ and that 
$\gamma_{i,i} \in \C^*$ (it cannot be $0$ since $\Gamma$ is invertible) so
at last $\Gamma \in \Dnc$. A similar argument works for $\Delta$ (the scalar
$x^\mu$ factor gets simplified at once). \\
Now, from $F M_0 = N_0 \Gamma$ and $F M_\infty = N_\infty \Delta$ we draw:
$$
M = M_0^{-1} M_\infty = (F M_0)^{-1} F M_\infty = 
(N_0 \Gamma)^{-1} N_\infty \Delta = \Gamma^{-1} N_0^{-1} N_\infty \Delta = 
\Gamma^{-1} N \Delta
$$
as expected.
\Finprcourt

\begin{thm}[``Riemann-Hilbert-Birkhoff correspondence'', first version]
\label{thm:RHBversion1}
The map \eqref{equation:RHBcorresp} defined in proposition \ref{prop:RHBcorresp}
is bijective.
\end{thm}
\Pr
The proof comes in two parts.

\paragraph{Injectivity.}

Using the usual notations, let $A \in E_{R,S,\x}$, resp. $B \in E_{R,S,\x}$,
have image the class of $M = M_0^{-1} M_\infty$, resp. the class of 
$N = N_0^{-1} N_\infty$ in $F_{R,S,\x}/\sim$ and assume these images are the
same, that is $M \sim N$, so that $M = \Gamma^{-1} N \Delta$ where 
$\Gamma, \Delta \in \GLn(\Ka)$. Then:
$$
M_0^{-1} M_\infty = \Gamma^{-1} N_0^{-1} N_\infty \Delta \Longrightarrow
N_0 \Gamma M_0^{-1} = N_\infty \Delta M_\infty^{-1}.
$$
Call the latter matrix $F$. Then: 
$$
F \in \GLn(\M(\C)) \cap \GLn(\M(\C_\infty)) = \GLn(\M(\Sr)) = \GLn(\Kg).
$$
On the other hand, we have a commutative diagram:
$$
\xymatrix{
R \ar@<0ex>[rr]^{M_0} \ar@<0ex>[d]^{\Gamma} 
& & A \ar@{.>}[d]^{F} & & 
S x^\mu \ar@<0ex>[ll]_{M_\infty} \ar@<0ex>[d]^{\Delta} \\
R \ar@<0ex>[rr]^{N_0}
& & B  & & 
S x^\mu  \ar@<0ex>[ll]_{N_\infty}
}
$$
in which $F = N_0 \Gamma M_0^{-1} = N_\infty \Delta M_\infty^{-1}$ is, by force,
an isomorphism, \ie\ $F[A] = B$, so that the classes of $A$ and $B$ in
$E_{R,S,\x}/\sim$ are the same, which concludes the proof of injectivity.
Note that the relation $F[A] = B$ can also be deduced by direct computation:
$$
F M_0 = N_0 \Gamma \Longrightarrow 
\sq (F M_0) R (F M_0)^{-1} = \sq (N_0 \Gamma) R (N_0 \Gamma)^{-1} 
\Longrightarrow \sq F A F^{-1} = B.
$$

\paragraph{Surjectivity.}

$R,S,\x$ being given (and satisfying Fuchs relation), let $M \in F_{R,S,\x}$. 
We draw on the Riemann sphere $\Sr$ a closed analytic curve separating $\x$ 
from $q \x$, so that, with the notations of subsection \ref{subsection:BFAN},
$\Sigma_0 = q^{\N^*} \x$ and $\Sigma_\infty = q^{-\N} \x$. Using Birkhoff
factorisation theorem and in particular corollary \ref{cor:BFAN}, we obtain
a decomposition $M = M_0^{-1} M_\infty$, where:
\begin{itemize}
\item{$M_0$ is regular (\ie\ holomorphic, invertible with holomorphic inverse)
on $\C \setminus q^{-\N} \x$, with simple poles on $q^{-\N} \x$.}
\item{$M_0^{-1}$ is holomorphic over $\C$ and $\det M_0$ has simple zeroes
$q^{-\N} \x$.}
\item{$M_\infty$ is regular on $\C^* \setminus q^{\N^*} \x$, it is holomorphic
over $q^{\N^*} \x$ and $\det M_\infty$ has simple zeroes there. It is also
holomorphic at $\infty$, although it cannot be required to be regular there.}
\item{$M_\infty^{-1}$ is holomorphic over $\C^*$ and meromorphic at $\infty$.}
\end{itemize}
From the condition $\sq M = R M (S x^\mu)^{-1}$, expressed by the arrow
$R \overset{M}{\leftarrow} S x^\mu$, we see that $M_0[R] = M_\infty[S x^\mu]$.
Call $A$ this matrix, whence a diagram:
$$
R \overset{M_0}{\longrightarrow} A \overset{M_\infty}{\longleftarrow} S x^\mu
$$
We want to show that $A \in E_{R,S,\x}$ and that its class in
$\E_{R,S,\x} = E_{R,S,\x}/\sim$
is the preimage of the class of $M$ in $\F_{R,S,\x} =F_{R,S,\x}/\sim$. Clearly:
$$
A \in \GLn(\M(\C)) \cap \GLn(\M(\C_\infty)) = \GLn(\M(\Sr)) = \GLn(\Kg).
$$
Actually, more can be said. From the listed properties of $M_0$ and $M_\infty$,
one gets that $A$ is holomorphic over $\C$ and meromorphic at $\infty$, whence
polynomial:
$$
A = A_0 + \cdots + x^d A_d, \quad A_0, \ldots, A_d \in \Matnc,\; A_d \neq 0.
$$
Here $A_0 = A(0) = C R C^{-1} \in \GLnc$, where $C = M_0(0) \in \GLnc$. \\

From the given relations, we also see that $a(x) := \det A(x)$ has simple 
zeroes at $\x$ and nowhere else. Thus $a(x) = s (x - x_1) \cdots (x - x_N)$
for some $s \in \C^*$. Then $a(0) = s (-1)^N x_1 \cdots x_N$, but also
$a(0) = \det A_0 = \det C R C^{-1} = \det R$, so by Fuchs relation
$s =\det S$. Now we use the relation $A = M_\infty[S x^\mu]$; taking the
determinant and setting $f := \det M_\infty \in \Kai^*$, we draw:
$$
\dfrac{\sq f}{f} = \dfrac{a}{s x^N} = \prod_{i=1}^N (1 - x_i/x) \Longrightarrow
f = \phi \prod_{i=1}^N \dfrac{1}{(x_i/x;q)_\infty},
$$
where $\phi$ is elliptic. But at the same time, $\phi \in \Kai^*$, so that
actually $\phi \in \C^*$. This implies that $f(\infty) = \phi$, so that
$M_\infty$ is regular at $\infty$ (which Birkhoff factorisation did not
automatically imply). Setting $D := M_\infty(\infty) \in \GLnc$, we see
that $A = M_\infty[S x^\mu]$ is asymptotic to $(D S D^{-1}) x^\mu$ at $\infty$.
Since it is also asymptotic to $x^d A_d$, we get that $d = \mu$ and
$A_d = D S D^{-1} \in \GLnc$. It is then immediate that $A \in E_{R,S,\x}$ and 
the fact that its class in $\E_{R,S,\x} = E_{R,S,\x}/\sim$ is the antecedent of
the class of $M$ in $\F_{R,S,\x} = F_{R,S,\x}/\sim$ follows from the various
equalities we found out
during the computation (\ie\ the construction of $M$ goes through the
$M_0$ and $M_\infty$ used in the proof).
\Finprcourt


\subsubsection{Comparison with Birkhoff classification}
\label{subsubsection:Comparison with Birkhoff classification}

Recall the notation $e_{q,c}$ at the end of \ref{subsection:somefunctions}.
Let $e_R := \Diag(e_{q,\rho_1},\ldots,e_{q,\rho_n})$, so that 
$\sq e_R = R e_R = e_R R$. Then $\X^{(0)} := M_0 e_{q,R}$ is a solution
of the system $\sq X = A X$ which maybe considered ``local at $0$''. \\

Let likewise $e_S := \Diag(e_{q,\sigma_1},\ldots,e_{q,\sigma_n})$, so that 
$\sq e_S = S e_S = e_S S$. Then $\X^{(\infty)} := M_\infty e_{q,S} \thq^{-\mu}$ 
is a solution of the system $\sq X = A X$ which maybe considered ``local 
at $\infty$''. \\

Birkhoff connection matrix\footnote{The original definition of Birkhoff
(taken up by Jimbo and Sakai) involved multivalued choices as solutions
of the elementary equations $\sq e_c = c e_c$, resp. $\sq f = x^{-1} f$,
such as $x^{\log_q c}$, resp. $q^{- \log_q x (\log_q x-1)/2}$. This does not
impact the present discussion.} is then defined as:
$$
P := \X_0^{-1} \X_\infty \in \GLn(\Kwg)^{\sq} = \GLn(\M(\Eq)).
$$
In our notations, $P = e_R^{-1} M e_{q,S} \thq^{-\mu}$. We think of $M$ as
freed from the local contributions at $\infty$ present in $P$. In the
interpretation of Birkhoff connection matrix as encoding the monodromy,
we think of $M$ as encoding specifically the ``intermediate'' monodromy
related to the singularities of $A$ on $\C^*$. As explained in
\ref{subsubsection:uniformsolutions}, this was shown to be
necessary in Galois theory \cite{JSGAL} and we think it to be useful
here as well.


\subsection{A more general classification theorem}

The previous result, theorem \ref{thm:RHBversion1}, is only valid for semi-simple
local data $R$ and $S$ (equivalently, $A(0)$ and $A(\infty)$). Here we relax this
assumption. In the next version, theorem \ref{thm:RHBversion2}, we take the following
local data:
\begin{itemize}
\item $R, S \in \GLnc$ such that $\Sp R, \Sp S \subset \Cq$;
\item $x_1,\ldots,x_N \in \C^*$ such that
$i \neq j \Rightarrow x_i \not\equiv x_j$,
\item and moreover subject to Fuchs relation: 
$x_1 \cdots x_N = (-1)^N \dfrac{\det R}{\det S}$.
\end{itemize}

We define the set $E_{R,S,\x}$ and its equivalence relation exactly like at the
beginning of subsection \ref{subsubsection:Local reductions}; and the set $F_{R,S,\x}$
in the same way as just after corollary \ref{cor:traficotedconnectionmatrix}.
But we take the equivalences $M \sim M^{(\Gamma, \Delta)}$ among those induced by
matrices $\Gamma, \Delta \in \GLnc$ such that $\Gamma$ commutes with $R$ and
$\Delta$ commutes with $S$:
$$
\forall M,N \in F_{R,S,\x} \;,\;
M \sim N \underset{def}{\Longleftrightarrow}
\exists \Gamma, \Delta \in \GLnc \;:\;
\begin{cases}
N = M^{(\Gamma, \Delta)} := \Gamma^{-1} M \Delta, \\
[\Gamma,R] = [\Delta,S] = 0.
\end{cases}
$$
(Under the assumptions of \ref{subsubsection:Local reductions}, this boils down
to the previous definition.)

\begin{thm}[``Riemann-Hilbert-Birkhoff correspondence'', second version]
\label{thm:RHBversion2}
There is a natural bijection:
$$
\E_{R,S,\x} = E_{R,S,\x}/\sim \;\; \longrightarrow \; \F_{R,S,\x} = F_{R,S,\x}/\sim.
$$
\end{thm}
\Pr
We only sketch the modifications to the proof of the first version. \\
According to subsection \ref{subsection:LRFE}, we can write
$A = M^{(0)}[R] = M^{(\infty)}[S x^\mu]$, but these gauge matrices are not
unique (and respectively have poles at $0$, $\infty$). \\
The polarity properties of $M^{(0)}$ and $M^{(\infty)}$ on $\C^*$ are exactly
the same as before, because we only used the fact that these matrices were
regular (\ie\ holomorphic with holomorphic inverse) in a punctured neighborhood
of $0$, resp. $\infty$. The matrix $M := (M^{(0)})^{-1} M^{(\infty)}$ belongs to
$F_{R,S,\x}$. We thus obtain a correspondence between $E_{R,S,\x}$ and $F_{R,S,\x}$,
but not a mapping in either direction. \\
Let $N^{(0)}$, $N^{(\infty)}$ another choice of gauge transformations realizing
the same reductions. Then $\Gamma := (M^{(0)})^{-1} N^{(0)}$ is such that
$\Gamma[R] = R$. Using subsection \ref{subsection:LRFE} one shows that
$\Gamma \in \GLnc$ and $[\Gamma,R] = 0$. Similarly for $\Delta$. This gives
an injective map
$\E_{R,S,\x} = E_{R,S,\x}/\sim \; \rightarrow \; \F_{R,S,\x} = F_{R,S,\x}/\sim$. \\
The rest of the argument does not change.
\Finprcourt

\begin{rmk}
Local data $R$ and $S$ play a symmetric role in the following sense:
$M \in F_{R,S,\x} \Leftrightarrow {}^t\! M \in F_{{}^t\! S,{}^t\! R,\x}$; also
the equivalence relations on these two sets correspond to each other.
This observation allows one to shorten some case studies.
\end{rmk}


\subsection{Reducibility criteria}
\label{subsection:reducibilitycriteria}

We shall have need for the possibility of determining if a system $\sq X = AX$
is reducible by looking at its image $(R,M,S)$ by the Riemann-Hilbert-Birkhoff 
correspondence. We give such a criterion in the generic case, the general case
can be tackled similarly but the formulation would be more complicated and we
do not need it here.

\begin{thm}
\label{thm:reducibilitycriterion}
We assume strong non resonancy as above, \ie\ the $n$ eigenvalues of $R$, resp.
of $S$, are pairwise non congruent modulo $q^{\Z}$. Then the system 
$\sq X = A X$ is reducible if, and only if some matrix $M'$ obtained from $M$
by permutation of lines and columns is block triangular.
\end{thm}

\begin{cor}
\label{cor:reducibilitycriterion}
If $n = 2$, under the same generic assumptions, reducibility is equivalent to:
$M$ has a zero coefficient.
\end{cor}

To prove the theorem we use a tannakian criterion based on the Galois theory
as expounded for instance in \cite{RS3}. Recall that to the system $\sq X = A X$
is attached a Galois \emph{groupoid} with base $\{0,\infty\}$ and its 
canonical representation. It can be realized as $(G(0),G(0,\infty),G(\infty))$
operating on $(V(0),V(\infty))$, \ie\ $G(0)$, resp. $G(\infty)$ is an algebraic
group of automorphisms of the linear space $V(0)$, resp. $V(\infty)$; and
$G(0,\infty)$ is a set of isomorphisms $V(0) \rightarrow V(\infty)$. By Galois
correspondence, reducibility of the system is equivalent to reducibility of the
representation, \ie\ to the existence of non trivial (non zero and non whole)
subspaces $V_0$ of $V(0)$ and $V_\infty$ of $V(\infty)$ such that $G(0)$ leaves
$V_0$ stable, $G(\infty)$ leaves $V_\infty$ stable and $G(0,\infty)$ sends $V_0$
to $V_\infty$. \\

We can and will take $V(0) = V(\infty) = \C^n$. After \cite{RS3} we introduce:
\begin{enumerate}
\item The subgroup $G_0$ of $G(0) \subset \GLnc$ made up of all diagonal matrices 
$\Diag(\phi(\rho_1),\ldots,\phi(\rho_n))$ where $\phi: \C^* \rightarrow \C^*$
is a group morphism such that $\phi(q) = 1$;
\item the subgroup $G_\infty$ of $G(\infty) \subset \GLnc$ built similarly, except 
that we also allow invertible scalar matrices\footnote{This accounts for the non
trivial ``slope'' $\mu$ in $S x^\mu$.} and of course all resulting products;
\item the subset $G_{0,\infty}$ of $G(0,\infty)$ made up of all values $M(x)^{-1}$
at regular points (\ie\ where $M(x)$ is invertible).
\end{enumerate}
A preliminary fact is that these three components are included in the corresponding
components of the Galois groupoid. This is analogous to the classical fact that the
monodromy group of a (complex, linear, analytic) differential equation is included
in its Galois group. Then a Schlesinger type density theorem states that the whole
Galois groupoid is the Zariski closure of the subgroupoid generated by these three
components (two local components and a connection component). \\

Non congruent elements of $\C^*$ can be separated by morphisms of the above
kind. Under our strong non resonancy assumption, this implies that the 
subspaces of $V(0) = \C^n$ stable under $G_0$ are exactly those generated by 
a subset of the canonical basis; and the same at $\infty$. It follows that the 
canonical representation is reducible if and only if there are two non trivial 
(non empty, non whole) subsets $B_0,B_\infty$ of the canonical basis such that 
all invertible values $M(x)$ send $\Vect(B_\infty)$ isomorphically to 
$\Vect(B_0)$. The criterion of the theorem is just a rephrasing of that fact.


\subsection{The hypergeometric class}
\label{subsection:thehypergeometricclass}

This is the case\footnote{We shall abusively call $q$-hypergeometric the 
systems classified herebelow, without checking if they really come from a 
$q$-hypergeometric equation. By \cite{JRRigidity} this is generically true
but certainly false if the system is reducible, \ie, from the above study,
if at least one coefficient of $M$ vanishes.} $n = 2$, $\mu = 1$, $N = 2$, 
$x_1 x_2 = \det R/\det S$, $x_1/x_2 \not\in q^\Z$. As in the case of 
ordinary differential equations, we shall (generically) find \emph{rigidity}, 
\ie\, these local data being fixed, there are no continuous moduli. \\

We may assume that each $R,S$ take one of the following forms, respectively called
\emph{generic}, \emph{trivial} and \emph{logarithmic}:
\begin{align*}
R &= \begin{pmatrix} \rho_1 & 0 \\ 0 & \rho_2 \end{pmatrix}, \rho_1 \neq \rho_2, 
\text{~or~} \begin{pmatrix} \rho & 0 \\ 0 & \rho \end{pmatrix} \text{~or~}
\begin{pmatrix} \rho & \rho \\ 0 & \rho \end{pmatrix}
\\
S &= \begin{pmatrix} \sigma_1 & 0 \\ 0 & \sigma_2 \end{pmatrix}, \sigma_1 \neq \sigma_2, 
\text{~or~} \begin{pmatrix} \sigma & 0 \\ 0 & \sigma \end{pmatrix} \text{~or~}
\begin{pmatrix} \sigma & \sigma \\ 0 & \sigma \end{pmatrix}
\end{align*}
where $\rho_1,\rho_2,\rho,\sigma_1,\sigma_2,\sigma \in \Cq$. Also we shall
write for short $E := E_{R,S,\x}$ and $F := F_{R,S,\x}$; and also $\E := E/\sim$
and $\F := F/\sim$. \\

Note that, from the relations $x_1 x_2 = \det R/ \det S$ and
$\sq (\det M) = (\det R/ \det S) \det M$, we draw, using subsection
\ref{subsection:someequations}, that $\det M$ vanishes at $x_1$ if,
and only if, it vanishes at $x_2$, so we need use only one of these conditions
to test whether $M \in F$. \\

Also note that cases 1, 2 and 3 herebelow are special in that $x_1,x_2$ are imposed
by $R,S$ (other values would mean that $E$ is empty); and that cases 1 and 3 do
not fall under our assumptions for $x_1,x_2$ (but we all the same describe
$\E$ and $\F$).

\paragraph{Case 1, trivial/trivial:}
$R = \begin{pmatrix} \rho & 0 \\ 0 & \rho \end{pmatrix}$,
$S = \begin{pmatrix} \sigma & 0 \\ 0 & \sigma \end{pmatrix}$. \\
Here $E = \{(\rho + \sigma x) I_2\}$ and \emph{a fortiori} $\E$ is a singleton.
Also $\det A(x) = (\rho + \sigma x)^2$ so $x_1 = x_2 = -\rho/\sigma$ so we are not
within the assumptions of our theorem. \\
Actually, one sees easily that matrices $M \in F$ have the form
$\thq\left(\frac{\sigma}{\rho}x\right) C$, $C \in \GL_2(\C)$ being arbitrary,
with equivalences $M \sim \Gamma^{-1} M \Delta$ for arbitrary
$\Gamma, \Delta \in \GL_2(\C)$ so $\F$ is also a singleton.

\paragraph{Case 2, trivial/generic:}
$R = \begin{pmatrix} \rho & 0 \\ 0 & \rho \end{pmatrix}$,
$S = \begin{pmatrix} \sigma_1 & 0 \\ 0 & \sigma_2 \end{pmatrix}$. \\
Here, $E = \left\{\text{conjugates of~} \rho I_2 + x S =
\begin{pmatrix} \rho + \sigma_1 x & 0 \\
0 & \rho + \sigma_2 x \end{pmatrix}\right\}$
(conjugacy by $\GL_2(\C)$), so $\E$ is a singleton. We have
$\det A(x) = (\rho + \sigma_1 x)(\rho + \sigma_2 x)$ so the only non void
possibility (up to reindexing) is $x_i = -\rho/\sigma_i$, $i = 1,2$; we are under
the assumptions of our theorem. \\
Matrices $M \in F$ have the form
$$
M = \begin{pmatrix}
\alpha_{1,1} \thq\left(\frac{\sigma_1}{\rho} x\right) &
\alpha_{1,2} \thq\left(\frac{\sigma_2}{\rho} x\right) \\
\alpha_{2,1} \thq\left(\frac{\sigma_1}{\rho} x\right) &
\alpha_{2,2} \thq\left(\frac{\sigma_2}{\rho} x\right)
\end{pmatrix} = C T(x) \text{~where~} C \in \Mat_2(\C) \text{~and~}
T(x) := \begin{pmatrix} \thq\left(\frac{\sigma_1}{\rho} x\right) & 0 \\
0 & \thq\left(\frac{\sigma_2}{\rho} x\right) \end{pmatrix}.
$$
Since $\det T$ vanishes at $x_1,x_2$ but not identically, we get that
$E = \{C T(x) \tq C \in \GL_2(\C)\}$. The action of $\Gamma,\Delta$
comes here with arbitrary $\Gamma$ so $\F$ is also a singleton.

\paragraph{Case 3, trivial/logarithmic:}
$R = \begin{pmatrix} \rho & 0 \\ 0 & \rho \end{pmatrix}$,
$S = \begin{pmatrix} \sigma & \sigma \\ 0 & \sigma \end{pmatrix}$. \\

Here, $E = \{\text{conjugates of~} \rho I_2 + x S\}$ (conjugacy by
$\GL_2(\C)$), so $\E$ is a singleton. We have
$\det A(x) = (\rho + \sigma x)^2$ so the only non void possibility
is $x_1 = x_2 = -\rho/\sigma$; we are not under
the assumptions of our theorem. Yet we go on ! \\
Coefficients of matrices $(m_{i,j}) \in F$ must satisfy the functional
equations
$\begin{cases} \sq m_{i,1} = \dfrac{\rho}{\sigma x} m_{i,1}, \\
\sq m_{i,2} = \dfrac{\rho}{\sigma x} (m_{i,2} - m_{i,1}), \end{cases}$ for $i = 1,2$,
which we solve using lemma \ref{lem:logarithmiccase} of subsection
\ref{subsubsection:qlog}. More precisely, holomorphic solutions
of $\sq f = \dfrac{\rho}{\sigma x} f$ have the form
$f = \alpha \thq\left(\frac{\sigma}{\rho} x\right)$ and then holomorphic solutions
of $\sq g = \dfrac{\rho}{\sigma x} (g-f)$ have the form
$g = \alpha \frac{\sigma}{\rho} x \thqp\left(\frac{\sigma}{\rho} x\right) +
\beta \thq\left(\frac{\sigma}{\rho} x\right)$. Therefore, matrices $M \in F$
have the form
$$
M = \begin{pmatrix} \alpha_1 \thq\left(\frac{\sigma}{\rho} x\right) &
\alpha_1 \frac{\sigma}{\rho} x \thqp\left(\frac{\sigma}{\rho} x\right) +
\beta_1 \thq\left(\frac{\sigma}{\rho} x\right) \\
\alpha_2 \thq\left(\frac{\sigma}{\rho} x\right) &
\alpha_2 \frac{\sigma}{\rho} x \thqp\left(\frac{\sigma}{\rho} x\right) +
\beta_2 \thq\left(\frac{\sigma}{\rho} x\right) \end{pmatrix} = 
\begin{pmatrix} \alpha_1 & \beta_1 \\ \alpha_2 & \beta_2 \end{pmatrix} T(x),
\text{~where~} T(x) := \begin{pmatrix} \thq\left(\frac{\sigma_1}{\rho} x\right) &
\frac{\sigma}{\rho} x \thqp\left(\frac{\sigma}{\rho} x\right) \\
0 & \thq\left(\frac{\sigma_2}{\rho} x\right) \end{pmatrix}.
$$
Since $\det T$ is non trivial but vanishes at $x_1 = x_2$, we get that
$F = \{C T(x) \tq C \in \GL_2(\C)\}$. The action of $\Gamma,\Delta$
comes here with arbitrary $\Gamma$ so $\F$ is again a singleton. \\

In cases 4,5 and 6 (those truly of interest), the space $E$ and its quotient
$\E$ are more complicated to study directly so our technology comes handy.

\paragraph{Case 4, generic/generic:}
$R = \begin{pmatrix} \rho_1 & 0 \\ 0 & \rho_2 \end{pmatrix}$,
$S = \begin{pmatrix} \sigma_1 & 0 \\ 0 & \sigma_2 \end{pmatrix}$. \\
Matrices $M \in F$ have the form
$$
M = \begin{pmatrix}
\alpha_{1,1} \thq\left(\frac{\sigma_1}{\rho_1} x\right) &
\alpha_{1,2} \thq\left(\frac{\sigma_2}{\rho_1} x\right) \\
\alpha_{2,1} \thq\left(\frac{\sigma_1}{\rho_2} x\right) &
\alpha_{2,2} \thq\left(\frac{\sigma_2}{\rho_2} x\right)
\end{pmatrix}.
$$
Such a matrix is completely determined by the quadruple
$(\alpha_{1,1}, \alpha_{1,2}, \alpha_{2,1}, \alpha_{2,2}) \in \C^4$.
Since the functions
$\thq\left(\frac{\sigma_1}{\rho_1} x\right) \thq\left(\frac{\sigma_2}{\rho_2} x\right)$
and
$\thq\left(\frac{\sigma_1}{\rho_2} x\right) \thq\left(\frac{\sigma_2}{\rho_1} x\right)$
are linearly independent, the condition that $\det M$ does not vanish identically is
equivalent to $(\alpha_{1,1} \alpha_{2,2},\alpha_{1,2} \alpha_{2,1}) \neq (0,0)$.
On the other hand, the gauge freedom on $F$ is expressed by the fact that arbitrary
invertible diagonal matrices $\Diag(\gamma_1,\gamma_2), \Diag(\delta_1,\delta_2)$ act
on $M$, so that:
$$
(\alpha_{1,1}, \alpha_{1,2}, \alpha_{2,1}, \alpha_{2,2}) \sim
\left(\dfrac{\delta_1}{\gamma_1} \alpha_{1,1},
\dfrac{\delta_2}{\gamma_1} \alpha_{1,2},
\dfrac{\delta_1}{\gamma_2} \alpha_{2,1},
\dfrac{\delta_2}{\gamma_2} \alpha_{2,2}\right)
$$
for all $\delta_1,\delta_2,\gamma_1,\gamma_2 \in \C^*$, which in turn implies that:
$$
(\alpha_{1,1} \alpha_{2,2},\alpha_{1,2} \alpha_{2,1}) \sim
\dfrac{\delta_1 \delta_2}{\gamma_1 \gamma_2}
(\alpha_{1,1} \alpha_{2,2},\alpha_{1,2} \alpha_{2,1}).
$$
We thus obtain a well defined map
$M \mapsto \alpha(M) := (\alpha_{1,1} \alpha_{2,2}:\alpha_{1,2} \alpha_{2,1})$
from $F/\sim$ to $\PC$. We shall see that $\alpha$ is ``almost injective''. \\
On the other hand, the condition $\det M(x_1) = 0$ (equivalently $\det M(x_2) = 0$)
says that the image of this map is reduced to a single point. To make this more
precise while legible, we identify $\PC$ with $\C \cup \{\infty\}$ and
$(a_1:a_2) \in \PC$ with $a_1/a_2$. Also we introduce the $q$-elliptic function:
$$
\Phi(x) := \dfrac
{\thq\left(\frac{\sigma_1}{\rho_1} x\right) \thq\left(\frac{\sigma_2}{\rho_2} x\right)}
{\thq\left(\frac{\sigma_1}{\rho_2} x\right) \thq\left(\frac{\sigma_2}{\rho_1} x\right)}
\cdot
$$
Then $\det M(x_1) = 0 \Leftrightarrow \alpha(M) = \dfrac{1}{\Phi(x_1)} \cdot$
The apparent dissymmetry of this condition with respect with $x_1,x_2$ disappears
if one notes that $\Phi$ admits an involution under which $x_1 \leftrightarrow x_2$:
$$
\Phi(x') = \Phi(x'') \Longleftrightarrow
\left(x' x'' \equiv \dfrac{\rho_1 \rho_2}{\sigma_1 \sigma_2}
\text{~~~~or~~~~} x' = x''\right).
$$
In the last step, we want to recover the quadruple $(\alpha_{i,j})$ (up to the
gauge action by $\Gamma,\Delta)$ from the point $\alpha(M) = \dfrac{1}{\Phi(x_1)}$
in $\PC$. A small computation shows that for $\alpha(M) \neq 0,\infty$, that
is for $\alpha(M) \in \C^*$, the preimage is unique; while for $\alpha(M) = 0$
there are three preimages: $(0,1,1,0)$, $(0,1,1,1)$ and $(1,1,1,0)$; and for
$\alpha(M) = \infty$, there are three preimages: $(1,0,0,1)$, $(1,1,0,1)$ and
$(1,0,1,1)$. \\
The conclusion is that the nature of $\F = F/\sim$ depends on the element
$\Phi(x_1) = \Phi(x_2) \in \PC$, which is totally determined by the local data.
If this element is $0$ or $\infty$, the space $\F$ has three elements;
otherwise (thus generically) it is a singleton.

\begin{rmk}
The special cases $\Phi(x_1) = \Phi(x_2) \in \{0,\infty\}$ correspond to those
when one of the coefficients $\alpha_{i,j}$ vanishes, \ie, by corollary
\ref{cor:reducibilitycriterion}, to the case of a reducible system.
\end{rmk}

\paragraph{Case 5, generic/logarithmic:}
$R = \begin{pmatrix} \rho_1 & 0 \\ 0 & \rho_2 \end{pmatrix}$,
$S = \begin{pmatrix} \sigma & \sigma \\ 0 & \sigma \end{pmatrix}$. \\
The functional equation for $M = (m_{i,j})$ is:
$\sq M = \begin{pmatrix}
\frac{\rho_1}{\sigma x} m_{1,1} & \frac{\rho_1}{\sigma x} (m_{1,2} - m_{1,1}) \\
\frac{\rho_2}{\sigma x} m_{2,1} & \frac{\rho_2}{\sigma x} (m_{2,2} - m_{2,1})
\end{pmatrix}$,
so using \ref{subsubsection:qlog} more or less as in case 3, we get:
$$
M = \begin{pmatrix} \alpha_1 \thq\left(\frac{\sigma}{\rho_1} x\right) &
\alpha_1 \frac{\sigma}{\rho_1} x \thqp\left(\frac{\sigma}{\rho_1} x\right) +
\beta_1 \thq\left(\frac{\sigma}{\rho_1} x\right) \\
\alpha_2 \thq\left(\frac{\sigma}{\rho_2} x\right) &
\alpha_2 \frac{\sigma}{\rho_2} x \thqp\left(\frac{\sigma}{\rho_2} x\right) +
\beta_2 \thq\left(\frac{\sigma}{\rho_2} x\right) \end{pmatrix} =
\begin{pmatrix} \alpha_1 \phi_1 & \alpha_1 x \phi'_1 + \beta_1 \phi_1 \\
\alpha_2 \phi_2 & \alpha_2 x \phi'_2 + \beta_2 \phi_2 \end{pmatrix},
$$
where $\phi_i(x) := \thq\left(\frac{\sigma}{\rho_i} x\right)$, $i = 1,2$.
The space of interest (parameterizing $F$) is the space of quadruples
$(\alpha_1,\beta_1,\alpha_2,\beta_2) \in \C^4$. \\
The gauge action is by diagonal matrices $\Gamma = \Diag(\gamma_1,\gamma_2)$
and by unipotent matrices $\Delta = \begin{pmatrix} 1 & \delta \\ 0 & 1 \end{pmatrix}$
(the omitted scalar factor in $\Delta$ can be accounted for in the action of $\Gamma$).
The corresponding allowed transformations for quadruples are best described seing
$L_i := (\alpha_i,\beta_i)$, $i = 1,2$ as lines; and the $\alpha$- and $\beta$-part
respectively as columns $C_1,C_2$. The operations then are dilatations of lines
$L_i \leftarrow \gamma_i L_i$, $\gamma_i \in \C^*$; and transvection
$C_2 \leftarrow C_2 + \delta C_1$. \\
So we find that:
$$
\det M = (\alpha_1 \beta_2 - \alpha_2 \beta_1) \phi_1 \phi_2 +
\alpha_1 \alpha_2 x (\phi_2' \phi_1 - \phi_1' \phi_2) =
x \phi_1 \phi_2 \left(\dfrac{\alpha_1 \beta_2 - \alpha_2 \beta_1}{x} +
\alpha_1 \alpha_2 \left(\frac{\phi_2'}{\phi_2} -\frac{\phi_1'}{\phi_1}\right)\right),
$$
\ie\ $\det M/(x \phi_1 \phi_2)$ is the logarithmic derivative of
$x^{\alpha_1 \beta_2 - \alpha_2 \beta_1} \left(\frac{\phi_2}{\phi_1}\right)^{\alpha_1 \alpha_2}$,
from which we draw that $\det M$ vanishes identically if, and only if
$\alpha_1 \beta_2 - \alpha_2 \beta_1 = \alpha_1 \alpha_2 = 0$. This bad set
within $\C^4$ has three components: $\alpha_1 = \alpha_2 = 0$,
$\alpha_1 = \beta_1 = 0$ and $\alpha_2 = \beta_2 = 0$. Each of these three
components is invariant under the gauge $(\Gamma,\Delta)$ action. \\
We must now check the condition $\det M(x_1) = 0$ while staying within the good
part of $\C^4$; the latter is the union of three disjoint components, each invariant
under the gauge action, and we discuss\footnote{In essence, we are looking for
normal forms.} the corresponding cases. \\
\underline{$\alpha_1 \alpha_2 \neq 0$:} Up to $\Gamma$-action, we may
assume that $\alpha_1 = \alpha_2 = 1$ and up to $\Delta$-action, we may
assume that $\beta_1 = 0$, whence (writing $\beta$ instead of $\beta_2$):
$$
M = \begin{pmatrix} \phi_1 & x \phi'_1 \\ \phi_2 & x \phi'_2 + \beta \phi_2 \end{pmatrix}
\Longrightarrow
\det M(x_1) = \beta (\phi_1 \phi_2)(x_1) + x_1 (\phi_2' \phi_1 - \phi_1' \phi_2)(x_1).
$$
We have three possibilities:
\begin{enumerate}
\item If $(\phi_1 \phi_2)(x_1) \neq 0$ there is a unique $\beta$ such that
$\det M(x_1) = 0$. There is exactly one corresponding class in $\F$.
\item If $\phi_1(x_1) = 0$, then $\phi'_1(x_1) \neq 0$ (because $\thq$ only
has simple zeroes) and $\phi_2(x_1) \neq 0$ (because of non resonancy). There
is no such corresponding class. 
\item If $\phi_2(x_1) = 0$, same conclusion for symmetric reasons.
\end{enumerate}
\underline{Subset $\alpha_1 = 0 \Rightarrow \alpha_2 \beta_1 \neq 0$:}
We may suppose $\alpha_2 = \beta_1 = 1$ and $\beta_2 = 0$, so
$M = \begin{pmatrix} 0 & x \phi_1 \\ \phi_2 & x \phi'_2\end{pmatrix}$.
If $(\phi_1 \phi_2)(x_1) = 0$, this yields one class; otherwise none. \\
\underline{Subset $\alpha_2 = 0 \Rightarrow \alpha_1 \beta_2 \neq 0$:}
Same conclusion for symmetric reasons.. \\

We now summarize the discussion in intelligible form:
\begin{itemize}
\item If $(\phi_1 \phi_2)(x_1) \neq 0$, then $\F$ is a singleton.
\item If $(\phi_1 \phi_2)(x_1) = 0$, that is, if
$\{x_1,x_2\} = \left\{-\frac{\rho_1}{\sigma},-\frac{\rho_2}{\sigma}\right\}$,
then $\F$ has two elements respectively corresponding to the classes
$\alpha_1 = 0$ and $\alpha_2 = 0$. Representatives have been described above.
\end{itemize}

\paragraph{Case 6, logarithmic/logarithmic:}
$R = \begin{pmatrix} \rho & \rho \\ 0 & \rho \end{pmatrix}$,
$S = \begin{pmatrix} \sigma & \sigma \\ 0 & \sigma \end{pmatrix}$. \\
Under our assumptions, $\F$ is isomorphic to $\C$. We just give normal
forms leaving to the interested reader to provide the necessary arguments. 
Let $\phi(x) := \thq\left(\frac{\sigma}{\rho}x\right)$ and $\psi := x \phi'$.
Then each class admits a unique representative of the form
$\begin{pmatrix} - \psi & m_0 + \lambda \phi \\ \phi & \psi \end{pmatrix}$,
$\lambda \in \C$, where $m_0 \in \Rwg$ is a particular solution of the
functional equation $\sq m = \dfrac{\rho}{\sigma x} (m - \phi + 2 \psi)$. 
We have no simple explicit formula for $m$ but it can be proven that for any
$c \in \C^*$ and $g \in \Rwg$ there exists $f \in \Rwg$ such that
$\sq f = \dfrac{c}{x} (f + g)$. One way of obtaining $f$ is to iterate the
operator $f \mapsto - g + \dfrac{x}{c} \sq f$.



\section{The Jimbo-Sakai family (I)}
\label{section:JimboSakaiFamilyI}

The Jimbo-Sakai family studied in \cite{JimboSakai} is a path inside a subspace of
the space:
$$
\left\{A_0 + x A_1 + x^2 A_2 \in \Mat_2(\Rg) \Tq A_0,A_2 \in \GL_2(\C)\right\}.
$$
The subspace is restrained by conditions on the local monodromy data at $0$ and $\infty$
and also by conditions on intermediate singularities. Sakai gave a direct description
of the space of equations $\sq X = \left(A_0 + x A_1 + x^2 A_2\right) X$ as an algebraic
surface (indeed, a rational surface); this is what we consider as the ``left hand side''
of the Riemann-Hilbert-Birkhoff correspondence. In this section, we introduce the
corresponding ``right hand side'', the space of monodromy data. \\

We shall, here and in all this paper, consider the local data at $0$ and $\infty$
(denoted $R$ and $S$) and the intermediate singularities (denoted $\x$) as fixed.


\subsection{Definitions and assumptions}
\label{subsection:Defandassump}

In this section we model the family studied by Jimbo and Sakai in \cite{JimboSakai}:
$n = 2$ and $\mu = 2$, whence $N = 4$. The local data are $R := \Diag(\rho_1,\rho_2)$,
$S := \Diag(\sigma_1,\sigma_2)$ (thus $R$ and $S$ are the exponents at $0$ and $\infty$);
and $\x := \{x_1,x_2,x_3,x_4\}$ (the so-called ``intermediate singularities'', \ie\
those in $\C^*$). We assume \emph{Fuchs relation} and \emph{strong non
resonancy}\footnote{The usual non resonancy condition would only require that
$\dfrac{\rho_1}{\rho_2}$ and $\dfrac{\sigma_1}{\sigma_2}$ do not belong to
$q^{\Z \setminus \{0\}}$, \ie\ equalities $\rho_1 = \rho_2$ or $\sigma_1 = \sigma_2$
would be allowed. Although life is simpler with strong non resonancy, our results
probably extend to the more general case.} in the following form:
\begin{itemize}
\item[(FR)] $x_1 x_2 x_3 x_4 \equiv \dfrac{\rho_1 \rho_2}{\sigma_1 \sigma_2}$;
\item[(NR)] $\dfrac{\rho_1}{\rho_2},\dfrac{\sigma_1}{\sigma_2} \not\in q^\Z$ and
for $1 \leq k < l \leq 4$, $\dfrac{x_k}{x_l} \not\in q^\Z$.
\end{itemize}


\subsubsection{Riemann-Hilbert-Birkhoff correspondence}
\label{subsubsection:RHB}

\paragraph{Left hand side of the RHB correspondence.}
Our primary objects of interest are the space of corresponding systems and its quotient
by the gauge equivalence relation:
$$
E_{R,S,\x} := \left\{A = A_0 + x A_1 + x^2 A_2 \in \Mat_2(\Rg) \Tq
\begin{array}{lr}
A_0 \text{~is conjugate to~} R, \\
A_2 \text{~is conjugate to~} S, \\
\det A(x) \text{~vanishes on~} \x
\end{array}
\right\}, \quad
\E_{R,S,\x} := \dfrac{E_{R,S,\x}}{\sim} \cdot
$$
Note that, for any $A \in E_{R,S,\x}$, one has
\begin{align*}
\det A(x) &= (\det A_2) (x-x_1) (x-x_2) (x-x_3) (x-x_4) \\
&= \sigma_1 \sigma_2 (x-x_1) (x-x_2) (x-x_3) (x-x_4) \\
&= (\det A_0) (1 - x/x_1) (1 - x/x_2) (1 - x/x_3) (1 - x/x_4) \\
&= \rho_1 \rho_2 (1 - x/x_1) (1 - x/x_2) (1 - x/x_3) (1 - x/x_4).
\end{align*}

\paragraph{Right hand side of the RHB correspondence.}
Our secondary objects of interest are dictated by the classification theorems stated
in the previous section, \ie the space of corresponding ``monodromy'' matrices $M$ and
its quotient by the equivalence relation defined there:
$$
F_{R,S,\x} := \left\{M \in \Mat_2(\Rwg) \Tq
\begin{array}{lr}
\sq M = R M (S x^2)^{-1}, \\
\det M \neq 0, \\
\det M \text{~vanishes on~} \x
\end{array}
\right\}, \quad
\F_{R,S,\x} := \dfrac{F_{R,S,\x}}{\sim} \cdot
$$
Note that, for any $M \in F_{R,S,\x}$, one has
$\sq(\det M) = \dfrac{\rho_1 \rho_2}{\sigma_1 \sigma_2} x^{-4} (\det M)$, so that
$$
\det M = C \, \thq(-x/x_1) \thq(-x/x_2) \thq(-x/x_3) \thq(-x/x_4)
\text{~for some~} C \in \C^*,
$$
and therefore it has simple zeroes at $[\x;q]$ and nowhere else, entailing that $M^{-1}$
has simple poles at $[\x;q]$ and nowhere else. \\

Recall the equivalence relation on $F_{R,S,\x}$ (here $\Ddc$ denotes the group of
invertible diagonal $2 \times 2$ matrices):
$$
\forall M,N \in F_{R,S,\x} \;,\; M \sim N \underset{def}{\Longleftrightarrow}
\exists \Gamma, \Delta \in \Ddc \;:\; N = \Gamma^{-1} M \Delta.
$$
We shall write for short $F := F_{R,S,\x}$ and $\F := \F_{R,S,\x}$. In section
\ref{section:A Birkhoff type classification theorem}, we have described a natural
bijective Riemann-Hilbert correspondence $\E \overset{\sim}{\longleftrightarrow} \F$.
In the rest of this section (from \ref{subsection:Falgsurf} on) we attempt at a
description of $\F$ as an \emph{algebraic surface}, at least under some generically
satisfied assumptions. More complete results will come in section
\ref{section:JimboSakaiFamilyII}.


\subsubsection{Reducibility in the Jimbo-Sakai family}
\label{subsubsection:reducibilityJSF}

First we discuss reducibility in the Jimbo-Sakai family. This will be helpful
later in order to understand many exceptional cases, singularities, etc. 
Applying corollary \ref{cor:reducibilitycriterion} of theorem
\ref{thm:reducibilitycriterion}, we see that matrices
$M = (m_{i,j})_{i,j = 1,2} \in F$ corresponding to reducible systems in $E$ are
those such that $m_{i,j} = 0$ for some $i,j \in \{1,2\}$. Note that this property
is invariant under the \emph{left} action by diagonal matrices
$M \mapsto \Gamma M \Delta^{-1}$ (from now and for simplicity, on this replaces
the previous \emph{right} action $M \mapsto \Gamma^{-1} M \Delta$ without any
unwanted logical consequence).

\begin{defn}
\label{defn:reduciblematrices}
We say that $M \in F$ is \emph{reducible} if it corresponds to a reducible
system. (Since the local data $R,S$ are fixed, this makes sense.) We then
say that the class of $M$ in $\F$ is reducible. (By the invariance stated
above, this also makes sense.)
\end{defn}

We discuss the case that $m_{1,1} = 0$, the other ones being entirely similar.
If $m_{1,1} = 0$, then $\det M = - m_{1,2} m_{2,1} \neq 0$ and it must vanish
over $\x$. But each of $m_{1,2},m_{2,1}$ has \emph{a priori} two $q$-spirals
of simple zeroes, the union of all four of them being the whole of $[\x;q]$.
This implies that $m_{1,2} = c x^r \thq(-x/x_k) \thq(-x/x_l)$ and
$m_{2,1} = d x^s \thq(-x/x_m) \thq(-x/x_n)$ for some $c,d \in \C^*$, some
$r,s \in \Z$ and some ``splitting'' $\{1,2,3,4\} = \{k,l\} \sqcup \{m,n\}$.
Since $\sq m_{i,j} = (\rho_i/\sigma_j) m_{i,j}$, this in turn implies that
$x_k x_l \equiv \rho_1/\sigma_2$ and $x_m x_n \equiv \rho_2/\sigma_1$. \\

Conversely, if these congruences are satisfied, we can produce matrices
$M \in F$ with $m_{1,1} = 0$, $m_{1,2} = c x^r \thq(-x/x_k) \thq(-x/x_l)$,
$m_{2,1} = d x^s \thq(-x/x_m) \thq(-x/x_n)$ and an arbitrary $m_{2,2}$ such
that $\sq m_{2,2} = (\rho_2/\sigma_2) m_{2,2}$. Alternatively, we can take
$m_{2,2} = 0$ and an arbitrary $m_{1,1}$ such that
$\sq m_{1,1} = (\rho_1/\sigma_1) m_{1,1}$.

\begin{defn}
\label{defn:splittingofFR}
We say that there is a \emph{splitting} of Fuchs relation (FR) if there is
a permutation $(i,j,k,l)$ of $(1,2,3,4)$ and a permutation $(m,n)$ of $(1,2)$
such that:
$$
x_i x_j \equiv \dfrac{\rho_m}{\sigma_n} \text{~and~}
x_k x_l \equiv \dfrac{\rho_n}{\sigma_m}.
$$
Note that, because of (FR), these two congruence relations are actually equivalent.
\end{defn}

The general conclusion is as follows:

\begin{thm}
\label{theorem:reducibilityJSF}
The space $E$ contains reducible systems (equivalently, the space $F$ contains
reducible matrices) if, and only if there exist a ``splitting'' of Fuchs relation.
\end{thm}
\Finprcourt


\subsection{$\F$ as an algebraic surface: heuristics}
\label{subsection:Falgsurf}

We begin by some heuristic considerations. Let
$M := (m_{i,j})_{i,j = 1,2} \in \Mat_2(\Rwg)$. Then:
$$
M \in F \Longleftrightarrow \begin{cases}
\sq m_{i,j} = \dfrac{\rho_i}{\sigma_j} x^{-2} m_{i,j}, i,j = 1,2, \\
m_{1,1} m_{2,2} \neq m_{1,2} m_{2,1}, \\
m_{1,1} m_{2,2} - m_{1,2} m_{2,1} \text{~vanishes on~} \x.
\end{cases}
$$
The space $\hat{F}$ of all matrices $M \in \Mat_2(\Rwg)$ such that
$\sq M = R M (S x^2)^{-1}$ is the product of the four spaces defined by
$\sq m_{i,j} = \dfrac{\rho_i}{\sigma_j} x^{-2} m_{i,j}$, $i,j = 1,2$ thus it
is a linear space of dimension $8$ (see section \ref{subsection:someequations}). \\

Condition $m_{1,1} m_{2,2} \neq m_{1,2} m_{2,1}$ defines a dense Zariski open subset
of $\hat{F}$. The four conditions $(m_{1,1} m_{2,2} - m_{1,2} m_{2,1})(x_i) = 0$
actually represent \emph{three} independent conditions (because of (FR)), so $F$
has dimension $5$. (Beware however that $F$ is \emph{not} dense or Zariski-dense in
$\hat{F}$.) \\

Now, if $M := (m_{i,j}), N := (n_{i,j}) \in \Mat_2(\Rwg)$:
\begin{align*}
M \sim N & \Longleftrightarrow
n_{i,j} = \dfrac{\delta_j}{\gamma_i} m_{i,j}, i,j = 1,2, 
\text{~for some~} \gamma_1, \gamma_2, \delta_1, \delta_2 \in \C^* \\
& \Longleftrightarrow
n_{i,j} = \lambda_{i,j} m_{i,j}, i,j = 1,2, 
\text{~for some~} \lambda_{i,j} \in \C^* \text{~such that~}
\lambda_{1,1} \lambda_{2,2} = \lambda_{1,2} \lambda_{2,1}.
\end{align*}
Thus $\F$ is actually the quotient of $F$ by the free action of a $3$ dimensional torus:
therefore $\F$ has dimension $2$, it is a surface. \\

What follows (in this section) rests on the observation that the above conditions
for $M \in F$ can be nicely expressed in terms of the diagonal and antidiagonal
products $m_{1,1} m_{2,2}$ and $m_{1,2} m_{2,1}$, and similarly for for $M \sim N$
where $N := (n_{i,j})$. So, for $M,N$ as above in $\Mat_2(\Rwg)$, we write:
$$
f_1 := m_{1,1} m_{2,2}, \quad f_2 := m_{1,2} m_{2,1} \quad
g_1 := n_{1,1} n_{2,2}, \quad g_2 := n_{1,2} n_{2,1}.
$$
Then we observe:
\begin{itemize}
\item If $M \in \hat{F}$, then $f_1,f_2$ are solutions of the same equation
$\sq f = \dfrac{\rho_1 \rho_2}{\sigma_1 \sigma_2} x^{-4} f$, which defines in $\Rwg$
a linear space $W$ of dimension $4$. This space will be studied in subsection
\ref{subsection:Spacesoffunctions again}.
\item Let $M \in \hat{F}$; then $M \in F$ if and only if$f_1 \neq f_2$ and
$f_1(x_i) = f_2(x_i)$, $i = 1,\ldots,4$; actually, under (FR), three of these four
conditions are enough.
\item If $M \sim N$, then $(g_1,g_2) = \lambda (f_1,f_2)$ for some
$\lambda \in \C^*$; actually, with the above notations,
$\lambda = \lambda_{1,1} \lambda_{2,2} = \lambda_{1,2} \lambda_{2,1}$.
\end{itemize}

So there is a natural map $M \mapsto (f_1,f_2)$ from $\hat{F}$ to $W \times W$
which goes to the quotient as:
\begin{align*}
\F & \rightarrow \dfrac{W \times W}{\C^*}, \\
\text{class of~} M & \mapsto \text{class of~} (m_{1,1} m_{2,2},m_{1,2} m_{2,1}).
\end{align*}

In the following subsections we shall find out what exactly the image $\G$ of
this map is and how far it is from injective. Then we shall describe $\G$ in
algebro-geometrical terms. Here, as in section \ref{section:JimboSakaiFamilyII},
the idea is that we have explicit finite dimensional spaces of functions and
that our (seemingly) transcendental conditions can be expressed in the language
of linear and multilinear algebra. In order to do so, we need again to discuss
some spaces of solutions of elementary $q$-difference equations.


\subsection{Spaces of functions again}
\label{subsection:Spacesoffunctions again}

For $k \in \N^*$ and $a \in \C^*$, let:
$$
V_{k,a} := \{f \in \Rwg \tq \sq f = a x^{-k} f\},
$$
a $\C$-linear space of dimension $k$ (see section \ref{subsection:someequations}).
An explicit basis is, for instance, the family of all $\thq(x/\alpha)^k$ where
$\alpha^k = a$. All elements of $V_{k,a}$ have the form
$\lambda \thq(x/\alpha_1) \cdots \thq(x/\alpha_k)$ where $\lambda \in \C$ and
$\alpha_1 \cdots \alpha_k = a$. We write $V_{k,a}^* := V_{k,a} \setminus \{0\}$. \\

Let $k,l \in \N^*$ and $a,b \in \C^*$. There is a natural map
\begin{align*}
V_{k,a} \times V_{l,b} &\rightarrow V_{k+l,ab}, \\
(f,g) &\mapsto fg.
\end{align*}
We study it in case $k = l = 2$, writing it:
$$
p_{a,b}: V_{2,a} \times V_{2,b} \to V_{4,ab}.
$$

\begin{prop}
\label{prop:XT-YZ,image}
(i) Let $a \not\equiv b$. Then the image of $p_{a,b}$ in $V_{4,ab}$ is a homogeneous
quadric hypersurface, its equation is $XT - YZ = 0$ in some coordinate system. \\
(ii) Let $b = q^k a$, $k \in \Z$. Then the image is a hyperplane.
\end{prop}
\Pr
(i) Let $\alpha,\beta \in \C^*$ such that $\alpha^2 = a$, $\beta^2 = b$, so that
$(u_1,u_2) := (\thq(x/\alpha)^2,\thq(-x/\alpha)^2)$ is a basis of $V_{2,a}$ and
$(v_1,v_2) := (\thq(x/\beta)^2,\thq(-x/\beta)^2)$ is a basis of $V_{2,b}$. Then
$(u_1 v_1,u_1 v_2,u_2 v_1,u_2 v_2)$ is a basis of $V_{4,ab}$ (this can be checked
either by arguing on the zeroes or by using the following proposition, which is
independent). Written in this basis, the image of $V_{2,a} \times V_{2,b}$ in $V_{4,ab}$
is:
$$
\{(x_1 y_1,x_1 y_2,x_2 y_1,x_2 y_2) \tq x_1,x_2,y_1,y_2 \in \C\} =
\{(X,Y,Z,T) \in \C^4 \tq XT - YZ = 0\}.
$$
(ii) Let $(u_1,u_2)$ any basis of $V_{2,a}$, so that $(x^k u_1,x^k u_2)$ is a basis
of $V_{2,b}$. Then $(x^k u_1^2, x^k u_1 u_2,x^k u_2^2)$ is a free system in $V_{4,ab}$.
Complete it into a basis. Written in this basis, the image of $V_{2,a} \times V_{2,b}$
in $V_{4,ab}$ is:
is:
$$
\{(x_1 y_1,x_1 y_2 + x_2 y_1,x_2 y_2,0) \tq x_1,x_2,y_1,y_2 \in \C\} = \C^3 \times \{0\}.
$$
\Finprcourt

\begin{prop}
\label{prop:XT-YZ,fibers}
Let $f \in V_{2,a}^*$, $g \in V_{2,b}^*$ so that $fg \in V_{4,ab}^*$. \\
If $a \not\equiv b$, the preimage of $fg$ is:
$p_{a,b}^{-1}(fg) = \{(\lambda^{-1} f,\lambda g) \tq \lambda \in \C^*\}$. \\
If $b = q^k a$, $k \in \Z$, the preimage is:
$p_{a,b}^{-1}(fg) =
\{(\lambda^{-1} f,\lambda g) \tq \lambda \in \C^*\} \cup
\{(\lambda^{-1} x^{-k} g,\lambda x^k f) \tq \lambda \in \C^*\}$.
\end{prop}
\Pr
Let $f,g$ as in the statement and $f_1 \in V_{2,a}$, $g_1 \in V_{2,b}$ such that
$f_1 g_1 = f g$. Then $f_1,g_1 \neq 0$. We have an equality of elliptic functions:
$\dfrac{f_1}{f} = \dfrac{g}{g_1} \cdot$ If there is cancellation of zeroes, our
elliptic function is constant (order $1$ is impossible), which yields the first
part of the preimage $\{(\lambda^{-1} f,\lambda g) \tq \lambda \in \C^*\}$. \\
If there is no cancellation of zeroes, these are elliptic functions with order $2$
and they have the same divisor of zeroes. So $h := g/f_1$ satisfies the relation
$\sq h = (b/a) h$ and it has neither zeroes nor poles on $\C^*$. Since $h \neq 0$,
this is only possible if $b/a = q^k$ for some $k \in \N^*$ and $h = \lambda x^k$ for
some $\lambda \in \C^*$, whence the conclusion.
\Finprcourt

\begin{rmk}
We see that in all cases there is a $\C^*$-action on $V_{2,a}^* \times V_{2,b}^*$,
given by $\left(\lambda,(f,g)\right) \mapsto (\lambda^{-1} f,\lambda g)$. In case
that $b = q^k a$, there is moreover an involution $(f,g) \mapsto (x^{-k} g, x^k f)$.
The latter does not commute to the former (the involution conjugates $\lambda$ with
$\lambda^{-1}$) and we eventually get the action of the corresponding semi-direct product
$\C^* \ltimes (\Z/2\Z)$.
\end{rmk}

\begin{cor}
\label{cor:planesinquadrics}
We suppose that $a \not\equiv b$. Then the planes included in the
quadric hypersurface $\Sigma := \Im\ p_{a,b} \subset V_{4,ab}$ are either of the
form $f V_{2,b}$ for some $f \in V_{2,a} \setminus \{0\}$, or of the form $g V_{2,a}$
for some $g \in V_{2,b} \setminus \{0\}$.
\end{cor}
\Pr
It is clear that such sets $f V_{2,b}$, $g V_{2,a}$ are indeed planes included
in $\Sigma$. \\
Conversely, let $f_1,f_2 \in V_{2,a} \setminus \{0\}$ and
$g_1,g_2 \in V_{2,b} \setminus \{0\}$ be such that
$\Vect(f_1 g_1,f_2 g_2) \subset \Sigma$. We shall prove that either $f_1,f_2$
are colinear, or $g_1,g_2$ are; this will entail our conclusion. So assume for
instance that $f_1,f_2$ are not colinear. By assumption:
\begin{align*}
f_1 g_1 + f_2 g_2 &= fg \qquad \qquad \qquad \quad
\text{for some~} f \in V_{2,a}, g \in V_{2,b} \\
&= (\alpha_1 f_1 + \alpha_2 f_2) g \qquad \text{for some~} \alpha_1,\alpha_2 \in \C,
\end{align*}
whence $f_1(g_1 - \alpha_1 g) = f_2 (\alpha_2 g - g_2)$, so that, by proposition
\ref{prop:XT-YZ,fibers} (and the assumption that $f_1,f_2$ are not colinear),
$g_1 = \alpha_1 g$ and $g_2 = \alpha_2 g$, \ie\ $g_1,g_2$ are colinear.
\Finprcourt

\begin{cor}
\label{cor:noplaneinintersectionofquadrics}
Let $c = a_1 b_1 = a_2 b_2$, so that $\Sigma_1 := \Im\ p_{a_1,b_1} \subset V_{4,c}$
and $\Sigma_2 := \Im\ p_{a_2,b_2} \subset V_{4,c}$. In the generic situation that
the two decompositions of $c$ are ``essentially inequivalent'', \ie\ none of
$a_1/a_2$, $a_1/b_2$, $b_1/a_2$, $b_1/b_2$ belongs to $q^\Z$, the intersection
$\Sigma_1 \cap \Sigma_2$ cannot contain a plane.
\end{cor}
\Pr
Indeed, by the previous corollary, this would imply the existence of a non
trivial common factor, \ie\ a non zero element in one of the spaces
$V_{2,a_1} \cap x^k V_{2,a_2}$, $V_{2,a_1} \cap x^k V_{2,b_2}$, $V_{2,b_1} \cap x^k V_{2,a_2}$, 
$V_{2,b_1} \cap x^k V_{2,b_2}$, for some $k \in \Z$.
\Finpr

Our intended application is with
$$
(a_1,b_1) := (\rho_1/\sigma_1,\rho_2/\sigma_2), \quad \text{resp.} \quad
(a_2,b_2) := (\rho_1/\sigma_2,\rho_2/\sigma_1),
$$
thus defining two quadric hypersurfaces
$$
\Sigma_1 := \Im\ p_{a_1,b_1} \subset W
\quad \text{and} \quad \Sigma_2 := \Im\ p_{a_2,b_2} \subset W
$$
of the \emph{same} linear space
$$
W := V_{4,c}, \quad \text{where} \quad
c := (\rho_1 \rho_2)/(\sigma_1 \sigma_2) = a_1 b_1 = a_2 b_2.
$$
We shall also write
$$
W_{i,j} := V_{2,\rho_i/\sigma_j}.
$$

By the non resonancy assumption (NR) we are then in the generic situation mentioned in
the above corollary (the two decompositions of $c$ are ``essentially inequivalent'').
As a consequence, the intersection surface can be explicitly described; we leave it to
the reader to check that:

\begin{prop}
Recall that $(a_1,b_1) := (\rho_1/\sigma_1,\rho_2/\sigma_2)$ and
$(a_2,b_2) := (\rho_1/\sigma_2,\rho_2/\sigma_1)$. Then:
$$
\Sigma_1 \cap \Sigma_2 =
\left\{\lambda \thq(x/\alpha) \thq(x/\alpha') \thq(x/\beta) \thq(x/\beta') \Tq
\lambda \in \C, \alpha, \alpha', \beta, \beta' \in \C^*,
\begin{array}{lr}
  \alpha \alpha' = a_1, \\
  \beta \beta' = b_1, \\
  \alpha \beta = a_2, \\
  \alpha' \beta' = b_2
\end{array}
\right\}.
$$
It can be parameterized by $\lambda,\alpha$ by taking:
$$
\alpha' := a_1/\alpha, \quad \beta := a_2/\alpha, \quad \text{and~}
\beta' := b_1 \alpha/a_2 = b_2 \alpha/a_1 = b_1 b_2 \alpha/c = c \alpha/(a_1 a_2).
$$
\end{prop}
\Finpr

In a second step, we shall consider the \emph{projective} quadric surfaces,
the respective images $\widetilde{\Sigma}_1$ of $\Sigma_1$ and $\widetilde{\Sigma}_2$
of $\Sigma_2$ in $\P(W)$.


\subsubsection{Incidence relations of planes of $\Sigma_1,\Sigma_2$}
\label{subsubsection:IncidenceinW}

The motivation for the following study was explained in \ref{subsection:contents}.
Recall that $(a_1,b_1) := (\rho_1/\sigma_1,\rho_2/\sigma_2)$ and
$(a_2,b_2) := (\rho_1/\sigma_2,\rho_2/\sigma_1)$. \\

First we note the following facts:
\begin{itemize}
\item Two distinct planes in $W$ intersect along a line if, and only if, they are
contained in a common hyperplane (equivalently: their sum is a hyperplane). This is
because the ambient space $W$ has dimension $4$. In this situation, we shall say that
the two planes are \emph{incident}.
\item A line of $W$ has the form 
$D_{\alpha,\beta,\gamma,\delta} :=
\C \thq(x/\alpha) \thq(x/\beta) \thq(x/\gamma) \thq(x/\delta)$ for some fixed
$\alpha,\beta,\gamma,\delta \in \C^*$ such that $\alpha \beta \gamma \delta = c$.
\item Write $H_x$ the hyperplane $\{f \in W \tq f(x) = 0\}$ (of course, $H_{qx} = H_x$).
The only such hyperplanes containing $D_{\alpha,\beta,\gamma,\delta}$ are $H_\alpha$,
$H_\beta$, $H_\gamma$ and $H_\delta$.
\end{itemize}

Combining these facts with corollary \ref{cor:planesinquadrics}, we get the following:

\begin{prop}
All pairs of incident planes $P_1 \subset \Sigma_1$, $P_2 \subset \Sigma_2$ are obtained
by the following process. Let $\alpha,\beta,\gamma,\delta \in \C^*$ such that:
$$
\alpha \beta = a_1, \quad \alpha \gamma = a_2, \quad
\gamma \delta = b_1, \quad \beta \delta = b_2.
$$
(These conditions of course imply $\alpha \beta \gamma \delta = c$.)
Then one of the following four possibilities holds:
\begin{itemize}
\item $P_1 = \thq(x/\alpha) \thq(x/\beta) W_{2,2}$ and
$P_2 = \thq(x/\alpha) \thq(x/\gamma) W_{2,1}$; then $P_1 + P_2 = H_\alpha$;
\item $P_1 = \thq(x/\alpha) \thq(x/\beta) W_{2,2}$ and
$P_2 = \thq(x/\beta) \thq(x/\delta) W_{1,2}$; then $P_1 + P_2 = H_\beta$;
\item $P_1 = \thq(x/\gamma) \thq(x/\delta) W_{1,1}$ and
$P_2 = \thq(x/\alpha) \thq(x/\gamma) W_{2,1}$; then $P_1 + P_2 = H_\gamma$;
\item $P_1 = \thq(x/\gamma) \thq(x/\delta) W_{1,1}$ and
$P_2 = \thq(x/\beta) \thq(x/\delta) W_{1,2}$; then $P_1 + P_2 = H_\delta$.
\end{itemize}
In all cases, $P_1 \cap P_2 = D_{\alpha,\beta,\gamma,\delta}$.
\end{prop}
\Finpr

We see that those hyperplanes of $W$ which cut $\Sigma_1$ and $\Sigma_2$ along
planes are of a very special nature: they have the form $H_x$ defined above.
Among their incidence properties, we see that through a generic line in $W$ pass
four such hyperplanes and each of them gives rise to two pairs of planes (one in
each of $\Sigma_1$, $\Sigma_2$) having a common intersection line; we say a bit more
about that at the end of \ref{subsubsection:IncidenceinPW} herebelow; and we return
to them in \ref{subsubsection:TheHyperplanesHx}.


\subsubsection{Incidence relations of lines of $\widetilde{\Sigma}_1,\widetilde{\Sigma}_2$}
\label{subsubsection:IncidenceinPW}

We shall write $\tilde{D}_{\alpha,\beta,\gamma,\delta} \in \P(W)$ the image of line
$D_{\alpha,\beta,\gamma,\delta}$ (thus a point), $\tilde{H}_x \subset \P(W)$ the image
of hyperplane $H_x$ (thus a projective plane) and $\tilde{P}_i \subset \P(W)$ the image
of plane $P_i$ (thus a projective line).

\begin{cor}
All pairs of intersecting lines $\Delta_1 \subset \widetilde{\Sigma}_1$, 
$\Delta_2 \subset \widetilde{\Sigma}_2$ are obtained by taking
$\Delta_i := \tilde{P}_i$, $i = 1,2$, with the above construction. \\
For a given $\alpha \in \C^*$, parameters $\beta,\gamma,\delta \in \C^*$ are uniquely
determined. The process defines four lines (two in each of
$\widetilde{\Sigma}_1,\widetilde{\Sigma}_2$). All four lines meet at the point
$\tilde{D}_{\alpha,\beta,\gamma,\delta}$.
\end{cor}
\Finpr

From what we said hereabove at the end of \ref{subsubsection:IncidenceinW},
we draw that through a generic point of $\P(W)$ pass four planes that cut
$\widetilde{\Sigma}_1$ and $\widetilde{\Sigma}_2$ along two lines each, all the four
lines corresponding to one such special plane having one common point. In the
application to our problem detailed in \ref{subsection:coveringquadric} and
\ref{rough4surf}, we consider the particular point $\tilde{L} \in \P(W)$ image of
the line
$$
L := \Vect\left(\thq(-x/x_1)\thq(-x/x_2)\thq(-x/x_3)\thq(-x/x_4)\right) \subset W,
$$
\ie\ the intersection of the four hyperplanes $H_{x_i}$. It seems likely that
the sixteen corresponding lines in $\widetilde{\Sigma}_1$, $\widetilde{\Sigma}_2$ are
related somehow to the ``special fibers'' we shall encounter in our analysis
in section \ref{section:JimboSakaiFamilyII}.


\subsubsection{The family of hyperplanes $(H_x)_{x \in \C^*}$}
\label{subsubsection:TheHyperplanesHx}

As noted at the end of \ref{subsubsection:IncidenceinW}, the hyperplanes $H_x$ of $W$,
$x \in \C^*$, seem to play a special role. They have interesting properties which we feel
to be relevant in understanding the geometry of $\F$. Although we have not been able
to exploit them completely, we summarize here some of these properties. \\

Recall $W$ is the set of holomorphic solutions of some equation $\sq f = c x^{-4} f$.
Let $W'$ the dual of the space $W$. Each $H_x$ can be seen as a point of $\P(W')$, which
is isomorphic to $\P^3(\C)$. For $f \in W$, the functional equation $\sq f = c x^{-4} f$
implies that $f(qx) = 0 \Leftrightarrow f(x) = 0$, \ie\ $H_{qx} = H_x$. Therefore we get
an holomorphic map $\overline{x} \mapsto H_x$ from $\Eq$ to $\P(W')$. That map is injective.
Indeed, if $x \not\equiv y$ one can easily build $f \in W$ such that $f(x) = 0 \neq f(y)$,
so that $H_x \neq H_y$. It follows that the image $\{H_x \tq x \in \C^*\}$ is a holomorphic
curve isomorphic to $\Eq$, hence an elliptic curve. We write it $\mathbf{E}_{q}'$. \\

To understand $\mathbf{E}_{q}'$ as an embedded algebraic curve, note that any three
distinct points on it generate a (projective) plane which cuts the curve at a fourth
point (the latter maybe not distinct from the three former ones). Let indeed
$a_1,a_2,a_3 \in \C^*$ pairwise non congruent defining $H_{a_1}$, $H_{a_2}$, $H_{a_3}$
our three distinct points on $\mathbf{E}_{q}'$. Let $x_4 \in \C^*$ be chosen such
that $a_1 a_2 a_3 a_4 \equiv c$. Then $H_{a_1} \cap H_{a_2} \cap H_{a_3} \cap H_{a_4}$
is the line generated by $\thq(-x/a_1)\thq(-x/a_2)\thq(-x/a_3)\thq(-x/a_4)$; while,
for any other choice of $a_4$ the intersection would be $\{0\}$. \\

To find equations for $\mathbf{E}_{q}'$ in $\P(W')$, we proceed as follows. Recall
that the affine algebra of $W'$ (the algebra of polynomials on $W'$) is the symmetric
algebra on its dual, \ie\ $\Sym(W)$. We look for a homogeneous ideal in $\Sym(W)$.
Selecting a basis $(f_1,f_2,f_3,f_4)$ of $W$ gives an identification of $\Sym(W)$
with $\C[X_1,X_2,X_3,X_4]$. We map the homogeneous component $\C[X_1,X_2,X_3,X_4]_d$
into $W^{(d)} := \{f \in \Rwg \tq \sq f = c^d x^{-4d} f\}$ by sending $F(X_1,X_2,X_3,X_4)$
to $F(f_1,f_2,f_3,f_4)$. Then $F$ is an equation of $\mathbf{E}_{q}'$ if, and only if,
$F(f_1,f_2,f_3,f_4)$ vanishes at all $x \in \C^*$, \ie\ if $F$ is in the kernel of
$\C[X_1,X_2,X_3,X_4]_d \rightarrow W^{(d)}$. The source of this map has dimension
$(d+1)(d+2)(d+3)/6$, while its target has dimension $4 d$. Therefore, for $d \geq 2$,
the kernel is non trivial. For instance, taking $d = 2$, we find two independent
quadratic forms vanishing on $\mathbf{E}_{q}'$. Therefore this curve is, at any rate,
a component of the intersection of two quadric surfaces. We do not know if it is
the only component. \\

One way to find the quadratic forms (\ie\ the degree $2$ component of the ideal
of $\mathbf{E}_{q}'$) is to choose pairwise non congruent $a_1,a_2,a_3,a_4 \in \C^*$
such that $a_1 a_2 a_3 a_4 \not\equiv c$; and then $f_1,f_2,f_3,f_4 \in W$ such that
$f_i(a_j) = \delta_{i,j}$ (easy using theta functions). Then we require that
$\sum\limits_{1 \leq i \leq j \leq 4} c_{i,j} f_i f_j$ vanish at all $a_i$ at order $2$
(in $W^{(2)}$ this implies that it is $0$). Evaluation at $a_i$ yields $c_{i,i} = 0$.
Then we are left with four linear conditions (vanishing of the derivatives) on the
six coordinates $c_{i,j}$, $1 \leq i < j \leq 4$.


\subsection{$\F$ as a degree $2$ covering of a quadric surface}
\label{subsection:coveringquadric}

Let $W_1 := V_{2,\frac{\rho_1}{\sigma_1}} \times V_{2,\frac{\rho_2}{\sigma_2}}$ and
$W_2 := V_{2,\frac{\rho_1}{\sigma_2}} \times V_{2,\frac{\rho_2}{\sigma_1}}$. As seen in
\ref{subsection:Falgsurf}, $\hat{F}$ is identified with $W_1 \times W_2$ by
$(m_{i,j}) \mapsto ((m_{1,1},m_{2,2}),(m_{1,2},m_{2,1}))$. \\

Let $W := V_{4,\frac{\rho_1 \rho_2}{\sigma_1 \sigma_2}}$. We have defined in
\ref{subsection:Spacesoffunctions again} multiplication maps $W_1 \rightarrow W$
and $W_2 \rightarrow W$. Composing, we get a map:
$$
\begin{cases}
\hat{F} \rightarrow W \times W, \\
M = (m_{i,j}) \mapsto (m_{1,1} m_{2,2},m_{1,2} m_{2,1}),
\end{cases}
$$
the image of which is (by proposition \ref{prop:XT-YZ,image}) $\Sigma_1 \times \Sigma_2$,
where $\Sigma_i$, image of $W_i \rightarrow W$ for $i = 1,2$, is a homogeneous
quadric hypersurface of the four dimensional $\C$-linear space $W$. Thus we have
a surjective mapping:
$$
\begin{cases}
\hat{F} \rightarrow \Sigma_1 \times \Sigma_2, \\
(m_{i,j}) \mapsto (f_1,f_2) := (m_{1,1}m_{2,2},m_{1,2}m_{2,1}).
\end{cases}
$$
The condition: $\det M \neq 0$ (on elements of $F$) translates to: $f_1 \neq f_2$.
The condition: $\det M$ vanishes on $\x$ translates to: $f_1 - f_2$ belongs to the
intersection of the four hyperplanes $\Ker(f \mapsto f(x_i))$. However, because of
$(FR)$, these are really three independent linear conditions and that intersection is
the line
$$
L := \Vect\left(\thq(-x/x_1)\thq(-x/x_2)\thq(-x/x_3)\thq(-x/x_4)\right) \subset W.
$$
We write $L^* := L \setminus \{0\}$ and deduce a surjective map:
$$
F \rightarrow G :=
\{(f_1,f_2) \in \Sigma_1 \times \Sigma_2 \tq f_1 - f_2 \in L^*\}.
$$
The torus action on $F$ (\ie\ the diagonal action of $\C^*$) corresponds in $G$ to
the obvious $\C^*$-action 
$\left(\lambda,(f_1,f_2)\right) \mapsto (\lambda f_1,\lambda f_2)$. Thus
we eventually get a surjective map:
$$
\F \rightarrow \G := \dfrac{G}{\C^*} \cdot
$$

We shall now formulate assumptions on the local data $R$, $S$ and $\x$ under
which this map is bijective. Note however that these assumptions are generically
satisfied\footnote{And that moreover we shall relax them in section
\ref{section:JimboSakaiFamilyII} where a quite different approach is followed.}.
There are actually two causes of non injectivity. One of them comes from
the second case of proposition \ref{prop:XT-YZ,fibers}; the other comes
from the fact that proposition \ref{prop:XT-YZ,fibers} adresses only the case
of non zero functions. \\

The former cause is taken care by the assumption (nicknamed ``special condition''):
$$
\text{(SC)} \quad 
\dfrac{\rho_1}{\sigma_1} \not\equiv \dfrac{\rho_2}{\sigma_2}
\quad \text{~and~} \quad 
\dfrac{\rho_2}{\sigma_1} \not\equiv \dfrac{\rho_1}{\sigma_2}.
$$
This assumption guarantees that the first part of proposition \ref{prop:XT-YZ,fibers}
can be applied. \\

As for the latter cause, let $f = m_{1,1} m_{2,2} \in W$ (the case of $m_{1,2} m_{2,1}$
is similar) and suppose that also $f = m'_{1,1} m'_{2,2}$ (with obvious notations).
If $f \neq 0$, we know from \ref{prop:XT-YZ,fibers} that, under assumption (SC), 
$(m'_{1,1},m'_{2,2}) = (c m_{1,1},c^{-1} m_{2,2})$ for some $c \in \C^*$. From this,
one easily deduces that, if $\prod m_{i,j} \neq 0$, then the class of $M = (m_{i,j})$
in $\F$ is alone in the preimage of its image. So defects of injectivity are only
possible if there are matrices such that at least one $m_{i,j} = 0$. But this is
the case of reducibility discussed in theorem \ref{theorem:reducibilityJSF}. Now
recall the terminology introduced in definition \ref{defn:splittingofFR}.

\begin{prop}
If (SC) holds and there is no splitting of (FR), the map $\F \rightarrow \G$
is bijective.
\end{prop}
\Pr
Immediate from the above discussion and theorem \ref{theorem:reducibilityJSF}.
\Finpr

Now we give a preliminary crude description of $\F$ and $\G$. This will be refined
in the next section \ref{rough4surf}. \\

For that, we project $G$ to $\Sigma_1$ by $(f_1,f_2) \mapsto f_1$. The preimage
of $f_1$ is $\{f_1\} \times \left((f_1 + L^*) \cap \Sigma_2\right)$. Assume that
$f_1 \not\in \Sigma_1 \cap \Sigma_2$. Then the punctured affine line $f_1 + L^*$
meets $\Sigma_2$ at two points, or at a double point in case of tangency (that
is, if $f_1$ is in the critical locus of the projection). We dare not call theorem
the following, because of the difficulty of stating precise assumptions; but it
is plainly true ``in general''. Note that, as a (non closed) algebraic surface,
$\G$ is endowed with the Zariski topology inherited from $\P(W^2)$.

\begin{fact}
Under the projection $(f_1,f_2) \mapsto f_1$, an open dense subset of the space
$\G$ is a degree $2$ ramified covering of a quasi-projective quadric surface
(an open dense subset of the image of $\Sigma_1$ in $\P(W)$).
\end{fact}

If we try to make more precise where the projection fails to be a covering,
we are led to define two loci: the set $\Sigma_1 \cap \Sigma_2$ on the one hand;
and the set of those $f_1 \in \Sigma_1$ such that $f_1 + L$ is tangent to $\Sigma_2$.
\emph{These two loci are generically not the same}. Indeed, it is a general
fact that the locus of points on $\Sigma_2$ where $L$ is a tangent direction is
included in a hyperplane\footnote{If $B$ is a nondegenerate symmetric bilinear
form on a space $V$ and $\Sigma$ is the quadric hypersurface $B(x,x) = 0$, then
the set of points of $\Sigma$ where the direction $\Vect(u)$ is tangent is
(except if $B(u,u) = 0$) the intersection of $\Sigma$ with the hyperplane
$B(x,u) = 0$. Here $B(u,u) \neq 0$ because a generator of the line $L$ is
$u := \thq(-x/x_1) \thq(-x/x_2) \thq(-x/x_3) \thq(-x/x_4) \not\in \Sigma_2$.}.
By symmetry reason, if the two loci were generically the same,
$\Sigma_1 \cap \Sigma_2$ would be included in two distinct hyperplanes,
so it would be plane, in contradiction with corollary
\ref{cor:noplaneinintersectionofquadrics}. \\

We shall not pursue here this \emph{first approach} to the geometry of $\F$ based
on deducing incidence relations from theta relations. Hereafter we propose a more
algebraic approach, but our main attack will come later in sections
\ref{section:JimboSakaiFamilyII} and \ref{section:Geometrysurgeryandpants} where
we can use our main tool, Mano decomposition.


\subsection{Embedding of $\F$ into $\left(\PC\right)^4$}
\label{rough4surf}

We will give a description of $\F$ as a (non closed) surface inside $\left(\PC\right)^4$.
More precisely, it is a \emph{first attempt} towards such a description. We will assume
some (imprecise \dots) genericity hypothesis and our description is partly heuristic, in
the ``old italian algebraic geometry" style. \\

We have an isomorphism $\F\rightarrow \G$, therefore we have a description of the
surface $\F$ as a (non closed) algebraic surface of the projective space associated
to $W^2$ (which is isomorphic to $\P^7(\C)$). \\

We set~:
\[
\overline G :=
\{(f_1,f_2) \in \Sigma_1 \times \Sigma_2 \tq f_1 - f_2 \in L\} ~~
\text{and} ~~ \overline{\G} := \dfrac{G}{\C^*} \cdot
\]
Thus $\G$ is open and Zariski dense in the projective surface $\overline{\G}$. \\

We denote by $\P(W)$ the projective space associated to $W$ and
$\widetilde \Sigma_1$,  $\widetilde \Sigma_2$, $\tilde L$ the respective
images of $\Sigma_1$,  $\Sigma_2$, $L$ in $\P(W)$;
$\widetilde \Sigma_1$ and $\widetilde \Sigma_2$ are quadric surfaces and
$\tilde L$ is a point. We consider the tangent cones to 
$\widetilde \Sigma_1$, resp.  $\widetilde \Sigma_2$, directed from the point
$\tilde L$,
that we denote by $C(\tilde L,\widetilde \Sigma_1)$, resp.
$C(\tilde L,\widetilde \Sigma_2)$; they are quadratic cones. The intersections 
$Q_1':=\widetilde \Sigma_1\cap C(\widetilde L,\widetilde \Sigma_2)$,
$Q_2':=\widetilde \Sigma_2\cap C(\widetilde L,\widetilde \Sigma_1)$
and  $Q'':=\widetilde \Sigma_1\cap \widetilde \Sigma_2$
are quadratic curves.  \\

We denote $\varpi_1: \G\rightarrow \tilde \Sigma_1$ the projection induced
by $(f_1,f_2)\rightarrow f_1$. We have the following description of the fibration
of $\G$ by $\varpi_1$.
\begin{itemize}
\item
The image of $\varpi_1$ is $\widetilde \Sigma_1 \setminus (Q_1'\cap Q'')$.
Generically $Q_1'\cap Q''$ is a finite set (at most $16$ points). 
\item 
Above $\widetilde\Sigma_1\setminus (Q_1'\cup Q'')$ there are exactly two points
of $\mathcal G$.
\item
Above  $Q_1'\setminus Q''$  there is exactly one point of $\G$. 
\item
Above $\widetilde\Sigma_1 \setminus Q''$ we have a degree two covering ramified
on $Q_1'\setminus Q''$. 
\end{itemize}

\smallskip

We have similar properties for $\varpi_2: \G\rightarrow \tilde \Sigma_1$,
the projection induced by $(f_1,f_2)\rightarrow f_2$. \\

The maps $\varpi_1,~\varpi_2$ extend repectively into  maps
$\overline{\varpi}_1: \overline{\G} \rightarrow \widetilde \Sigma_1$
and $\overline{\varpi}_2: \overline{\G} \rightarrow \widetilde \Sigma_2$.

\bigskip

The quadric surface $\widetilde \Sigma_1$ is a bi-ruled surface:
$\widetilde \Sigma_1\sim \PC\times \PC$.
If $D$ is a line of one of the two families, then 
$D$ intersects the quadric $C(\tilde L,\widetilde \Sigma_1)$
at at most two points and generically at two points.
Then $Q_1'$ can be considered as a curve of bidegree $(2,2)$ in 
$\PC\times \PC$. In order to simplify the notations
we will identify $\widetilde \Sigma_1$ and 
$\PC\times \PC$. In particular we will interpret $\overline{\varpi}_1$
as a map from $\overline{\G}$ to $\PC\times \PC$. We have a similar
description for $\widetilde \Sigma_2$.  \\

\begin{prop}
\begin{itemize}
\item[(i)]
The map $\overline{\varpi}_1$ is a double ramified covering of $\left(\PC\right)^2$
branched over $Q_1'$.
\item[(ii)]
There exists a projective curve $\Gamma\subset \overline{\G}$ such that
$\overline{\varpi}_1 (\Gamma)=Q''$ and $\overline{\G}\setminus \Gamma=\G$.
\end{itemize}
\end{prop}
  
\Pr
\begin{itemize}
\item[(i)]
Follows easily from the description at the beginning of this paragraph.
\item[(ii)]
We denote $\iota:\overline{\G}\rightarrow \overline{\G}$ the involution
associated to the covering.

We verify that the closure in $\overline{\G}$ of the intersection 
$\G\cap \varpi^{-1}\left(Q''\setminus (Q'\cap Q'')\right)$ is
a curve $\Gamma$. Then $\varpi^{-1}(Q'')=\Gamma \cup \iota(\Gamma)$.
We verify that $\varpi(\Gamma)=Q''$. We have $\G=\overline{\G}\setminus \Gamma$.
\end{itemize}
\Finprcourt


\smallskip

Let $D_0$ be a line of one of the families \big($D_0=\{x\}\times \PC$\big).
Then $\overline{\varpi}_1$ induces a double covering 
$\overline{\varpi}_1^{-1}(D_0)\rightarrow D_0$ ramified above
$D_0\cap Q_1'$, that is generically above $2$ points. If $D_0$ is tangent
to $Q'_1$, then we get a ramification above a unique point. 

\smallskip

The map~:
$
(\overline{\varpi}_1,\overline{\varpi}_2): \overline{\G} \rightarrow
\left(\PC\right)^4
$
is clearly regular (\ie\ a morphism of algebraic varieties) and injective.

\begin{prop}\footnote{Stricly speaking it is partly conjectural: \cf\ the introduction
of this paragraph.}
\label{prop1descp4}
The map~:
\[
(\overline{\varpi}_1,\overline{\varpi}_2): \overline{\G} \rightarrow
\left(\PC\right)^4
\] 
realizes an embedding of $\overline{\G}$ into $\left(\PC\right)^4$.
\end{prop}

Then $\F$ can be interpreted as a closed algebraic surface of 
$\left(\PC\right)^4$ minus a closed curve.

\bigskip

Later on, using Mano decomposition, we will give a more precise description
of the surface $\F$. The problem of an explicit relation between our two
descriptions will be tackled in \ref{subsubqpants}.



\section{Mano decomposition}
\label{section:ManoDecomposition}

This extremely useful process was inspired to us by the paper \cite{Mano} of
Toshiyuki Mano. However, we shall give a more precise statement and a direct
proof of the particular property of interest here. We show that degree $2$
order $2$ equations (the ones that appear in the Jimbo-Sakai family, that is
in the linear isomonodromic model of the $q$-Painlev\'e equation) can, in some
sense, be \emph{split} into $q$-hypergeometric components\footnote{As noted 
in a footnote at the beginning of \ref{subsection:thehypergeometricclass}, we 
abusively term $q$-hypergeometric all order $2$ degree $1$ systems, although 
the reducible ones cannot be such by \cite{JRHG1,JRHG2}.}. \\

We saw at the end of section \ref{subsection:thehypergeometricclass} that
the set of classes of systems of $q$-hypergeometric type, as seen through
their monodromy data, admitted a nice geometric classification. So Mano
decomposition allows for an enhanced geometric classification of the Jimbo-Sakai
family. This will be done in section \ref{section:JimboSakaiFamilyII}. \\

Mano decomposition can also be understood as providing a splitting of the
\emph{global} monodromy around the four intermediate singularities into
\emph{local} monodromies around two pairs of singularities\footnote{We
expressed in \ref{subsubsection:intermediatesingularities} our opinion
that defining local monodromies and local Galois groups at intermediate
singularities is one of the most important open problems in modern
$q$-difference theory. This will require a more general version of Mano
decomposition. Extension to higher degrees should be easy along the same
lines, but extension to higher orders (polynomial matrices with coefficients
in $\Matnc$) seems more difficult.}. A discussion of what we accomplish here
appears in our concluding section \ref{section:Geometrysurgeryandpants}, where
a geometrical interpretation (``surgery of pants'') is provided. \\

Since the objects and processes here seem to be new, and since their study is
full of special cases, we tried to present it as clearly and cleanly as possible;
all the more since these special cases seem to have some geometric meaning.
The main result, synthetizing our rather lengthy discussion, is theorem
\ref{thm:existenceMD}, stated at the beginning of \ref{subsection:existenceMD}.


\subsection{Setting, general facts}
\label{subsection:Setting, general facts}

Remember our general setting with generic local data\footnote{Note however
that, these generic local data being given, \emph{all} equations of the
corresponding class will be shown to admit a Mano decomposition.}
$\rho_1,\rho_2 \in \C^*$ (exponents at $0$), $\sigma_1,\sigma_2 \in \C^*$
(exponents at $\infty$), $x_1,x_2,x_3,x_4 \in \C^*$ (``intermediate''
singularities) subject to the following conditions:
\begin{description}
\item[(FR)] Fuchs relation:
$x_1 x_2 x_3 x_4 \equiv \rho_1 \rho_2/\sigma_1 \sigma_2$.
\item[(NR)] Non resonancy: $\rho_1/\rho_2,\sigma_1/\sigma_2 \not\in q^\Z$ and
for $1 \leq k < l \leq 4$, $x_k/x_l \not\in q^\Z$.
\end{description}

Now we select a particular pair of singularities among $x_1,x_2,x_3,x_4$, with
the idea of partially ``localize'' the monodromy around that pair. Let the
indexing be chosen such that $\{x_1,x_2\}$ is the selected pair. We shall
need the following supplementary condition:

\begin{description}
\item[(NS)] Non splitting: for all $i,j = 1,2$,
$\rho_i/\sigma_j \not\equiv x_1 x_2$.
\end{description}

\begin{lem}
Assuming (FR) and (NR), there always is an indexing of $x_1,x_2,x_3,x_4$ such
that (NS) holds.
\end{lem}
\Pr
We leave to the reader the simple combinatorial proof of this fact (hint:
there are three splittings of $\{x_1,x_2,x_3,x_4\}$ in two pairs, while there
are only two splittings of $\rho_1 \rho_2/\sigma_1 \sigma_2$ as a product of
two fractions).
\Finpr

Note that, because of (FR), condition (NS) is equivalent to: for all $i,j = 1,2$,
$\rho_i/\sigma_j \not\equiv x_3 x_4$. However, beyond the partition
of $\{x_1,x_2,x_3,x_4\}$ into $\{x_1,x_2\}$ and $\{x_3,x_4\}$, the two components
do not play a symmetric role: $\{x_1,x_2\}$ will be related to the left factor
and $\{x_3,x_4\}$ to the right one. (Thus there is a kind of ``dual'' decomposition
obtained by permuting the roles of these two pairs.) \\

We let as usual $R := \begin{pmatrix} \rho_1 & 0 \\ 0 & \rho_2 \end{pmatrix}$
and $S := \begin{pmatrix} \sigma_1 & 0 \\ 0 & \sigma_2 \end{pmatrix}$. Also
we write $\x := \{x_1,x_2,x_3,x_4\}$. Then we consider a matrix
$M \in \Mat_2(\Rwg)$ such that:
\begin{enumerate}
\item $\sq M = R M (S x^2)^{-1}$,
\item $\det M \neq 0$,
\item $\det M$ vanishes at $\x$.
\end{enumerate}

Clearly we have $\det M \in \Kwg$ and
$$
\dfrac{\sq(\det M)}{\det M} = \dfrac{\rho_1 \rho_2}{\sigma_1 \sigma_2 x^4} =
\dfrac{x_1 x_2 x_3 x_4}{x^4},
$$
so that by \ref{subsection:someequations}, we deduce from the third condition
that $\det M$ vanishes at all points of $[x_1,x_2,x_3,x_4;q]$ with multiplicity
one, and nowhere else.

\begin{rmk}
\label{rmk:Birkhoffmu=n=2}
Write $M =: \begin{pmatrix} m_{1,1} & m_{1,2} \\ m_{2,1} & m_{2,2} \end{pmatrix}$,
so that $m_{i,j} \in \Rwg$ and $\sq m_{i,j} = (\rho_i/\sigma_j) x^{-2} m_{i,j}$.
From \ref{subsection:someequations} we thus have for $i,j = 1,2$
$m_{i,j} = \lambda_{i,j} \thq(x/\alpha_{i,j}) \thq(x/\beta_{i,j})$, where
$\lambda_{i,j} \in \C$ and $\alpha_{i,j} \beta_{i,j} = \rho_i/\sigma_j$.
This is the ``encoding'' used by Birkhoff in \cite{Birkhoff1} (actually,
the ``characteristic constants'' alluded to in \ref{subsubsection:RHCforqDE}).
However, we shall not use that form in the present section.
\end{rmk}


\subsubsection{A projective invariant}
\label{subsubsection:A projective invariant}

For $k = 1,2,3,4$, $M(x_k) \neq 0$ (otherwise, $\det M$ would have a multiple
zero at $x_k$). Let $\begin{pmatrix} f_k \\ g_k \end{pmatrix}$ a non zero
column of $M(x_k)$ and $\begin{pmatrix} f'_k \\ g'_k \end{pmatrix}$ the other
column, so that 
$\begin{pmatrix} f'_k \\ g'_k \end{pmatrix} \in
\C \begin{pmatrix} f_k \\ g_k \end{pmatrix}$. In particular, $f_k = 0$ means
that $M(x_k) = \begin{pmatrix} 0 & 0 \\ \ast & \ast \end{pmatrix}$ and $g_k = 0$
means that $M(x_k) = \begin{pmatrix} \ast & \ast \\ 0 & 0 \end{pmatrix}$.

\begin{lem}
\label{lem:invariantdeM}
(i) One cannot have $f_1 g_2 = f_2 g_1 = 0$. \\
(ii) The ratio $(f_1 g_2 : f_2 g_1) \in \PC$ is well defined from $M$,
independently of the particular choices of non-zero columns. \\
(iii) This ratio is invariant by the group action of diagonal matrices
$M \mapsto \Gamma M \Delta^{-1}$, $\Gamma,\Delta \in \Ddc$.
\end{lem}
\Pr
(i) We first prove that one cannot have $f_1 = f_2 = 0$ nor $g_1 = g_2 = 0$.
We prove only the first impossibility, the second one is similar. 
So assume that $f_1 = f_2 = 0$. Then $M(x_1)$ and $M(x_2)$ have the form
$\begin{pmatrix} 0 & 0 \\ \ast & \ast \end{pmatrix}$, \ie\ both $m_{1,1}$
and $m_{1,2}$ vanish at $x_1$ and $x_2$. Since $\det M \neq 0$, they cannot
both be null, say $m_{1,j} \neq 0$. Then by \ref{subsection:someequations},
the fact that $\sq m_{1,j} = (\rho_1/\sigma_j) x^{-2} m_{1,j}$ implies that
$x_1 x_2 \equiv \rho_1/\sigma_j$, contradicting (NS). \\
Now to the point: assume for instance that $f_1 = 0$ (the case $g_2 = 0$
being similar). Since by the previous statement $f_2 \neq 0$, the assumption
$f_2 g_1 = 0$ would entail $g_1 = 0$, contradicting the definition of
$\begin{pmatrix} f_k \\ g_k \end{pmatrix}$ as a non zero column. \\
(ii) If for instance $\begin{pmatrix} f'_1 \\ g'_1 \end{pmatrix} \neq 0$, then
$\begin{pmatrix} f'_1 \\ g'_1 \end{pmatrix} = \lambda
\begin{pmatrix} f_1 \\ g_1 \end{pmatrix}$ with $\lambda \neq 0$, whence
$(f'_1 g_2,f_2 g'_1) = \lambda (f_1 g_2,f_2 g_1)$, so that
$(f'_1 g_2 : f_2 g'_1) = (f_1 g_2 : f_2 g_1)$. \\
(iii) Let $M' := \Gamma M \Delta^{-1}$ and let $j_1,j_2$ the indexes of the selected
nonzero columns. With obvious notations, $f'_i = \dfrac{\gamma_1}{\delta_{j_i}} f_i$
and $g'_i = \dfrac{\gamma_2}{\delta_{j_i}} g_i$, whence
$f'_1 g'_2 = \dfrac{\gamma_1 \gamma_2}{\delta_{j_1} \delta_{j_2}} f_1 g_2$ and
$f'_2 g'_1 = \dfrac{\gamma_1 \gamma_2}{\delta_{j_1} \delta_{j_2}} f_2 g_1$.
The conclusion follows.
\Finpr

We thus obtain a well defined map\footnote{\label{footnote:Pi12etnonPi}This map
is obviously related to our selection of the pair $\{x_1,x_2\}$ and it should
properly be denoted $\Pi_{1,2}$ but we don't need to do that for the moment; see 
section \ref{section:Geometrysurgeryandpants}.}
$$
\Pi: F \rightarrow \PC,
$$
which goes to the quotient as
$$
\Pi: \F \rightarrow \PC
$$
(the shorthand notations $F := F_{R,S,\x}$ and $\F := \F_{R,S,\x}$ were introduced
in part \ref{subsubsection:RHB} of subsection \ref{subsection:Defandassump}).
These maps (denoted by the same letter, which, hopefully, will cause no confusion)
will come back in full glory when we describe $F$ and $\F$ as fibered spaces. \\

We leave to the reader the proof of the following, which gives in generic
cases a necessary and sufficient condition for two ordered pairs to give
the same point in $\PC$.

\begin{lem}
\label{lem:mysterieusefractioncroisee}
Let $f_1, f_2, g_1, g_2 \in \C^*$ and $p_1, p_2, q_1, q_2 \in \C^*$. Then:
$$
\left\{\dfrac{p_1 q_2}{p_2 q_1} = \dfrac{f_1 g_2}{f_2 g_1}\right\}
\Longleftrightarrow
\left\{\exists\ \lambda,\mu \in \C^* \;:\;
\begin{pmatrix} f_i \\ g_i \end{pmatrix} =
\begin{pmatrix} \lambda & 0 \\ 0 & \mu \end{pmatrix}
\begin{pmatrix} p_i \\ q_i \end{pmatrix}, i = 1,2
\right\}.
$$
\end{lem}
\Finprcourt


\subsubsection{Two special fibers of the projective invariant $\Pi$}
\label{subsubsection:specialfibersofPi}

The fibers\footnote{We shall find in \ref{subsection:DescriptionF_C} that there
are two more special fibers.} $\Pi^{-1}(0),\Pi^{-1}(\infty) \subset F$ have special
significance. For example, they contain all reducible matrices. Indeed, if $M \in F$
is reducible, \ie\ (definition \ref{defn:reduciblematrices}) if it corresponds
to a reducible system in $E$, then one of its coefficients vanishes (\cf\
\ref{subsubsection:reducibilityJSF}) and $M(x_1)$, $M(x_2)$ acquire one of
the forms $\begin{pmatrix} 0 & 0 \\ \ast & \ast \end{pmatrix}$,
$\begin{pmatrix} \ast & \ast \\ 0 & 0 \end{pmatrix}$. Then $f_1 f_2 g_1 g_2 = 0$,
so that $\Pi(M) = 0$ or $\Pi(M) = \infty$. \\

However, the converse is not true. We briefly describe irreducible matrices $M$
in $\Pi^{-1}(0) \cup \Pi^{-1}(\infty)$. This means that $f_1 f_2 g_1 g_2 = 0$:
$$
\Pi^{-1}(0)  = (f_1 = 0) \cup (g_2 = 0) \text{~and~}
\Pi^{-1}(\infty)  = (f_2 = 0) \cup (g_1 = 0).
$$
We treat the case $f_1 = 0$, the other ones being entirely similar. To simplify
the discussion, we shall assume that $F$ contains no reducible matrices, \ie\
that there is no splitting of (FR) (theorem \ref{theorem:reducibilityJSF}).
A more general situation is tackled by other means in section
\ref{subsection:DescriptionF_C}. \\

If $f_1 = 0$, then $M(x_1) = \begin{pmatrix} 0 & 0 \\ \ast & \ast \end{pmatrix}$.
Since by assumption of irreducibility $m_{1,1}, m_{1,2} \neq 0$, we see that
these coefficients have the forms $m_{1,1} = c_1 \thq(-x/x_1) \thq(-x/x'_1)$,
$m_{1,2} = c_2 \thq(-x/x_1) \thq(-x/x''_1)$, with arbitrary $c_1,c_2 \in \C^*$
and $x'_1,x''_1 \in \C^*$ such that $x_1 x'_1 = \rho_1/\sigma_1$,
$x_1 x''_1 = \rho_1/\sigma_2$. \\

Up to the action $M \mapsto \Gamma M \Delta^{-1}$ by diagonal matrices
$\Gamma := \Diag(\gamma_1,\gamma_2)$, $\Delta := \Diag(\delta_1,\delta_2)$,
we may and shall take $c_1 = c_2 = 1$. Then $m_{2,1} \in V_{2,\frac{\rho_2}{\sigma_1}}$
and $m_{2,2} \in V_{2,\frac{\rho_2}{\sigma_2}}$ are arbitrary non zero elements,
except for the conditions on $\det M$ (vanishing at $\x$ but not everywhere).
The remaining gauge freedom (while retaining the form $c_1 = c_2 = 1$) requires
$\gamma_1 = \delta_1 = \delta_2$, \ie\ only $\C^*$-action on $(m_{2,1},m_{2,2})$
is allowed. So we shall identify the image of $f_1 = 0$ in $\F$ to a subset
of $\dfrac{V_{2,\frac{\rho_2}{\sigma_1}} \times V_{2,\frac{\rho_2}{\sigma_2}}}{\C^*} \cdot$

\paragraph{Non triviality of $\det M$.}

We have equivalences:
\begin{align*}
\det M = 0 &\Longleftrightarrow
m_{2,1} \thq(-x/x_1) \thq(-x/x''_1) = m_{2,2} \thq(-x/x_1) \thq(-x/x'_1) \\
& \Longleftrightarrow m_{2,1} \thq(-x/x''_1) = m_{2,2} \thq(-x/x'_1).
\end{align*}
Since $x'_1 \not\equiv x''_1$ (this follows from (NR)), we conclude
that $m_{2,1} = c \thq(-x/x'_1) \thq(-x/x'_2)$ and
$m_{2,2} = d \thq(-x/x''_1) \thq(-x/x''_2)$, where $c,d \in \C^*$ and
$x'_1 x'_2 = \rho_2/\sigma_1$, $x'_1 x''_2 = \rho_2/\sigma_2$. Actually,
$x'_2 = x_1 \rho_2/\rho_1 = x''_2$; and $\det M = 0 \Leftrightarrow c = d$.
So we see that matrices $M \in (f_1 = 0)$ such that $\det M = 0$ correspond
to a $\C^*$-line in $V_{2,\frac{\rho_2}{\sigma_1}} \times V_{2,\frac{\rho_2}{\sigma_2}}$
and their classes in $\F$ to a unique point in the associated projective space.

\paragraph{Vanishing of $\det M$ at $\x$.}

Vanishing of $\det M$ at $x_1$ is ensured by the above chosen form. Vanishing
at $x_2,x_3,x_4$ then amounts to \emph{two} linear conditions (this is because
$\sq(\det M)/\det M = x_1 x_2 x_3 x_4$):
$$
\begin{cases}
m_{2,1}(x_2) \thq(-x_2/x_1) \thq(-x_2/x''_1) =
m_{2,2}(x_2) \thq(-x_2/x_1) \thq(-x_2/x'_1), \\
m_{2,1}(x_3) \thq(-x_3/x_1) \thq(-x_3/x''_1) =
m_{2,2}(x_3) \thq(-x_3/x_1) \thq(-x_3/x'_1).
\end{cases}
$$
Simplifications are allowed by (NR) and yield:
$$
\begin{cases}
m_{2,1}(x_2) \thq(-x_2/x''_1) = m_{2,2}(x_2) \thq(-x_2/x'_1), \\
m_{2,1}(x_3) \thq(-x_3/x''_1) = m_{2,2}(x_3) \thq(-x_3/x'_1).
\end{cases}
$$
Each of the two linear conditions is non trivial (because, again by (NR),
$x'_1 \not\equiv x''_1$), whence defines a hyperplane in the
$4$-dimensional product space
$V_{2,\frac{\rho_2}{\sigma_1}} \times V_{2,\frac{\rho_2}{\sigma_2}}$. These
hyperplanes are distinct: indeed, using non splitting, one can find
$f \in V_{2,\frac{\rho_2}{\sigma_1}}$ and $g \in V_{2,\frac{\rho_2}{\sigma_2}}$,
each of them vanishing at $x_2$ but not at $x_3$; then, for a proper
choice of $c,d \in \C^2$, the pair $(m_{2,1},m_{2,2}) := (cf,dg)$ belongs
to the first hyperplane but not to the second one. Therefore their intersection
is a plane containing the line $\det M = 0$. Going to the associated projective
space yields a line. So we conclude that the image of the component $(f_1 = 0)$
is a projective line minus a point, \ie\ an \emph{affine} line $\C$. \\

The fiber $\Pi^{-1}(0)$ has the two components $f_1 = 0$ and $g_2 = 0$ each
isomorphic to the line $\C$ (this meaning that the obvious bijections are
biregular in the algebro-geometric sense). Matrices in the intersection of
these components have the form:
$$
M = \begin{pmatrix}
c_{1,1} \thq(-x/x_1) \thq(-x/x'_1) & c_{1,2} \thq(-x/x_1) \thq(-x/x''_1) \\
c_{2,1} \thq(-x/x_2) \thq(-x/x'_2) & c_{2,2} \thq(-x/x_2) \thq(-x/x''_2)
\end{pmatrix}
$$
for some arbitrary $c_{i,j} \in \C^*$ and $x'_i,x''_i \in \C^*$ determined by
obvious conditions. The $\Gamma,\Delta$-action allows one to take three of
the $c_{i,j}$ with value $1$ and the fourth is then determined by the vanishing
of $\det M$ at $x_3$, so our two projective lines intersect at exacly one point,
corresponding to the double degeneracy $f_1 = g_2 = 0$.

\begin{thm}
\label{thm:specialfibersofPi}
Assume non splitting of (FR) (\ie\ all classes in $\F$ are irreducible). The
special fibers $\Pi^{-1}(0)$ and $\Pi^{-1}(\infty)$ of $\Pi: \F \rightarrow \PC$
are both made of two affine lines intersecting at one point.
\end{thm}
\Finprcourt

The study of the other fibers will be done in section
\ref{section:JimboSakaiFamilyII}.


\subsubsection{An elliptic function}
\label{subsubsection:anellipticfunction}

Along with the projective invariant $\Pi$, a central role will be played by the
following auxiliary function\footnote{This should really be denoted $\Phi_{1,2}$,
see footnote \ref{footnote:Pi12etnonPi} page \pageref{footnote:Pi12etnonPi}.}:
$$
\Phi(\xi) := \dfrac
{\thq\left(\frac{x_1}{\rho_1}\xi\right)\thq\left(\frac{x_2}{\rho_2}\xi\right)}
{\thq\left(\frac{x_1}{\rho_2}\xi\right)\thq\left(\frac{x_2}{\rho_1}\xi\right)}
\cdot
$$
It is readily verified that $\Phi$ is an elliptic function on $\C^*$, \ie\ that
$\Phi(q \xi) = \Phi(\xi)$, so, after the conventions of subsection
\ref{subsection:q-notations}, we may and will see it as a mapping (denoted the
same) $\Phi: \Eq \rightarrow \P^1(\C)$. \\

Under assumption (NR) there is no cancellation of zeroes between the numerator
and denominator (which are both holomorphic over $\C^*$), so the elliptic function
$\Phi$ has degree $2$. Also:
\begin{equation}
\label{eqn:involution1}
\forall \xi_1,\xi_2 \in \C^* \;,\;
\xi_1 \xi_2 \equiv \dfrac{\rho_1 \rho_2}{x_1 x_2}
\Longrightarrow \Phi(\xi_1) = \Phi(\xi_2).
\end{equation}
(Argument: by ellipticity of $\Phi$, one may assume equality in the premise; and
then direct calculation yields the result.) More precisely, we get:

\begin{prop}
(i) The elliptic function $\Phi$ realizes a degree $2$ ramified covering 
$\Eq \rightarrow \PC$ with $4$ critical values. \\
(ii) Generic fibers have the form $\{\overline{\xi_1},\overline{\xi_1}\}$,
where $\xi_1,\xi_2 \in \Cq$,
$\xi_1 \xi_2 \equiv \dfrac{\rho_1 \rho_2}{x_1 x_2}$,
$\xi_1 \not\equiv \xi_2$. \\
(iii) The four singular fibers have the form $\{\overline{\xi}\}$, where
$\xi \in \Cq$, $\xi^2 \equiv \dfrac{\rho_1 \rho_2}{x_1 x_2}$.
\end{prop}
\Finpr

Also note that the divisor of $\Phi$ depends only on the local data:
$$
\div_\Eq(\Phi) =
\left[\pi\left(-\frac{\rho_1}{x_1}\right)\right] +
\left[\pi\left(-\frac{\rho_2}{x_2}\right)\right] -
\left[\pi\left(-\frac{\rho_1}{x_2}\right)\right] -
\left[\pi\left(-\frac{\rho_2}{x_1}\right)\right].
$$
(Recall that $\pi: \C^* \rightarrow \Eq$ is the canonical projection.)
Up to a non zero constant, $\Phi$ is determined by this divisor.


\subsubsection{An involution of $\Eq$}
\label{subsubsection:involutionofEq}

Here is the geometrical meaning of $\Phi$. The relation
$\xi_1 \leftrightarrow \xi_2$ whenever
$\xi_1 \xi_2 \equiv \dfrac{\rho_1 \rho_2}{x_1 x_2}$ defines an
involution on $\Eq$, and the mapping $\Phi: \Eq \rightarrow \PC$ defines
a quotient of $\Eq$ by this involution. Up to a dilatation in $\PC$, this
is the only realisation of this quotient with the following special fibers:
\begin{align*}
\text{Fiber at~} 0 &= \left\{
\pi\left(-\frac{\rho_1}{x_1}\right),\pi\left(-\frac{\rho_2}{x_2}\right)
\right\}, \\
\text{Fiber at~} \infty &= \left\{
\pi\left(-\frac{\rho_1}{x_2}\right),\pi\left(-\frac{\rho_2}{x_1}\right)
\right\}.
\end{align*}


\subsection{Definition, gauge freedom}
\label{subsection:defngaugefreedom}

\begin{defn}
A \emph{Mano decomposition $M = PQ$ with factor $C \in \GL_2(\C)$} is given by
$P,Q \in \Mat_2(\Rwg)$ such that $\sq P = R P (Cx)^{-1}$ and $\det P$ vanishes
at $x_1,x_2$. (It is understood that $x_1,x_2$ have previously been chosen and
we do not mention them in the terminology.)
\end{defn}

It follows immediately from the definition, the assumptions on $M$ and the
properties recalled in \ref{subsection:someequations} that:
\begin{itemize}
\item $\det P \neq 0$ and $\det P$ vanishes at $[x_1,x_2;q]$ with simple zeroes
  and nowhere else,
\item $\det Q \neq 0$ and $\det Q$ vanishes at $[x_3,x_4;q]$ with simple zeroes
  and nowhere else,
\item $\sq Q = C Q (S x)^{-1}$,
\item $\dfrac{\sq(\det P)}{\det P} = \dfrac{\rho_1 \rho_2}{(\det C) x^2}$,
\item $\dfrac{\sq(\det Q)}{\det Q} = \dfrac{\det C}{\sigma_1 \sigma_2 x^2}$,  
\item $\det C \equiv (\rho_1 \rho_2)/(x_1 x_2)$.
\end{itemize}


\subsubsection{First consequences}
\label{subsubsection:First consequences}

Since we expect $\det P$ to have simple zeroes at $x_1, x_2$, we can,
for $i = 1,2$, choose a non zero column $\begin{pmatrix} p_i \\ q_i \end{pmatrix}$
of $P(x_i)$; and the other column $\begin{pmatrix} p'_i \\ q'_i \end{pmatrix}$ then
necessarily belongs to $\C \begin{pmatrix} p_i \\ q_i \end{pmatrix}$. Arguments
similar to those used in the proof of lemma \ref{lem:invariantdeM} yield the
following.

\begin{lem}
\label{lem:invariantdeP}
(i) One cannot have $p_1 q_2 = p_2 q_1 = 0$. \\
(ii) The ratio $(p_1 q_2 : p_2 q_1) \in \PC$ is well defined from $P$,
independently of the particular choices of non-zero columns.
\end{lem}
\Finpr

The ``invariant'' $(p_1 q_2 : p_2 q_1)$ will be related to $\Pi(M)$ in proposition
\ref{prop:lesdeuxfractionscroisees} and to the values of the elliptic function $\Phi$
in lemma \ref{lem:autrefractioncroisee}.


\subsubsection{Gauge freedom}
\label{subsubsection:gaugefreedom}

\begin{prop}[Gauge freedom]
\label{prop:gaugefreedom}
Let $M = P Q$ a Mano decomposition with factor $C$. \\
(i) Let $\Lambda \in \GL_2(\Rwg)$ be such that $\Lambda[C] \in \GL_2(\C)$.
Let $P' := P \Lambda^{-1}$, $Q' := \Lambda Q$ and $C' := \Lambda[C]$. Then
$M = P' Q'$ is a Mano decomposition with factor $C'$. \\
(ii) All Mano decompositions of $M$ are obtained that way.
\end{prop}
\Pr
(i) comes by a mechanical calculation. \\
(ii) Let $M = P' Q'$ a Mano decomposition with factor $C'$. Set
$\Lambda := {P'}^{-1} P = Q' Q^{-1}$. Since ${P'}^{-1} P$ has poles only at
$[x_1,x_2;q]$ and $Q' Q^{-1}$ at $[x_3,x_4;q]$ and since by (NR) these sets
do not meet, $\Lambda \in \GL_2(\Rwg)$. From
$(\sq P) C P^{-1} = \dfrac{1}{x} R = (\sq P') C' {P'}^{-1}$ we draw that
$C' = \Lambda[C]$.
\Finpr

Actually, $\Lambda$ is a Laurent polynomial (with matrix coefficients). Indeed,
write $\Lambda = \sum x^n \Lambda_n$. Then relation $C' = \Lambda[C]$ implies
that $C' \Lambda_n = q^n \Lambda_n C$. But equation with matricial unknown $X$
$C' X - X (q^n C) = 0$ has non trivial solutions if and only if $\Sp C'$ and
$\Sp(q^n C)$ intersect, which is possible only for a finite number of values of $n$.
This statement can be made more precise using proposition \ref{prop:NormalformsforC}
and its corollary.


\subsubsection{Normal forms for $C$}
\label{subsubsection:NormalformsforC}

\begin{prop}[Normal forms for $C$]
\label{prop:NormalformsforC}
The central factor $C$ can be taken in one and only one of the following forms:
\begin{align*}
C &= \begin{pmatrix} \xi & 0 \\ 0 & \xi \end{pmatrix} \quad &\text{(trivial form)} \\
C &= \begin{pmatrix} \xi_1 & 0 \\ 0 & \xi_2 \end{pmatrix},\;\;\; \xi_1 \neq \xi_2 \quad
&\text{(generic form)}, \\
C &= \begin{pmatrix} \xi & \xi \\ 0 & \xi \end{pmatrix} \quad &\text{(logarithmic form)},
\end{align*}
with $\xi_1,\xi_2$ or $\xi$ in the fundamental annulus $\Cq: |q| < |z| \leq 1$. \\
Moreover, the form is unique except that in generic form $\xi_1,\xi_2$ can be permuted.
\end{prop}
\Pr
Among gauge transforms are conjugacies (by $\GL_2(\C)$) so, by standard reduction
theory, $C$ can be taken diagonal or in the form
$\begin{pmatrix} \xi & \xi \\ 0 & \xi \end{pmatrix}$. Then gauge transformation
by so called ``shearing matrices'' $\Diag(x^\mu,x^\nu)$, $\mu,\nu \in \Z$, allows
one to bring $\xi_1,\xi_2$ or $\xi$ into $\Cq$. The possibility to permute
$\xi_1,\xi_2$ comes from the equality:
$$
\begin{pmatrix} \xi_2 & 0 \\ 0 & \xi_1 \end{pmatrix} =
J\left[\begin{pmatrix} \xi_1 & 0 \\ 0 & \xi_2 \end{pmatrix}\right],
\text{~where~}
J := \begin{pmatrix} 0 & 1 \\ 1 & 0 \end{pmatrix}.
$$
We now show that this is the only defect of uniqueness. So let $C,C'$ in one of
the quoted forms and let $\Lambda = (\lambda_{i,j})_{1 \leq i,j \leq 2} \in \GL_2(\Rwg)$
such that $\Lambda[C] = C'$, or, equivalently, $\sq \Lambda = C' \Lambda C^{-1}$. \\

\underline{Case 1, $C$ and $C'$ generic or trivial:} Write $C = \Diag(\xi_1,\xi_2)$
and $C' = \Diag(\xi'_1,\xi'_2)$, whence
$\sq \lambda_{i,j} = \dfrac{\xi'_i}{\xi_j} \lambda_{i,j}$. By section
\ref{subsection:someequations} we know that this is possible with
$\lambda_{i,j} \in \Rwg \setminus \{0\}$ only if $\dfrac{\xi'_i}{\xi_j} \in q^\Z$;
since $\xi'_i,\xi_j \in \Cq$, this would mean $\xi'_i = \xi_j$ and
$\lambda_{i,j} \in \C$. The end of the proof is standard linear algebra. \\

\underline{Case 2, $C$ and $C'$ logarithmic:}
Write $C = \begin{pmatrix} \xi & \xi \\ 0 & \xi \end{pmatrix}$
$C' = \begin{pmatrix} \xi' & \xi' \\ 0 & \xi' \end{pmatrix}$, whence:
$$
\sq \Lambda = \dfrac{\xi'}{\xi}
\begin{pmatrix}
\lambda_{1,1} + \lambda_{2,1} & \lambda_{1,2} + \lambda_{2,2} - (\lambda_{1,1} + \lambda_{2,1}) \\
\lambda_{2,1} & \lambda_{2,2} - \lambda_{2,1}
\end{pmatrix}.
$$
If $\xi' \neq \xi$, equation $\sq \lambda_{2,1} = \dfrac{\xi'}{\xi} \lambda_{2,1}$
implies $\lambda_{2,1} = 0$; then equation
$\sq \lambda_{1,1} = \dfrac{\xi'}{\xi} \lambda_{1,1}$ implies $\lambda_{1,1} = 0$;
but then $\Lambda$ is not invertible, contradiction. \\
So $\xi' = \xi$ and $\lambda_{2,1} \in \C$. But then we know from subsection
\ref{subsubsection:qlog} that $\sq \lambda_{1,1} = \lambda_{1,1} + \lambda_{2,1} $
is possible with $\lambda_{1,1} \in \Rwg$ only if $\lambda_{2,1} = 0$. The end
of the proof along the same lines is easy. \\

\underline{Case 3, mixed case:} similar calculations left to the reader show
this case to be impossible because there is no $q$-logarithm in $\Rwg$.
\Finpr

From the above proof, one can also draw:

\begin{cor}
\label{cor:NormalformsforC}
If $C$ is in normal form, then
$\Lambda[C] = C \Longleftrightarrow
\left(\Lambda \in \GL_2(\C) \text{~and~} [\Lambda,C] = 0\right)$.
\end{cor}
\Finprcourt


\subsection{Necessary conditions for the equality $M = P Q$}

As a preliminary observation, note that if $M = P Q$ is a Mano decomposition with
factor $C$ in normal form, then we can still replace $P$ by any $P' = P \Lambda^{-1}$
where $\Lambda \in \GL_2(\C)$ commutes with $C$. The apparently more general case
$\Lambda \in \GL_2(\Rwg)$, $\Lambda[C] = C$ boils down to this one by 
corollary \ref{cor:NormalformsforC}.


\subsubsection{Possible forms of the factor $P$}
\label{subsubsection:possibleforms}

We can (and will) search the factor $C$ in one of the forms shown in proposition
\ref{prop:NormalformsforC}. We also know from the property of $\det C$ stated in 
\ref{subsection:defngaugefreedom} that, according to the case,
$\xi_1\xi_2 \equiv (\rho_1 \rho_2)/(x_1 x_2)$ or
$\xi^2 \equiv (\rho_1 \rho_2)/(x_1 x_2)$.

\paragraph{Trivial form.}
In the case $C = \begin{pmatrix} \xi & 0 \\ 0 & \xi \end{pmatrix}$, we must have:
$$
P = \begin{pmatrix}
\alpha_{1,1} \thq\left(\frac{\xi}{\rho_1}x\right) &
\alpha_{1,2} \thq\left(\frac{\xi}{\rho_1}x\right) \\
\alpha_{2,1} \thq\left(\frac{\xi}{\rho_2}x\right) &
\alpha_{2,2} \thq\left(\frac{\xi}{\rho_2}x\right)
\end{pmatrix} =
\begin{pmatrix} \thq\left(\frac{\xi}{\rho_1}x\right) & 0 \\
0 & \thq\left(\frac{\xi}{\rho_2}x\right) \end{pmatrix}
\begin{pmatrix} \alpha_{1,1} & \alpha_{1,2} \\
\alpha_{2,1} & \alpha_{2,2} \end{pmatrix},
$$
for some $\alpha_{i,j} \in \C$, $i,j = 1,2$. Thus
$\det P = 0 \Leftrightarrow \alpha_{1,1} \alpha_{2,2} - \alpha_{1,2} \alpha_{2,1} = 0$.
Actually, the constant right factor $(\alpha_{i,j})$ can be taken rid of by the
observation at the beginning of this section, \ie\ we can take:
$$
P = \begin{pmatrix} \thq\left(\frac{\xi}{\rho_1}x\right) & 0 \\
0 & \thq\left(\frac{\xi}{\rho_2}x\right) \end{pmatrix}.
$$

\paragraph{Generic form.}
In the case $C = \begin{pmatrix} \xi_1 & 0 \\ 0 & \xi_2 \end{pmatrix}$,
$\xi_1 \not\equiv \xi_2$, we must have:
$$
P = \begin{pmatrix}
\alpha_{1,1} \thq\left(\frac{\xi_1}{\rho_1}x\right) &
\alpha_{1,2} \thq\left(\frac{\xi_2}{\rho_1}x\right) \\
\alpha_{2,1} \thq\left(\frac{\xi_1}{\rho_2}x\right) &
\alpha_{2,2} \thq\left(\frac{\xi_2}{\rho_2}x\right)
\end{pmatrix},
$$
for some $\alpha_{i,j} \in \C$, $i,j = 1,2$. Then:
$$
\det P = \alpha_{1,1} \alpha_{2,2}
\thq\left(\frac{\xi_1}{\rho_1}x\right) \thq\left(\frac{\xi_2}{\rho_2}x\right) -
\alpha_{1,2} \alpha_{2,1}
\thq\left(\frac{\xi_2}{\rho_1}x\right) \thq\left(\frac{\xi_1}{\rho_2}x\right),
$$
so that:
$$
\det P = 0 \Longleftrightarrow
\alpha_{1,1} \alpha_{2,2} = \alpha_{1,2} \alpha_{2,1} = 0.
$$
Indeed, since $\xi_1 \not\equiv \xi_2$, because of (NR) the functions
$\thq\left(\frac{\xi_1}{\rho_1}x\right) \thq\left(\frac{\xi_2}{\rho_2}x\right)$
and
$\thq\left(\frac{\xi_2}{\rho_1}x\right) \thq\left(\frac{\xi_1}{\rho_2}x\right)$
have no common zero.

\paragraph{Logarithmic form.}
In the case $C = \begin{pmatrix} \xi & \xi \\ 0 & \xi \end{pmatrix}$,
the coefficients $p_{i,j}$ of $P$ satisfy
$\begin{cases}
\sq p_{i,1} = \frac{\rho_i}{\xi} p_{i,1}, \\
\sq p_{i,2} = \frac{\rho_i}{\xi} (p_{i,2} - p_{i,1}),
\end{cases}$ $i = 1,2$. 
According to  lemma \ref{lem:logarithmiccase} in subsection \ref{subsubsection:qlog},
setting $\phi_i(x) := \thq\left(\frac{\xi}{\rho_i}x\right)$ and
$\psi_i(x) := x \phi_i'(x) =
\frac{\xi}{\rho_i} x \thqp\left(\frac{\xi}{\rho_i}x\right)$,
we must have:
$P = \begin{pmatrix}
\alpha_{1,1} \phi_1 & \alpha_{1,1} \psi_1 + \alpha_{1,2} \phi_1 \\
\alpha_{2,1} \phi_2 & \alpha_{2,1} \psi_2 + \alpha_{2,2} \phi_2
\end{pmatrix}$
for some $\alpha_{i,j} \in \C$, $i,j = 1,2$, so that
\begin{align*}
\det P &= (\alpha_{1,1} \alpha_{2,2} - \alpha_{2,1} \alpha_{1,2}) \phi_1 \phi_2 +
\alpha_{1,1} \alpha_{2,1} x (\phi_1 \phi_2' - \phi_2 \phi_1') \\
&=
x \phi_1 \phi_2 \times \text{~logarithmic derivative of~}
\left(x^{\alpha_{1,1} \alpha_{2,2} - \alpha_{2,1} \alpha_{1,2}}
\left(\frac{\phi_2}{\phi_1}\right)^{\alpha_{1,1} \alpha_{2,1}}\right),
\end{align*}
whence
$$
\det P = 0 \Longleftrightarrow
\begin{cases}
\alpha_{1,1} \alpha_{2,1} = 0 \text{~and~} \\
\alpha_{1,1} \alpha_{2,2} - \alpha_{2,1} \alpha_{1,2} = 0
\end{cases}
\Longleftrightarrow
\begin{cases}
(\alpha_{1,1} = \alpha_{2,1} = 0) \text{~or~} \\
(\alpha_{1,1} = \alpha_{1,2} = 0) \text{~or~} \\
(\alpha_{2,1} = \alpha_{2,2} = 0).
\end{cases}
$$


\subsubsection{Further necessary conditions}

From lemma \ref{lem:invariantdeP}, we know that, $P$ being given, 
$(p_1 q_2 : p_2 q_1) \in \PC$ is well defined. We identify
this point of the projective line with
$\dfrac{p_1 q_2}{p_2 q_1} \in \Sr = \C \cup \{\infty\}$, considered as
the target space of elliptic functions, in particular of the function
$\Phi$ defined in \ref{subsubsection:anellipticfunction}. \\

In what follows, we write $\xi_1,\xi_2$ the eigenvalues of $C$, with
maybe (in the trivial or in the logarithmic case) $\xi_1 = \xi_2$.
Their images $\overline{\xi_1},\overline{\xi_2} \in \Eq$ are the
\emph{exponents} of $C$.

\begin{lem}
\label{lem:autrefractioncroisee}
In all cases above, the exponents $\overline{\xi_1},\overline{\xi_2}$
constitute the fiber $\Phi^{-1}\left(\dfrac{p_1 q_2}{p_2 q_1}\right)$,
\ie:
$$
\Phi^{-1}\left(\dfrac{p_1 q_2}{p_2 q_1}\right) =
\left\{\overline{\xi_1},\overline{\xi_2}\right\}.
$$
\end{lem}
\Pr
Since $\xi_1\xi_2 \equiv (\rho_1 \rho_2)/(x_1 x_2)$, it will be
enough to prove that:
$$
\Phi(\xi_1) = \dfrac{p_1 q_2}{p_2 q_1} \cdot
$$
In all cases above, the first column of $P(x)$ is
$A(x) := \begin{pmatrix} \alpha_{1,1} \thq\left(\frac{\xi_1}{\rho_1}x\right) \\
\alpha_{2,1} \thq\left(\frac{\xi_1}{\rho_2}x\right) \end{pmatrix}$. \\
If $A(x_1)$ and $A(x_2)$ are both non zero, we can take them as column
vectors $\begin{pmatrix} p_i \\ q_i \end{pmatrix}$. An immediate calculation
then gives $\dfrac{p_1 q_2}{p_2 q_1} = \Phi(\xi_1)$ as wanted. \\
Assume that for instance $A(x_1) = 0$ (the case $A(x_2) = 0$ being entirely
similar). We have
$\alpha_{1,1} \thq\left(\frac{\xi_1}{\rho_1}x\right) =
\alpha_{2,1} \thq\left(\frac{\xi_1}{\rho_2}x\right) = 0$.
We cannot have $\alpha_{1,1} = \alpha_{2,1} = 0$ because then $\det P$ would
vanish identically; neither can we have
$\thq\left(\frac{\xi_1}{\rho_1}x\right) =
\thq\left(\frac{\xi_1}{\rho_2}x\right) = 0$, because of condition (NR).
Therefore, we have two cases to consider:
\begin{enumerate}
\item $\alpha_{1,1} \neq 0$, $\alpha_{2,1} = 0$ and
$\thq\left(\frac{\xi_1}{\rho_1}x\right) = 0$. Then $\Phi(\xi_1) = 0$.
Also $\thq\left(\frac{\xi_1}{\rho_2}x\right) \neq 0$, so $A(x_2) \neq 0$,
so we can take it as $\begin{pmatrix} p_2 \\ q_2 \end{pmatrix}$, so
$q_2 = \alpha_{2,1} \thq\left(\frac{\xi_1}{\rho_2}x\right) = 0$, so
$\dfrac{p_1 q_2}{p_2 q_1} = 0 = \Phi(\xi_1)$ as wanted.
\item $\alpha_{1,1} = 0$, $\alpha_{2,1} \neq 0$ and
$\thq\left(\frac{\xi_1}{\rho_2}x\right) = 0$. A similar calculation
yields $\Phi(\xi_1) = \dfrac{p_1 q_2}{p_2 q_1} = \infty$.
\end{enumerate}
\Finprcourt


\subsection{Existence of Mano decomposition}
\label{subsection:existenceMD}

Let $R,S,\x$ the local data described in subsection
\ref{subsection:Setting, general facts}, subject to the assumptions (FR), (NR)
and (NS) stated there. \\

Let $\Phi(\xi) := \dfrac
{\thq\left(\frac{x_1}{\rho_1}\xi\right)\thq\left(\frac{x_2}{\rho_2}\xi\right)}
{\thq\left(\frac{x_1}{\rho_2}\xi\right)\thq\left(\frac{x_2}{\rho_1}\xi\right)}$
the elliptic function defined and studied in subsection
\ref{subsubsection:anellipticfunction}. \\

For every matrix $M \in F_{R,S,\x}$, recall $\Pi(M) := (f_1 g_2 : f_2 g_1) \in \PC$
the projective invariant defined and studied in subsection
\ref{subsubsection:A projective invariant}.

\begin{thm}[Existence of Mano decomposition]
\label{thm:existenceMD}
Every matrix $M \in F_{R,S,\x}$ admits a Mano decomposition $M = PQ$ with (some)
factor $C$. More precisely:
\begin{enumerate}
\item If $\Pi(M)$ is not a critical value of $\Phi$, the decomposition is
generic (\ie\ the factor $C$ is in generic form).
\item If $\Pi(M)$ is a critical value of $\Phi$ and $\Pi(M) \neq 0,\infty$,
the decomposition is logarithmic.
\item If $\Pi(M)$ is a critical value of $\Phi$ and $\Pi(M) = 0$ or $\infty$,
the decomposition is logarithmic, except for the respective degenerate cases:
$f_1 = g_2 = 0$ and $f_2 = g_1 = 0$; in these degenerate cases, the decomposition
is trivial.
\end{enumerate}
\end{thm}

The rest of this section is devoted to the proof. We shall keep the notations
$\begin{pmatrix} f_i \\ g_i \end{pmatrix}$, $i = 1,2$, for a non zero column
of $M(x_i)$ (\ref{subsubsection:A projective invariant}) and similarly
$\begin{pmatrix} p_i \\ q_i \end{pmatrix}$, $i = 1,2$, for a non zero column
of the left factor $P(x_i)$ (\ref{subsubsection:First consequences}). \\

We write the fiber of $\Phi$ at $\Pi(M) \in \PC$ as:
$$
\Phi^{-1}\left(\Pi(M)\right) = \Phi^{-1}\left(\dfrac{f_1 g_2}{f_2 g_1}\right) =
\left\{\overline{\xi_1},\overline{\xi_2}\right\},
$$
where $\xi_1,\xi_2 \in \Cq$ (subsection \ref{subsubsection:anellipticfunction}).
Thus we have:
\begin{equation}
\label{eqn:involution2}
\xi_1 \xi_2 \equiv \dfrac{\rho_1 \rho_2}{x_1 x_2}.
\end{equation}


\subsubsection{Preliminary reductions}

\begin{prop}
\label{prop:reducedconditions}
In order to prove theorem \ref{thm:existenceMD}, it is enough to find a factor $C$
and $P \in \Mat_2(\Rwg)$ such that:
\begin{enumerate}
\item $\sq P = R P (Cx)^{-1}$,
\item $\det P \neq 0$ and $\det P$ vanishes at $x_1$,
\item $P^{-1} M$ is well defined (\ie\ has no pole) at $x_1,x_2$.
\end{enumerate}
\end{prop}
\Pr
By \ref{subsection:someequations}, $\det P$ vanishes at $[x_1,x_2;q]$
and $Q := P^{-1} M$ vanishes at $[x_3,x_4;q]$ under the same conditions
as usual and immediate calculation shows that we got a Mano decomposition.
\Finpr

Now, from Cramer's rule $P^{-1} = (\det P)^{-1} \tilde{P}$, we see that the
last condition can be replaced by: $(\tilde{P} M)(x_i) = 0$, $i = 1,2$. This
in turn is equivalent to: $\begin{pmatrix} f_i \\ g_i \end{pmatrix}$ is a
linear combination of the columns of $P(x_i)$ for $i = 1,2$, \ie\ that it
is proportional to the selected non zero columns.

\begin{cor}
\label{cor:reducedconditions}
Condition 3 above can be replaced by:
\begin{enumerate}
\item[3'.] $\begin{pmatrix} f_i \\ g_i \end{pmatrix} \in
\C \begin{pmatrix} p_i \\ q_i \end{pmatrix}, i = 1,2$.
\end{enumerate}
\end{cor}

Using lemmas \ref{lem:mysterieusefractioncroisee} and \ref{lem:invariantdeP},
we complete the above as follows:

\begin{prop}
\label{prop:lesdeuxfractionscroisees}
Let $P$ a left factor of $M$ in a Mano decomposition and keep previous
notations $\begin{pmatrix} f_i \\ g_i \end{pmatrix}$ and
$\begin{pmatrix} p_i \\ q_i \end{pmatrix}$, $i = 1,2$. Then one cannot have
$p_1 q_2 = p_2 q_1 = 0$, the ratio $(p_1 q_2 : p_2 q_1) \in \PC$
is well defined from $P$ and $(p_1 q_2 : p_2 q_1) = (f_1 g_2 : f_2 g_1)$.
\end{prop}


\subsubsection{Proof of existence, case I: $f_1f_2g_1g_2 \neq 0$}

From the proposition above, $p_1p_2q_1q_2 \neq 0$. We shall use lemma
\ref{lem:mysterieusefractioncroisee}.

\paragraph{Subcase Ia: $\xi_1 \neq \xi_2$.}

For some $\alpha_1,\alpha_2 \in \C$ to be determined (not both zero), we set:
$$
P_0 = \begin{pmatrix}
\thq\left(\frac{\xi_1}{\rho_1}x\right) &
\alpha_{1} \thq\left(\frac{\xi_2}{\rho_1}x\right) \\
\thq\left(\frac{\xi_1}{\rho_2}x\right) &
\alpha_{2} \thq\left(\frac{\xi_2}{\rho_2}x\right)
\end{pmatrix}
\quad \text{~and~} C := \begin{pmatrix} \xi_1 & 0 \\ 0 & \xi_2 \end{pmatrix}.
$$
Then, writing
$T_1(x) :=
\thq\left(\frac{\xi_1}{\rho_1}x\right) \thq\left(\frac{\xi_2}{\rho_2}x\right)$
and
$T_2(x) :=
\thq\left(\frac{\xi_2}{\rho_1}x\right) \thq\left(\frac{\xi_1}{\rho_2}x\right)$,
we have:
$$
\det P_0 = \alpha_{2} T_1 - \alpha_{1} T_2,
$$
so taking $\alpha_1 := T_2(x_1)$ and $\alpha_2 := T_1(x_1)$ (which, by (NR),
cannot both be $0$), we see that $P_0$ satisfies the first two conditions of
proposition \ref{prop:reducedconditions}. Writing 
$\begin{pmatrix} p_i \\ q_i \end{pmatrix}$ the first column of $P_0(x_i)$,
$i = 1,2$, we fall, by choice of $\xi_1$, under the assumptions of lemma
\ref{lem:mysterieusefractioncroisee}. We set
$\Lambda := \begin{pmatrix} \lambda & 0 \\ 0 & \mu \end{pmatrix}$ with
$\lambda,\mu$ as provided by lemma \ref{lem:mysterieusefractioncroisee}
and then $P := \Lambda P_0$. Using the fact that $\Lambda$ and $C$
commute with each other, we easily conclude that $P$ satisfies all three
conditions of proposition \ref{prop:reducedconditions}, so theorem
\ref{thm:existenceMD} is proved in this case.

\paragraph{Subcase Ib: $\xi_1 = \xi_2 =: \xi$.}

From the definition of $\Phi$, we see that
$\thq\left(\frac{x_j}{\rho_i} \xi\right) \neq 0$ for $i,j = 1,2$.
Writing $\phi_i(x) := \thq\left(\frac{\xi x}{\rho_i}\right)$, $i = 1,2$,
it follows that $(\phi_1 \phi_2)(x_1) \neq 0$. \\
We take $C := \begin{pmatrix} \xi & \xi \\ 0 & \xi \end{pmatrix}$
and set $P_0 := \begin{pmatrix} \phi_1 & x \phi'_1 + \alpha_1 \phi_1 \\
\phi_2 & x \phi'_2 + \alpha_2 \phi_2 \end{pmatrix}$,
so that $\sq P_0 = R P (Cx)^{-1}$. Also 
$\det P_0 = (\alpha_2 - \alpha_1) \phi_1 \phi_2 + x(\phi_1 \phi'_2 - \phi'_1 \phi_2)$
is equal to $x \phi_1 \phi_2$ times the logarithmic derivative of
$x^{\alpha_2 - \alpha_1} \frac{\phi_2}{\phi_1}$ and certainly does not vanish identically. \\
Since $(\phi_1 \phi_2)(x_1) \neq 0$, a proper choice of $\alpha_1, \alpha_2$
yields $\det P_0(x_1) = 0$ and the argument can be completed exactly as in
subcase Ia to obtain a left factor $P := \Lambda P_0$.


\subsubsection{Proof of existence, case II: $f_1f_2g_1g_2 = 0$}

Here we have $(f_1g_2)(f_2g_1) = 0$ while, by lemma \ref{lem:invariantdeM},
$f_1g_2$ and $f_2g_1$ cannot both be zero. If $f_1g_2 = 0 \neq f_2g_1$, we
must have $\Phi(\xi_1) = \Phi(\xi_2) = 0$; if $f_1g_2 \neq 0 = f_2g_1$, we must
have $\Phi(\xi_1) = \Phi(\xi_2) = \infty$. We consider the former case only, the
latter being entirely similar.
Thus we have $\Phi^{-1}(0) = \left\{\overline{\xi_1},\overline{\xi_2}\right\}$.

\paragraph{Subcase IIa: $\xi_1 \neq \xi_2$.}

If $f_1 = 0$, we take, for some $\alpha \in \C$ to be determined:
$$
P = \begin{pmatrix}
\thq\left(\frac{\xi_1}{\rho_1}x\right) & 0 \\
\alpha \thq\left(\frac{\xi_1}{\rho_2}x\right) &
\thq\left(\frac{\xi_2}{\rho_2}x\right)
\end{pmatrix}.
$$
Then $\det P =
\thq\left(\frac{\xi_1}{\rho_1}x\right) \thq\left(\frac{\xi_2}{\rho_2}x\right)$
does not vanish identically; moreover:
$$
P(x_1) = \begin{pmatrix} 0 & 0 \\ \alpha \ast & \ast \end{pmatrix}
\quad \text{~and~} \quad
P(x_2) = \begin{pmatrix} \ast & 0 \\ \alpha \ast & 0 \end{pmatrix},
$$
where each $\ast$ stands for some non zero complex number. Both determinants
vanish, as required. We can take as $\begin{pmatrix} p_1 \\ q_1 \end{pmatrix}$
the right column of $P(x_1)$, which is indeed non zero and proportional to
$\begin{pmatrix} f_1 \\ g_1 \end{pmatrix} =
\begin{pmatrix} 0 \\ \ast \end{pmatrix}$.
Likewise, we can take as
$\begin{pmatrix} p_2 \\ q_2 \end{pmatrix}$ the left column of $P(x_2)$, which
is indeed non zero and can be made proportional to
$\begin{pmatrix} f_2 \\ g_2 \end{pmatrix} =
\begin{pmatrix} \ast \\ g_2 \end{pmatrix}$ by an appropriate choice of
$\alpha$. This terminates the proof in this case. \\

If $g_2 = 0$, we take, for some $\alpha \in \C$ to be determined:
$$
P = \begin{pmatrix}
\thq\left(\frac{\xi_1}{\rho_1}x\right) & 
\alpha \thq\left(\frac{\xi_2}{\rho_1}x\right) \\
0 & \thq\left(\frac{\xi_2}{\rho_2}x\right)
\end{pmatrix}.
$$
We leave to the reader to complete the argument in this case.

\paragraph{Subcase IIb: $\xi_1 = \xi_2 =: \xi$.}

Then $\xi \equiv -\dfrac{\rho_1}{\sigma_1} \equiv -\dfrac{\rho_2}{\sigma_2}$
and $0$ is a critical value\footnote{In the case $f_2 g_1 = 0$ and $\xi_1 = \xi_2$,
we would have
$\xi \equiv -\dfrac{\rho_1}{\sigma_2} \equiv -\dfrac{\rho_2}{\sigma_1}$ and
the critical value would be $\infty$.} of $\Phi$. If $f_1 = 0$, one has
$M(x_1) = \begin{pmatrix} 0 & 0 \\ \ast & \ast \end{pmatrix}$ and we want $p_1 = 0$.
If $g_2 = 0$, one has
$M(x_2) = \begin{pmatrix} \ast & \ast \\ 0 & 0 \end{pmatrix}$ and we want $q_2 = 0$. \\

\underline{IIb (i): $f_1 = g_2 = 0$.}
We take $C := \Diag(\xi,\xi)$ and
$P := \Diag\left(\thq\left(\frac{\xi}{\rho_1}x\right),
\thq\left(\frac{\xi}{\rho_2}x\right)\right)$.
The right column of $P(x_1)$ is $\begin{pmatrix} p_1 \\ q_1 \end{pmatrix}$ and
the left column of $P(x_2)$ is $\begin{pmatrix} p_2 \\ q_2 \end{pmatrix}$. \\

\underline{IIb (ii): $f_1 = 0$, $g_2 \neq 0$.}
We take $C := \begin{pmatrix} \xi & \xi \\ 0 & \xi \end{pmatrix}$. We take
$P := \begin{pmatrix}
\alpha_{1,1} \phi_1 & \alpha_{1,1} \psi_1 + \alpha_{1,2} \phi_1 \\
\alpha_{2,1} \phi_2 & \alpha_{2,1} \psi_2 + \alpha_{2,2} \phi_2
\end{pmatrix}$
for some $\alpha_{i,j} \in \C$, $i,j = 1,2$ with the usual notations for $\phi_i$
(see subsection \ref{subsubsection:possibleforms}). We have here
$\phi_1(x_1) = \phi_2(x_2) = 0$ and, as a consequence,
$\phi_1'(x_1), \phi_2'(x_2), \phi_2(x_1), \phi_1(x_2) \neq 0$.
Then $P(x_1) = \begin{pmatrix} 0 & 0 \\ \ast & ? \end{pmatrix}$, its left column
is $\begin{pmatrix} p_1 \\ q_1 \end{pmatrix}$, indeed colinear with 
$\begin{pmatrix} f_1 \\ g_1 \end{pmatrix}$. \\
Also, setting $\alpha_{1,1} := 0$ and $\alpha_{2,1} := 1$, one has
$P(x_2) = \begin{pmatrix} 0 & \alpha_{1,2} \phi_1(x_2) \\ 0 & \ast \end{pmatrix}$, its
right column can be chosen as $\begin{pmatrix} p_2 \\ q_2 \end{pmatrix}$ and made
colinear to $\begin{pmatrix} f_2 \\ g_2 \end{pmatrix}$ by an appropriate choice
of $\alpha_{1,2}$. \\

\underline{IIb (iii): $f_1 \neq 0$, $g_2 = 0$.} We leave to the reader to find
the argument in this case (symmetric of the previous one). \\

This ends the proof of the theorem.



\section{The Jimbo-Sakai family (II)}
\label{section:JimboSakaiFamilyII}

We first recall our general assumptions on the Jimbo-Sakai family, as described
in section \ref{subsection:Defandassump} and completed in section
\ref{subsection:Setting, general facts}. The local data are:
$$
R := \Diag(\rho_1,\rho_2), \quad S := \Diag(\sigma_1,\sigma_2) \quad
\text{~and~} \quad \x := \{x_1,x_2,x_3,x_4\}.
$$
We assume \emph{Fuchs relation} (FR), \emph{strong non resonancy} (NR) and add
\emph{non splitting} (NS) with respect to the selected pair 
$\x' := \{x_1,x_2\}$. We also write $\x'' := \{x_3,x_4\}$. \\

Our goal here is to give a geometric description of the \emph{monodromy data 
space} $\F := \F_{R,s,\x}$ underlying the Jimbo-Sakai approach to the study of 
the discrete Painlev\'e equation $q$-PVI. We gave such a (crude) description
in section \ref{section:JimboSakaiFamilyI}. In this section, we obtain a more
general and more precise description, using for that the Mano decomposition 
studied in section \ref{section:ManoDecomposition}. \\

In hope that the reader doesn't get lost in the maze of computations of cases
and subcases, here is a brief summary of the process. In short, $\F$ will be
considered as fibered over the base consisting of all possible factors $C$ in
the Mano decomposition:
\begin{enumerate}
\item{The space $X$ of possible matrices $C$ is that of all those $C \in \GL_2(\C)$
such that $\det C$ is compatible with the prescribed local data:
$\det C \equiv \dfrac{\rho_1 \rho_2}{x_1 x_2} \equiv \sigma_1 \sigma_2 x_3 x_4$.}
\item{Since $\F_C$, the fiber over $C$, is really determined by the class of $C$ under
gauge equivalence,
the true base space $B$ of our fibration is a quotient of $X$. We describe it
with the help of normal forms for $C$; generically they are diagonal and unique
up to permutation of the diagonal terms, so we should take in account an involution.}
\item{The space $\F$ is the quotient of $F := F_{R,s,\x}$ by the equivalence relation
$\sim$, which was defined in \ref{subsubsection:CMRHBC}.
Let $F_C$ the subspace of $F$ made up of those matrices admitting Mano
decomposition with factor $C$ and $\F_C := \frac{F_C}{\sim}$. We parameterize
each $F_C$ by some space of complex matrices\footnote{As in all the paper, this
is possible because the spaces of solutions of ``elementary'' $q$-difference
equations come equipped with explicit finite bases.} and then each $\F_C$ can
be described with the help of some (multi-)linear algebra.}
\item{Of course the description of $\F_C$ is not the same in the \emph{generic},
\emph{trivial}\footnote{Actually the trivial case will be excluded \emph{de facto},
see herebelow \ref{subsubsection:notrivialcase}.} and \emph{logarithmic} cases.
But moreover in the generic case there are some \emph{special} values for which
the fiber is not the same as the \emph{general} fiber (for instance fibers related
to those we encountered in \ref{subsubsection:specialfibersofPi}).}
\item{Last we rebuild $\F$ as he union of the fibers $\F_C$.}
\end{enumerate}


\subsection{An assumption on the local data and a preliminary consequence}
\label{subsection:assumptiononlocaldata}

In \ref{subsubsection:specialcases}, we shall be led to introduce two subsets of
the fundamental annulus $\Cq$, related to the special fibers mentioned hereabove.
For $x \in \C^*$, recall (from the $q$-notations in \ref{subsection:q-notations})
that we write $R(x)$ the unique representative of $x$ in $\Cq$ (\ie\ $R(x) \in \Cq$
and $R(x) \equiv x$). Let:
\begin{align*}
\Xi' &:= \{R(-\rho_1/x_1),R(-\rho_1/x_2),R(-\rho_2/x_1),R(-\rho_2/x_2)\}, \\
\Xi'' &:= \{R(-\sigma_1 x_3),R(-\sigma_1 x_4),R(-\sigma_2 x_3),R(-\sigma_2 x_4)\}.
\end{align*}
Then, in order to simplify the exposition, we shall assume\footnote{Our analysis
can easily be extended to other cases, yielding similar though slightly different
geometries.} (see equation \eqref{eqn:Hyp8}) that:
$$
\text{Assumption~} \mathbf{Hyp}_{8}: \quad
\Xi := \Xi' \cup \Xi'' \text{~has eight (pairwise distinct) elements.}
$$


\subsubsection{Our assumption excludes the trivial case}
\label{subsubsection:notrivialcase}

We draw at once a consequence of this assumption: \emph{the trivial case for $C$
(see the third situation described in theorem \ref{thm:existenceMD}, page
\pageref{thm:existenceMD}) cannot occur.}
Indeed, by \ref{subsubsection:possibleforms}, in the case
$C = \begin{pmatrix} \xi & 0 \\ 0 & \xi \end{pmatrix}$, we must have:
$$
P = \begin{pmatrix}
\alpha_{1,1} \thq\left(\frac{\xi}{\rho_1}x\right) &
\alpha_{1,2} \thq\left(\frac{\xi}{\rho_1}x\right) \\
\alpha_{2,1} \thq\left(\frac{\xi}{\rho_2}x\right) &
\alpha_{2,2} \thq\left(\frac{\xi}{\rho_2}x\right)
\end{pmatrix} =
\begin{pmatrix} \thq\left(\frac{\xi}{\rho_1}x\right) & 0 \\
0 & \thq\left(\frac{\xi}{\rho_2}x\right) \end{pmatrix} A, \text{~where~}
A := \begin{pmatrix} \alpha_{1,1} & \alpha_{1,2} \\
\alpha_{2,1} & \alpha_{2,2} \end{pmatrix},
$$
In the same way:
$$
Q = \begin{pmatrix}
\beta_{1,1} \thq\left(\frac{\sigma_1}{\xi}x\right) &
\beta_{1,2} \thq\left(\frac{\sigma_2}{\xi}x\right) \\
\beta_{2,1} \thq\left(\frac{\sigma_1}{\xi}x\right) &
\beta_{2,2} \thq\left(\frac{\sigma_2}{\xi}x\right)
\end{pmatrix} = B
\begin{pmatrix} \thq\left(\frac{\sigma_1}{\xi}x\right) & 0 \\
0 & \thq\left(\frac{\sigma_2}{\xi}x\right) \end{pmatrix},
\text{~where~}
B := \begin{pmatrix} \beta_{1,1} & \beta_{1,2} \\
\beta_{2,1} & \beta_{2,2} \end{pmatrix},
$$
Usual conditions on $P,Q$ imply first that $\det A, \det B \neq 0$; and then,
since $\det P(x_1) = \det Q(x_3) = 0$, that:
$$
\big(\xi \equiv -\rho_1/x_1 \text{~or~}
\xi \equiv -\rho_2/x_1 \big)
\text{~and~}
\big(\xi \equiv -\sigma_1 x_3 \text{~or~}
\xi \equiv -\sigma_2 x_3 \big).
$$
This would imply $\xi \in \Xi' \cap \Xi''$, which, by assumption, is impossible.


\subsubsection{Under our assumption, $0$ and $\infty$ cannot be critical values
  of $\Phi$}

This means that case 3 in theorem \ref{thm:existenceMD} cannot occur at all,
so this statement actually subsumes the statement hereabove about the
impossibility of the trivial case. \\

The argument is the following: if for instance $0$ is a critical value, then
from the definition of $\Phi$ (recalled at the beginning of subsection
\ref{subsection:existenceMD}, just before theorem \ref{thm:existenceMD}),
$-\rho_1/x_1 \equiv - \rho_2/x_2$. Similarly, if $\infty$ is a critical value,
then $-\rho_1/x_2 \equiv - \rho_2/x_1$. Both congruences are excluded by the
fact that all elements of $\Xi'$ are distinct.


\subsection{Fibration of $\F$}
\label{subsection:FibrationdeF}

Recall that $\F = F/\sim$ where $F := F_{R,s,\x}$ and the equivalence relation
$\sim$ is defined by $M \sim M^{(\Gamma, \Delta)} := \Gamma M \Delta^{-1}$, where
$\Gamma$ and $\Delta$ are diagonal matrices. It is induced by the group action
of the product of the group of $2 \times 2$ invertible diagonal matrices by
itself (actually, by the free action of the $3$-dimensional torus
$\dfrac{\Ddc \times \Ddc}{\C^*(I_2,I_2)}$, see \ref{subsection:Falgsurf}).


\subsubsection{Partition of $F$}
\label{subsubsection:partitionF}

Let 
$$
X := \left\{C \in \GL_2(\C) \tq 
\det C \equiv \dfrac{\rho_1 \rho_2}{x_1 x_2}\right\},
$$
the set of possible factors $C$ for Mano decompositions. For all $C \in X$, 
set:
\begin{align*}
F_C' & := \left\{P \in \Mat_2(\Rwg) \TQ 
\begin{cases} \sq P = R P (C x)^{-1}, \\ \det P \neq 0, \\
\det P \text{~vanishes on~} \x' \end{cases} \right\} \\
F_C'' & := \left\{Q \in \Mat_2(\Rwg) \TQ 
\begin{cases} \sq Q = C Q (S x)^{-1}, \\ \det Q \neq 0, \\
\det Q \text{~vanishes on~} \x'' \end{cases} \right\}
\end{align*}
Thus, as usual, all zeroes of $\det P$ and $\det Q$ are simple and located
on $[\x';q]$ and $[\x'';q]$ respectively. \\

For every $C \in X$, we have a well defined product map:
\begin{align*}
F_C' \times F_C'' & \rightarrow F, \\
(P,Q) & \mapsto PQ.
\end{align*}
We call $F_C$ its image: it is the set of those $M \in F$ that admit
a Mano decomposition with factor $C$. By theorem \ref{thm:existenceMD}:
$$
F = \bigcup_{C \in X} F_C.
$$
Also, by proposition \ref{prop:gaugefreedom}, $F_{C_1} = F_{C_2}$ if, and only
if, $C_1$ and $C_2$ are gauge equivalent; otherwise $F_{C_1}$ and $F_{C_2}$
are disjoint. So, writing $\sim$ the restriction of the gauge equivalence
relation to $X$, with a slight abuse of notation (confusing $C$ with its class
in $X/\sim$) we have a partition:
$$
F = \bigsqcup_{C \in \frac{X}{\sim}} F_C.
$$


\subsubsection{Partition of $\F$}

In order to apply the partition of $F$ to the quotient $\F = F/\sim$, we need
to complete proposition \ref{prop:gaugefreedom} (gauge freedom for $C$).

\begin{prop}
Let $C \in N$ and let $(P_1,Q_1), (P_2,Q_2) \in F_C' \times F_C''$. Then:
$$
P_1 Q_1 \sim P_2 Q_2 \Longleftrightarrow
(P_2,Q_2) = (\Gamma P_1 \Lambda^{-1},\Lambda Q_1 \Delta^{-1})
\text{~for some~} \Gamma, \Delta \text{~diagonal and~} \Lambda \in \GL_2(\C)
\text{~such that~} \Lambda C = C \Lambda.
$$
\end{prop}
\Pr
We have $P_2 Q_2 = \Gamma P_1 Q_1 \Delta^{-1}$ for some $\Gamma, \Delta$
diagonal and $(\Gamma P_1) (Q_1 \Delta^{-1})$ is another Mano decomposition
with factor $C$ for $P_2 Q_2$, so proposition \ref{prop:gaugefreedom}
yields some $\Lambda$ such that $\Lambda[C] = C$. Since $C$ is normalized, 
$\Lambda \in \GL_2(\C)$ and $\Lambda C = C \Lambda$.
\Finpr

As a consequence, we define on each $F_C' \times F_C''$, $C \in N$, an
equivalence relation by:
$$
(P_1,Q_1) \sim (P_2,Q_2) \underset{def}{\Longleftrightarrow}
(P_2,Q_2) = (\Gamma P_1 \Lambda^{-1},\Lambda Q_1 \Delta^{-1})
\text{~for some~} \Gamma, \Delta \text{~diagonal and~} \Lambda \in \GL_2(\C)
\text{~such that~} \Lambda C = C \Lambda.
$$

We then deduce from the discussion in \ref{subsubsection:partitionF} that:
$$
\F = \bigsqcup_{C \in \frac{X}{\sim}} \F_C,
$$
where we have a well defined bijection:
$$
\dfrac{F_C' \times F_C''}{\sim} \overset{\sim}{\longrightarrow} \F_C.
$$


\subsubsection{The base space of the fibration}
\label{subsubsection:Basespace}

\emph{A priori}, we should take as base space of our fibration the quotient
$\dfrac{X}{\sim}$ of the set $X \subset \GL_2(\C)$ by the gauge equivalence
relation. However, under the assumption formulated in subsection
\ref{subsection:assumptiononlocaldata} we saw in \ref{subsubsection:notrivialcase}
that the \emph{trivial case} is impossible and reduced scalar matrices may be
excluded. \\

For a more precise description, we shall use the normal forms for $C$ found
in \ref{subsubsection:NormalformsforC}, proposition \ref{prop:NormalformsforC}; 
we define the following subsets of $X$:
\begin{align*}
N_g &:= 
\left\{\begin{pmatrix} \xi_1 & 0 \\ 0 & \xi_2 \end{pmatrix} \TQ
\xi_1,\xi_2 \in \Cq \text{~and~} \xi_1 \neq \xi_2 \text{~and~}
\xi_1 \xi_2 \equiv \dfrac{\rho_1 \rho_2}{x_1 x_2}\right\}, \\
N_t &:= 
\left\{\begin{pmatrix} \xi & 0 \\ 0 & \xi \end{pmatrix} \TQ
\xi \in \Cq \text{~and~}
\xi^2 \equiv \dfrac{\rho_1 \rho_2}{x_1 x_2}\right\}, \\
N_u &:= 
\left\{\begin{pmatrix} \xi & \xi \\ 0 & \xi \end{pmatrix} \TQ
\xi \in \Cq \text{~and~}
\xi^2 \equiv \dfrac{\rho_1 \rho_2}{x_1 x_2}\right\}, \\
N &:= N_g \sqcup N_t \sqcup N_u, \\
N^* &:= N_g \sqcup N_u.
\end{align*}
Every element of $X$ is equivalent to a element of $N$ and the only possible
non trivial equivalences within $N$ are relations within $N_g$ of the form
$\Diag(\xi_1,\xi_2) \sim \Diag(\xi_2,\xi_1)$. Factors $C \in N_g$, resp.
$C \in N_t$, resp. $C \in N_u$ correspond to what we called the generic, resp.
the trivial, resp. the logarithmic case. Since any element in the group $\Eq$
has exactly $4$ square roots (recall that we use the multiplicative notation),
$\card\ N_t = \card\ N_u = 4$. \\

Our true base space ($N_t$ being excluded by \ref{subsubsection:notrivialcase})
is therefore the quotient:
$$
B := \dfrac{N^*}{\sim},
$$
where $\sim$ is the relation induced by the involution trivial on $N_t$ and
defined on $N_g$ as $\Diag(\xi_1,\xi_2) \mapsto \Diag(\xi_2,\xi_1)$. Sending
$\Diag(\xi_1,\xi_2)$ to $\xi_1$ and
$\begin{pmatrix} \xi & \xi \\ 0 & \xi \end{pmatrix}$ to $\xi$ defines a bijection
of $N^*$ to $\Cq$, hence to $\Eq$. The corresponding involution on $\Eq$ is the map:
$$
\alpha \mapsto \dfrac{a}{\alpha},
$$
where $a$ is the class of $\dfrac{\rho_1 \rho_2}{x_1 x_2}$. We already met
this involution in \ref{subsubsection:involutionofEq}.

\begin{prop}
\label{prop:Basespace}
As a quotient holomorphic curve, the base space $B$ is isomorphic to the projective
line $\PC$. The quotient map is realized by $\Phi$, \ie\ it is the mapping:
$$
N^* \rightarrow \PC, \quad
\begin{pmatrix} \xi & \ast \\ \ast & \ast \end{pmatrix} \mapsto \Phi(\xi).
$$
\end{prop}
\Finpr

Let $M \in F$, let $\xi_1,\xi_2$ the eigenvalues of the factor $C$ in the Mano
decomposition of $M$ and $P$ the left factor. The projective invariant $\Pi(M)$
was defined after lemma \ref{lem:invariantdeM}. \\

In the discussion at the beginning of \ref{subsection:existenceMD}, we used lemma
\ref{lem:autrefractioncroisee}, according to which:
$$
\Phi^{-1}\left(\dfrac{p_1 q_2}{p_2 q_1}\right) =
\left\{\overline{\xi_1},\overline{\xi_2}\right\},
$$
where the $p_i,q_i$ are related to $P$ as explained there. Then, in proposition
\ref{prop:lesdeuxfractionscroisees}, we found that:
$$
(p_1 q_2 : p_2 q_1) = (f_1 g_2 : f_2 g_1),
$$
where $(f_1 g_2 : f_2 g_1) = \Pi(M)$. Combining these facts, we find that:
$$
\Pi(M) = \Phi(\xi_1) = \Phi(\xi_2).
$$
Using proposition \ref{prop:Basespace} above, we can now recognize the true role
of the projective invariant $\Pi$:

\begin{thm}
\label{thm:fibrationavecPi}
The fibration obtained from $F \rightarrow X$ when going to the quotient is
(up to natural identifications):
$$
\F \overset{\Pi}{\longrightarrow} \PC.
$$
\end{thm}
\Pr
Recall from \ref{subsubsection:partitionF} the partition
$F = \bigsqcup\limits_{C \in \frac{X}{\sim}} F_C$, which allows us to define a map
$F \rightarrow \dfrac{X}{\sim} \cdot$ The identification of $\dfrac{X}{\sim}$
with $\PC$ is provided by the map $\Phi$. We get a commutative diagram:
$$
\xymatrix{
F \ar@<0ex>[rr] \ar@<0ex>[d] & & X/\sim \ar@<0ex>[d]^\Phi \\
\F \ar@<0ex>[rr]_\Pi & & \PC
}
$$
The lower horizontal line is thus identified with the fibration.
\Finprcourt


\subsection{Description of $\F_C$ in the generic case}
\label{subsection:DescriptionF_C}

As can be guessed, we have to distinguish the generic and logarithmic cases for $C$
(the trivial case has been excluded, see \ref{subsection:assumptiononlocaldata}). \\

Let $C = \Diag(\xi_1,\xi_2)$, where $\xi_1,\xi_2 \in \Cq$, $\xi_1 \neq \xi_2$
and $\det C = \xi_1 \xi_2 \equiv \dfrac{\rho_1 \rho_2}{x_1 x_2}$.
The matrices $\Lambda \in \GL_2(\C)$ commuting with $C$ are the diagonal
matrices.


\subsubsection{Spaces of matrices}

The elements of $F_C'$ are the matrices:
$$
P = \begin{pmatrix}
\alpha_{1,1} \thq\left(\frac{\xi_1}{\rho_1}x\right) &
\alpha_{1,2} \thq\left(\frac{\xi_2}{\rho_1}x\right) \\
\alpha_{2,1} \thq\left(\frac{\xi_1}{\rho_2}x\right) &
\alpha_{2,2} \thq\left(\frac{\xi_2}{\rho_2}x\right)
\end{pmatrix} \text{~for some~} \alpha_{i,j} \in \C, i,j = 1,2,
\text{~such that~} \begin{cases}
(\alpha_{1,1} \alpha_{2,2},\alpha_{1,2} \alpha_{2,1}) \neq (0,0), \\
\det P(x_1) = \det P(x_2) = 0. \end{cases}
$$
Condition $\det P(x_1) = \det P(x_2) = 0$ leads us to introduce the elliptic
function:
$$
\Phi_C'(x) := \dfrac
{\thq\left(\frac{\xi_1}{\rho_1}x\right) \thq\left(\frac{\xi_2}{\rho_2}x\right)}
{\thq\left(\frac{\xi_2}{\rho_1}x\right) \thq\left(\frac{\xi_1}{\rho_2}x\right)},
$$
which is of degree $2$ since there is no cancellation of zeroes (this follows
from (NR) and the fact that we are in the generic case); also $\Phi_C'$ is such 
that $\Phi_C'(x_1) = \Phi_C'(x_2)$ (because of the condition on $\det C$). Then:
$$
\left(\det P(x_1) = \det P(x_2) = 0\right) \Longleftrightarrow
\dfrac{\alpha_{1,1} \alpha_{2,2}}{\alpha_{1,2} \alpha_{2,1}} = 
\dfrac{1}{\Phi_C'(x_i)}, i = 1,2.
$$
Said equality is understood to hold in $\PC = \C \cup \{\infty\}$. We are thus
led to set:
$$
s := \dfrac{1}{\Phi_C'(x_1)} = \dfrac{1}{\Phi_C'(x_2)} \in \C \cup \{\infty\}
$$
and to define the following spaces of matrices; first:
$$
\Mat_2(\C)^* := \{A := (\alpha_{i,j}) \in \Mat_2(\C) \tq
(\alpha_{1,1} \alpha_{2,2},\alpha_{1,2} \alpha_{2,1}) \neq (0,0)\},
$$
which we endow with a mapping:
\begin{align*}
\alpha: \Mat_2(\C)^* &\rightarrow \PC = \C \cup \{\infty\}, \\
A := (\alpha_{i,j}) &\mapsto 
\dfrac{\alpha_{1,1} \alpha_{2,2}}{\alpha_{1,2} \alpha_{2,1}} \cdot
\end{align*}
Second:
$$
\Mat_2(\C)_C' := \{A := (\alpha_{i,j}) \in \Mat_2(\C)^* \tq \alpha(A) = s\}.
$$
We then have a bijection:
\begin{align*}
\Mat_2(\C)_C' & \rightarrow F_C', \\
A := (\alpha_{i,j}) & \mapsto 
P := \left(\alpha_{i,j} \thq\left(\frac{\xi_j}{\rho_i}x\right)\right).
\end{align*}
We now observe by direct computation (or reasoning on line and columns)
that, if $A \mapsto P$ and if $\Gamma,\Lambda \in \GL_2(\C)$ are diagonal, 
then $\Gamma A \Lambda^{-1} \mapsto \Gamma P \Lambda^{-1}$. \\

A similar study about the right factor $Q \in F_C''$ in the Mano decomposition
leads us to introduce the elliptic function:
$$
\Phi_C''(x) := \dfrac
{\thq\left(\frac{\sigma_1}{\xi_1}x\right) 
\thq\left(\frac{\sigma_2}{\xi_2}x\right)}
{\thq\left(\frac{\sigma_2}{\xi_1}x\right) 
\thq\left(\frac{\sigma_1}{\xi_2}x\right)},
$$
which is of degree $2$ and such that $\Phi_C''(x_3) = \Phi_C''(x_4)$; and to set:
$$
t := \dfrac{1}{\Phi_C''(x_3)} = \dfrac{1}{\Phi_C''(x_4)} \in \C \cup \{\infty\}.
$$
We then define the space of matrices:
$$
\Mat_2(\C)_C'' := \{B := (\beta_{i,j}) \in \Mat_2(\C)^* \tq \alpha(B) = t\}.
$$
This yields a bijection:
\begin{align*}
\Mat_2(\C)_C'' & \rightarrow F_C'', \\
B := (\beta_{i,j}) & \mapsto 
Q := \left(\beta_{i,j} \thq\left(\frac{\sigma_j}{\xi_i}x\right)\right).
\end{align*}
Again we observe that, if $B \mapsto Q$ and if $\Lambda,\Delta \in \GL_2(\C)$ 
are diagonal, 
then $\Lambda B \Delta^{-1} \mapsto \Lambda Q \Delta^{-1}$.

\begin{prop}
We thereby define a bijection:
$$
\dfrac{\Mat_2(\C)_C' \times \Mat_2(\C)_C''}{\sim} 
\overset{\sim}{\longrightarrow} \dfrac{F_C' \times F_C''}{\sim}
\overset{\sim}{\longrightarrow} \F_C,
$$
where the equivalence relation in the left hand side is defined by 
$(A,B) \sim (\Gamma A \Lambda^{-1},\Lambda B \Delta^{-1})$ for all 
$(A,B) \in \Mat_2(\C)_C' \times \Mat_2(\C)_C''$ and for all invertible
diagonal matrices $\Gamma, \Lambda,\Delta$. The rightmost arrow
was previously defined: it is induced by the multiplication map
$F_C' \times F_C'' \rightarrow \F_C$, $(P,Q) \mapsto PQ$.
\end{prop}
\Finprcourt


\subsubsection{Keeping track of the involution
$C = \Diag(\xi_1,\xi_2) \mapsto \tilde{C} := JCJ = \Diag(\xi_2,\xi_1)$}

Recall that $J = \begin{pmatrix} 0 & 1 \\ 1 & 0 \end{pmatrix}$ and that,
writing $\tilde{C} := J C J$, we have $F_C = F_{\tilde{C}}$ by proposition
\ref{prop:gaugefreedom} and therefore $\F_C = \F_{\tilde{C}}$. We intend to 
produce a ``coordinate'' on $\F_C$ by using the bijection described in
the above proposition. This is not the same if we use $C$ or $\tilde{C}$:
we have to make explicit this involutive relationship. \\

More precisely, we have a bijection:
\begin{align*}
F_C' \times F_C'' &\rightarrow F_{\tilde{C}}' \times F_{\tilde{C}}'', \\
(P,Q) &\mapsto (PJ,JQ)
\end{align*}
yielding the left vertical map of a commutative diagram:
$$
\xymatrix{
F_C' \times F_C'' \ar@<0ex>[d] \ar@<0ex>[r] & F_C \ar@<0ex>[d]^= \\
F_{\tilde{C}}' \times F_{\tilde{C}}'' \ar@<0ex>[r] & F_{\tilde{C}}.
}
$$
Intepreting the multiplications by $J$ in terms of exchanging lines or columns, 
we get another commutative diagram of bijections:
$$
\xymatrix{
\Mat_2(\C)_C' \times \Mat_2(\C)_C'' \ar@<0ex>[d] \ar@<0ex>[r] & 
F_C' \times F_C''  \ar@<0ex>[d] \\
\Mat_2(\C)_C' \times \Mat_2(\C)_C'' \ar@<0ex>[r] & 
F_{\tilde{C}}' \times F_{\tilde{C}}'',
}
$$
where the left vertical map is $(A,B) \mapsto (AJ,JB)$. \\

Write $\tilde{\Lambda} := J \Lambda J$. From the obvious relations:
$$
(\Gamma A \Lambda^{-1}) J = \Gamma (AJ) \tilde{\Lambda}^{-1}
\text{~and~}
J (\Lambda B \Delta^{-1}) = \tilde{\Lambda} (JB) \Delta^{-1},
$$
we deduce that said map $(A,B) \mapsto (AJ,JB)$ is compatible with the 
equivalence relation on $\Mat_2(\C)_C' \times \Mat_2(\C)_C''$. So in the end
we get a commutative diagram:
$$
\xymatrix{
\dfrac{\Mat_2(\C)_C' \times \Mat_2(\C)_C''}{\sim} \ar@<0ex>[d] \ar@<0ex>[r] & 
\dfrac{F_C' \times F_C''}{\sim} \ar@<0ex>[d] \ar@<0ex>[r] & 
\F_C \ar@<0ex>[d]^= \\
\dfrac{\Mat_2(\C)_C' \times \Mat_2(\C)_C''}{\sim} \ar@<0ex>[r] & 
\dfrac{F_{\tilde{C}}' \times F_{\tilde{C}''}}{\sim} \ar@<0ex>[r] & 
\F_{\tilde{C}}
}
$$


\subsubsection{A list of special cases}
\label{subsubsection:specialcases}

We shall describe
$$
\dfrac{\Mat_2(\C)_C' \times \Mat_2(\C)_C''}{\sim}
$$
by looking for normal forms for pairs $(A,B) \in \Mat_2(\C)_C' \times \Mat_2(\C)_C''$.
When all coefficients of $A$ and $B$ are non zero, the answer is simpler, so we 
first discuss here the possibility that some coefficients vanish. We shall only
analyse the case that $\alpha_{1,1} = 0$; the other cases will be expounded
dogmatically, the arguments being entirely similar. \\

If $\alpha_{1,1} = 0$, then $\alpha_{1,2} \alpha_{2,1} \neq 0$. Let $P \in F_C'$ 
the corresponding matrix; then
$$
\det P = -\alpha_{1,2} \alpha_{2,1}
\thq\left(\frac{\xi_2}{\rho_1}x\right) \thq\left(\frac{\xi_1}{\rho_2}x\right)
$$
must vanish at $x_1$ and $x_2$, which means that one of its theta factors
vanishes at $x_1$ and the other at $x_2$. This implies that either
$\xi_1 \equiv -\dfrac{\rho_2}{x_1}$ and $\xi_2 \equiv -\dfrac{\rho_1}{x_2}$ or
$\xi_1 \equiv -\dfrac{\rho_1}{x_2}$ and $\xi_2 \equiv -\dfrac{\rho_2}{x_1}$ .
(We could alternatively use the fact that
$\dfrac{\alpha_{1,1} \alpha_{2,2}}{\alpha_{1,2} \alpha_{2,1}} = 
\dfrac{1}{\Phi_C'(x_1)} = \dfrac{1}{\Phi_C'(x_2)}$.) With the notations of 
\ref{subsection:Setting, general facts} and \ref{subsection:defngaugefreedom},
we also see that either $f_2 = 0$ and $\Phi(\xi_1) = \Phi(\xi_2) = \infty$,
or $f_1 = 0$ and $\Phi(\xi_1) = \Phi(\xi_2) = 0$. \\

Conversely, if for instance
$\xi_1 \equiv -\dfrac{\rho_2}{x_1}$ (and therefore
$\xi_2 \equiv -\dfrac{\rho_1}{x_2}$), we have 
$\thq\left(\frac{\xi_1}{\rho_2} x_1\right) = 0$, whence 
$\det P(x_1) = \alpha_{1,1} \alpha_{2,2}
\thq\left(\frac{\xi_1}{\rho_1}x_1\right) 
\thq\left(\frac{\xi_2}{\rho_2}x_1\right) = 0$.
Under the assumption (NR) and due to the fact that we are in the generic case,
$\thq\left(\frac{\xi_1}{\rho_1}x\right) \thq\left(\frac{\xi_2}{\rho_2}x\right)$
has no common zero with $\thq\left(\frac{\xi_1}{\rho_2}x\right)$. So the above
in turn implies that $\alpha_{1,1} \alpha_{2,2} = 0$. \\

We summarize in the following table the list of all possible cases of 
vanishing of some $\alpha_{i,j}$: \\

\begin{tabular}{|l|l|l|l|l|}
  \hline
  vanishing   & class of $\xi_1$ & class of $\xi_2$ & value of one & 
  $\Phi(\xi_1) = $ \\
  coefficient & $\pmod{q^\Z}$  & $\pmod{q^\Z}$ & of the $f_i$, $g_i$ &
  $\Phi(\xi_2) =$ \\
  \hline
  $\alpha_{1,1} = 0$ & 
  $\xi_1 \equiv -\rho_2/x_1$ & $\xi_2 \equiv -\rho_1/x_2$ &
  $f_2 = 0$ &
  $\infty$ \\
  $\alpha_{1,1} = 0$ & 
  $\xi_1 \equiv -\rho_2/x_2$ & $\xi_2 \equiv -\rho_1/x_1$ &
  $f_1 = 0$ &
  $0$ \\
  \hline
  $\alpha_{2,2} = 0$ & 
  $\xi_1 \equiv -\rho_2/x_1$ & $\xi_2 \equiv -\rho_1/x_2$ &
  $g_1 = 0$ & $\infty$ \\
  $\alpha_{2,2} = 0$ & 
  $\xi_1 \equiv -\rho_2/x_2$ & $\xi_2 \equiv -\rho_1/x_1$ &
  $g_2 = 0$ & $0$ \\
  \hline
  $\alpha_{1,2} = 0$ & 
  $\xi_1 \equiv -\rho_1/x_1$ & $\xi_2 \equiv -\rho_2/x_2$ &
  $f_1 = 0$ & $0$ \\
  $\alpha_{1,2} = 0$ & 
  $\xi_1 \equiv -\rho_1/x_2$ & $\xi_2 \equiv -\rho_2/x_1$ &
  $f_2 = 0$ & $\infty$ \\
  \hline
  $\alpha_{2,1} = 0$ & 
  $\xi_1 \equiv -\rho_1/x_1$ & $\xi_2 \equiv -\rho_2/x_2$ &
  $g_2 = 0$ & $0$ \\
  $\alpha_{2,1} = 0$ & 
  $\xi_1 \equiv -\rho_1/x_2$ & $\xi_2 \equiv -\rho_2/x_1$ &
  $g_1 = 0$ & $\infty$ \\
  \hline
\end{tabular}

\medskip

Just for the few following definitions, we shall, for every $a \in \C^*$, write
$R(a)$ its unique representative in $\Cq$, \ie\ $\{R(a)\} = [a;q] \cap \Cq$. We
introduce the sets of special values:
$$
\Xi'_1 := \{R(-\rho_1/x_1),R(-\rho_1/x_2)\}, \quad
\Xi'_2 := \{R(-\rho_2/x_1),R(-\rho_2/x_2)\}, \quad
\Xi' := \Xi'_1 \cup \Xi'_2 \subset \C^*,
$$
Then we can also express our conditions as:
$$
\xi_1 \in \Xi'_1 \Longleftrightarrow \xi_2 \in \Xi'_2
\Longleftrightarrow \alpha_{1,1} \alpha_{2,2} = 0, \quad
\xi_1 \in \Xi'_2 \Longleftrightarrow \xi_2 \in \Xi'_1
\Longleftrightarrow \alpha_{1,2} \alpha_{2,1} = 0.
$$

There are similar results for the vanishing of the $\beta_{i,j}$, but we do not
tabulate them (although they will be used when necessary). In short, all the
$\beta_{i,j}$ are non zero except maybe if $\xi_1$ is congruent modulo $q^\Z$
to one of the $-\sigma_j x_i$, $i = 3,4$, $j = 1,2$. (Here, one must take in 
account vanishing of such expressions as 
$\thq\left(\frac{\sigma_j}{\xi_1} x_i\right)$.) So we define new sets of special
values:
$$
\Xi''_1 := \{R(-\sigma_1 x_3),R(-\sigma_1 x_4)\}, \quad
\Xi''_2 := \{R(-\sigma_2 x_3),R(-\sigma_2 x_4)\}, \quad
\Xi'' := \Xi''_1 \cup \Xi''_2  \subset \C^*.
$$
We have:
$$
\xi_1 \in \Xi''_2 \Longleftrightarrow \xi_2 \in \Xi''_1
\Longleftrightarrow \beta_{1,1} \beta_{2,2} = 0, \quad
\xi_1 \in \Xi''_1 \Longleftrightarrow \xi_2 \in \Xi''_2
\Longleftrightarrow \beta_{1,2} \beta_{2,1} = 0.
$$

Altogether we get the following set of special values:
$$
\Xi := \Xi' \cup \Xi''.
$$
From now on, we shall assume for simplicity that these eight special values (\ie\
their classes in $\Eq$) are pairwise distinct:
\begin{equation}
\label{eqn:Hyp8}
\text{Assumption~} \mathbf{Hyp}_{8}: \quad
\card\ \Xi = 8.
\end{equation}
Besides, this implies that $\xi_1 \not\equiv \xi_2$ for all $\xi_1 \in \Xi_0$, 
so all the ``general values'' in $\C^* \setminus \Xi_0$ actually fall under 
the generic case being presently studied.


\subsubsection{General fiber}

We assume here that $\xi_1, \xi_2 \not\in \Xi$. Then, for all pairs
$(A,B) \in \Mat_2(\C)_C' \times \Mat_2(\C)_C''$, all the coefficients
$\alpha_{i,j}, \beta_{i,j}$, $i,j = 1,2$ are non zero. \\

Straightforward use of the action of $\Gamma,\Lambda$ allows one to bring $A$ 
to the form $\begin{pmatrix} 1 & 1 \\ 1 & ? \end{pmatrix}$ and the missing
coefficient is necessarily $s$. Then the only possible actions of $\Gamma$,
$\Lambda$ preserving that form are those such that
$\gamma_1/\lambda_1 = \gamma_1/\lambda_2 = \gamma_2/\lambda_1 = 1$, whence
$\Gamma = \Lambda = \lambda I_2$ for some $\lambda \in \C^*$.This means that
we can only act on $B$, while preserving the normal form for $A$, by maps
$B \mapsto \lambda B \Delta^{-1}$. This allows one to bring $B$ to the form
$\begin{pmatrix} 1 & 1 \\ \eta & y \end{pmatrix}$ with the relation $y = t \eta$.
In the end, we obtain a normal form for $(A,B)$:
$$
(A,B) \sim (A_0,D_\eta B_0), \text{~where~}
A_0 := \begin{pmatrix} 1 & 1 \\ 1 & s \end{pmatrix}, \quad
B_0 := \begin{pmatrix} 1 & 1 \\ 1 & t \end{pmatrix}, \quad
D_\eta := \begin{pmatrix} 1 & 0 \\ 0 & \eta \end{pmatrix}.
$$
Here $s,t$ are fixed by the value of the base point $(\xi_1,\xi_2) \in X$ but
$x \in \C^*$ is the one free parameter (indeed, a coordinate) characterizing
the class of $(A,B)$. \\

The action of the involution permutes the columns of $A$ and the lines of $B$.
It sends $(A_0,D_\eta B_0)$ to:
$$
(A_0J,JD_\eta B_0) =
\left(\begin{pmatrix} 1 & 1 \\ s & 1 \end{pmatrix},
\begin{pmatrix} \eta & \eta t \\ 1 & 1 \end{pmatrix}\right) \sim
\left(\begin{pmatrix} 1 & 1 \\ 1 & s^{-1} \end{pmatrix},
\begin{pmatrix} 1 & 1 \\ \eta^{-1} & \eta^{-1} t^{-1} \end{pmatrix}\right) ,
$$
the equivalence being induced by the action of:
$$
(\Gamma,\Lambda,\Delta) := \left(\Diag(1,s^{-1}),I_2,\Diag(\eta,\eta t)\right).
$$
Note that here we consider the fiber over $\tilde{C} = \Diag(\xi_2,\xi_1)$
and it is readily checked that $s$ and $t$ must respectively be replaced
by $s^{-1}$ and $t^{-1}$. We conclude that the coordinate $\eta$ must correspondingly
be replaced by $\eta^{-1}$.


\subsubsection{Special fibers}
\label{subsubsection:specialfibers}

Assume for instance that $\xi_1 \equiv -\rho_1/x_1$ (the other possibilities
are entirely similar). Let $(A,B) \in \Mat_2(\C)_C' \times \Mat_2(\C)_C''$. Then
all the coefficients $\beta_{i,j}$, $i,j = 1,2$ as well as $\alpha_{1,1}$ and
$\alpha_{2,2}$ are non zero; but $\alpha_{1,2} \alpha_{2,1} = 0$ (see the table
in \ref{subsubsection:specialcases}). \\

If $\alpha_{1,2} = 0 \neq \alpha_{2,1}$, use of $\Gamma,\Lambda$ allows one to
reduce $A$ to the form $\begin{pmatrix} 1 & 0 \\ 1 & 1 \end{pmatrix}$. This
pattern can only be preserved by further transformations such that
$\Gamma = \Lambda = \lambda I_2$, so possible transformations of $B$ must have
the form $B \mapsto \lambda B \Delta^{-1}$. This allows one to bring $B$ to the
form $D_\eta B_0$ as before. So pairs $(A,B)$ of this type give rise to one line
in $\F_C$, parameterized by $\C^*$. The case $\alpha_{1,2} \neq 0 = \alpha_{2,1}$
is obviously similar and leads to the same conclusion. \\

Now, if $\alpha_{1,2} = \alpha_{2,1} = 0$, $A$ can be brought to the form $I_2$ but
all pairs $(\Gamma,\Lambda)$ such that $\Gamma = \Lambda$ preserve that form, so
that all actions $B \mapsto \Lambda B \Delta^{-1}$ are allowed and $B$ can be
brought to the form $B_0$. This means that all these pairs define a unique point
in $\F_C$. Since this point is a degeneracy of each of the two punctured $\C^*$
lines found above, the resulting figure (for the special fiber over
$\xi_1 \equiv -\rho_1/x_1$ or for any of the other possibilities) is: two $\C$
lines intersecting at a point. \\

Now we discuss the involution. It replaces $\xi_1 \equiv -\rho_1/x_1$ by
$\xi_1 \equiv -\rho_2/x_2$ and the reduced pair
$\left(\begin{pmatrix} 1 & 0 \\ 1 & 1 \end{pmatrix}, D_\eta B_0\right)$ by
$\left(\begin{pmatrix} 0 & 1 \\ 1 & 1 \end{pmatrix}, J D_\eta B_0\right)$,
which is equivalent to 
$\left(\begin{pmatrix} 0 & 1 \\ 1 & 1 \end{pmatrix},
\begin{pmatrix} 1 & 1 \\ \eta^{-1} & \eta^{-1} t^{-1} \end{pmatrix}\right)$.
So the two $\C^*$-lines above $-\rho_1/x_1$ go isomorphically to the two special
lines above $-\rho_2/x_2$; and obviously, the degenerate point to the degenerate
point. Therefore, after going to the quotient by the involution, each special
fiber still consists in two $\C$ lines intersecting at a point.


\subsection{Description of $\F_C$ in the logarithmic case}
\label{subsection:FCinlogarithmic case}

Let $C := \begin{pmatrix} \xi & \xi \\ 0 & \xi \end{pmatrix}$, where
$\xi \in \Cq$ and $\xi^2 \equiv \dfrac{\rho_1 \rho_2}{x_1 x_2} \cdot$
Matrices $\Lambda \in \GL_2(\C)$ commuting with $C$ are those of the form
$\begin{pmatrix} \lambda & \mu \\ 0 & \lambda \end{pmatrix}$,
where $\lambda \in \C^*$ and $\mu \in \C$. \\

From \ref{subsubsection:possibleforms}, we know that corresponding left
factors in Mano decomposition have the form
$P = \begin{pmatrix}
\alpha_{1,1} \phi_1 & \alpha_{1,1} \psi_1 + \alpha_{1,2} \phi_1 \\
\alpha_{2,1} \phi_2 & \alpha_{2,1} \psi_2 + \alpha_{2,2} \phi_2
\end{pmatrix}$
for some $\alpha_{i,j} \in \C$, $i,j = 1,2$, where
$\phi_i(x) := \thq\left(\frac{\xi}{\rho_i}x\right)$ and
$\psi_i(x) := x \phi_i'(x) =
\frac{\xi}{\rho_i} x \thqp\left(\frac{\xi}{\rho_i}x\right)$,
so that:
$$
\det P = (\alpha_{1,1} \alpha_{2,2} - \alpha_{2,1} \alpha_{1,2}) \phi_1 \phi_2 +
\alpha_{1,1} \alpha_{2,1} x (\phi_1 \phi_2' - \phi_2 \phi_1').
$$
Also condition $\det P \neq 0$ is equivalent to
$(\alpha_{1,1} \alpha_{2,1},\alpha_{1,1} \alpha_{2,2} - \alpha_{2,1} \alpha_{1,2})
\neq (0,0)$. \\

We note that $\phi_1 \phi_2(x_1) \neq 0$. Indeed, if for instance
$\phi_1(x_1) = 0$, then $\xi \equiv -\rho_1/x_1$, so $\xi \in \Xi'$
which we saw was impossible. (Direct argument: the congruence property
on $\xi^2$ would imply $\xi \equiv -\rho_2/x_2$, whence 
$\rho_1/x_1 \equiv \rho_2/x_2$, contrary to the assumption that $\Xi$
has eight pairwise distinct elements.) \\

We also deduce that $\alpha_{1,1} \alpha_{2,1} \neq 0$ because otherwise
the condition $\det\ P(x_1) = 0$ would imply 
$\alpha_{1,1} \alpha_{2,2} - \alpha_{2,1} =0$, then
$(\alpha_{1,1} \alpha_{2,1},\alpha_{1,1} \alpha_{2,2} - \alpha_{2,1} \alpha_{1,2})
= (0,0)$. So in fact a small calculation yields:
$$
\det\ P(x_1) = 0 \Longleftrightarrow
\dfrac{\alpha_{2,2}}{\alpha_{2,1}} - \dfrac{\alpha_{1,2}}{\alpha_{1,1}} = u,
\text{~where~}
u := x_1 \left(\dfrac{\phi_1'}{\phi_1} - \dfrac{\phi_2'}{\phi_2}\right)(x_1).
$$

So here we introduce the space of matrices:
$$
\Mat_2(\C)_C' := \left\{A := (\alpha_{i,j}) \in \Mat_2(\C) \TQ
\begin{cases}
\alpha_{1,1}, \alpha_{2,2} \neq 0, \\
\dfrac{\alpha_{2,2}}{\alpha_{2,1}} - \dfrac{\alpha_{1,2}}{\alpha_{1,1}} = u.
\end{cases} \right\}.
$$
An easy calculation shows that action $P \mapsto \Gamma P \Lambda^{-1}$,
with diagonal $\Gamma$ and $\Lambda$ as shown above, translates to the
similar action on $A$. \\

We now go into the corresponding calculations for the right factors $Q$
in Mano decomposition; they have the form
$Q = \begin{pmatrix}
  \beta_{1,1} \overline{\phi}_1 - \beta_{2,1} \overline{\psi}_1 &
  \beta_{1,2} \overline{\phi}_2 - \beta_{2,2} \overline{\psi}_2 \\
\beta_{2,1} \overline{\phi}_1 & \beta_{2,2} \overline{\phi}_2
\end{pmatrix}$
for some $\beta_{i,j} \in \C$, $i,j = 1,2$, where
$\overline{\phi}_j(x) := \thq\left(\frac{\sigma_j}{\xi}x\right)$ and
$\overline{\psi}_j(x) := x \overline{\phi}_j'(x)$,
so that:
$$
\det Q = (\beta_{1,1} \beta_{2,2} - \beta_{2,1} \beta_{1,2})
\overline{\phi}_1 \overline{\phi}_2 +
\beta_{2,1} \beta_{2,2} \, x
(\overline{\phi}_1 \overline{\phi}_2' - \overline{\phi}_2 \overline{\phi}_1').
$$
In the same way as before, we are led to set:
$$
v := x_3 \left(\dfrac{\overline{\phi}_1'}{\overline{\phi}_1} -
\dfrac{\overline{\phi}_2'}{\overline{\phi}_2}\right)(x_3)
$$
and to define:
$$
\Mat_2(\C)_C'' := \left\{B := (\beta_{i,j}) \in \Mat_2(\C) \TQ
\begin{cases}
\beta_{1,1}, \beta_{2,2} \neq 0, \\
\dfrac{\beta_{2,2}}{\beta_{1,2}} - \dfrac{\beta_{2,1}}{\beta_{1,1}} = v.
\end{cases} \right\}.
$$
Again we find that action $Q \mapsto \Lambda Q \Delta^{-1}$, with diagonal
$\Delta$ and $\Lambda$ as shown above, translates to the similar action
on $B$. So we get a bijection:
$$
\dfrac{\Mat_2(\C)_C' \times \Mat_2(\C)_C''}{\sim}
\overset{\sim}{\longrightarrow} \F_C,
$$
where relation $\sim$ on $\Mat_2(\C)_C' \times \Mat_2(\C)_C''$ is defined by
the action
$(A,B) \mapsto \left(\Gamma A \Lambda^{-1},\Lambda B \Delta^{-1}\right)$
of triples $(\Gamma,\Lambda,\Delta)$, where $\Gamma,\Delta$ are diagonal
invertible and
$\Lambda := 
\begin{pmatrix} \lambda & \mu \\ 0 & \lambda \end{pmatrix}$,
for some $\lambda \in \C^*$ and $\mu \in \C$. It is easily checked that
this action indeed sends $\Mat_2(\C)_C' \times \Mat_2(\C)_C''$ to itself. \\

Action of $\Gamma$ can be used to reduce $A$ to the form
$\begin{pmatrix} 1 & ? \\ 1 & ? \end{pmatrix}$, then action of $\Lambda$
to the form $\begin{pmatrix} 1 & 0 \\ 1 & u \end{pmatrix}$, where the
down right coefficient $u$ is forced upon us by the condition defining
$\Mat_2(\C)_C'$. Then the only possibility to preserve this form is to have
$\Gamma = \Lambda$, a scalar matrix. So the remaining possible actions
on $B$ are $B \mapsto \lambda B \Delta^{-1}$, $\lambda \in \C^*$ and
$\Delta$ diagonal invertible. This can be used to reduce $B$ to the form
$\begin{pmatrix} 1 & 1 \\ x & y \end{pmatrix}$ with condition $y - x = v$.
So the mapping:
$$
x \mapsto \text{~the class of~} \begin{pmatrix} 1 & 1 \\ x & x+v \end{pmatrix}
$$
induces a bijective parameterisation of $\F_C$ by $\C$.


\subsection{Putting it all together}

Recall our assumptions from the beginning of section \ref{section:JimboSakaiFamilyII}:
Fuchs relation (FR), strong non resonancy (NR) and non splitting (NS); to which
we added in \ref{subsubsection:specialcases} assumption $\mathbf{Hyp}_{8}$. \\

There are three components in $\F$, two of which project to finite subsets
of the base:
\begin{enumerate}
\item The logarithmic part $\bigcup\limits_{C \in N_u} \F_C$ is in bijection with
$\Upsilon \times \C$, where $\Upsilon$ is the set of square roots of
$\pi\left(\dfrac{\rho_1 \rho_2}{x_1 x_2}\right) = \pi(\sigma_1 \sigma_2 x_3 x_4)$
in $\Eq$ (so $\card\ \Upsilon = 4$). We simplify the formulation by saying that
``the logarithmic part is $\Upsilon \times \C$'', and similarly for the following.
\item Putting together the special fibers in the generic part and identifying
$\Xi$ by its image in $\Eq$, we have (set theoretically) the quotient of
$\Xi \times (\C^* \sqcup \C^* \sqcup \{.\})$ by the involution. Choosing
a representative subset $\Xi_0$ of $\Xi$ for the involution
$\xi_1 \leftrightarrow \xi_2$ (thus $\card\ \Xi_0 = 4$), we see that this
quotient can be identified with $\Xi_0 \times (\C^* \sqcup \C^* \sqcup \{.\})$.
\item Putting together the general fibers, we have (again set theoretically)
the quotient of $\left(\Eq \setminus (\Xi \cup \Upsilon)\right) \times \C^*$
by the involution. We shall give a closer look at this component in
\ref{subsubsection:algebrogeometricdescription}.
\end{enumerate}


\subsubsection{The fibering: generic part}
\label{subsubsection:fiberinggenericpart}

We must justify our contention that $\Pi: \F \rightarrow \PC$ is a fibration
and that its general (non logarithmic) part has exactly four special fibers.
So we set $\Eqg := \Eq \setminus \Upsilon$ and $\PCg$ its image under $\Phi$,
a projective line minus four points (actually the critical values of $\Phi$),
so that the restriction $\Phi: \Eqg \rightarrow \PCg$ is an \emph{unramified}
degree $2$ covering. We write $\Cqg$ the subset of $\Cq$ corresponding to
$\Eqg$, so that $\pi$ induces a bijection $\Cqg \rightarrow \Eqg$. \\

For $\xi \in \Cqg$, we write $\tilde{\xi} \in \Cqg$ its image under the
involution, \ie\ the unique element of $\Cqg$ such that
$\xi \tilde{\xi} \equiv \rho_1 \rho_2/(x_1 x_2) = \sigma_1 \sigma_2 x_3 x_4$;
and $C(\xi) := \Diag\left(\xi,\tilde{\xi}\right)$. Accordingly, we write
$F_\xi := F_{C(\xi)}$ and $\F_\xi := \F_{C(\xi)}$. \\

Last, in order to have a unified picture, we set:
$$
\Phi'(\xi) := \Phi'_{C(\xi)}(x_1) = \Phi'_{C(\xi)}(x_2) =
\dfrac
{\thq\left(\frac{\xi}{\rho_1}x_1\right)
\thq\left(\frac{\tilde{\xi}}{\rho_2}x_1\right)}
{\thq\left(\frac{\tilde{\xi}}{\rho_1}x_1\right)
\thq\left(\frac{\xi}{\rho_2}x_1\right)}
=
\dfrac
{\thq\left(\frac{\xi}{\rho_1}x_1\right)
\thq\left(\frac{\rho_1}{x_2 \xi}\right)}
{\thq\left(\frac{\xi}{\rho_2}x_1\right)
\thq\left(\frac{\rho_2}{x_2 \xi}\right)}
= \dfrac{\rho_1}{\rho_2} \;
\dfrac
{\thq\left(\frac{x_1}{\rho_1}\xi\right)
\thq\left(\frac{x_2}{\rho_1}\xi\right)}
{\thq\left(\frac{x_1}{\rho_2}\xi\right)
\thq\left(\frac{x_2}{\rho_2}\xi\right)}
$$
and
$$
\Phi''(\xi) := \Phi_{C(\xi)}''(x_3) = \Phi_{C(\xi)}''(x_4) =
\dfrac
{\thq\left(\frac{\sigma_1}{\xi}x_3\right) 
\thq\left(\frac{\sigma_2}{\tilde{\xi}}x_3\right)}
{\thq\left(\frac{\sigma_2}{\xi}x_3\right) 
\thq\left(\frac{\sigma_1}{\tilde{\xi}}x_3\right)}
=
\dfrac
{\thq\left(\frac{\sigma_1}{\xi}x_3\right) 
\thq\left(\frac{\xi}{\sigma_1 x_4}\right)}
{\thq\left(\frac{\sigma_2}{\xi}x_3\right) 
\thq\left(\frac{\xi}{\sigma_2 x_4}\right)}
= \dfrac{\sigma_1}{\sigma_2} \;
\dfrac
{\thq\left(\frac{\xi}{\sigma_1 x_3}\right) 
\thq\left(\frac{\xi}{\sigma_1 x_4}\right)}
{\thq\left(\frac{\xi}{\sigma_2 x_3}\right) 
\thq\left(\frac{\xi}{\sigma_2 x_4}\right)}
\cdot
$$
(We used the Fuchs relation and the functional equation $\thq(1/x) =(1/x) \thq(x)$,
see \ref{subsection:somefunctions}, page \pageref{subsection:somefunctions}.) \\

After \ref{subsection:DescriptionF_C}, we have a cartesian square:
$$
\xymatrix{
\Mat_2(\C)^* \times \Mat_2(\C)^* \ar@<0ex>[d]_{\alpha \times \alpha} &
\underset{\xi \in \Cqg}{\bigsqcup} F_\xi \ar@<0ex>[d] \ar@<0ex>[l] \\
\PCg \times \PCg & 
\Cqg \ar@<0ex>[l]^{\qquad (1/\Phi',1/\Phi'')}
}
$$
Taking in account the $(\Gamma,\Lambda,\Delta)$ action, this gives rise
to a bigger commutative diagram:
$$
\xymatrix{
\dfrac{\Mat_2(\C)^* \times \Mat_2(\C)^*}{(\Gamma,\Lambda,\Delta)-\text{action}}
\ar@<0ex>[d]_{\alpha \times \alpha} & &
\underset{\xi \in \Cqg}{\bigsqcup} \F_\xi
\ar@<0ex>[d] \ar@<0ex>[ll] \ar@<0ex>[rrr]^{\text{quotient by involution}} & & &
\F \ar@<0ex>[d]^\Pi \\
\PCg \times \PCg & &
\Cqg \ar@<0ex>[ll]^{\qquad (1/\Phi',1/\Phi'')} \ar@<0ex>[rrr]_\Phi & & &
\PC
}
$$
The left hand square is cartesian, the right hand square is only commutative. \\

Now it follows that the fibers in the generic part have the following form:
$$
\Pi^{-1}(\Phi(\xi)) = \dfrac{\F_\xi \sqcup \F_{\tilde{\xi}}}{\text{involution}},
$$
and the case-by-case computation in \ref{subsection:DescriptionF_C} says that
the fiber is degenerate (``special case'') if, and only if
$\Phi'(\xi) \in \{0,\infty\}$ or $\Phi''(\xi) \in \{0,\infty\}$, that is, after
the above computations, if $\xi \in \Xi'$ or if $\xi \in \Xi''$ respectively.
Taking in account the involution $\xi \leftrightarrow \tilde{\xi}$, for which
each of $\Xi'$, $\Xi''$ is invariant, this means that there are \emph{four}
critical values in $\PCg$ giving rise to special fibers $\Pi^{-1}(-)$ of the
form $\C^* \sqcup \C^* \sqcup \{.\}$:
\begin{align*}
\Phi\left(-\rho_1/x_1\right) = \Phi\left(-\rho_2/x_2\right) &= 0, \\
\Phi\left(-\rho_1/x_2\right) = \Phi\left(-\rho_2/x_1\right) &= \infty, \\
\Phi\left(-\sigma_1 x_3\right) = \Phi\left(-\sigma_2 x_4\right) &=
\dfrac
{\thq\left(\frac{\sigma_1}{\rho_1} x_1 x_3\right)
\thq\left(\frac{\sigma_1}{\rho_2} x_2 x_3\right)}
{\thq\left(\frac{\sigma_1}{\rho_1} x_2 x_3\right)
\thq\left(\frac{\sigma_1}{\rho_2} x_1 x_3\right)}, \\
\Phi\left(-\sigma_1 x_4\right) = \Phi\left(-\sigma_2 x_3\right) &=
\dfrac
{\thq\left(\frac{\sigma_2}{\rho_1} x_1 x_3\right)
\thq\left(\frac{\sigma_2}{\rho_2} x_2 x_3\right)}
{\thq\left(\frac{\sigma_2}{\rho_1} x_2 x_3\right)
\thq\left(\frac{\sigma_2}{\rho_2} x_1 x_3\right)}
\cdot
\end{align*}
Other expressions are possible for the last two critical values, but we
found no simple ones.


\subsubsection{More about special fibers}
\label{subsubsection:Moreaboutspecialfibers}

So special fibers correspond to the vanishing of one of the $\alpha_{i,j}$ or one
of the $\beta_{i,j}$; with for instance the conditions $\alpha_{1,1} = 0$ and
$\alpha_{2,2} = 0$ each providing a line of the same special fiber (and of course
these two lines intersect where $\alpha_{1,1} = \alpha_{2,2} = 0$). \\

The first two of the above four special fibers were already met when analyzing
the behaviour of $\Pi$ in \ref{subsubsection:specialfibersofPi}. Indeed, when
encoding the factors $P$ and $Q$ of the Mano decomposition $M = P Q$ by matrices
$A := (\alpha_{i,j})$ and $B := (\beta_{i,j})$, we find that one of the $\alpha_{i,j}$
is zero if, and only if one coefficient of $P$ vanishes and this, by already
explained arguments, is equivalent to: $M(x_1)$ or $M(x_2)$ has a null line.
We already saw in \ref{subsubsection:A projective invariant} that this is
equivalent to one (at least) of $f_1,f_2,g_1,g_2$ vanishes, \ie\ to
$\Pi(M) \in \{0,\infty\}$. And actually each of the four lines found in
\ref{subsubsection:specialfibers} corresponds to one of those conditions. \\

This means that $M$ belonging to one of the lines of one of the first two special
fibers can be read either on $M(x_1)$ or on $M(x_2)$, independently of each other.
So if we refine our notation and write $\Pi_{1,2}$ for the above $\Pi$ and more
generally $\Pi_{i,j}$ for the one obtaining by using $x_i,x_j$ instead of $x_1,x_2$,
we see that each of these four lines is a line of one of the first two special
fibers of $\Pi_{1,3}$ or $\Pi_{1,4}$ or $\Pi_{2,3}$ or $\Pi_{2,4}$. \\

Now we are going to see that the last two of the four special fibers can be read
on $M(x_3)$ and $M(x_4)$, although not on $\Pi_{3,4}$. Actually, the same line of
argument as above leads to: one of the $\beta_{i,j}$ is zero if, and only if
$M(x_3)$ or $M(x_4)$ has a zero \emph{column}. This leads us to introduce for each
of them (both have rank $1$) one non zero \emph{line}, say $(u_3,v_3)$ and
$(u_4,v_4)$; and to define:
$$
\Pi'_{3,4}(M) := \dfrac{u_3 v_4}{u_4 v_3} \in \PC.
$$
Then we conclude from the above discussion that $M$ is in one of the last two
special fibers if, and only if $\Pi'_{3,4}(M) \in \{0,\infty\}$.

\begin{rmk}
The relation of the $\Pi_{i,j}$ and $\Pi'_{i,j}$ projective invariants is subtle
and interesting in its own right. Matrices $M(x_i)$ have rank one, thus can be
written $C_i \times L_i$, a product of a column by a line matrix, both non zero
and defined up to a non zero scalar factor. We saw in lemma
\ref{lem:mysterieusefractioncroisee} of \ref{subsubsection:A projective invariant}
that generically $\Pi_{i,j}$ is a complete invariant for the left action of
diagonal matrices on pairs $(C_i,C_j)$. Clearly, $\Pi'_{i,j}$ is a complete
invariant for the right action of diagonal matrices on pairs $(L_i,L_j)$.
\end{rmk}


\subsubsection{Algebro-geometric description of the general component}
\label{subsubsection:algebrogeometricdescription}

The space of interest is the quotient of
$\left(\Eq \setminus (\Xi \cup \Upsilon)\right) \times \C^*$ by the
the involution $(\overline{\xi_1},\eta) \mapsto (\overline{\xi_2},\eta^{-1})$,
where $\overline{\xi_1} \overline{\xi_2} = a$, the particular class
written above (recall that we write multiplicatively the group law
in $\Eq$). We extend this involution to $\Eq \times \C^*$. Also after
choosing a particular square root $\alpha$ of $a$, we can instead
use a parameter $t \in \Eq$ such that $\overline{\xi_1} = \alpha t$
and $\overline{\xi_2} = \alpha t^{-1}$. So we must find the quotient
of $\Eq \times \C^*$ by the involution $\tau: (t,\eta) \mapsto (t^{-1},\eta^{-1})$. \\

In the usual projective model of $\Eq$, the point at infinity is its
own inverse. So it makes sense to restrict the involution to the affine
algebraic set $\Eqp \times \C^*$, where $\Eqp := \Eq \setminus \{\infty\}$.
For the latter we have an algebro-geometric model:
$$
\Eqp = \text{Spec}\ \C[x,y] \text{~where~} \C[x,y] := \dfrac{\C[X,Y]}{Y^2 - f(X)},
$$
for some separable cubic polynomial $f(X)$. The inversion map on $\Eqp$
is dual to the automorphism of $\C[x,y]$ defined by $y \mapsto -y$. In
this model (the indeterminate $z$ denoting a ``coordinate'' for $\eta$):
$$
\Eqp \times \C^* = \text{Spec}\ \C[x,y][z,1/z]
$$
and the involution is dual to the automorphism of $\C[x,y][z,1/z]$
defined by $y \mapsto -y$, $z \mapsto 1/z$. \\

The quotient of a complex affine algebraic variety by a finite group
(here the group generated by the involution) is obtained by computing
its affine algebra as the subalgebra fixed by the dual action of the
group. So:
$$
\dfrac{\Eqp \times \C^*}{\text{involution~} \tau} =
\text{Spec}\ \C[x,y][z,1/z]^\tau.
$$
We proceed to compute the invariant subalgebra $\C[x,y][z,1/z]^\tau$.
An element of $\C[x,y][z,1/z]$ can be uniquely written as
$$
g = \sum_{n \in \Z} a_n(x) z^n + y \sum_{n \in \Z} b_n(x) z^n,
\text{~where all the~} a_n,b_n \in \C[X].
$$
Invariance by $\tau$ translates into:
$$
g = a_0(x) + \sum_{n \geq 1} a_n(x) (z^n + z^{-n}) + 
y \sum_{n \geq 1} b_n(x) (z^n - z^{-n}).
$$
Setting $w := \dfrac{z + z^{-1}}{2}$ and $v := \dfrac{z - z^{-1}}{2} y$,
we get the form:
$$
g = A(x,w) + v B(x,w).
$$
Clearly, $x$ and $w$ are algebraically independent and $v \not\in \C[x,w]$.
Also:
$$
v^2 = \dfrac{(z - z^{-1})^2}{4} y^2 = (w^2 - 1) f(x).
$$
This describes an algebraic surface, a degree $2$ covering of the plane
$\text{Spec}\ \C[x,w]$ ramified over the set of equation $(w^2 - 1) f(x) = 0$
(a union of six lines). \\

Since for any fixed $v_0 \neq \pm 1$ the equation $w^2 = f(x) (v_0^2 - 1)$
defines an affine elliptic curve, our surface is also a pencil of elliptic
curves parameterized by $\C \setminus \{+1,-1\}$. \\

To recover the \emph{projective} picture from the above \emph{affine} one,
just note that each fixed point $\alpha$ of the involution, \ie\ each
$\alpha = \overline{\xi}$ with $\xi^2 \equiv \dfrac{\rho_1 \rho_2}{x_1 x_2}$,
gives rise to such an affine chart. This will be done in some detail in section
\ref{section:Geometrysurgeryandpants}. Herebelow we attempt at a geometric
description of the projective surface.


\subsection{Geometric description of the whole of $\F$}


\subsubsection{Definitions}

We recall the classical definitions of an \emph{elliptic surface}.
The first one \cite{Beau} is the following.

\begin{defn}
\label{defellfib1}
Let $S$ be a complex projective surface. We will say that it is an elliptic surface
if there exists a smooth curve $B$ and a surjective morphism $p:S\rightarrow B$ whose
generic fiber is an elliptic curve.
\end{defn}

In fact a variant of this definition is better for our purposes \cite{SchSh}
(definition 3.1, page 7).

\begin{defn}
\label{defellfib2}
An elliptic surface $S$ over $B$ is a smooth projective surface S 
with an elliptic fibration over $B$ , i.e. a surjective morphism 
$p: S \rightarrow B$, such that~:
\begin{itemize}
\item[(i)]
almost all fibers are smooth curves of genus $1$; 
\item[(ii)]
no fiber contains an exceptional curve of the first kind. 
\end{itemize}
\end{defn}

It is better to say ``smooth curves of genus $1$" than ``elliptic curve". An elliptic
curve is a smooth curve of genus $1$ with a marked point and definition
\ref{defellfib1} could suggest that there exists a section. \\

We recall that an exceptional curve of the first kind is a smooth rational curve of
self-intersection $-1$ \big(also called $(-1)$-curve\big). Naturally, $(-1)$-curves
occur as exceptional divisors of blow-ups of surfaces at smooth points. One can
always successively blow-down $(-1)$-curves to reduce to a smooth minimal model.
Therefore an elliptic fibration in the sense of the first definition \ref{defellfib1}
can be transformed by a succession of blow-downs into an elliptic fibration in the
sense of the second definition \ref{defellfib2}. \\

A \emph{section} of an elliptic surface $p: S \rightarrow B$ is a morphism 
$s: B\rightarrow S$  such that $p \circ s= \Id_B$. An elliptic surface
does not necessarily admit a section.

    
\subsubsection{An elliptic fibration. Algebraic charts on $\F$}
\label{subsubellipticfib}

We set $\mathcal{Y}:=\dfrac{\Eq\times \C^*}{\text{involution~} \tau}$.
(recall that the involution $\tau$ was defined at the beginning of
\ref{subsubsection:algebrogeometricdescription}). 
We have two maps $p: \mathcal{Y} \rightarrow \C$ and
$\Psi: \mathcal{Y} \rightarrow \P^1(\C)$~: the maps induced respectively by
$$
p: (\xi,\eta) \mapsto w := \frac12(\eta+1/\eta) \quad \text{and} \quad
\Psi: (\xi,\eta)\mapsto \Phi(\xi).
$$

\begin{itemize}
\item 
The map $p$ gives an elliptic fibration of $\mathcal{Y}$ with two exceptionnal
fibers above $\pm 1$.
\item
For $\eta\neq \pm 1$, the canonical map $\Eq\times \{\eta\}\rightarrow \mathcal{Y}$
induces an isomorphism of $\Eq$ onto the generic fiber
$p^{-1}(w)$ \big($w=\frac12(\eta+\eta^{-1})$\big).
\item
The map $\Psi$ induces an isomorphism between the exceptional fiber $p^{-1}(1)$
\big(resp. $p^{-1}(-1)$\big) and $\PC$.
\item
For $u\in \Phi(\Upsilon)\subset \PC$, the fiber $\Psi^{-1}(u)$ is parameterized
bijectively by $w\in \C$. We get $4$ sections of the elliptic fibration $p$.
\item
For $u\in \PC\setminus \Phi(\Upsilon)$ we can describe the fiber $\Psi^{-1}(u)$
as two copies of $\C$ glued at the points $1$ and $-1$ on each copy.
\end{itemize}

The surface $\Eq\times \C^*$ is smooth, therefore if $a\in \mathcal{Y}$ is not the image
of a \emph{fixed point} of the involution $\tau$, then $\mathcal{Y}$ is smooth at $a$.
There are $8$ fixed points: $\Upsilon\times \{\pm 1\}$. Their images belong to the union
of the $4$ logarithmic fibers and to the union of the two exceptional fibers
$p^{-1}(\pm 1)$. One verifies that these images are isolated singular
points\footnote{They are rational double points.} of $\mathcal{Y}$. \\

We can extend $\tau$ into an involution $\tilde\tau$ on $\Eq\times \PC$.
We set $\widetilde{\mathcal{Y}}:=\dfrac{\Eq\times \PC}{\text{involution~} \tilde\tau}$
and we extend the elliptic fibration $p$ into an elliptic fibration 
$\tilde p: \widetilde{\mathcal{Y}}\rightarrow \PC$. We extend $\Psi$ into
$\widetilde \Psi$. \\

If we remove from $\Eq$ the set $\Upsilon$, that is the $4$ fixed points of the
involution $\tau:\xi\mapsto \rho_1\rho_2/x_1x_2\xi$, we get 
$\Eqg=\Eq \setminus \Upsilon$. If we remove from $\Eq$ the $8$ points of $\Xi$,
we get $\Eq^{\dag}:=\Eq\setminus \Xi$. We set $\Eqg^\dag:=\Eqg\cap \Eq^{\dag}$.
We recall that the points of $\Xi$ are not fixed by $\tau$, therefore we have
removed $12$ points. The surfaces $\Eqg\times \C^*$, $\Eq^\dag\times \C^*$ and
$\Eq^{\bullet\dag}\times \C^*$ are invariant by $\tau$. \\

We set~:
\[
\mathcal{Y}^{\bullet}:=\dfrac{\Eqg \times \C^*}{\text{involution~} \tau}, ~~
\mathcal{Y}^{\dag}:=\dfrac{\Eq^\dag \times \C^*}{\text{involution~} \tau}, ~
 \mathcal{Y}^{\bullet\dag}:=\dfrac{\Eqg^\dag \times \C^*}{\text{involution~} \tau}\cdot 
 \]
 We have some ``punctured elliptic fibrations"~:
 \[
 p^{\bullet}: \mathcal{Y}^{\bullet}\rightarrow \C, ~~
 p^{\dag}: \mathcal{Y}^{\dag}\rightarrow \C, ~~
p^{\bullet\dag}: \mathcal{Y}^{\bullet\dag}\rightarrow \C
\]
induced (by restriction) by the elliptic fibration 
$p:\mathcal{Y}\rightarrow \C$. \\

We recall that we have an injective analytic map
$\psi:\mathcal{Y}^{\bullet\dag}\rightarrow \F$. Considering the algebraic structure
on $\Eq$ and the $3$ affine charts described in the preceding section, we can
interpret this map (in $3$ different ways) as an \emph{algebraic chart} of the
surface $\F$. The image of this chart misses $12$ lines of $\F$.


\subsubsection{Description of some fibers of $\Pi$}

If we remove from $\PC$ the set $\Pi(\Upsilon)$ (that is the images of the $4$
``logarithmic fibers") we get $\PC^{\bullet}=\PC\setminus \Pi(\Upsilon)$. If we
remove from $\PC$ the $4$ points of $\Pi(\Xi)$ (that is $0$, $\infty$ and two
other points), we get $\PC^{\dag}:=\PC\setminus \Pi(\Xi)$. We set 
$\PC^{\bullet \dag}:=\PC^\bullet\cap \PC^\dag$: we have removed $8$ points. \\

Using the elliptic fibration $p$ we can describe the fibers of $\Pi$ above
$\PC^{\bullet\dag}$ (the generic fibers). \\

Let $u\in \PC^{\bullet\dag}$. Then $\psi$ is defined and we have $\Pi\circ \psi=\Psi$,
hence $\psi$ induces an isomorphism of $\Psi^{-1}(u)$ onto $\Pi^{-1}(u)$. Therefore
we can describe the algebraic curve $\Pi^{-1}(u)$ as two copies of $\C$ glued at
the points $1$ and $-1$ on each copy\footnote{A model is given by two parabolas
in $\C^2$ in general position.}. \\

The  open set $\Psi^{-1}\left(\P^1(\C)^{\bullet}\right)=\mathcal{Y}^{\bullet}$ is smooth
(it does not contain singular points). A fortiori
$\Psi^{-1}\left(\P^1(\C)^{\bullet \dag}\right) =\mathcal{Y}^{\bullet \dag}$ is smooth and 
$\psi(\mathcal{Y}^{\bullet \dag})=\Pi^{-1}\left(\P^1(\C)^{\bullet \dag}\right)\subset \F$
is also smooth. We conjecture that $\F$ is also smooth in a convenient neighborhood
of each logarithmic fiber. (We will return to this question later, \cf\ \ref{subsubqpants}).

\begin{conj}
\label{conjlisslog}
The inverse image $\Pi^{-1}\left(\P^1(\C)^{\dag}\right)\subset \F$ is smooth.
\end{conj}

Otherwise speaking we can add to $\psi(\mathcal{Y}^{\bullet \dag})$
the $4$ logarithmic fibers and we get a smooth open subset of $\F$. \\

When $u\in \P^1(\C)^{\bullet\dag}$ tends to
$u_0\in \Pi((\Upsilon))$, the generic fiber tends is some sense towards a logarithmic
fiber (\cf\ \ref{subsubqpants}). In order to prove the conjecture it seems necessary to
understand more precisely what happens. Blow ups of the $8$ singular points of
$\mathcal{Y}$ could be useful (\cf\ \ref{subsubqpants}). \\

The fibers above $u\in \Pi(\Xi)$ are made of two affine lines intersecting at
one point (for $u=0$ and $u=\infty$ it is the theorem \ref{thm:specialfibersofPi},
for the two other values , see \ref{subsubsection:specialfibers}). \\

It is difficult to describe what happens to the generic fibers $\Pi^{-1}(u)$ when
$u\in \PC^{\bullet\dag}$ tends to $u_1\in \Pi(\Xi)$. It seems that the two points
$\Pi^{-1}(u)\cap p^{-1}(1)$ and the two points $\Pi^{-1}(u)\cap p^{-1}(-1)$ glue into
an unique point: the special point of $\Pi^{-1}(u_1)$ which is the intersection of
the two affine lines. \\

A possible approach for a description is to reparameterize $\C\setminus \{\pm 1\}$
using an affine transform of $\C$ sending $\pm 1$ to $\pm\, \epsilon$. We get a
family of elliptic fibrations above the family of punctured lines
$\left\{\C\setminus \{\pm\, \epsilon\}\right\}_\epsilon$. Then we can try to describe
what happens when $\epsilon\rightarrow 0$. A similar method works perfectly for a
description of the fibration of the cubic surface $\mathcal{S}_{VI}$ analog to the
fibration by $\Pi$: \cf\ below \ref{simplemodelblowup} \emph{``A simple model''}.


\subsubsection{An heuristic description of the fibration by $\Pi$}

For a generic $u\in \PC$, we can describe the algebraic curve
$\widetilde{\Psi}^{-1}(u)$ as two copies of the projective line $\PC$ glued at
the points $1$ and $-1$ on each copy. We will give another description of the
abstract algebraic curve $\widetilde{\Psi}^{-1}(u)$ (with its fibration induced
by $p$) as an algebraic curve of $\left(\PC\right)^2$ (with the fibration induced
by the projection on the first factor). Using this picture we will give herebelow
a simple \emph{heuristic} description of the fibration of $\F$ by $\Pi$ as a family
of curves. \\

Let $\widetilde W$ be a $(2,2)$ curve of $\left(\PC\right)^2$ decomposed into
two $(1,1)$ curves bitangent at two distinct points $A^+$ and $A^{-}$. We denote
$\widetilde{\mathbf{p}}$ the projection of $\widetilde W$ on $\PC$ induced by
restriction of the projection on the first factor $\left(\PC\right)^2\rightarrow \PC$.
Up to a M\"{o}bius transform on the first factor, we can suppose that
$\widetilde{\mathbf{p}}(A^{\pm})=\pm \epsilon$ ($\epsilon\in \C^*$). The fibers of
$\widetilde{\mathbf{p}}$ above a generic point 
$\mathbf{w}\in \PC\setminus \{\pm\, \epsilon\}$ are sets made of two points.
The special fibers above $\pm\, \epsilon$ are one point sets. \\

For $u\in \PC^{\bullet\dag}$, the two fibrations of abstract algebraic curves
$\tilde p:\overline{\Pi^{-1}(u)}\rightarrow \PC$ and
$\widetilde{\mathbf{p}}: \widetilde W\rightarrow \PC$ are isomorphic:
$\C \setminus \{\pm 1\}$ is send bijectively to 
$\PC\setminus \{\pm\, \epsilon\}$. (We leave the verification to the reader.) \\

We set $W:=\widetilde W\setminus \widetilde{\mathbf{p}}^{-1}(\infty)$.
Then the two algebraic fibrations $p:\Pi^{-1}(u)\rightarrow \C$ and
$\mathbf{p}: W\rightarrow \C$ are isomorphic. \\

We consider an algebraic family of fibered curves
$\left(\widetilde W_{\lambda},\widetilde{\mathbf{p}}_{\lambda},\PC\right)_{\lambda\in \PC}$
such that for a generic value of $\lambda$ the pair
$(\widetilde W_{\lambda},\widetilde{\mathbf{p}}_{\lambda})$  is of the type
$(\widetilde W,p)$. We allow for the curve $\widetilde W_{\lambda}$, as an algebraic
curve of $\left(\PC\right)^2$, some degeneracies of the two following types.

\begin{itemize}
\item 
The two $(1,1)$ curves of $\widetilde W_\lambda$ degenerate into a double $(1,1)$
curve. Then the projection  $\widetilde{\mathbf{p}}_{\lambda}$ becomes an isomorphism.
\item
The two $(1,1)$ curves of $\widetilde W_\lambda$ degenerate into two double lines
intersecting at only one point $B$ ($\epsilon \rightarrow 0$, $A^+$ and $A^-$ glue
together into the point $B$).
\end{itemize}



\section{Geometry, surgery and pants}
\label{section:Geometrysurgeryandpants}

We will freely use some notations introduced in section
\ref{subsection:IsomonodromyandPainleve}.


\subsection{The classical geometry of a smooth cubic complex surface and
  the representations of a free group of rank $3$}

There are strong relations between the classical geometry of a smooth complex
cubic surface ($27$ lines, $45$ tritangent planes \dots, \cf\ \cite{Cay}) and some
properties of the representations into ${\text{SL}}_{2}(\C)$ of a free group of
rank $3$. As far as we know these (simple but important) relations remained unnoticed
until recently\footnote{They are due to Martin Klimes, Emmanuel Paul and the second
author, and they are studied in a work in progress on the confluence of the Painlev\'e
equations \cite{Kli,KliPR}.}.


\subsubsection{The geometry of the cubic surfaces $\mathcal{S}(a)$}
\label{subsubgeomcubic}

By the classical theory, if  $\Sigma$ is a smooth projective cubic surface:
\begin{itemize}
\item
$\Sigma$ admits $27$ lines and each line is a $(-1)$-line;
\item 
$\Sigma$ admits $45$ tritangent planes;
\item
the intersection of each tri-tangent plane with $\Sigma$ is the union of $3$ lines
forming a triangle;
\item
each line of $\Sigma$ belongs exactly to $5$ tri-tangent planes.
\end{itemize}

\smallskip

The equation of the projective surface 
$\overline{\mathcal{S}(a)}\subset \mathbf{P}^3(\C)$ in projective coordinates
$(\tilde X_0,\tilde X_t,\tilde X_1,\tilde T)$ is~:
\[
\tilde X_0\tilde X_t\tilde X_1+\tilde X_0^2\tilde T
+\tilde X_t^2\tilde T+\tilde X_1^2\tilde T -A_0\tilde X_0\tilde T^2
-A_t\tilde X_t\tilde T^2-A_1\tilde X_1\tilde T^2+A_\infty\tilde T^3=0.
\]
The plane at infinity $\tilde T=0$ is a tri-tangent plane and its intersection
with the surface is the triangle $\tilde X_0\tilde X_t\tilde X_1=0$. Therefore
the affine cubic surface $\mathcal{S}(a)$ contains exactly $24$ lines. We have
the following description of these lines. \\

Each line at infinity is contained in $4$ tri-tangent planes different from
the plane at infinity. The intersection of such a tri-tangent plane and
$\overline{\mathcal{S}(a)}$ is a triangle, therefore the intersection with
$\mathcal{S}(a)$ is the union of $2$ affine lines with a common point.
Therefore for each line at infinity we get $8$ affine lines on $\mathcal{S}(a)$.
Using the coordinates $X_0,X_t,X_1$ we see that for each $l=0,t,1$
there exists $4$ exceptional values of $X_l$ such that 
$\{X_l=0\}\cap \mathcal{S}(a)$ is the union of $2$ affine lines. \\

Below we will interpret the $24$ lines on $\mathcal{S}(a)$ in terms of representations.

\begin{prop}
\label{martinlines}
Let $a_0,a_t,a_1,a_\infty\in \C$ \emph{arbitrary}. The $24$ lines \emph{distinct
or not} defined in $\C^3$ by the following equations are contained in the cubic
surface $\mathcal{S}(a)\subset \C^3$:
\begin{equation}
\label{martinequa}
X_k=e_ie_j^{-1}+e_je_i^{-1}, \quad e_iX_i+e_jX_j
=a_\infty+e_ie_ja_k,
\end{equation}
\[
X_k=e_ie_j^{-1}+e_je_i^{-1}, \quad e_iX_j+e_jX_i
=a_k+e_ie_ja_\infty,
\]
\[
X_k=e_ie_j+e_i^{-1}e_j^{-1}, \quad X_i+e_ie_jX_j
=e_ja_k+e_ia_\infty,
\]
\[
X_k=e_ie_j+e_i^{-1}e_j^{-1}, \quad X_j+e_ie_jX_i
=e_ja_\infty+e_ia_k,
\]
\[
X_k=e_ke_\infty^{-1}+e_\infty e_k^{-1}, \quad e_\infty X_i+e_kX_j
=a_i+e_ke_\infty a_j,
\]
\[
X_k=e_ke_\infty^{-1}+e_\infty e_k^{-1}, \quad e_kX_i+e_\infty X_j
=a_j+e_ke_\infty a_i,
\]
\[
X_k=e_ke_\infty+e_k^{-1}e_\infty^{-1}, \quad X_i+e_ke_\infty X_j
=e_ka_j+e_\infty a_i,
\]
\[
X_k=e_ke_\infty+e_k^{-1}e_\infty^{-1}, \quad X_j+e_ke_\infty X_i
=e_ka_i+e_\infty a_j.
\]
\end{prop}

\Pr
The result follows immediately from some decompositions of $F(X,a)$ (\cf\ \cite{Kli},
proposition 4.5). We give only one of these decompositions~:
\begin{align*}
F(X:a) &= (X_k-e_ie_j^{-1}-e_je_i^{-1})
(F_{X_k}-X_k+e_ie_j^{-1}+e_je_i^{-1}) \\
& -e_i^{-1}e_j^{-1}(e_iX_i+e_jX_j-a_\infty-e_ie_ja_k)
(e_iX_j+e_jX_i-a_k-e_ie_ja_\infty).
\end{align*}
\Finprcourt


\subsubsection{Reducibility of representations of the free group of rank $2$}

We recall the well known conditions of reducibility for the representations
of a free group of rank $2$ into $\SLdc$ and some classical
results (\cf\ \cite{IKSY}). We denote $\Gamma_2 := \langle u,v \rangle$ the free group of
rank $2$ generated by the letters $u,~v$.

\begin{defn}
A pair of matrices $(M',M'')\in \left(\SLdc\right)^2$ is said
\emph{reducible} if there exists a common (non trivial, non total) invariant
subspace.
\end{defn}

It is equivalent to say that the corresponding representation
$\omega: \Gamma_2 \rightarrow \SLdc$ defined by
$\omega(u) := M'$ and $\omega(v) := M''$ is reducible. \\

Let $\omega:~ \Gamma_2 \rightarrow \SLdc$ be a linear representation.
We set $M':=\omega(u)$ and $M'':=\omega(v)$. We denote  $e'$ and $(e')^{-1}$
(resp. $e''$ and $(e'')^{-1}$) the eigenvalues of $M'$ (resp. $M''$). We denote
$e$ and $e^{-1}$ the eigenvalues of $M:=M'M''$.

\begin{prop}
\label{repr2}
The following assertions are equivalent
\begin{itemize}
\item[(i)]
The representation $\omega$ is reducible.
\item[(ii)]
The pair  $(M',M'')$ is reducible.
\item[(iii)]
We have~: $e=e'e''$ or $e=e'(e'')^{-1}$ or $e=(e')^{-1}e''$ or
$e=(e')^{-1}(e'')^{-1}$.
\item[(iv)]
We have $\Tr~M=e'e''+(e'e'')^{-1}$ or 
$\Tr~M=e'(e'')^{-1}+(e')^{-1}e''$. 
\end{itemize}
\end{prop}
\Pr
\begin{itemize}
\item 
The assertions (i) and (ii) are evidently equivalent. 
\item
If the pair 
$(M',M'')$ is reducible, then there exists a common eigenvector $v\in \C^2$.
Then $M'v=\lambda' v$,  $M''v=\lambda'' v$ and $Mv=\lambda v$, where $\lambda'$,
$\lambda''$ and $\lambda$ are respectively eigenvalues of $M'$, $M''$ and $M$.
Therefore one of the conditions of (iii) is satisfied.
\item
The assertion (iii) implies clearly the assertion (iv).
\item
If $\omega$ is irreducible, then one can prove that~:
\[
\Tr M\neq e'e''+(e'e'')^{-1} ~~ \text{\emph{and}} ~~
\Tr~M\neq e'(e'')^{-1}+(e')^{-1}e''
\]
(\cf\ \cite{IKSY}, (4.2.9), page 83). Therefore if (iv) is satisfied, then $\omega$
is necessarily reducible. This proves the assertion (i).
\end{itemize}
\Finpr

Let $\omega$ be a representation and $M'$, $M''$ as above. We suppose that
$\Tr~M'\neq \pm 2$ and $\Tr~M''\neq \pm 2$. Then $M'$ and $M''$ are diagonalisable.
There allways exists a \emph{mixed basis} $\{v',v''\}$ of $\C^2$ formed by an eigenvector
of $M'$ and an eigenvector of $M''$. In general there are (up to rescaling of the
eigenvectors) $4$ ways one can form such a basis. The $4$ cases of
reducibility of $\omega$ correspond to the degeneracy of one of these $4$ basis. \\

We recall that if $\omega$ is irreducible, then it is determined, up to equivalence,
by the traces of $M'$, $M''$ and $M$ (\cf\ \cite{IKSY}  Theorem 4.2.1, page 80).


\subsubsection{Partial reducibility and lines on $\mathcal{S}(a)$}
\label{subsubpartialreduclines}

We describe a relation between a notion of partial reducibility of a representation
and the lines on the cubic surface $\mathcal{S}(a)$. This relation is apparently
new\footnote{If we replace representations by wild representations, then this relation
can be extended to all the Painlev\'e equations.}, it has been found recently by
M. Klimes, E. Paul and the second author \cite{Kli,KliPR}. \\

As in \ref{subsubiso}, we denote $\Gamma_3:= \langle u_0,u_t,u_1 \rangle$ the free group
of rank $3$
generated by the letters $u_0,~u_t,~u_1$ and we set $u_\infty=u_1^{-1}u_t^{-1}u_0^{-1}$.
Let $\rho: \Gamma_3\rightarrow \SLdc$ be a linear representation.
We set $M_l:=\rho(u_i)$ ($i=0,t,1,\infty$). We denote 
$e_l$ and $e_l^{-1}$ ($i=0,t,1,\infty$) the eigenvalues of $M_i$. \\

We have the following characterizations of smoothness (some are classical and some
are apparently new).

\begin{thm}
\label{lissites}
Let $a\in \C^4$. We suppose $a_l\neq \pm 2$ ($l=0,t,1,\infty$) (non resonance).
The following conditions are equivalent~:
\begin{itemize}
\item[(i)]
The affine cubic surface $\mathcal{S}(a)$ is smooth.
\item[(ii)]
The projective cubic surface $\overline{\mathcal{S}(a)}$ is smooth.
\item[(iii)]
The $24$ lines (\ref{martinequa}) are pairwise distinct.
\item[(iv)]
The $3$ following conditions are satisfied:
\begin{itemize}
\item 
the $4$ numbers built from the $e_l$ ($l=0,t,1,\infty$) 
\begin{equation}
\label{specialdecvaluesexcl0}
e_te_1^{-1}+e_1e_t^{-1}, ~~
e_te_1+e_t^{-1}e_1^{-1},  ~~ 
e_0e_\infty^{-1}+e_\infty e_0^{-1},  ~~ 
e_0e_\infty+e_0^{-1}e_\infty^{-1}
\end{equation}
are pairwise distinct;
\item
the $4$ numbers~:
\begin{equation}
\label{specialdecvaluesexclt}
e_1e_\infty^{-1}+e_\infty e_1^{-1}, ~~
e_1e_\infty+e_1^{-1}e_\infty^{-1},  ~~ 
e_te_0^{-1}+e_0 e_t^{-1},  ~~ 
e_te_0+e_t^{-1}e_0^{-1}
\end{equation}
are pairwise distinct,
\item
the $4$ numbers~:
\begin{equation}
\label{specialdecvaluesexcl1}
e_1 e_\infty^{-1}+e_\infty e_1^{-1}, ~~
e_1e_\infty+e_1^{-1}e_\infty^{-1},  ~~ 
e_te_0^{-1}+e_0 e_t^{-1},  ~~ 
e_te_0+e_t^{-1}e_0^{-1}
\end{equation}
are pairwise distinct.
\end{itemize}
\item[(v)]
We have the $8$ conditions: $e_0\, e_t^{\pm 1}e_1^{\pm 1}e_\infty^{\pm 1}\neq1$
(the $3$ signs are chosen independantly).
\item[(vi)]
If $\rho$ is a representation such that $\Tr~ \rho(u_l)=a_l$ for all $l=0,t,1,\infty$,
then it is irreducible.
\end{itemize}
\end{thm}
\Pr
\begin{itemize}
\item 
A singular point of $\overline{\mathcal{S}(a)}$ is allways contained into $\mathcal{S}(a)$, 
therefore (i) $\Leftrightarrow$ (ii).
\item
The conditions (iv) and (v) are clearly equivalent.
\item
We have (iii) $\Leftrightarrow$ (i).
If the $24$ lines are distinct, then $\overline{\mathcal{S}(a)}$ contains $27$ distinct
lines and therefore it is smooth.
\item
We have (iv) $\Rightarrow$ (iii). If we suppose (iv), then we get $3$ sets of $8$ two
by two distinct lines. Each set corresponds to $4$ distinct values of $X_0$ or $X_t$
or $X_1$. It is easy to check that a line cannot belong to $2$ different such sets.
\item
We have (i) $\Rightarrow$ (v) $\Rightarrow$ (iii). 
If $\mathcal{S}(a)$ is smooth, then (\ref{reduct3}) is impossible, therefore (v)
is true; (iv) and (iii) follows.
\item
We have (ii) $\Leftrightarrow$ (vi). Easy.
\end{itemize}
\Finpr

If we use the parameters $\theta_l$, then the conditions (v) are translated into:
\[
\theta_0\pm\theta_l\pm\theta_1\pm\theta_\infty\in \Z
\]
(\cf\ \cite{Mazz} Theorem 4.1). \\

We suppose that the surface $\mathcal{S}(a)$ is smooth.
For each pair $(l,m)$ of elements of $\{0,t,1,\infty\}$, each of the $2$ planes~:
\begin{equation}
\label{equapartialreduct}
\begin{cases}
X_n=e_le_m+e_l^{-1}e_m^{-1}\\
X_n=e_le_m^{-1}+e_l^{-1}e_m
\end{cases} \text{~~with} ~~ (l,m,n)=
\begin{cases}
(i,j,k) \\
(k,\infty,k)
\end{cases}
\end{equation}
intersects $\mathcal{S}(a)$ at $2$ lines. The resulting $4$ lines correspond to
the reducibility of the pair of matrices $(M_l,M_m)$. \\

More precisely if two matrices $M_l$ and $M_m$ are diagonalizable, for each of them
there exists a pair of invariant subspaces giving rise to a basis. Then there are
in general $4$ possibilities of pairing of invariant subspaces out of which one can
form a \emph{mixed} basis. The cases of reducibility of the pair $(M_l,M_m)$ corresponds
to the degeneracy of (at least) one of these mixed bases. Each of the $4$ lines
corresponds to such a case of degeneracy.

\begin{defn}
  Let $\rho:\Gamma_3\rightarrow \SLdc$. We will say that $\rho$ is
  \emph{partially reducible} if there exists $i,j\in \{0,t,1,\infty\}$, $i\neq j$,
  such that the pair of matrices $\left(\rho(u_i),\rho(u_j)\right)$ is reducible.
\end{defn}

\begin{prop}
\label{propprreprandlines}
Let $a\in \C^4$ arbitrary and $\rho:\Gamma_3\rightarrow \text{SL}_2(\C)$ a representation
such that $\Tr\, \rho(u_l)=a_l$ ($l=0,t,1,\infty$).
\begin{itemize}
\item[(i)]
If the representation $\rho$ is \emph{partially reducible}, then 
its equivalence class belongs to one of the $24$ lines (distinct or not) defined in
proposition \ref{martinlines}.
\item[(ii)]
We suppose that $\mathcal{S}(a)$ is smooth. Then $\rho$ is partially reducible if
and only if its equivalence class belongs to one of the $24$ lines of 
$\mathcal{S}(a)$.
\end{itemize}
\end{prop}
\Pr
The representation $\rho$ is partially reducible if and only if there exists $k=1,2,3$
such that~:
\begin{equation}
\label{specialdecvalues}
X_k=e_ie_j^{-1}+e_ie_j^{-1} ~~ \text{or} ~~
X_k=e_ie_j+e_i^{-1}e_j^{-1}  ~~ \text{or} ~~
X_k=e_ke_\infty^{-1}+e_\infty e_k^{-1}  ~~ \text{or} ~~
X_k=e_ke_\infty+e_k^{-1}e_\infty^{-1}.
\end{equation}
Then assertion (i) follows from proposition \ref{martinlines} and assertion (ii)
follows from theorem \ref{lissites}.
\Finpr

The condition (\ref{specialdecvalues}) appear in various papers in ``$theta$ notation":
\[
\sigma_k\pm \theta_i \pm \theta_j \in 2 \Z, \quad
\sigma_k\pm \theta_k \pm \theta_\infty \in 2 \Z,
\] 
where $X_k=e^{i\pi \sigma_k}+e^{-i\pi \sigma_k}$.
\cf\ \cite{JimboMonodromy} p 1141: condition $(A.3)_{PVI}$ and
\cite{Gu}, footnote $9$, page $85$.

\paragraph{A fibration.} 
\label{fibrationdiff}
We suppose that we are in the "generic case" (i. e. ${\mathcal S}_{VI}(a)$ is smooth). \\

Let $\Pi_0:{\mathcal S}_{VI}(a)\rightarrow \C$, $\Pi_0:(X_0,X_t,X_\infty)\mapsto X_0$.
We recall~:
\[
  {\mathcal S}_{VI}(a)=\left\{(X_0,X_t,X_1)\in \C^3 \tq
  X_0X_tX_1+X_0^2+X_t^2+X_1^2-A_0X_0-A_tX_t-A_1X_1
+A_\infty=0\right\}.
\]
For $c\in \C$, $\Pi_0^{-1}(c)$ is interpreted as an affine conic in the 
$(X_t,X_1)$-plane:
\[
X_t^2+X_1^2+cX_tX_1-A_tX_t-A_1X_1-cA_0+A_\infty=0.
\]
The generic fiber is isomorphic to $\C^*$. The exceptional fibers are of of two types~:
\begin{itemize}
\item 
either $X_0=\pm 2$, then 
$X_t^2+X_1^2\pm 2X_tX_1=(X_t\pm X_1)^2$,
 the fiber is a \emph{parabola} and it is isomorphic to $\C$\, ;
\item
either we are in a partially reductible case, that is in one of the $4$ cases~:
\begin{equation}
X_0=e_te_1^{-1}+e_1e_t^{-1} ~~ \text{or} ~~
X_0=e_te_1+e_t^{-1}e_1^{-1}  ~~ \text{or} ~~
X_0=e_0e_\infty^{-1}+e_\infty e_0^{-1}  ~~ \text{or} ~~
X_0=e_0e_\infty+e_0^{-1}e_\infty^{-1}, 
\end{equation}
then the fiber is \emph{degenerated into two lines}. The intersection of these
two lines is a critical point of $\Pi_0$. Its image is a critical value of $\Pi_0$. \\

If we remove the $6$ exceptional fibers (that is $8$ lines and two curves)
from ${\mathcal S}_{VI}(a)$ then we can parameterize the remaining set by
a Zariski open set of $\C\times \C^*$. Such parameterizations appear in
many papers \cite{JimboMonodromy}, \cite{ILT} \dots We will return to this
question below (\cf\ \ref{subsubpantparclassical}).
In the $q$-case, we will deduce later from the Mano decomposition a similar
parameterization (\cf\ \ref{subsubqpants}). \\

If $c_0$ is a critical value of $\Pi_0$, then we can describe
$\Pi_0^{-1}(\Delta)$ for a small open disc $\Delta\subset \C$ centered at $c_0$.
We blow down one of the lines of
$\Pi^{-1}(c_0)$ into the smooth  surface $\Pi_0^{-1}(\Delta)$ 
(the closure of such a line in $\overline{\mathcal{S}(a)}$ is a $(-1)$-line).
We get a surface analytically diffeomorphic to 
$\Delta\times \C^*$. Therefore $\Pi_0^{-1}(\Delta)$ is 
analytically diffeomorphic to a blow-up of $\Delta\times \C^*$. \\

The projective version of the above description is the following.
We choose one of the lines $\Delta$ of the triangle at infinity of
the surface $\overline{\mathcal{S}_{VI}(a)}$. The family of planes passing
by this line cut the surface along a linear family of cubic curves.
Each one is decomposed into the union of $\Delta$ and a conic curve.
There are $5$ conic curves degenerated into $2$ lines:
they correspond to the $5$ tritangent planes, the plane at infinity
and $4$ other planes. There are two other exceptional curves: when
the conic is tangent to $\Delta$.

\paragraph{A simple model.}
\label{simplemodelblowup}

We will give two descriptions of a hyperbolic paraboloid of $\C^3$
(an affine quadric surface)~:
\begin{itemize}
\item 
by a blow-up of a point in $\C^2$;
\item
by a singular fibration.
\end{itemize}

Let $\mathcal{Q}\subset \C^3$ defined by the equation $Y=XZ$ (hyperbolic
paraboloid). Let $f: \C^3\rightarrow \C$ defined by $f:(X,Y,Z)\mapsto Y$
and $g:\C^3\rightarrow \C^2$ (the $(X,Y)$ plane) defined by
$g:(X,Y,Z)\mapsto \left(X,f(X,Y,Z)=XZ\right)$. We denote
$\Phi$ (resp. $\phi$) the restriction of $f$ (resp. $g$) to
$\mathcal{Q}$. 

\begin{itemize}
\item 
The map $\phi$ is a bijection of 
$\mathcal{Q}\setminus \phi^{-1}(0,0)$ onto $\C^2\setminus \{(0,0)\}$
and $\phi^{-1}(0,0)$ is the line of $\mathcal{Q}$ defined by
$\{(X,Y,Z) \in C^3\vert~X=Y=0\}$. If we parameterize $\mathcal{Q}$
by $(X,Z)$ then $\phi$ is expressed by $(X,T)\mapsto (X,Y=XT)$.
The quadric $\mathcal{Q}$ is a blow-up of $\C^2$ at the point $(0,0)$.
If we blow down into $\mathcal{Q}$ the line $(X=Z=0)$ by $g$ we get $\C^2$.
\item
We can describle $\mathcal{Q}$ by the fibers $\Phi^{-1}(c)$ of $\Phi$.
If $c\neq 0$, then the fiber is the affine conic
$\{(X,c,Z)\in \C^3\vert~X\neq 0,~Z=c/X\}$. If $c=0$, then the fiber is
the union of the two lines $(X=0)$ and $(Z=0)$.
\end{itemize}
\end{itemize}


\subsubsection{Dynamics on ${\mathcal S}_{VI}(a)$}
\label{subsubcubicsurfdyn}

The cubic surface of $P_{VI}$ admits  $3$ polynomial involutions $s_0,~s_t,~s_1$.
These involutions are \emph{anti-symplectic} and they generate a subgroup of
the group of algebraic automorphisms of $\mathcal{S}_{VI}$ which is isomorphic
to $\Z_2*\Z_2*\Z_2$ \cite[Theorem 3.1, page 2948]{CL}. The products
$g_{i,j}:=s_i\circ s_j$ ($i\neq j$) are \emph{symplectic polynomial automorphisms}.
The corresponding automorphisms of the Okamoto variety of initial conditions
obtained by conjugation by $RH$ are the \emph{non-linear monodromies} of PVI
around the singular points $0,1,\infty$ \cite{CL}.


\subsection{Pant decompositions, surgery, parametrizations}


\subsubsection{Pant decompositions of representations of a free group of rank $3$.}

Let $\rho: \Gamma_3\rightarrow \SLdc$ be a linear representation. We use the same
notations as before. We suppose $a$ fixed such that $a_l=\Tr~M_l\neq \pm 2$,
$l=0,t,1,\infty$. We associate to $\rho$ the two representations of $\Gamma_2$~:
\[
\omega_{0,t}:\Gamma_2\rightarrow \text{SL}_2(\C) ~~\text{and} ~~
\omega_{1,\infty}:\Gamma_2\rightarrow \text{SL}_2(\C) 
\]
We will say that $(\omega_{0,t},\omega_{1,\infty})$ is the \emph{pant decomposition}
of $\rho$ associated to the partition $\{0,t,1,\infty\}=\{0,t\}\cup \{1,\infty\}$.
There are two others pant decompositions associated to the two other partitions. \\

We denote $\tilde{\rho}$ and $\tilde{\omega}$ the equivalence classes of representations.
The pair $(\tilde{\omega}_{0,t},\tilde{\omega}_{1,\infty})$ depends only on $\tilde{\rho}$. \\

We suppose that $\omega_{0,t}$ and $\omega_{1,\infty}$ are irreducible, equivalently
$\tilde \rho$ belongs to $\mathcal{S}(a)$ minus the $8$ critical lines. Then the
knowledge of $(\tilde{\omega}_{0,t},\tilde{\omega}_{1,\infty})$ is equivalent to the
knowledge of $X_0=\Tr~M_0M_t=\Tr~M_1M_\infty$ ($a$ is fixed). \\

In order to recover $\tilde \rho$ from $X_0$, we need another parameter.
We set $X_0=e+e^{-1}$. There are two different cases.
\begin{itemize}
\item 
We suppose $X_0\neq \pm 2$ (equivalently $e\neq e^{-1}$). 
We can choose representations $\omega'$ and $\omega''$ of $\Gamma_2$ such that~:
\[
\Tr~\omega'(u)=a_0, ~~ \Tr~\omega'(v)=a_t, ~~
\omega'(uv)=\Diag(e,e^{-1}),
\]
\[
\Tr~\omega''(u)=a_\infty, ~~ \Tr~\omega''(v)=a_1, ~~
\omega''(uv)=\Diag(e^{-1},e).
\]
There is a freeness in the choice: we can replace $\omega'$ (resp. $\omega''$)
by an overall conjugate by an arbitrary matrix commuting with $\Diag(e,e^{-1})$,
that is of the form $\Diag (t,t^{-1})$ with $t\in \C^*$
\item
We suppose $X_0=\pm 2$. We verify that the ``trivial case" $\omega'(uv)=\pm I_2$
is impossible. Then we can choose representations $\omega'$ and $\omega''$ of
$\Gamma_2$ such that~:
\[
\Tr~\omega'(u)=a_0, ~~ \Tr~\omega'(v)=a_t, ~~
\omega'(uv)=\begin{pmatrix}
e & 1 \\
0 & e
\end{pmatrix},
\]
\[
\Tr~\omega''(u)=a_\infty, ~~ \Tr~\omega''(v)=a_1, ~~
\omega''(uv)=\begin{pmatrix}
e^{-1} & -1 \\
0 & e^{-1}
\end{pmatrix}.
\]
There is a freeness in the choice: we can replace $\omega'$ (resp. $\omega''$) by
an overall conjugate by an arbitrary matrix commuting with
$\begin{pmatrix}
e & 1 \\
0 & e
\end{pmatrix}$, that is of the form 
$\begin{pmatrix}
1 & t \\
0 & 1
\end{pmatrix}$ with $t\in \C$.
\end{itemize}


\subsubsection{Pant decompositions and pant parametrizations}
\label{subsubpantparclassical}

\paragraph{Pant decompositions of a $n$-punctured sphere. Trace coordinates.}

Let $n \in \N$, $n\geq 4$. We denote $S_n^2$ the $n$-punctured sphere. We replace
the $n$ punctures by little holes obtained by cutting along non-intersecting simple
closed curves surrounding the punctures, we get a \emph{$n$-holed sphere} that we
also denote $S_n^2$. 

\begin{defn}
A \emph{pant decomposition} is defined by cutting $S_n^2$ along $n-3$ simple closed
curves $\gamma_r$, $r=1,\ldots,n-3$, on $S_n^2$ in such a way that this will decompose
$S_n^2$ into a disjoint union of $n-2$ three-holed spheres $S_{3,t}^2$, $t=1,\ldots,n-2$.
The collection $\{\gamma_1,\ldots ,\gamma_{n-3}\}$ of curves is called the \emph{cut
system}. 
\end{defn}

The origin of the terminology is clear: a pair of pants is homeomorphic to a $3$
holed sphere. \\

If $n=4$, then the cut system if a set of $3$ curves. Each curve separates the set
of punctures $\{0,t,1,\infty\}$ into two unordered pairs of unordered sets of two
elements:
\[
\left((0,t),(1,\infty)\right), ~~ \left((0,1),(t,\infty)\right), ~~ 
\left((0,\infty),(t,1)\right).
\]
We will also call (abusively \dots) such a pair \emph{a pant decomposition} of
the $4$-punctured sphere. In that sense we get $3$ pant decompositions. \\

We can interpret the character variety $\mathcal{S}$ associated to the free
group $\Gamma_{n-1}$ (of rank $n-1$) as the set of representations $\rho$ of
$\pi_1(S_n^2)$ into $\SLdc$ modulo equivalence of representations.
Then useful sets of coordinates on $\mathcal{S}$ are given by the trace functions 
$\Tr\, \rho(\gamma)$ associated to any simple closed curve $\gamma$ on $S_4^2$.
It is a classical fact that minimal sets of trace functions that can be used to
parameterize $\mathcal{S}$ can be identified using \emph{pant decompositions}.
We will give only the basic idea. For more details see for example \cite{ILT} (that
we follow in our description). In the next part we will detail the case $n=4$. To a
pant decomposition we can associate decompositions of a representation $\rho$. More
abstractly we get a notion of pant decomposition of a representation of a free group
of rank $n$. This generalizes the definition introduced in the preceding paragraph. \\

To each curve $\gamma_r$ we can associate the union of the two $3$-holed spheres
which have $\gamma_r$ in their boundary (one of the $3$ components). We get a $4$-holed
sphere $S_{n,r}^2$. We choose an orientation on each curve $\gamma_r$. This allows us
to introduce a numbering of the $4$ boundary components of $S_{n,r}^2$. Then we can
consider the curves $\gamma_s^r$ and $\gamma_t^r$ which encircle respectively the
pair of component $(1,2)$ and $(2,3)$. The collection of pairs of trace functions 
$\left(\Tr\, \rho(\gamma_s^r),\Tr\, \rho(\gamma_t^r)\right)$, $r=1,\ldots n-3$, can
be used to parameterize the character variety.
 
\paragraph{Pant parametrization of a $4$ holed sphere. Jimbo formulae.}

Let $S_4^2$ be the four punctured sphere. Its fundamental group $\pi_1(S_4^2)$ is
isomorphic to a free group of rank $3$: we can choose as generators the homotopy
classes of three simple loops turning around $3$ punctures. We choose simple loops 
$\gamma_i$, $i=1,2,3,4$ turning respectively around the $4$ punctures and based at
a point $z_0$ of the punctured sphere. \\

As above, up to a M\"{o}bius transformation,  we can choose as punctures $0,t,1,\infty$
for some value of $t\in \C\setminus\{0,1\}$. Then we denote $\gamma_i$, $i=0,t,1,\infty$,
the simple loops and $M_i=\rho(\gamma_i)$. We suppose the $M_i$ semi-simple with
eigenvalues $e_i,e_i^{-1}$ and we set, as above, 
$a_i:=\Tr\, M_i=e_i+e_i^{-1}$, $i=0,t,1,\infty$, and 
$$
X_0=\Tr~M_1M_t, \quad  X_t=\Tr~M_1M_0, \quad X_1=\Tr~M_tM_0.
$$

We can apply this to the monodromy representation of a system \ref{systemp6}, then
we have local monodromy exponents $\theta_i\in \C$ and $a_i=2\cos \theta_i$. We suppose
that the non resonance conditions are satisfied: $a_i\neq \pm 2$ or equivalently
$\theta_i\notin \Z$. \\

We can interpret $M_tM_0$ as the monodromy associated to an oriented curve separating
the singularities in two packs $(0,t)$ and $(1,\infty)$, and therefore the $4$ punctured
sphere $S_4^2$ into two $3$ punctured sphere $S_3^2$. The corresponding monodromy exponent
is denoted $\sigma_1$: $X_1=2\cos 2\pi\sigma_1$. We define similarly $\sigma_0$ and
$\sigma_t$: $X_0=2\cos 2\pi\sigma_0$ and $X_t=2\cos 2\pi\sigma_t$. \\
   
We will recall some Jimbo formulae \cite{JimboMonodromy} and interpret it in
relation with (a variant of\footnote{We will fix $X_1$ in place of $X_0$.}) the
fibration of the cubic surface ${\mathcal S}_{VI}(a)$ described page
\pageref{fibrationdiff}. We use the presentation of \cite{ILT}, \cf\ 6.1, page 19,
with a change of notations. \\

We fix $a$. If we further fix $X_1$, then the equation
\[
X_0X_tX_1+X_0^2+X_t^2+X_1^2-A_0X_0-A_tX_t-A_1X_1
+A_\infty=0
\]
(where the $A_l$, which depends only on $A$ are fixed)
defines a \emph{conic} in the variables $X_0,X_t$. This conic admits a \emph{rational
parameterization} \cite{JimboMonodromy},
\cite{ILT} (\cf\ (6.67a), (6.67b), (6.68a), (6.68b), (6.68c), (6.67d))~:
\begin{equation}
\label{equajimboform}
\begin{split}
(X_1^2-4)X_0=D_{0,+} s + D_{0,-} s^{-1}+D_{0,0}  \\
(X_1^2-4)X_t=D_{t,+} s + D_{t,-} s^{-1}+D_{t,0},
\end{split}
\end{equation}
with coefficients given by~:
\begin{equation}
\label{equajimboform2}
\begin{split}
D_{0,0} &:=X_1A_t-2A_0, ~~ D_{t,0}=X_1A_0-2A_t, \\
D_{0,\pm} &:=16 \prod_{\epsilon=\pm 1} \sin \pi(\theta_t\mp \sigma_1+\epsilon \theta_0) +
\sin \pi(\theta_1\mp \sigma_1+\epsilon \theta_\infty), \\
D_{t,\pm} &:=-D_{0,\pm}\, e^{\mp 2i\pi \sigma_1}.
\end{split}
\end{equation}

More precisely the above formulae give a rational parametrization of the conic if
we suppose that he following conditions are satisfied~:
\begin{enumerate}
\label{classicalexceptfibers}
\item 
$X_1=\pm 2$, or equivalently $\sigma_1\notin \Z$;
\item
the $4$ conditions \ref{specialdecvalues} for $X_1$ are excluded, or equivalently~:
\[
\sigma_1\pm \theta_i \pm \theta_j \in 2 \Z, \quad
\sigma_1\pm \theta_k \pm \theta_\infty \in 2 \Z.
\] 
\end{enumerate}
The first case correspond to the $2$ parabolic fibers. The second case to the $4$ cases
of decomposition of the conic into two lines. \\

We can compare with the fibration of $\F$ by $\Pi$ described in
\ref{section:JimboSakaiFamilyII}. The first case correspond to the logarithmic fibers.
The second to the exceptional non logarithmic fibers. \\

Formulae \ref{equajimboform} define a parametrization of the surface
${\mathcal S}_{VI}(a)$ by  the $(X_1,s)$ (resp. $(\sigma_1,s)$) satisfying the above
restrictions. It is a \emph{pants parametrization}. The image misses $8$ lines and
$2$ parabolas. There are two others similar pants parametrizations (we replace $X_1$
by $X_0$ or $X_t$). 


\subsubsection{$q$-pants parametrizations}
\label{subsubqpants}

\paragraph{$q$-pants decompositions.}

At the beginning of \ref{subsection:assumptiononlocaldata} we introduced special
values associated to the decomposition $\{1,2,3,4\}=\{1,2)\}\cup \{3,4\}$. We recall
these values (adding indices)~:
\begin{align*}
\Xi'_{1,2} &:= \{R(-\rho_1/x_1),R(-\rho_1/x_2),R(-\rho_2/x_1),R(-\rho_2/x_2)\}, \\
\Xi''_{1,2} &:= \{R(-\sigma_1 x_3),R(-\sigma_1 x_4),R(-\sigma_2 x_3),R(-\sigma_2 x_4)\}.
\end{align*}
We assumed (this was $\mathbf{Hyp}_{8}$, see equation \eqref{eqn:Hyp8} at the end of
\ref{subsubsection:specialcases})
that $\Xi_{1,2} := \Xi'_{1,2} \cup \Xi''_{1,2}$ has eight (pairwise distinct)
elements. We can consider similar conditions for the five other decompositions. We can
assume the six conditions, then we will say that $\mathbf{Hyp}_{48}$ is satisfied. \\

In all this part we suppose that (FR), (NR), (NS) and $\mathbf{Hyp}_{48}$ are satisfied. \\
 
We consider the six decompositions of the set $\{1,2,3,4\}$ (indexing the intermediate
singularities $x_1,x_2,x_3,x_4$) into two \emph{ordered packs} of \emph{unordered}
elements~:
\[
\left((1,2),(3,4)\right), ~ \left((1,3),(2,4)\right), ~
\left((1,4),(2,3)\right), ~ \left((2,3),(1,4)\right), 
\left((2,4),(1,3)\right), ~ \left((3,4),(1,2)\right);
\]
\[
\left((i,j),(k,l)\right)=\left((j,i),(k,l)\right)=\left((i,j),(l,k)\right)
=\left((j,i),(l,k)\right).
\]
A decomposition $\left((i,j),(k,l)\right)$ is indexed\footnote{It can be convenient
to allow also indexation by $(j,i)$, with $j>i$, when it simplifies some notations.}
by $(i,j)$ such that $i<j$.  The corresponding decomposition of the set of intermediate
singularities $\left((x_i,x_j),(x_k,x_l)\right)$ is called a \emph{$q$-pants decomposition}.
As we explained before, the heuristic idea is to select a particular pair of singularities
$x_i,x_j$ among $x_1,x_2,x_3,x_4$, with the idea of ``localize'' the ``$q$-monodromy" around
that pair. \\

Be careful, in the classical case of representations of the free group $\Gamma_3$ generated
by $u_0,~u_t,~u_1,~u_\infty$ up to the the relation $u_0u_tu_1u_\infty=1$ (or equivalently of
the fundamental group of the $4$-punctured sphere $\P^1(\C)\setminus \{0,t,1,\infty\}$),
the pant-decompositions are indexed by the $3$ decompositions of $\{0,t,1,\infty\}$ into
two \emph{unordered packs} of unordered elements~:
\[
\left((0,t),(1,\infty)\right), ~~ \left((0,1),(t,\infty)\right), ~~
\left((0,\infty),(t,1)\right).
\]
This is an important difference between the representations and the ``$q$-representations''.
For the Fricke coordinates we have $\Tr\, M_0M_t=\Tr\, M_1M_\infty$ but for the $q$-analogs
we have $\Pi_{1,2}\neq \Pi_{3,4}$. \\

To each indices decomposition is associated a Mano decomposition. If necessary we will
index the objects appearing in the study of this Mano decomposition by the corresponding
$(i,j)$: $\Pi_{i,j}$, $\Xi_{1,2}, \ldots$.

\paragraph{$q$-pants parameterizations and $q$-pants charts.}

We remove from $\Eq$ the $4$ fixed points of the involution 
$\xi\mapsto \rho_1\rho_2/x_1x_2\xi$ and the $8$ points in $\Xi_{1,2}$. We get a punctured
elliptic curve denoted by $\Eq^{\bullet,\dag;1,2}$. We denote $U^{\bullet,\dag;1,2}\subset \C^*$
the inverse image of the punctured elliptic curve by the canonical map
$\C^*\rightarrow \Eq$. We set\footnote{If necessary one can precise $s_{1,2}$
and $t_{1,2}$.}~:
\begin{equation}
\begin{split}
s(\xi_1,\xi_2) := \dfrac
{\thq\left(\frac{\xi_2}{\rho_1}x_1\right) \thq\left(\frac{\xi_1}{\rho_2}x_1\right)}
{\thq\left(\frac{\xi_1}{\rho_1}x_1\right) \thq\left(\frac{\xi_2}{\rho_2}x_1\right)} ~~
\text{and} ~~
s(\xi):=s(\xi,\rho_1\rho_2/x_1x_2\xi), \\
t(\xi_1,\xi_2) := \ \dfrac
{\thq\left(\frac{\sigma_2}{\xi_1}x\right) 
\thq\left(\frac{\sigma_1}{\xi_2}x\right)}
{\thq\left(\frac{\sigma_1}{\xi_1}x\right) 
\thq\left(\frac{\sigma_2}{\xi_2}x\right)}, ~~ \text{and} ~~
t(\xi):=t(\xi,\rho_1\rho_2/x_1x_2\xi),
\end{split}
\end{equation}

\begin{equation}
\begin{split}
  \underline{P}(\xi_1,\xi_2,x;\alpha_{11},\alpha_{12},\alpha_{21},\alpha_{22}):=
  P\left(\alpha_{ij}\, \theta_q\left(\frac{\xi_j}{\rho_i}x\right)\right) ~~ \text{and} ~~
P_{1,2}(\xi,x):=\underline{P}\left(\xi,\rho_1\rho_2/x_1x_2\xi,x;1,1,1,s(\xi)\right), \\
\underline{Q}(\xi_1,\xi_2,x;\beta_{11},\beta_{12},\beta_{21},\beta_{22}):=Q\left(\beta_{ij}\,
\theta_q\left(\frac{\sigma_j}{\xi_i}x\right)\right) ~~ \text{and} ~~
Q_{1,2}(\xi,x):=\underline{P}\left(\xi,\rho_1\rho_2/x_1x_2\xi,x;1,1,1,t(\xi)\right).
\end{split}
\end{equation}

\begin{lem}
We suppose the $\rho_i$ and $\sigma_j$ ($i,j=1,2$) fixed (satisfying the ``good
conditions"). We have~:
\[
\theta_q\left(\frac{\xi_j}{\rho_i}x\right) \in V_{1,\frac{\rho_i}{\xi_j}}, ~~
\theta_q\left(\frac{\sigma_j}{\xi_i}x\right) \in V_{1,\frac{\xi_i}{\sigma_j}}, ~~
\theta_q\left(\frac{\xi_h}{\rho_i}x\right)\theta_q\left(\frac{\sigma_j}{\xi_h}x\right)
\in V_{2,\frac{\rho_i}{\sigma_j}}.
\]
The maps $\zeta_{i,j,h}:\C^*\mapsto V_{2,\frac{\rho_i}{\sigma_j}}$ ($i,j,h=1,2$) defined by
$\zeta_{i,j,h}: \xi_h \mapsto
\theta_q\left(\frac{\xi_h}{\rho_i}x\right)\theta_q\left(\frac{\sigma_j}{\xi_h}x\right)$
are analytic on $\C^*$.
\end{lem}
\Pr
We fix $(i,j)$. Let $(e_1,e_2)$ be a basis of $V_{2,\frac{\rho_i}{\sigma_j}}$. There exist
two functions $C_1(\xi_h)$ and $C_2(\xi_h)$ of $\xi_h$, uniquely determined, such that~:
\[
\zeta_{i,j,h}(\xi_h)=C_1(\xi_h)e_1+C_1(\xi_h)e_2.
\]
We have~:
\[
\zeta_{i,j,h}(\xi_h)(qx)=C_1(\xi_h)e_1(qx)+C_1(\xi_h)e_2(qx).
\]
For $x\in \C^*$ fixed, the two functions~:
$\xi_h\rightarrow \zeta_{i,j,h}(\xi_h)(x)$ and
$\xi_h\rightarrow \zeta_{i,j,h}(\xi_h)(qx)$ are analytic on $\C^*$, therefore $C_1$
and $C_2$ are analytic on $\C^*$.
\Finpr

The functions $s$ and $t$ are meromorphic on $\C^*$, they are analytic on the inverse
image $U^{\bullet,\dag;1,2}$ of $\Eq^{\bullet,\dag;1,2}$. Identifying
$V:=V_{2,\frac{\rho_1}{\sigma_1}} \times V_{2,\frac{\rho_1}{\sigma_2}} \times
V_{2,\frac{\rho_2}{\sigma_1}}\times V_{2,\frac{\rho_2}{\sigma_2}}$
with a set of matrices, we define a map 
$\underline M_{1,2}:U^{\bullet,\dag;1,2}\times \C^*\rightarrow V$ by~:
\begin{equation}
\underline M_{1,2}(\xi,\eta) := P_{1,2}(\xi,x)\, \Diag(1,\eta)\, Q_{1,2}(\xi,x).
\end{equation}
If we fix $\xi$, then $M_{1,2}(\xi,w)$ is \emph{linear} in $\eta$ (in the trivial sense).
The map $\underline M_{1,2}$ is analytic in the variable $\xi$ and it extends uniquely in
a map meromorphic on $\C^*\times\C^*$. We have~: $M_{1,2}(q\xi,w)=M_{1,2}(\xi,w)$. Therefore
$M_{1,2}$ induces an analytic map 
$M_{1,2}:\Eq^{\bullet,\dag;1,2}\times \C^*\rightarrow V$ and this map extends uniquely into
a map meromorphic on $\Eq\times \C^*$. Therefore it can be interpreted as a \emph{rational
map} from $\Eq\times \C^*$ to the linear space $V$. This map is \emph{regular} on
$\Eq^{\bullet,\dag;1,2}\times \C^*$. The image of $M_{1,2}$ is contained in $F$, therefore
we get by corestriction a map (abuse of notations \dots)
$M_{1,2}:\Eq^{\bullet,\dag;1,2}\times \C^*\rightarrow F$. \\

If we compose by the quotient map $F\mapsto \F$, then we get
$\overline M_{1,2}: \Eq^{\bullet,\dag;1,2}\times \C^*\rightarrow \F$.
This map is called the \emph{$q$-pant parameterization} associated to the $q$-pant
decomposition $(1,2)$. \\

We have\footnote{We denote abusively $\Pi\circ M=\Pi\circ \overline M$.}~ 
\begin{equation}
\label{equam12}
\Pi_{12}\circ M_{1,2}(\xi,\eta)=\Phi_{1,2}(\xi), ~~
M_{1,2}(q\xi,\eta)\sim M_{1,2}(\xi,\eta), ~~
M_{1,2}(\xi,\eta)= M_{1,2}(\rho_1\rho_2/x_1x_2\xi,\eta^{-1}).
\end{equation}

The map $\overline M_{1,2}$ is not injective: the fiber is 
$\left((\xi,\eta),(\rho_1\rho_2/x_1x_2\xi,\eta^{-1})\right)$. Therefore
$\overline{M}_{1,2}$ induces an \emph{injective} map~:
 \[
 \psi_{1,2}:\mathcal{Y}^{\bullet,\dag,1,2}\rightarrow \mathcal{F},  
 \]
where  
$$
\mathcal{Y}^{\bullet\dag;1,2} :=
\dfrac{\Eq^{\bullet\dag;1,2} \times \C^*}{\text{involution~} \tau_{12}} \cdot
$$

This map is called the \emph{$q$-pant chart} associated to the $q$-pant decomposition
$(1,2)$. We can interpret $\mathcal{Y}^{\bullet\dag;1,2}$ as an algebraic variety. We
denote it $\mathcal{Y}_{alg}^{\bullet\dag;1,2}$ and we denote~:
\[
\psi_{alg;1,2}: \mathcal{Y}_{alg}^{\bullet\dag;1,2}\rightarrow
\mathcal{F}
\]
the corresponding \emph{regular} map. It is called the algebraic $q$-pant chart
associated to the $q$-pant decomposition $(1,2)$. There are respectively $6$ similar
$q$-pant parameterizations, $q$-pant charts, algebraic $q$-pant charts. \\

We can consider the maps~:
\[
\Pi_{i',j'}\circ \psi_{alg;i,j}: \mathcal{Y}^{\bullet\dag;i,j}\rightarrow 
\P^1(\C).
\]
They are \emph{regular} maps and can be computed \emph{explicitly}. At the level
of $q$-pant parameterization we can compute explicitly $\Pi_{i',j'}\circ M_{i,j}$ using
$M_{i,j}(x_{i'})$ and $M_{i,j}(x_{j'})$. \\

If we interpret the six $\Pi_{i',j'}$ as ``coordinates'', then we get \emph{$q$-coordinate
charts from the pant decomposition $(i,j)$}, the $q$-analogs of the coordinate charts
from a pant decomposition of the classical case given by \ref{equajimboform}. 

\begin{rmk}
\label{remextendparam}
In \ref{subsubellipticfib}
we extended the involution $\tau$ into an involution $\tilde\tau$ on 
$\Eq\times \P^1(\C)$ and set 
$\widetilde{\mathcal{Y}}:=\dfrac{\Eq\times \P^1(\C)}{\text{involution~} \tilde\tau}$.
Afterwards we extended the elliptic fibration $p$ into an elliptic fibration 
$\tilde p: \widetilde{\mathcal{Y}}\rightarrow \P^1(\C)$ and the application 
$\Psi$ into an application $\widetilde \Psi$. \\
At the end of \ref{subsubsection:fiberinggenericpart}
we remarked that the surface $\left\{v^2=(w^2-1)f(x)\right\}$ is a double covering
of the affine plane $\C^2$ (the $(x,w)$ plane) ramified above $6$ lines. \\
We consider the map 
$\pi:=(\tilde p,\widetilde \Psi): \widetilde{\mathcal{Y}}\rightarrow 
\left(\P^1(\C)\right)^2$. It is a double covering of $\left(\P^1(\C)\right)^2$
ramified along $6$ lines \big(a $(2,3)$sextic of $\left(\P^1(\C)\right)^2$\big).
The surface $\widetilde{\mathcal{Y}}$ has $8$ singular points (rational double points),
above the pairwise intersections of the $6$ lines. Blowing up these $8$ points, we get
a smooth surface $\mathcal{X}$ (one can compare with example \ref{exadoubleplane}).
It is possible to compute explicitely a system of algebraic charts for $\mathcal{X}$.
Then the surface $\mathcal{X}$ minus the $12$ lines above
$\widetilde \Psi^{-1}\left(\P^1(\C)^{\bullet \dag}\right)$, that we denote $\mathcal{X}^{\dag}$,
could perhaps be used for a parameterization of $\F$ by explicit Zariski open sets. \\
More precisely, we conjecture that it is possible to extend the $q$-pant chart $\psi$
into a regular injective map $\mathcal{X}\rightarrow \F$ (\cf\ \ref{remextendparam})
such that its image is a smooth open set of $\F$ containing the $4$ logarithmic fibers.
Conjecture \ref{conjlisslog} would follow.
\label{extendparam}
\end{rmk}

\paragraph{A smoothness conjecture.}

We end this subsection with a conjecture\footnote{It is a $q$-analog of theorem
\ref{lissites}.}. This conjecture is strongly related to the configuration of the
lines on the surface $\F$. We will return to this question in the next paragraphs.\\

We set $\mathcal{U}_{i,j}=\F\setminus\left(\Pi_{i,j}^{-1}(0)\cup
\Pi_{i,j}^{-1}(\infty)\cup (\Pi')_{i,j}^{-1}(0)\cup 
(\Pi')_{i,j}^{-1}(\infty) \right)$; it is $\F$ minus the $4$ exceptional non logarithmic
fibers. The image of $\psi_{i,j}$ is $\F$ minus all the exceptional fibers, therefore
it is $\mathcal{U}_{i,j}$ minus the logarithmic fibers.

\begin{conj}
\label{conjsmooth}
We suppose that (FR), (NR), (NS) and 
$\mathbf{Hyp}_{48}$ are satisfied. Then $\F$ is \emph{smooth}.
\end{conj}

We will prove below (\cf\ proposition \ref{propcoverqpants}) that $\mathcal{F}$ minus
all the logarithmic fibers is covered by the union of the images\footnote{The image of
each $q$-pants chart is $\F$ minus $8$ lines (depending on the chart).} of the six
$q$-pant charts $\psi_{i,j}$ ($i,j\in \{1,2,3,4\}, i<j$). Then the above
conjecture will follow immediately from conjecture \ref{conjlisslog} \\

The rational functions $\Pi_{i,j}$ on the surface $\F$ can be interpreted as $q$-analogs
of the Fricke coordinates. Then there are natural questions:
\begin{itemize}
\item 
What are the algebraic relations between the $\Pi_{i,j}$ ? How to compute them ?
\item
Is it possible to use the $\Pi_{i,j}$ to build an embedding of
$\F$ into some $\left(\P^1(\C)\right)^m$ ($m\in \N$, $m\geq 3$) ?
\item 
If we denote $X$ the closure of the image of $\F$ by such an embedding, what can we say
of the surface $X$ ?
\end{itemize}

We will return later to these questions, \cf\ \ref{subsubmainconj1}.


\subsubsection{$q$-pants decompositions and partial reducibility}

We consider the quotient of $\P^1(\C)$ by the action of $q^\Z$. We write it~: 
\[
\left[\P^1(\C);q\right]:=\{[0;q]\}\cup \Eq \cup \{[\infty;q]\}. 
\]
We have a ``cotangent bundle" (defined by the first projection)~:
\[
\left[\P^1(\C);q\right]\times \C^*\rightarrow \left[\P^1(\C);q\right].
\]
It is the $q$-analog of the cotangent bundle of $\P^1(\C)$ (of more generally a
Riemann surface) in the differential case; we will call it the \emph{$q$-cotangent
bundle}. \\

The fiber above $[a;q]$ is $[a;q]\times \C^*$. When the ``point'' $[a;q]$ is a
singularity, we will consider $\C^*$ as the space of \emph{possible monodromy
exponents}. More precisely above $[a;q]=[0;q]$ or $[a;q]=[\infty;q]$ we can choose
arbitrarily a monodromy exponent into $\C^*$ but if $[a;q]\in \Eq$, the only possible
choice above $[a;q]$ is $\zeta\in \C^*$ such that $[\zeta:q]=[-a^{-1};q]$. \\

In the context of our description of $q$-PVI we consider the following list of pairs
of points of the total space of the $q$-cotangent bundle\footnote{We skip the problem
of ordering or not of such pair.}~:
\begin{equation}
\label{listreduct}
\left(\left([0;q];\rho_1\right),\left([0;q];\rho_2\right)\right), ~~ 
\left(\left([\infty;q];\sigma_1\right),\left([\infty;q];\sigma_2\right)\right), ~~ 
\left(\left([x_i;q];-x_i^{-1}\right),\left([x_j;q];-x_j^{-1}\right)\right).
\end{equation}
($i=1,2,3,4$, $i\neq j$). Be careful~:
 \[
 \left(\left([x_i;q];-x_i^{-1}\right),\left([x_j;q];-x_j^{-1}\right)\right)
 \] 
is written in the chart of the $q$-cotangent bundle coming from $\C$. In the chart
coming from $\P^1(\C)\setminus \{0\}$, we write the same element
\[
\left(\left([x_i^{-1};q];-x_i\right),\left([x_j^{-1};q];-x_j\right)\right).
\]

The Mano decomposition allows us to decompose the \emph{global} monodromy
around $0$, $\infty$ and the four intermediate singularities into a pair
of \emph{local} monodromies~: around $0$ and one pair of singularities on
one side and one pair of singularities and $\infty$ on the other side. As we
will see it is better to interpret $0$, $\infty$ and the pair $(x_i,x_j)$ as
the corresponding elements in the list (\ref{listreduct}). \\

We have a criterion of reducibility for each local monodromy. Using the above
list (\ref{listreduct}), we will see that it is a perfect $q$-analog of the
criterium in the differential case (\cf\ proposition \ref{propprreprandlines}).\\

Looking on the left hand side of RH, the Mano decomposition can be interpreted as a
decomposion of the system into two hypergeometric systems. We recall that it is
the beginning of our story: in \cite{Mano} Mano gave a \emph{direct} method (based
on isomonodromy and a $q$-analogy with Jimbo decompoition \cite{JimboMonodromy}) in
order to decompose the original system. \\

The reducibility of a local monodromy on the right hand side of RH is equivalent to
the reducibility of the corresponding hypergeometric system on the left hand side (\cf\
subsection \ref{subsection:reducibilitycriteria}). \\

The reducibility of the local monodromy around $0$ and the pair $(x_i,x_j)$ is coded by~:
\[
 \left(\left([0;q];\rho_1\right),\left([0;q];\rho_2\right)\right)~~\text{and} ~~ \left(\left([x_i;q];-x_i^{-1}\right),\left([x_j;q];-x_j^{-1}\right)\right).
\]
The $4$ conditions of reducibility are~:
\[
\xi_1=-\rho_1/x_i, ~~ \xi_1=-\rho_1/x_j, ~~
\xi_1=-\rho_2/x_i, ~~ \xi_1=-\rho_2/x_j .
\]

They correspond to the following pairings~:
\begin{align*}
\left([0;q];\rho_1\right) &\longleftrightarrow \left([x_i;q];-x_i^{-1}\right),
\quad 
\left([0;q];\rho_1\right) \longleftrightarrow \left([x_j;q];-x_j^{-1}\right),
\\
\left([0;q];\rho_2\right) &\longleftrightarrow \left([x_i;q];-x_i^{-1}\right),
\quad 
\left([0;q];\rho_2\right) \longleftrightarrow \left([x_j;q];-x_j^{-1}\right).
\end{align*}
For simplicity we will denote these pairings
$(\rho_h\leftrightarrow x_i)$ ($h=1,2;~i=1,2,3,4$). \\

The reducibility of the local monodromy around the pair
$(x_k,x_l)$ and $\infty$ is coded by~:
\[
 \left(\left([\infty;q];\sigma_1\right),\left([\infty;q];\sigma_2\right)\right)~~\text{and} ~~ 
 \left(\left([x_k^{-1};q];-x_k\right),\left([x_l^{-1};q];-x_l\right)\right).
\]
The $4$ conditions of reducibility are~:
\[
\xi_1=-\sigma_1x_k, ~~ \xi_1=-\sigma_1x_l, ~~
\xi_1=-\sigma_2x_k, ~~ \xi_1=-\sigma_2x_l .
\]

They correspond to the following pairings~:
\begin{align*}
\left([\infty;q];\sigma_1\right) &\longleftrightarrow \left([x_k^{-1};q];-x_k\right),
\quad 
\left([\infty;q];\sigma_1\right) \longleftrightarrow \left([x_l^{-1};q];-x_l\right),
\\
\left([\infty;q];\sigma_2\right) &\longleftrightarrow \left([x_k^{-1};q];-x_k\right),
\quad 
\left([\infty;q];\sigma_2\right) \longleftrightarrow \left([x_l^{-1};q];-x_l\right).
\end{align*}
For simplicity we will denote these pairings
$(\sigma_h\leftrightarrow x_i)$ ($h=1,2;~i=1,2,3,4$). \\

We will see in the next subsection that each pairing $\longleftrightarrow$ corresponds
to a \emph{line} on the surface $\F$. We have $8$ pairings and therefore $8$ lines.\\

The pairing $(\rho_1 \leftrightarrow x_i)$ (resp. $(\rho_2 \leftrightarrow x_i)$)
appears in the $3$ Mano decompositions $(i,j)$, $(i,k)$ and $(i,l)$ (where $(i,j,k,l)$
is a permutation of $(1,2,3,4)$). In the next paragraph we will prove that the
corresponding $3$ lines \emph{coincide}. \\

Similarly the pairing $(\sigma_1\leftrightarrow  x_{i'})$ (resp.
$(\sigma_2\leftrightarrow  x_{i'})$) appears in the $3$ Mano decompositions
$(i,j)$, $(i,k)$ and $(i,l)$ ($(i,j,k,l)$ (where $(i,j)=\left((i,j),(i',j')\right)$ \dots).
We will also prove that the corresponding $3$ lines \emph{coincide}. 


\subsubsection{Description of some lines on the surface $\F$}

In all this part we suppose that (FR), (NR), (NS) and $\mathbf{Hyp}_{48}$ are satisfied. \\

We will give a global description of the special non-logarithmic lines on the surface
$\F$ and observe that this description is a translation by $q$-analogies of the
dictionnary between the set of lines on the cubic surface on one side and the partial
reducibility of representations of the other side that we described in the classical
cases (\cf\ \ref{subsubpartialreduclines}). \\

We will derive from this description a notion of \emph{reducibility of the local
monodromy around $0$ and $x_i$} (resp. \emph{$\infty$ and $x_i$}). \\

We recall~:
\begin{align*}
\Phi_{1,2}\left(-\rho_1/x_1\right) = \Phi_{1,2}\left(-\rho_2/x_2\right) &= 0, \\
\Phi_{1,2}\left(-\rho_1/x_2\right) = \Phi_{1,2}\left(-\rho_2/x_1\right) &= \infty, 
\end{align*}
and we set~:
\begin{align*}
\mathbf{e}_q^{1;1,2;3}(\underline{\rho},\underline{\sigma},\underline{x}):=\Phi_{1,2}\left(-\sigma_1 x_3\right) = \Phi_{1,2}\left(-\sigma_2 x_4\right) &=
\dfrac
{\thq\left(\frac{\sigma_1}{\rho_1} x_1 x_3\right)
\thq\left(\frac{\sigma_1}{\rho_2} x_2 x_3\right)}
{\thq\left(\frac{\sigma_1}{\rho_1} x_2 x_3\right)
\thq\left(\frac{\sigma_1}{\rho_2} x_1 x_3\right)}, \\
\mathbf{e}_q^{2;1,2;3}(\underline{\rho},\underline{\sigma},\underline{x}) :=
\Phi_{1,2}\left(-\sigma_1 x_4\right) = \Phi_{1,2}\left(-\sigma_2 x_3\right) &=
\dfrac
{\thq\left(\frac{\sigma_2}{\rho_1} x_1 x_3\right)
\thq\left(\frac{\sigma_2}{\rho_2} x_2 x_3\right)}
{\thq\left(\frac{\sigma_2}{\rho_1} x_2 x_3\right)
\thq\left(\frac{\sigma_2}{\rho_2} x_1 x_3\right)}
\cdot
\end{align*}

We verify~:
\[
\mathbf{e}_q^{1;1,2;3}(\underline{\rho},\underline{\sigma},\underline{x})=
\mathbf{e}_q^{2;1,2;4}(\underline{\rho},\underline{\sigma},\underline{x}) ~~
\text{and} ~~ \mathbf{e}_q^{2;1,2;3}(\underline{\rho},\underline{\sigma},\underline{x})=
\mathbf{e}_q^{1;1,2;4}(\underline{\rho},\underline{\sigma},\underline{x}).\]
 Using the other decompositions we define similarly 
$\mathbf{e}_q^{h;i,j;k}(\underline{\rho},\underline{\sigma},\underline{x})$
(where $h=1,2$ and $(i,j,k)$ is a set of $3$ distinct elements of $\{1,2,3,4\}$).
We have~: $\mathbf{e}_q^{h;i,j;k}=\mathbf{e}_q^{h;j,i;k}$ and we verify~:
\[
\mathbf{e}_q^{1;i,j;k}(\underline{\rho},\underline{\sigma},\underline{x})=
\mathbf{e}_q^{2;i,j;l}(\underline{\rho},\underline{\sigma},\underline{x}) ~~
\text{and} ~~ \mathbf{e}_q^{2;i,j;k}(\underline{\rho},\underline{\sigma},\underline{x})=
\mathbf{e}_q^{1;i,j;l}(\underline{\rho},\underline{\sigma},\underline{x}).\]

The function $\mathbf{e}_q^{h;i,j;k}$ is elliptic in $\sigma_h$. We interpret the $12$
functions as $q$-analogs of the $a_l=2\cos \theta_l$ ($l=0,t,1,\infty$) of the
differential case (the traces of the local monodromies). A big difference is that
the $\mathbf{e}_q^{h;i,j;k}$ involve \emph{all} the local data. \\

We will call the $\mathbf{e}_q^{h;i,j;k}$ the \emph{$q$-local monodromy invariants} 
\label{qmonoinvariants}. We conjecture that when the monodromy exponents
$\underline{\rho},~\underline{\sigma},~\underline{x}$ move these $q$-local monodromy
invariants ``parameterize algebraically" the variation of
$\F(\underline{\rho},\underline{\sigma},\underline{x})$.\\

We have~:
\begin{align*}
(\Pi')_{1,2}^{-1}(0) &=\Pi_{1,2}^{-1}\left(\mathbf{e}_q^{1;1,2;3}\right), ~~
(\Pi')_{1,2}^{-1}(\infty)=\Pi_{1,2}^{-1}\left(\mathbf{e}_q^{2;1,2;3}\right), ~~
\end{align*}

To the pairing $(\rho_1\leftrightarrow x_i)$ \big(resp.
$(\rho_2\leftrightarrow x_i)\big)$ we associate a line $L_{\rho_1,\sigma_i}$ (resp.
$L_{\rho_2,\sigma_i}$) of $\F$: it is the set of the classes of the matrices $M$
such that the first (resp. second) \emph{line} of $M(x_i)$ is null.
Similarly, to the pairing $(\sigma_1\leftrightarrow x_i)$ 
(resp. $(\sigma_2\leftrightarrow x_i))$ we associate a line $L_{\sigma_1,x_i}$
(resp. $L_{\sigma_2,x_i}$) of $\F$: it is the set defined by the classes of
the matrices $M$ such that the first (resp. second) \emph{column} of $M(x_i)$
is null. Using these lines we can describe the exceptional non logarithmic fibers
of $\Pi_{i,j}$ (and $\Pi'_{i,j}$). We detail the case $(i,j)=(1,2)$; the others
are similar. \\

We denote $\overline M\in \F$ the equivalence class of $M \in F$. We have~:
\[
f_i=0 \Leftrightarrow \overline M\in L_{\rho_1,x_i}
\Leftrightarrow  M(x_i) = \begin{pmatrix} 0 & 0 \\ \ast & \ast \end{pmatrix}
\quad \text{and} \quad
g_i=0 \Leftrightarrow \overline M\in L_{\rho_2,x_i}
\Leftrightarrow  M(x_i) = \begin{pmatrix}  \ast & \ast
\\ 0 & 0 \end{pmatrix}.
\]
Then (\cf\ theorem \ref{thm:specialfibersofPi})~:
\[
\Pi_{1,2}^{-1}(0)=\{f_1=0\}\cup \{g_2=0\}
=L_{\rho_1,x_1}\cup L_{\rho_2,x_2}
     ~~ \text{and} ~~ \Pi_{1,2}^{-1}(\infty)=\{f_2=0\}\cup \{g_1=0\}
     =L_{\rho_1,x_2}\cup L_{\rho_2,x_1}.
\]
\[
f_i=0 \Leftrightarrow \overline M\in L_{\rho_1,x_i}
\Leftrightarrow  M(x_i) = \begin{pmatrix} 0 & 0 \\ \ast & \ast \end{pmatrix}
\quad \text{and} \quad
g_i=0 \Leftrightarrow \overline M\in L_{\rho_2,x_i}
\Leftrightarrow  M(x_i) = \begin{pmatrix}  \ast & \ast
\\ 0 & 0 \end{pmatrix}.
\]
Similarly we have~:
\[
f'_i=0 \Leftrightarrow \overline M\in L_{\sigma_1,x_i}
\Leftrightarrow  M(x_i) = \begin{pmatrix} 0 & \ast \\ 0 & \ast \end{pmatrix}
\quad \text{and} \quad
g'_i=0 \Leftrightarrow \overline M\in L_{\sigma_2,x_i}
\Leftrightarrow  M(x_i) = \begin{pmatrix}  \ast & 0
\\ \ast & 0 \end{pmatrix}.
\]

Then~:
\[
({\Pi'})_{1,2}^{-1}(0)=\{f'_3=0\}\cup \{g'_4=0\}
=L_{\sigma_1,x_3}\cup L_{\sigma_2,x_4}
     ~~ \text{and} ~~ ({\Pi'})_{1,2}^{-1}(\infty)=\{f'_4=0\}\cup \{g'_3=0\}
     =L_{\sigma_1,x_4}\cup L_{\sigma_2,x_3}.
\]

The two by two intersections of the $4$ special fibers are empty, therefore
the $8$ lines are distinct and we have for these lines the following incidence
relations ~:
\begin{equation}
\begin{split}
L_{\rho_1,x_1} \cap L_{\rho_1,x_2}=\varnothing, ~
L_{\rho_2,x_1} &\cap L_{\rho_2,x_2}=\varnothing, ~~
L_{\sigma_1,x_3}\cap L_{\sigma_1,x_4}=\varnothing, ~
L_{\sigma_2,x_3}\cap L_{\sigma_2,x_4}=\varnothing,\\
 \forall \, h,h'=1,2, ~ & \forall \, i=1,2, ~ \forall \, j=3,4, ~~   
L_{\rho_h,x_i} \cap L_{\sigma_{h'},x_j}=\varnothing \\
 L_{\rho_1,x_1}\cap L_{\rho_2,x_2} &=\{\text{one~point}\}, ~~
L_{\rho_1,x_2}\cap  L_{\rho_2,x_1}=\{\text{one~point}\}, \\
L_{\sigma_1,x_3}\cap  L_{\sigma_2,x_4} &=\{\text{one~point}\}, ~~
L_{\sigma_2,x_3}\cap  L_{\sigma_1,x_4} =\{\text{one~point}\}.
\end{split}
\end{equation}
Replacing the decomposition $(1,2)$ by another decomposition, we get similar results. \\

We cannot have 
$M(x_i) = \begin{pmatrix} 0 & 0 \\ \ast & 0 \end{pmatrix}$,
$M(x_i) = \begin{pmatrix} \ast & 0 \\ 0 & 0 \end{pmatrix}$,
$M(x_i) = \begin{pmatrix} 0 & \ast \\ 0 & 0 \end{pmatrix}$, 
$M(x_i) = \begin{pmatrix} 0 & 0 \\ 0 & \ast \end{pmatrix}$,
therefore~:
\[
\forall \, h,h'=1,2,~i=1,2,3,4, \quad 
L_{\rho_h,x_i}\cap L_{\sigma_{h'},x_i}=\varnothing.
\]
For $h,h'=1,2$ and $i,j$ fixed, $i\neq j$, the lines $L_{\rho_h,x_i}$ and
$L_{\sigma_h',x_j}$ appear in two different exceptionnal fibers
$\Pi_{i,k}^{-1}$ ($k\neq i$ and $k\neq j$), therefore
$L_{\rho_h,x_i}\cap L_{\sigma_h',x_j}=\varnothing$. \\

Putting things together, we verify that we have $16$ \emph{different} lines
with the following incidence relations~:
\begin{itemize}
\item
for $h,h'=1,2$, $i,j=1,2,3,4$, $L_{\rho_h,x_i}\cap 
L_{\sigma_{h'},x_j}=\varnothing$.
\item 
for $h=1,2$, if $i\neq j$ then 
$L_{\rho_h,x_i}\cap L_{\rho_h,x_j}=\varnothing$ \, ;
\item 
for $h=1,2$, if $i\neq j$ then 
$L_{\sigma_h,x_i}\cap L_{\sigma_h,x_j}=\varnothing$ \, ;
\item
if $(i,j,k,l)$ is a permutation of $(1,2,3,4)$, then~: 
\begin{itemize}
\item 
$L_{\rho_1,x_i}$ meets $L_{\rho_2,x_j}$,
$L_{\rho_2,x_k}$, $L_{\rho_2,x_l}$ and the $3$ intersection points are \emph{distinct};
\item
$L_{\rho_2,x_i}$ meets $L_{\rho_1,x_j}$,
$L_{\rho_1,x_k}$, $L_{\rho_1,x_l}$ and the $3$ intersection points are \emph{distinct}
\end{itemize}
\item
if $(i,j,k,l)$ is a permutation of $(1,2,3,4)$, then~: 
\begin{itemize}
\item 
$L_{\sigma_1,x_i}$ meets $L_{\sigma_2,x_j}$,
$L_{\sigma_2,x_k}$, $L_{\sigma_2,x_l}$ and the $3$ intersection points are \emph{distinct};
\item
$L_{\sigma_2,x_i}$ meets $L_{\sigma_1,x_j}$,
$L_{\sigma_1,x_k}$, $L_{\sigma_1,x_l}$ and the $3$ intersection points are \emph{distinct}.
\end{itemize}
\end{itemize}

Each line is contained into exactly $3$ special fibers~:
\begin{equation}
\begin{split}
  L_{\rho_1,x_i} &\subset \Pi_{i,j}^{-1}(0), ~~ L_{\rho_1,x_i} \subset\Pi_{i,k}^{-1}(0),
  ~~L_{\rho_1,x_i}\subset\Pi_{i,l}^{-1}(0)  \\
L_{\rho_2,x_i} &\subset \Pi_{i,j}^{-1}(\infty), ~~  
L_{\rho_2,x_i} \subset \Pi_{i,k}^{-1}(\infty), ~~
L_{\rho_2,x_i}\subset \Pi_{i,l}^{-1}(\infty). \\
 L_{\sigma_1,x_i} &\subset \Pi_{j,k}^{-1}(\mathbf{e}_q^{1;j,k;i}), ~~ 
 L_{\sigma_1,x_i} \subset\Pi_{k,l}^{-1}(\mathbf{e}_q^{1;k,l;i}), ~~
 L_{\sigma_1,x_i}\subset\Pi_{l,j}^{-1}(\mathbf{e}_q^{1;l,j;i})  \\
 L_{\sigma_2,x_i} &\subset \Pi_{j,k}^{-1}(\mathbf{e}_q^{2;j,k;i}), ~~ 
 L_{\sigma_2,x_i} \subset\Pi_{k,l}^{-1}(\mathbf{e}_q^{2;k,l;i}), ~~
 L_{\sigma_2,x_i}\subset\Pi_{l,j}^{-1}(\mathbf{e}_q^{2;l,j;i}).
\end{split}
\end{equation}
Moreover this line is equal to each pairwise intersection of the $3$ special fibers. \\

If $(i,j,k,l)$ is a permutation of $(1,2,3,4)$, then we have~:
\begin{equation}
\begin{split}
L_{\rho_1,x_i} &=\Pi_{i,j}^{-1}(0)\cap \Pi_{i,k}^{-1}(0)
=\Pi_{i,k}^{-1}(0)\cap \Pi_{i,l}^{-1}(0)
=\Pi_{i,l}^{-1}(0)\cap \Pi_{i,j}^{-1}(0)  \\
L_{\rho_2,x_i} &=\Pi_{i,j}^{-1}(\infty)\cap \Pi_{i,k}^{-1}(\infty)
=\Pi_{i,k}^{-1}(\infty)\cap \Pi_{i,l}^{-1}(\infty)
=\Pi_{i,l}^{-1}(\infty)\cap \Pi_{i,j}^{-1}(\infty)  \\
L_{\sigma_1,x_i} &=\Pi_{j,k}^{-1}(\mathbf{e}_q^{1;j,k;i})\cap \Pi_{k,l}^{-1}(\mathbf{e}_q^{1;k,l;i})
=\Pi_{k,l}^{-1}(\mathbf{e}_q^{1;k,l;i})\cap \Pi_{l,j}^{-1}(\mathbf{e}_q^{1;l,j;i}) \\
&=\Pi_{l,j}^{-1}(\mathbf{e}_q^{1;l,j;i})\cap \Pi_{j,k}^{-1} (\mathbf{e}_q^{1;j,k;i})  \\
L_{\sigma_2,x_i} &=\Pi_{j,k}^{-1}(\mathbf{e}_q^{2;j,k;i})\cap \Pi_{k,l}^{-1}(\mathbf{e}_q^{2;k,l;i})
=\Pi_{k,l}^{-1}(\mathbf{e}_q^{2;k,l;i})\cap \Pi_{l,j}^{-1}(\mathbf{e}_q^{2;l,j;i}) \\
&=\Pi_{l,j}^{-1}(\mathbf{e}_q^{2;l,j;i})\cap \Pi_{j,k}^{-1} (\mathbf{e}_q^{2;j,k;i})  
\end{split}
\end{equation}
 
We recall the notation~:
$\mathcal{U}_{i,j}=\F\setminus\left(\Pi_{i,j}^{-1}(0)\cup
\Pi_{i,j}^{-1}(\infty)\cup (\Pi')_{i,j}^{-1}(0)\cup 
(\Pi')_{i,j}^{-1}(\infty) \right)$.

\begin{prop}
\label{propcoverqpants}
\begin{itemize}
\item[(i)]
Each line intersects $3$ other lines at $3$ different points. 
\item[(ii)]
If $a\in \F$, there are at most two different lines passing by $a$.
\item[(iii)]
The set $\bigcup \mathcal{U}_{i,j}$ contains $\F$ less all the logarithmic fibers.
\item[(iv)]
The set of $16$ lines has two connected components: the set 
$\{L_{\rho_h,x_i}\}_{h=1,2,;i=1,2,3,4}$ and the set 
$\{L_{\sigma_h,x_i}\}_{h=1,2,;i=1,2,3,4}$.
\end{itemize}
\end{prop}
\Pr
\begin{itemize}
\item[(i)] and\,  (ii) follows easily from the above relations
\item[(iii)]
Let $a\in \F$ which does not belong to a logarithmic fiber. We suppose that 
$a\notin \mathcal{U}_{1,2}\cup \mathcal{U}_{1,3}\cup \mathcal{U}_{1,4}$. Then we have
in particular $a\notin \mathcal{U}_{1,2}$ and $a$ must belong to one of the $4$
exceptional fibers and therefore to one of the eight lines. We have a similar
result for $(1,3)$ and $(1,4)$, therefore $a$ belongs to one of the four lines
$L_{\rho_h,x_1}$, $L_{\sigma_{h'},x_1}$. \\
We prove similarly that $a$ belongs to one of the four lines $L_{\rho_h,x_2}$,
$L_{\sigma_h,x_2}$, to one of the four lines $L_{\rho_h,x_3}$, $L_{\sigma_{h'},x_3}$ and
to one of the four lines $L_{\rho_h,x_4}$, $L_{\sigma_{h'},x_4}$. Then $a$ belongs to
$4$ distinct lines. This contradicts (ii).
\item[(iv)]
The intersection of the two sets is empty and, using (i), we can verify that each
set is connected.
\end{itemize}
\Finpr

If conjecture \ref{conjlisslog} is true, then the above proposition
implies that $\F$ is smooth, that is conjecture \ref{conjsmooth}. 


\subsubsection{An image of $\F$ in $\left(\PC\right)^3$}

\paragraph{The image and what we know about it.}

We consider the map $T_{1,2}:\F\rightarrow \left(\PC\right)^3$ defined by~:
\[
T_{1,2} := (\Pi_{1,2},\Pi_{2,3},\Pi_{3,4}):
\overline M\mapsto \left(u(\overline M)=\Pi_{1,2}(\overline M),
v(\overline M)=\Pi_{2,3}(\overline M),u'(\overline M)=\Pi_{3,4}(\overline M)\right).
\]
A better notation would be $T_{1,2,3,4}$: note that $T_{1,2,3,4} \neq T_{1,2,4,3}$ because
there is an exchange between $v$ and $u'$. We will use (carefully \dots) $T_{1,2}$ for
simplicity. \\

We would like to describe the Zariski closure\footnote{It is possible that $T_{1,2}(\F)$
is closed \dots} $Y$ in $\left(\PC\right)^3$ of the image $T_{1,2}(\F)$ and in particular
the closures of the images of the $16$ lines on $\F$. 

\begin{defn}
The \emph{$(1,2)$-skeleton} of $\F$ is the closure in $\left(\PC\right)^3$ of the image 
by $T_{1,2}$ of the set of the $16$ lines. We denote it
$$\text{Sk}_{1,2}(\F) \subset T_{1,2}(\F) \subset  \left(\PC\right)^3.
$$
\end{defn}

If $(i,j,k,l)$ is  a \emph{circular} permutation of $(1,2,3,4)$ we can define similarly
the  $(i,j)$-skeleton $\text{Sk}_{(i,j)}(\F)$. \\

It is important to understand the skeleton structure and to describe the inclusion
$\text{Sk}_{1,2}(\F) \subset T_{1,2}(\F)$. In the next paragraph we will explain how
(under some ``reasonable" conjectures) it is possible to ``sew the surface onto the
skeleton bones". \\

We have $L_{\rho_1,x_2} =\Pi_{1,2}^{-1}(0)\cap \Pi_{2,3}^{-1}(0)$, therefore the image of
this line in $\left(\PC\right)^3$ is contained in the line $\{u=v=0\}$, the intersection
of the planes $\{u=0\}$ and $\{v=0\}$. As $u'$ is not constant on the image, this image
is an open connected affine subset of the projective line $\{u=v=0\}$ and therefore it is
equal to $\{u=v=0\}$ punctured at a point. Its closure is the projective line. \\

We have similar results for the lines~:
\begin{equation}
\label{equagoodlines}
L_{\rho_2,x_2}, ~
L_{\rho_1,x_3}, ~ L_{\rho_2,x_3}, ~
L_{\sigma_1,x_1}, ~ L_{\sigma_2,x_1},
L_{\sigma_1,x_4}, ~ L_{\sigma_2,x_4}. 
\end{equation}
We will call \emph{half-skeleton} the set $\text{Sk}'_{1,2}(\F)$ formed by the closures
of the images of the $8$ lines by $T_{1,2}$. The half-skeleton is the closure of the union
of the images of the $4$ special fibers of $\Pi_{2,3}$. The half-skeleton is made of $4$
pairs of projective lines respectively contained in some planes $\{v=\beta\}$. Each line
is defined by $\{u=\alpha,v=\beta\}$ or $\{v=\beta,u'=\gamma\}$. \\

An important point is that it is possible to \emph{describe explicitly this half-skeleton}
using \emph{only some $q$-local monodromy invariants}.\\

We use $u,v,u'\in \C\cup \{\infty\}$ as coordinates on $\left(\PC\right)^3$.
We can describe the closed half-skeleton using some pairs of plane equations
taken from two different lines of the following table~:
\begin{equation}
\label{equaspecialplanes}
\begin{split}
&\{u=0\}, \quad \{u=\infty\}, \quad \{u=\mathbf{e}_q^{1;1,2;3}\}, \quad
\{u=\mathbf{e}_q^{2;1,2;3}\} \\
&\{v=0\}, \quad \{v=\infty\}, \quad \{v=\mathbf{e}_q^{1;2,3;4}\}, \quad
\{v=\mathbf{e}_q^{2;2,3;4}\} \\
&\{u'=0\}, \quad \{u'=\infty\}, \quad \{u'=\mathbf{e}_q^{1;3,4;1}\}, \quad
\{u'=\mathbf{e}_q^{2;3,4;1}\} 
\end{split}
\end{equation}

\begin{prop}
\label{propprojhalfskel}
The half-skeleton $\text{Sk}'_{1,2}(\F)$
is the union the $8$ lines~:
\begin{equation}
\label{equaspecialplanes2}
\begin{split}
\{u=v=0\} , ~~ \{v=u'=0\} &, ~~ \{u=v=\infty\}, ~~ \{v=u'=\infty\}, \\
\{u=\mathbf{e}_q^{1;1,2;3} ,v=\mathbf{e}_q^{2;2,3;4}\} &, ~~ 
\{u=\mathbf{e}_q^{2;1,2;3},v=\mathbf{e}_q^{1;2,3;4}\}, \\
 \{v=\mathbf{e}_q^{2;2,3;4},u'=\mathbf{e}_q^{1;3,4;1}\} &, ~~ 
 \{v=\mathbf{e}_q^{1;2,3;4},u'=\mathbf{e}_q^{2;3,4;1}\}.
\end{split}
\end{equation}
\end{prop}

The index $2,3$ and $1,4$ do not appear symmetrically in the definition of $T_{1,2}$.
It is more difficult to understand the images of the $8$ lines~:
\begin{equation}
\label{equabadlines}
L_{\rho_1,x_1}, ~ L_{\rho_2,x_1}, ~
L_{\rho_1,x_4}, ~ L_{\rho_2,x_4}, ~
L_{\sigma_1,x_2}, ~ L_{\sigma_2,x_2},
L_{\sigma_1,x_3}, ~ L_{\sigma_2,x_3}. 
\end{equation}
For example, we have $L_{\sigma_1,x_1}\subset \Pi_{1,2}^{-1}(0)$, therefore its image
by $T_{1,2}$ is contained into the plane $\{u=0\}=\{0\}\times \left(\PC\right)^2$.
But we have only one exceptional fiber in the picture and therefore we know a priori
\emph{only one plane} into $\left(\PC\right)^3$ containing the images
$T_{1,2}(L_{\sigma_1,x_1})$, \dots  \\

\paragraph{A first heuristic description of the image.}
\label{heuristicdescriptionp3}

It seems difficult to say more about the image $Y=T_{1,2}(\F)$ and in particular about
the skeleton $\text{Sk}_{1,2}(\F)$ without heavy computations involving the $q$-pants
charts. We plan to return to the question in a future work. Here we will only \emph{try
some guesses}.

\begin{conj}
The lines of the half-skeleton are double curves.
\end{conj}

More precisely, in a neighborhood of a smooth point of the double curve, we have in local
coordinates $Y = \{xy = 0\}$ (mild singularity). There could also exist pinch points (in
a neighborhood of a pinch point, $Y = \{x^2 - y z^2 =0\}$) and a finite number of more
complicated singular points. \\

The image of the half-skeleton $\text{Sk}'_{1,2}(\F)$ by the projection
$\left(\PC\right)^3\rightarrow \left(\PC\right)^2$ defined by $(u,v,u')\rightarrow (u,u')$
is a union of $8$ lines~:
\begin{equation}
\label{equaeightlinesbranch}
\begin{split}
& \{u=0\}, \quad \{u=\infty\}, \quad \{u=\mathbf{e}_q^{1;3,4;1}\}, \quad
\{u=\mathbf{e}_q^{2;1,3;4}\} \\
& \{u'=0\}, \quad \{u'=\infty\}, \quad \{u'=\mathbf{e}_q^{2;3,4;1}\}, \quad
\{u'=\mathbf{e}_q^{1;1,3;4}\} \\
\end{split}
\end{equation}
This image depends only on the four complex numbers $\mathbf{e}_q^{1;3,4;1}$,
$\mathbf{e}_q^{2;1,3;4}$, $\mathbf{e}_q^{2;3,4;1}$ and $\mathbf{e}_q^{1;1,3;4}$.
It is also clearly the image of the skeleton by the projection (the lines
(\ref{equabadlines}) are contained in special planes $\{u=\alpha\}$ or $\{u'=\beta\}$). \\

We recall the following definition.

\begin{defn}
Let $V$, $W$ two complex algebraic varieties. A morphism $f:V\rightarrow W$ is
a \emph{branched covering} if the two dimensions are the same and if the typical
fiber of $f$ is of dimension $0$.
\end{defn}

There is a Zariski dense open set $W'\subset W$ such that $f$ is unramified above $W'$
(a classical covering space). The complement of the largest possible $W'$ is called
the \emph{branching locus}. If $W'$ is connected then the cardinal of the fiber is
constant, it is the \emph{degree} of the branched covering. \\

Be careful, the classical notion of ramified covering is more restrictive: \emph{all
the fibers} are finite sets.

\begin{conj}
\label{conjk3first}
The map $\pi_2$ induced by the restriction to $Y$ of the projection 
$\left(\PC\right)^3\rightarrow \left(\PC\right)^2$ defined
by $(u,v,u')\rightarrow (u,u')$ is a branched double covering of
$\left(\PC\right)^2$. The branching set is
the set of eight lines defined by (\ref{equaeightlinesbranch}).
\end{conj}

We will see later (\cf\ page \pageref{fandkummer}) that there exists a K3 surface
(of Kummer type) $X'$ which is a double covering of $\left(\PC\right)^2$ branched
along the same set of $8$ lines. This suggest that $X'$ could be, up an isomorphism,
a projective completion of $\F$.

\paragraph{Another heuristic description of the image.}
\label{heuristicdescriptionp4}

There is another, more precise, heuristic description of the image. It is based on
the following conjecture.

\begin{conj}
\begin{itemize}
\item[(i)]
The surface $Y=\overline{T_{1,2}(\F)}$ is a $(2,2,2)$ surface of $\left(\PC\right)^3$. 
\item[(ii)]
The restriction of $T_{1,2}$ to $\F$ minus the $16$ lines is \emph{injective}.
\end{itemize}
\end{conj}

We return to the $12$ plane equations of (\ref{equaspecialplanes}). Each plane
$\{u=\alpha\}$, resp. $\{v=\beta\}$, resp. $\{u'=\gamma\}$ cuts $Y$ along a $(2,2)$ curve. \\

When $\beta$ moves there are four values such that the $(2,2)$ curve is decomposed into
two lines (with a common point) and these lines are necessary \emph{double} lines. The
four values are~: $\beta=0,~\infty~, ~ \mathbf{e}_q^{1;2,3;4}, ~\mathbf{e}_q^{2;2,3;4}$.
Each pair of lines is picked up in the list formed by the line $L_{\rho_1,x_2}$ and the lines
(\ref{equagoodlines}). More precisely the union of the $4$ pairs of lines is the
half-skeleton. \\

We have a similar situation when $\alpha$ and $\gamma$ move but we do not know \emph{a priori}
the four exceptional values of decomposition. When the $(2,2)$-curve is decomposed into two
lines, it is decomposed into a double line of the closed half-skeleton and \emph{another
double line}. According to proposition \ref{propcoverqpants}, the only possibilities for
such a double line seem to be~:
\begin{equation}
\label{4extralines}
\{u=0,u'=\infty\}, ~~ \{u=\infty,u'=0\}, ~~
\{u=\mathbf{e}_q^{1;3,4;1},u'=\mathbf{e}_q^{1;1,3;4}\}, ~~
\{u=\mathbf{e}_q^{2;3,4;1},u'=\mathbf{e}_q^{2;1,3;4}\},
\end{equation}
We are thus led to the following conjecture.

\begin{conj}
\label{skeleton12}
\begin{itemize}
\item[(i)]
The skeleton $\text{Sk}_{1,2}(\F)$ is a union of $12$ lines: the $8$ lines of the
half-skeleton $\text{Sk}'_{1,2}(\F)$ and the $4$ lines (\ref{4extralines}).
\item[(ii)]
The surface $Y$ is mildly singular along each line of the skeleton. 
\end{itemize}
\end{conj}

In order to put the set of the $16$ lines into $\left(\PC\right)^3$, it is necessay
to ``fold it". The skeleton has two connected components, it can be described as a
``split parallellepipedal structure". Each connected component is a deformed hexagonal
structure. \\
 
Moreover we can conjecture that there are $4$ pinch points on each line and that there
exist $12$ exceptional planes $\{u=\alpha\}$, $\{v=\beta\}$, $u'=\gamma\}$, such that
each one contains $4$ pinch points. This configuration seems to be related to the
logarithmic fibers. \\

We have a description of the fibrations of the surface by the coordinates. We detail it
for the coordinate $v$. When $\beta$ moves we have $3$ types of fiber~:
\begin{itemize}
\item 
the generic fiber is a $(2,2)$ curve with two nodes; it is decomposed into two $(1,1)$ curves;
\item
there are $4$ fibers decomposed into two double lines;
\item
there are $4$ fibers corresponding to the planes containing pinch points, they are
double $(1,1)$ curves.
\end{itemize}

According to the above description, we can verify that the images of the skeleton by
the projections $(u,v,u')\rightarrow (u,u')$, $(u,v,u')\rightarrow (u,v)$ and
$(u,v,u')\rightarrow (v,u')$ are sets of $8$ lines that we can explicit using only
the $q$-local monodromy invariants. \\

Then it is easy to prove that the restriction of each projection to $Y$ is a double
covering of $\left(\PC\right)^2$ branched \emph{exactly} along $8$ lines, a degenerated
$(4,4)$-curve (\cf\ conjecture\footnote{Now this conjecture follows from conjecture
\ref{skeleton12}} \ref{conjk3first}). We consider for example the first projection,
the ramification set is a $(4,4)$ curve of $\left(\PC\right)^2$ and it contains the $8$
lines defined by (\ref{equaeightlinesbranch}). The union of these $8$ lines is also a
$(4,4)$ curve, therefore we have equality.


\subsubsection{Conjectural embeddings of $\F$ into $\left(\PC\right)^6$ and into
$\left(\PC\right)^4$}
\label{subsubmainconj1}

It seems difficult to get an \emph{embedding} of $\F$ into $\left(\PC\right)^3$.
We can try to do better\footnote{We already know that there exists a embedding
of $\F$ in $\left(\PC\right)^4$, \cf\ part \ref{rough4surf}.} with maps which
involve more symmetrically the $\Pi_{i,j}$.

\begin{conj}
\label{conjembed}
We suppose that (FR), (NR), (NS) and $\mathbf{Hyp}_{48}$ are satisfied.
\begin{itemize}
\item[(i)]
The regular map~:
\[
T:=(\Pi_{1,2},\Pi_{2,3},\Pi_{3,4},\Pi_{1,4},\Pi_{1,3},\Pi_{2,4}):
\mathcal{F}\mapsto \left(\PC\right)^6
\]
defined by~:
\[
u:=\Pi_{1,2}(\bar M), ~~ v:=\Pi_{2,3}(\bar M), ~~
u':=\Pi_{3,4}(\bar M), ~~ v':=\Pi_{1,4}(\bar M), ~~
w:=\Pi_{1,3}(\bar M), ~~w':=\Pi_{2,4}(\bar M)
\]
is a regular embedding\footnote{That is the image is an affine surface into
$\left(\PC\right)^6$ and if we endow it with the induced Zariski topology, then
it is isomorphic to $\F$.}.
\item[(ii)]
The regular map~:
\[
T'_{1,2}:=(T_{1,2},\Pi_{1,4})=(\Pi_{1,2},\Pi_{2,3},\Pi_{3,4},\Pi_{1,4}): \mathcal{F}
\mapsto \left(\PC\right)^4
\]
defined by~:
\[
u:=\Pi_{1,2}(\bar M), ~~ v:=\Pi_{2,3}(\bar M), ~~
u':=\Pi_{3,4}(\bar M), ~~ v':=\Pi_{1,4}(\bar M)
\]
is a regular embedding. Similarly the maps $(\Pi_{1,2},\Pi_{1,3},\Pi_{3,4},\Pi_{2,4})$ and
$\Pi_{2,3},\Pi_{1,4},\Pi_{1,3},\Pi_{2,4}$ are regular embeddings.
\end{itemize}
\end{conj}

We end with a conjectural picture. 

\begin{conj}
\label{mainconj1}
Let $(\underline{\rho},\underline{\sigma},\underline{x})$ such that (FR), (NR), (NS) and 
$\mathbf{Hyp}_{48}$ are satisfied. We denote $u,v,u',v'\in \C\cup \infty$ coordinates
on $\left(\PC\right)^4$. Then there exist three polynomials $f_1,~f_2,~f_3$ ``on''
$\left(\PC\right)^3$, of tri-degree $(2,2,2)$, such that~:
\begin{itemize}
\item[(i)]
the $3$ equations~:
\[ 
\tilde f_1(u,v,u')=0, ~~ \tilde f_2(u,v,v')=0, ~~ \tilde f_3(u',v,v')=0
\]
define a \emph{smooth} surface $X$ of $\left(\PC\right)^4$;
\item[(ii)]
$\F(\underline{\rho},\underline{\sigma},\underline{x})$ is isomorphic to a Zariski open
subset of $X$.
\item[(iii)]
We can choose the coefficients of the polynomials $f_i$ as functions of
$(\underline{\rho},\underline{\sigma},\underline{x})$ in such a way that they depend
only on the $q$-local monodromy invariants, this dependance being rational.
\end{itemize}
\end{conj}

Each equation $\tilde f_i=0$ defines a $(2,2,2)$ surface of $\left(\PC\right)^3$. As
explained before, we think that this surface is \emph{singular}. We will suggest below
a method of computation of $\tilde f_i$ (\cf\ page \pageref{enriquesstylesurf}).


\subsection{K3 surfaces and conjectural description of $\F$}
\label{subsecfamilyk3}

This part is a stub and it contains mainly heuristics. However it could open some pathes
towards a clear synthesis of all the rigorous (but complicated \dots) informations that
we got on the surface $\F$. We plan to return to these questions in a future work.


\subsubsection{Definitions and exemples}

We recall the following definitions  \cite{HS}.

\begin{defn}
\begin{itemize}
\item[(i)]
A complex smooth projective surface $X$ is called \emph{K3 surface} if $X$ is simply
connected with trivial canonical bundle $\omega_X\approx \mathcal{O}_X$.
\item[(ii)]
An Enriques surface is a quotient of a K3 surface $X$ by a fixed point free involution
$\iota$ (called an Enriques involution).
\end{itemize}
\end{defn}

There exists a symplectic $2$-form on the K3 surface $X$ (unique up to a multiplicative
constant). We have\footnote{Another equivalent definition of a K3 surface, due to Andr\'e
Weil, around 1948, is $\omega_X\approx \mathcal{O}_X$ and $H^1(X;\mathcal{O}_X)=0$.} 
$H^1(X;\mathcal{O}_X)=0$, $H^1(X;\Z)=0$ and the rank of $H^2(X;\Z)$ is $22$. We recall
that the \emph{Betti numbers} $b_r(X)$ of a surface $X$  are the integers defined by
$b_r(X):=\dim_{\Q} H^r(X;\Q)$ (they are topological invariants). A surface $X$ is K3 if
and only if its canonical bundle is trivial and if $b_1(X)=0$. If $X$ is K3, then we
have $\dim H^2(X;\Z)=b_2(X)=22$. \\

Enriques surfaces and K3 surface have a null Kodaira dimension. They are not rational surfaces.\\

We recall the adjunction formula.

\begin{prop}[Adjunction formula]
If\, $Y\subset X$ is an hypersurface, with $X$ and $Y$ smooth, 
then we have the two equivalent equalities~:
\begin{equation}
\label{formadj}
K_Y=(K_X+Y)_{\vert Y}  \quad \text{and} \quad \omega_Y
=\left(\mathcal{O}(Y)\otimes \omega_X   \right)_{\vert Y} 
\end{equation}
\end{prop}

We recall the following result.

\begin{prop}
\label{-2line}
Let $L$ be a smooth curve of genus $g$ contained in a K3 surface, then $(L,L)=2g-2$;
$L$ is rational if and only if $(L,L)=-2$. 
\end{prop}
\Pr
This result follows from the adjunction formula:
$\omega_L=\mathcal{O}_L(L)$, therefore $2g-2=(L,L)$.
\Finpr

One can prove that if $L$ is a $(-2)$ curve in a K3 surface, then $L$ is smooth. \\

It follows from the above proposition that there does not exist $(-1)$ smooth curves on a K3
surface.

\begin{exa}
Every smooth quartic surface in $\P^3(\C)$ is a K3 surface. \\
Let $m\in N^*$. We have $\omega_{\P^m(\C)}=\mathcal{O}_{\P^m(\C)}(-m-1)$. Let 
$X\subset \P^m(\C)$ be a smooth hypersurface defined by an homogeneous polynomial
of degree $d$. By the adjunction formula~:
\[
\omega_X=\left(\omega_{\P^m(\C)}\otimes 
\mathcal{O}_{\P^m(\C)}(d)\right)_{\vert X}=\mathcal{O}_X(-m-1+d).
\]
Therefore, if $m=3$ and $d=4$, then $\omega_X=\mathcal{O}_X$. \\
The projective space $\P^3(\C)$ is simply connected, therefore, by Lefschetz theorem,
$X$ is also simply connected. \\
We can also prove that $H^1(X;\mathcal{O}_X)=0$ using the short exact sequence~:
\[
0\rightarrow \mathcal{O}_{\P^3(\C)}(-4)\rightarrow 
\mathcal{O}_{\P^3(\C)}\rightarrow \mathcal{O}_X\rightarrow 0
\]
and 
$H^1(\P^3(\C);\mathcal{O}_{\P^3(\C)})=H^2(\P^3(\C);\mathcal{O}_{\P^3(\C)}(-4))=0$. \\
An interesting example of a smooth quartic hypersurface in $\P^3(\C)$ is the Fermat
quartic~: $X^4+Y^4+Z^4+T^4=0$
\end{exa}

\begin{exa} [Double plane]
\label{exadoubleplane}
For this example \cf\ \cite{Huyb}. \\
Consider a double covering $\pi:X\rightarrow P^2(\C)$ branched along a sextic curve
$C\subset P^2(\C)$. Then
$\pi_*(\mathcal{O}_X)\approx   \mathcal{O}_{P^2(\C)}\oplus \mathcal{O}(-3)$ and therefore
$H^1(X;\mathcal{O}_X)=0$. \\  
We suppose that the branching curve $C$ is non-singular, then $X$ is non singular and
the canonical bundle formula for branched coverings shows that
$\omega_X=\pi^*(\omega_{P^2(\C)}\oplus \mathcal{O}(3))\approx \mathcal{O}_X$.
Therefore $X$ is a K3 surface, called a \emph{double plane}. \\
If the sextic $C$ is the union of $6$ generic lines in $P^2(\C)$, the double cover $X$
has $15$ rational double points. These $15$ points correspond to the pairwise intersections
of the $6$ lines. Blowing-up these $15$ singular points produces a K3 surface.
\end{exa}

\begin{prop}
The smooth $(2,2,2)$ surfaces are the K3 surfaces embedded in $\left(\PC\right)^3$.
\end{prop}
\Pr
Let $X\subset \left(\PC\right)^3$ be a smooth irreducible hypersurface of tri-degree
$(a,b,c)$. \\
Using the adjunction formula, we get $K_X=\mathcal{O}_X(a-2,b-2,c-2)$ and $K_X$ is trivial
if and only if $(a,b,c)=(2,2,2)$. \\
We suppose $(a,b,c)=(2,2,2)$. The fiber bundle $[X]$ is positive. Using the Lefschetz theorem
on hyperplane sections we get an isomorphism
$H^1\left(\left(\PC\right)^3;\Q\right)\rightarrow H^1(X;\Q)$
\big(induced by the canonical injection $X\rightarrow \left(\PC\right)^3$\big).
As $H^1\left(\left(\PC\right)^3;\Q\right)=0$, we have also $H^1(X;\Q)=0$ and $b_1(X)=0$.
\Finprcourt


\subsubsection{The Enriques surface and some surfaces in the same style}
\label{subsubenriquessurf}

\paragraph{The Enriques surface.}

Around 1895, after many discussions with G. Castelnuovo under the arcades of the city of
Bologna, F. Enriques discovered a very interesting surface \cite{Enr1,Enr2,Cast}. We quote
\cite{Enr2}~:

\emph{Nel 1896 mi si \`e presentata la superficie del $6^{\circ}$ ordine passante
  doppiamente per gli spigoli di un tetraedro come primo esempio di superficie di
  genere $p_g=p_a=0$, non razionale.}

The example of Enriques \cite{Dol} is a smooth normalization of a 
non-normal sextic surface Y in $\mathbf{P}^3(\C)$ that passes with multiplicity $2$ through
the edges of the coordinate tetrahedron. Its equation (in projective coordinates) is~:
\begin{equation}
\label{equaenriquessurf}
F := x^2_1x^2_2x^2_3+x^2_0x^2_2x^2_3+x^2_0x^2_1x^2_3
+x^2_0x^2_1x^2_2+x_0x_1x_2x_3\, q(x_0,x_1,x_2,x_3)=0,
\end{equation}
where $q$ is a non-degenerate quadratic form. \\

The surface $Y$ has the following singularities: a double curve $\Gamma$ with ordinary
triple points which are also triple points of the surface and some pinch points ($4$ on
each edge of the tetrahedron). \\

We choose an edge $\Delta$ of the tetrahedron. The family of planes passing by $\Delta$
cut the surface along a sextic curve. This sextic curve is decomposed into the double
line $\Delta$ and a quartic curve. There are the following exceptional cases~:
\begin{enumerate}
\item 
the plane is a face of the tetrahedron, then the quartic is decomposed into two double
lines (there are two such cases);
\item
the quartic pass by a pinch point on $\Delta$, then there appears a cusp (there are four
such cases).
\end{enumerate}
The quadric curves form a pencil. A base locus of this pencil is the union of $4$ edges,
excluding the union of $2$ opposite edges. \\

It is interesting to compare the above fibrations with the following fibrations~:
\begin{itemize}
\item 
the fibration of the cubic surface $\overline{\mathcal{S}_{VI}(a)}$ described in paragraph
\emph{A fibration}, page \pageref{fibrationdiff};
\item
the fibration of the image $Y=T_{1,2}(\F)$ of $\F$ into $\left(\PC\right)^3$ by the planes
$\{v=\text{cste}\}$;
\item
the fibration of $\F$ by each $\Pi_{ij}$.
\end{itemize}

More generally we can consider \emph{all} the sextic surfaces $Z$ in $\mathbf{P}^3(\C)$
mildly singular along the $6$ edges of the tetrahedron $x_0x_1x_2x_3=0$ (\cf\ \cite{Muk}).
An equation of such a surface is~:
\begin{equation}
\label{equaenriquessurf2}
F_{a,b,c,d;q} := ax^2_1x^2_2x^2_3+bx^2_0x^2_2x^2_3+cx^2_0x^2_1x^2_3
+x^2_0x^2_1x^2_2+x_0x_1x_2x_3\, q(x_0,x_1,x_2,x_3)=0,
\end{equation}
where $q$ is a non-degenerate quadratic form. There are four sextic monomials and the $10$
quadratic monomials of $q$. Quotienting by an action of $\left(\C^*\right)^4$, we get a $10$
parameters family. \\

Let $S$ be the normalization of $Z$. We have the following result (\cf\ \cite{Muk},
proposition 4.1, page 5).

\begin{prop}
The surface $S$ is an Enriques and the covering K3 surface $X$ is a $(2,2,2)$ surface in
$\left(\mathbf{P}^1(\C)\right)^3$ which is invariant by the involution
$(u,v,w)\rightarrow (-u,-v,-w)$.
\end{prop}
More precisely an equation of the surface $X$ is~:
\[
au^2v^2w^2+bu^2=cv^2+dw^2+uvw\, q(1,vw,uw,uv).
\]

\paragraph{Some $(2,2,2)$-surfaces in Enriques style.}
\label{enriquesstylesurf}

We consider a generalized version $\text{Sk}'$ of the half-skeleton, the set of $8$ lines
in $\left(\mathbf{P}^1(\C)\right)^3$~:
\begin{equation}
\label{equaspecialplanes4}
\begin{split}
\{u=\alpha_1,v=\beta_1\} &, ~~ \{v=\beta_1,u'=\gamma_1\}, ~~ 
\{u=\alpha_2,v=\beta_2\}, ~~ \{v=\beta_2,u'=\gamma_2\} \\
\{u=\alpha_3,v=\beta_4\} &, ~~ \{v=\beta_4,u'=\gamma_3\}, ~~
\{u=\alpha_4,v=\beta_3\}, ~~ \{v=\beta_3,u'=\gamma_4\}.
\end{split}
\end{equation}
parameterized by the $12$ ``numbers" 
$\underline{\alpha},~\underline{\beta},~\underline{\gamma}$~: each of
$\underline{\alpha}$ (resp. $\underline{\beta}$, resp. $\underline{\gamma}$) is an
arbitrary triple of \emph{distinct} elements of $\PC$. Up to M\"obius transformations
on each factor of $\left(\mathbf{P}^1(\C)\right)^3$ it is sufficient to consider the case~:
$\underline{\alpha}:=(\alpha_1,1,0,\infty)$,
$\underline{\beta}:=(\beta_1,1,0,\infty)$,
$\underline{\gamma}:=(\gamma_1,1,0,\infty)$. Then it remain only $3$ parameters. 
We will write $\text{Sk}'_{\underline{\alpha},~\underline{\beta},~\underline{\gamma}}$ if we want
to precise the parameters. \\

The half-skeleton corresponds to 
\begin{equation}
\label{equahskpecial}
\underline{\alpha}=(0,\infty,\mathbf{e}_q^{1;3,4;1},\mathbf{e}_q^{2;1,3;4}),~
\underline{\beta}=(0,\infty,\mathbf{e}_q^{1;2,3;4},\mathbf{e}_q^{2;2,3;4}) ~
\underline{\gamma}=(0,\infty,\mathbf{e}_q^{2;3,4;1},\mathbf{e}_q^{1;1,3;4}).
\end{equation} 

The projection of $\text{Sk}'_{\underline{\alpha},~\underline{\beta},~\underline{\gamma}}$ on
the $(u,u')$-plane is the following set of $8$ coordinates lines:

\begin{equation}
\label{equaspecialplanes3}
\begin{split}
& \{u=\alpha_1\}, \quad \{u=\alpha_2\}, \quad 
\{u=\alpha_3\}, \quad
\{u=\alpha_4\} \\
& \{u'=\gamma_1\}, \quad \{u'=\gamma_2\}, \quad 
\{u'=\gamma_3\}, \quad
\{u'=\gamma_4\}.
\end{split}
\end{equation}
 
There is exactly one line of $\text{Sk}'_{\underline{\alpha},~\underline{\beta},~\underline{\gamma}}$
above each line of (\ref{equaspecialplanes3}). \\

We can also generalise the skeleton: it is an union of $12$ lines denoted
$\text{Sk}_{\underline{\alpha},~\underline{\beta},~\underline{\gamma}}$. It is easy to write the lines
equations using $\underline{\alpha},~\underline{\beta},~\underline{\gamma}$
(\cf\ (\ref{4extralines}). \\

Enriques considered the family of sextic surfaces of $\P^3(\C)$ mildly singular along
the six edges of the tetrahedron $xyzt=0$. Similarly we will consider the (possibly
empty) family 
$\left\{\mathcal{S}_{\underline{\alpha},~\underline{\beta},~\underline{\gamma}}\right\}_{\underline{\alpha},~\underline{\beta},~\underline{\gamma}}$
of $(2,2,2)$ surfaces of $\left(\PC\right)^3$ mildly singular along the $12$ lines of 
$\text{Sk}_{\underline{\alpha},~\underline{\beta},~\underline{\gamma}}$. \\

Let $h_{\underline{\alpha},~\underline{\beta},~\underline{\gamma}}$ be a polynomial of tri-degree
$(2,2,2)$ such that
$\mathcal{S}_{\underline{\alpha},~\underline{\beta},~\underline{\gamma}}=V(h_{\underline{\alpha},~\underline{\beta},~\underline{\gamma}})$. \\

We can write some \emph{necessary conditions} on
$h_{\underline{\alpha},~\underline{\beta},~\underline{\gamma}}$. 

\begin{equation}
\label{linearsystemforequation}
\begin{split}
h_{\underline{\alpha},~\underline{\beta},~\underline{\gamma}} &=
a_1(u-\alpha_1)^2(u'-\gamma_1)^2+(v-\beta_1)^2(\cdots)
=a_2(u-\alpha_2)^2(u'-\gamma_2)^2+(v-\beta_2)^2(\cdots) \\
&=a_3(u-\alpha_3)^2(u'-\gamma_3)^2+(v-\beta_3)^2(\cdots),
\end{split}
\end{equation}
where each $(\cdots)$ is a $(2,2)$ polynomial in $(u,u')$ ($6$ monomials). There are similar
conditions replacing $(u,u')$ by $(u,v)$ or $(v,u')$. \\

Using the above conditions we can write a linear system where the unknown are the
coefficients of $h_{\underline{\alpha},~\underline{\beta},~\underline{\gamma}}$. Solving this
system one will obtain the family\footnote{It could be empty for generic values of
$\underline{\alpha},~\underline{\beta},~\underline{\gamma}$. For the values associated
to the $q$-monodromy invariants, we can conjecture that there exists an unique solution
up to scaling.} $\left\{\cal{HS}_{\underline{\alpha},~\underline{\beta},~\underline{\gamma}}\right\}$.
We can use M\"obius transformations in order to simplify the system. \\

Using resultants, we can write another linear system. We denote $R_t(P,Q)$ the resultant
of two polynomials in $t$. We write for simplicity
$h=h_{\underline{\alpha},~\underline{\beta},~\underline{\gamma},~\underline{\alpha}',~\underline{\gamma}'}$
Then we can consider the three resultants~:
\[
R_2:=R_v\left(h,\frac{\partial h}{\partial v}\right), ~~
R_1:=R_u\left(h,\frac{\partial h}{\partial u}\right)=0, ~~
R_3:=R_{u'}\left(h,\frac{\partial h}{\partial u'}\right)=0
\]
and write that they vanish respectively on three systems of $8$ lines. We get a linear
system: the unknown are the coefficients of $h$. A solution of the first system is clearly
a solution of the second. We do not know if they are equivalent. \\

We return to conjecture \ref{mainconj1}. 
If conjecture \ref{mainconj1} is true, then $Y$ is a $(2,2,2)$-surface mildly
singular along the $12$ lines of the skeleton and one can use the above method to get
an equation of $Y$ into $\left(\PC\right)^3$, or equivalently an algebraic relation
between $\Pi_{1,2}$, $\Pi_{2,3}$, and $\Pi_{3,4}$. 

\paragraph{The surface $\F$ and the Kummer surfaces.}
\label{fandkummer}

The $(2,2,2)$ surface $Y$ is a double covering of $\left(\PC\right)^2$ \emph{branched}
along the $8$ lines. We can conjecture that there exists a smooth projective completion
$X$ of $\F$ which is also a double covering of $\left(\PC\right)^2$ \emph{branched} along
the $8$ lines. \\

We will explain how to \emph{compute} a double covering of $\left(\PC\right)^2$ ramified
along the $8$ lines. It is a K3 surface, more precisely a K3 surface of Kummer type. It
is a good candidate for a projective completion of $\F$ (up to an isomorphism). \\

For simplicity we denote~:
\begin{equation}
\label{equaeightlinesbranch2}
\begin{split}
& \{u=\alpha_1\}, \quad \{u=\alpha_2\}, \quad \{u=\alpha_3\}, \quad
\{u=\alpha_4\} \\
& \{u'=\gamma_1\}, \quad \{u'=\gamma_2\}, \quad \{u'=\gamma_3\}, \quad
\{u'=\gamma_4\} \\
\end{split}
\end{equation}
the equations of the $8$ lines. \\

Let $p:Z\rightarrow \left(\PC\right)^2$ be a double \emph{ramified} covering, ramified
on the $8$ lines. The $16$ double points of the ramification locus are $(\alpha_i,\gamma_j)$.
We consider double ramified coverings $A\rightarrow \PC$ and $B\rightarrow \PC$ respectively
ramified above $(\alpha_1,\alpha_2,\alpha_3,\alpha_4)$ and
$(\gamma_1,\gamma_2,\gamma_3,\gamma_4)$, $A$ and $B$ being \emph{elliptic curves}. Then $Z$
is isomorphic to the quotient of $A\times B$ by the canonical involution $x \mapsto -x$.
It is a \emph{Kummer surface}. It is isomorphic to a nodal quartic surface into $\P^3(\C)$
with $16$ nodal points. Blowing up at the $16$ nodes, we get a K3 surface $X'$. We conjecture
that $X'$ is isomorphic to a projective completion of $\F$. \\

The K3 surface $X'$ is \emph{uniquely determined}, up to an isomorphism, by the two
cross-ratios $(\alpha_1,\alpha_2,\alpha_3,\alpha_4)$ and 
$(\gamma_1,\gamma_2,\gamma_3,\gamma_4)$ and a fortiori by the $4$ complex numbers~:
\[
\mathbf{e}_q^{1;3,4;1}(\underline{\rho},\underline{\sigma},\underline{x}), ~~ \mathbf{e}_q^{2;1,3;4}(\underline{\rho},\underline{\sigma},\underline{x}), ~~
\mathbf{e}_q^{2;3,4;1}(\underline{\rho},\underline{\sigma},\underline{x}) ~~ \text{and} ~~ 
\mathbf{e}_q^{1;1,3;4}(\underline{\rho},\underline{\sigma},\underline{x}). 
\]

This model could be a $q$-analog of the "algebraic dependence" of the smooth projective
cubic surface $\tilde{\mathcal{S}}_{A_0,A_t,A_1,A_\infty}$ on $(A_0,A_t,A_1,A_\infty)\in \C^4$.


 
\section{Conclusion: open questions and perspectives}
\label{section:conclusion}

We have nearly achieved our initial aim. We built the character variety of $q$-PVI
and gave a quite precise description of this variety. There remain some open problems.
More generally there is a lot of related questions and possible generalizations. We
will give a (non exhaustive\footnote{In particular we will not discuss the important
problems of symplectic structures.} \dots) list.

 
\subsection{Generalized versions of Riemann-Hilbert map}

In the differential case the Riemann-Hilbert map is a complex analytic morphism 
$RH: \mathcal{M} \rightarrow \widetilde{\text{Rep}}$ from a moduli space $\mathcal{M}$
of connections to a (categorical) moduli space of (generalized) monodromy data. More
precisely from a \emph{family} of moduli space of connections to a \emph{family} of
moduli space of monodromy data. In the fuchsian case (i. e. PVI) the parameters on
the left hand side are $(t,\theta)$ and $(t,a)$ on the right hand side 
($a_l=2\cos 2\pi \theta_l$). \\

In the irregular case it is necessary to add some generalized exponents into the
parameters and Stokes multipliers into the monodromy data \cite{MRAcceleration},
\cite{vdPSaito}. \\

Painlev\'e equations are derived from holomorphic flows on 
$\mathcal{M}$. The flows are transversal to the parameter fibration. The fibers are
the Okamoto spaces of initial conditions. \\

The above picture works perfectly in the PVI case for \emph{generic values} of the
parameters and we have generalized it to the $q$-PVI case for \emph{fixed generic values
of the parameters}. In the differential PVI case if one wants to allow the exceptional
values of the parameters, then it is necessary to replace connections by \emph{parabolic
connections}\footnote{Intuitively one ``adds a line".} \cite{IIS,Inaba}. There are
possible generalizations of our work (or of part of our work).
\begin{itemize}
\item 
The case of a fixed exceptional parameter. It will be necessary
to use the \emph{parabolic $q$-difference modules} of
Mochizuki \cite{Mochi}.
\item
The case of $q$-PVI \emph{with parameters}. We remark that even the simpler hypergeometric
case with parameters is not known.
\item
The case of the equations of the Murata's list with or without parameters.
\end{itemize}

We conjecture that it is possible to extend \emph{part of our results} to all the
equations of Murata's list (taking account of $q$-Stokes phenomena). In sharp contrast,
there is no hope to extend the Mano decompositions for \emph{all} the equations. \\

In \cite{JimboMonodromy} Jimbo give a splitting of PIII and PV respectively into~:
\begin{itemize}
\item 
two confluent hypergeometric equations,
\item
an hypergeometric equation and a confluent hypergeometric equation.
\end{itemize}
For a more detailed description of splittings of Painlev\'e equations, cf. \cite{GL},
\emph{Figure 3: CMR confluence diagram for Painlev\'e equations}. \\

It is possible to extend Mano result for PVI=P(A3) to the equations \text{P(A4)},
$\text{P(A5)}^\sharp$, $\text{P(A6)}^\sharp$. One gets respective splittings
into\footnote{\cf\ \cite{Ohyama} for basics on irregular $q$-hypergeometric equations.}~:
\begin{itemize}
\item
a $q$-Kummer and a Heine $q$-hypergeometric equation,
\item
two $q$-Kummer equations,
\item
a $q$-Kummer and a Hahn Exton $q$-Bessel equation. 
\end{itemize}

We conjecture that it is possible to extend our Mano decomposition for these three cases.


\subsection{Relations with $q$-difference Galois groups}

As we said above, defining local monodromies and local Galois groups at intermediate
singularities for $q$-difference equations is one of the most important open problems in
modern $q$-difference theory. In some sense it would close the problem of ``localisation''
of Galois groups: the problem of ``localisation'' of the Galois groups a $0$ and $\infty$
(\ie\ the description of the corresponding $q$-wild groups) was solved in full generality
in a series of papers of the two last authors \cite{RS3,JSTALPA}. \\

This will require a more general version of Mano decomposition. Extension to higher degrees
should be easy along the same lines, but extension to higher orders (polynomial matrices
with coefficients in $\Matnc$) seems more difficult. Moreover the ``basic bricks'' are not
clear.


\subsection{Exceptional lines and points on $q$-character varieties and exceptional
solutions of $q$-Painlev\'e equations}

In the differential case there is a fundamental heuristic principle: there is a dictionary
between the asymptotics of a solution of a Painlev\'e equation at the singular points and
some ``natural coordinates" of the corresponding point on the character variety. \\

This principle is illustrated for PVI by Jimbo formula \cite{JimboMonodromy}, \cite{Boa5}
(Appendix B). For the others Painlev\'e equations there is a lot of precise results in
this direction into the book \cite{FIKN}. \\

In \cite{KliPR}, M. Klimes, E. Paul and the second author propose another principle in
the same direction: the lines on the character variety (an affine cubic surface) correspond
to one parameter families of ``special'' solutions of the Painlev\'e equation, an
intersection of two such lines corresponds to a ``very special'' solution. A good
illustration is PII: there are $9$ lines, they correspond to Boutroux tronqu\'ees
solutions. The intersections of two lines correspond to tritronqu\'ees solutions or
to bitronqu\'ees solutions\footnote{Hahn-Exton style solutions.}. \\

It could be interesting to look at $q$-analogs, in particular about the $16$ exceptional
lines that we exhibited on $\F$.


\subsection{$q$-deformations of CFT}

During the last years appeared some papers about possible $q$-deformations of Conformal
Field Theory (CFT). In this context ordinary differential equations are replaced by
$q$-difference equations. We think that $q$-characters varieties, Mano-decompositions
and $q$-pants parametrizations could be useful (cf. in particular \cite{JNS}).
The irregular case ($q$-Stokes phenomena and $q$-sommations) seems interesting in such
approaches \cite{ABBDJJM}. \\

We quote \cite{GL} (cf. page 7).

\emph{Perhaps the most intriguing perspective is to extend our setup to $q$-isomonodromy
problems, in particular $q$-difference Painlev\'e equations, presumably related to the
deformed Virasoro algebra \cite{SKAO} and 5D gauge theories. Among the results pointing
in this direction, let us mention a study of the connection problem for $q$-Painlev\'e VI
\cite{Mano}
based on asymptotic factorization of the associated linear problem into two systems
solved by the Heine basic hypergeometric series ${}_2 \phi_1$ , and critical expansions
for sollutions of $q$-P( A1 ) equation recently obtained in \cite{JR}.} \\

We quote \cite{Tab}.

\emph{Localization techniques for supersymmetric quantum field theories allow one to
produce non-perturbative results such as computing partition functions exactly, in stark
contrast to general field theories. In many two-dimensional examples of supersymmetric
theories, the path integral or partition function is related to geometric invariants and
appears as a solution to certain differential equations with geometric and physical 
interpretation. Recently a program has been initiated to lift these constructions
from two- to three-dimensional theories. Beem, Dimofte and Pasquetti argued that the
natural 3D analogue of the differential equations whose solutions determine the partition
function in two-dimensions are $q$-difference equations, \dots}


\bibliographystyle{plain}

\bibliography{ORSRevise}

\printindex


\end{document}